\documentclass[twoside, 10pt]{article}

\usepackage[hmargin= 0.75in, height=9.4in]{geometry}
\usepackage{amssymb,amsthm, amsfonts,amsmath,mathrsfs, amsxtra, url, hyperref}

\usepackage{fancyhdr}
\fancypagestyle{plain}{
 \fancyhf{}

}

\renewcommand{\theequation}{\thesection.\arabic{equation}}
 \numberwithin{equation}{section}

\newtheorem {thm}{Theorem}[section]
\newtheorem {prop}{Proposition}[section]
\newtheorem {lemm}{Lemma}[section]

\newtheorem {rem}{Remark}[section]

\newtheorem {eg}{Example}[section]

\makeatletter
\newenvironment{highlightequation}{%
  \def\tagform@##1{\maketag@@@{(\ignorespaces##1\unskip\@@italiccorr*)}}%
  \ignorespaces
}{%
  \def\tagform@##1{\maketag@@@{(\ignorespaces##1\unskip\@@italiccorr)}}%
  \ignorespacesafterend
}
\makeatother

\def\ba{\begin{array}}
\def\ea{\end{array}}
\def\bea{\begin{eqnarray}}
\def\eea{\end{eqnarray}}
\def\beas{\begin{eqnarray*}}
\def\eeas{\end{eqnarray*}}
\def\bi{\begin{itemize}}
\def\ei{\end{itemize}}
\def\bc{\begin{cases}}
\def\ec{\end{cases}}

\def\bhe{\begin{highlightequation}  }
\def\ehe{\end{highlightequation}  }

\def\a{\alpha}
\def\ga{\gamma}
\def\d{\delta}
\def\e{\varepsilon}
\def\z{\zeta}
\def\k{\kappa}
\def\l{\lambda}
\def\vr{\varrho}
\def\si{\sigma}

\def\t{\tau}
\def\th{\theta}
\def\o{\omega}

\def\vth{\vartheta}

\def\D{\Delta}
\def\G{\Gamma}
\def\L{\Lambda}
\def\O{\Omega}

\def\Th{\Theta}
\def\U{\Upsilon}

\def\bF{{\bf F}}

\def\bd{{\bf d}}

\def\cA{{\cal A}}

\def\cE{{\cal E}}
\def\cF{{\cal F}}
\def\cG{{\cal G}}

\def\cI{{\cal I}}
\def\cJ{{\cal J}}

\def\cN{{\cal N}}

\def\cP{{\cal P}}

\def\cT{{\cal T}}

\def\cY{{\cal Y}}
\def\cZ{{\cal Z}}

\def\hC{\mathbb{C}}

\def\hE{\mathbb{E}}

\def\hN{\mathbb{N}}

\def\hP{\mathbb{P}}
\def\hQ{\mathbb{Q}}
\def\hR{\mathbb{R}}

\def\sE{\mathscr{E}}

\def\sX{\mathscr{X}}
\def\sY{\mathscr{Y}}
\def\sZ{\mathscr{Z}}

\def\fC{\mathfrak{C}}

\def\fI{\mathfrak{I}}

\def\fL{\mathfrak{L}}
\def\fM{\mathfrak{M}}

\def\fP{\mathfrak{P}}

\def\fR{\mathfrak{R}}
\def\fS{\mathfrak{S}}

\def\fU{\mathfrak{U}}

\def\fp{\mathfrak{p}}

\def\fl{\mathfrak{l}}

\def\fg{\mathfrak{g}}

\def\fn{\mathfrak{n}}

\def\ft{\mathfrak{t}}
\def\fs{\mathfrak{s}}
\def\fb{\mathfrak{b}}

\def\ti{\n \times \n}
\def\oti{\n \otimes \n}
\def\df{\n := \n}
\def\ls{\n \le \n}
\def\gs{\n \ge \n}
\def\={\n = \n}
\def\+{\n + \n}
\def\-{\n - \n}
\def\ins{\n \in \n}
\def\ld{\n \land \n}
\def\ve{\n \vee \n}
\def\sb{\n \subset \n}
\def\>{\n > \n}
\def\<{\n < \n}

\def\rY{\rho_{\raisebox{-1.5pt}{\scriptsize $Y$}}}
\def\wrY{\wh{\rho}_{\raisebox{-1.5pt}{\scriptsize $Y$}}}

\def\rX{\rho_{\n \raisebox{-2pt}{\scriptsize $\sX$}}}

\def\({\textnormal{(}}
\def\){\textnormal{)}}
\def\[{[\n[}
\def\]{]\n]}
\def\lan{\langle}
\def\ran{\rangle}

\def\no{\noindent}

\def\q{\quad}
\def\qq{\qquad}

\def\n{\negthinspace}
\def\dn{\n \n}
\def\tn{\n \n \n}

\def\ol{\overline}
\def\ul{\underline}
\def\ua{\mathop{\uparrow}}
\def\da{\mathop{\downarrow}}

\def\wt{\widetilde}
\def\wh{\widehat}
\def\fra{\mathfrak{a}}

\def\dtp{{\hbox{$dt \times d P-$a.s.}}}

\def\pas{{\hbox{$\hP-$a.s.}}}

\def\hb{\hbox}
\def\dis{\displaystyle}
\def\cd{\cdot}
\def\cds{\cdots}

\def\fa{\,\forall \,}

\def\es{\emptyset}

\def\b1{{\bf 1}}
\def\qed{\hfill $\Box$ \medskip}

\def\liminf{\mathop{\ul{\rm lim}}}
\def\limsup{\mathop{\ol{\rm lim}}}

\newcommand{\lsup}[1]{ \underset{#1}{\limsup}}
\newcommand{\linf}[1]{ \underset{#1}{\liminf}}
\newcommand{\lmt}[1]{ \underset{#1}{\lim}}
\newcommand{\lmtu}[1]{ \underset{#1}{\lim} \n \ua \,}
\newcommand{\lmtd}[1]{ \underset{#1}{\lim} \n \da \,}

\begin{document}

 \title{\bf      Optimal Stopping  with Random  Maturity \\ under Nonlinear Expectations
 }

\author{
  Erhan Bayraktar\thanks{ \noindent Department of
  Mathematics, University of Michigan, Ann Arbor, MI 48109; email:
{\tt erhan@umich.edu}.}  \thanks{E. Bayraktar is supported in part by the National Science Foundation,
 and in part by the Susan M. Smith Professorship.  Any opinions, findings, and conclusions
 or recommendations expressed in this material are
those of the authors and do not necessarily reflect the views of the National Science Foundation.} $\,\,$,
$~~$Song Yao\thanks{
\noindent Department of
  Mathematics, University of Pittsburgh, Pittsburgh, PA 15260; email: {\tt songyao@pitt.edu}. } }
\date{}

\maketitle

 \begin{abstract}

  We analyze   an optimal stopping problem
  $ \underset{\ga \in \cT}{\sup} \, \ol{\sE}_0   \big[\sY_{\ga \land \tau_0} \big] $ with random maturity $\tau_0$   under
  a nonlinear expectation $\ol{\sE}_0 [\cd] : = \underset{\hP \in \cP}{\sup} \, \hE_\hP  [\cd]$,
  where  $\cP$ is  a weakly compact  set   of mutually singular probabilities.
  The maturity $\tau_0$ is specified as  the hitting time to level $0$ of some  continuous index process $\sX$
  at which the payoff process $\sY$ is even allowed to have a positive jump.
  When $\cP$  collects a variety of semimartingale measures,
   the optimal stopping problem can be viewed as a {\it discretionary} stopping problem for a
   player who can influence both drift and   volatility of the dynamic of underlying stochastic flow.

    We utilize a martingale approach to construct an optimal pair $(\hP_*,\ga_*)$ for
  $   \underset{(\hP, \ga) \in \cP \times \cT}{\sup} \hE_\hP \big[ \sY_{\ga \land \tau_0} \big] $,
   in which   $\ga_*$ is  the first time  $\sY$ meets the limit $\sZ$
  of its  approximating  $\ol{\sE}-$Snell envelopes.  To overcome the technical  subtleties caused by
  the mutual singularity of  probabilities in $\cP$ and the discontinuity of the payoff process $\sY$,
  we approximate   $\tau_0$ by an increasing sequence of Lipschitz continuous stopping times
  and approximate   $\sY$ by a sequence of uniformly continuous processes.


\end{abstract}

 \smallskip


  \section{Introduction}

 We solve a continuous-time optimal  stopping problem with random maturity $\tau_0$ under
 an nonlinear expectation $\ol{\sE}_0 [\cd] : = \underset{\hP \in \cP}{\sup} \, \hE_\hP   [\cd]$,
 where $\cP$ is  a weakly compact   set   of mutually singular  probabilities
 on the canonical space $\O$  of continuous paths.
 More precisely, letting $\cT$ collect all stopping times with respect to the natural filtration $\bF$
 of the canonical  process $B$ on $\O$,
 we construct in Theorem \ref{thm_RDOSRT} an optimal pair
 $(\hP_*,\ga_*) \n \in \n  \cP  \n \times \n  \cT$ such that
  \bea \label{eq:RDOSRT}
  \underset{(\hP, \ga) \in \cP \times \cT}{\sup}  \,  \hE_\hP [\sY_{\ga \land \tau_0} ]
  =   \underset{\hP \in \cP }{\sup} \,  \hE_\hP [\sY_{\ga_* \land \tau_0} ]
  =  \hE_{\hP_*} \big[ \sY_{\ga_* \land \tau_0} \big] .
  \eea
 Here   the payoff process takes form of
 $\sY_t \n : = \n  \b1_{\{ t < \tau_0\}} L_t  \n + \n  \b1_{\{t \ge \tau_0\}} U_t $, $t  \n \in \n  [0,T]$
 for two bounded processes $L  \n \le \n  U$ that are uniformly continuous in sense of \eqref{eq:aa211},
 and the random maturity $\tau_0$ is the hitting time to level $0$
 of some continuous index process $\sX$ adapted to $\bF$. 
 Writing \eqref{eq:RDOSRT} alternatively as
 \bea   \label{eq:eu021}
  \underset{\ga  \in \cT}{\sup} \, \ol{\sE}_0  [\sY_{\ga \land \tau_0} ]
  = \ol{\sE}_0  [\sY_{\ga_* \land \tau_0} ]   ,
  \eea
 we see that  $\ga_*$ is an optimal stopping time for the optimal stopping with random maturity $\tau_0$ under
 nonlinear expectation $\ol{\sE}_0$.  When $\cP$ collects  measures under which
   $B$ is a semimartingale with uniformly bounded drift and diffusion coefficients
   (in this case, the nonlinear expectation $\ol{\sE}_0$ is the $G$-expectation
   in  sense of Peng \cite{Peng_G_2007b}),
   the optimal stopping problem can be viewed as a {\it discretionary} stopping problem for a
   player who can control both drift and volatility of $B$'s dynamic.

 The optimal  stopping problem with random maturity   under
 the  nonlinear expectation $\ol{\sE}_0$ was first studied by Ekren, Touzi and Zhang \cite{ETZ_2014}
 who took the random maturity to be the first exit time $\hb{\footnotesize H}$ of
   $B$  from some convex open domain $O$ and considered reward processes to have positive that have positive jumps but they do not allow for jumps at
$\hb{\footnotesize H}$, which is the case of interest for us.
Moreover, the convexity of $O$ is a restrictive assumption for the applications we have in mind
   in particular for finding an optimal triplet
   for robust Dynkin game in \cite{RDG}\footnote{The authors would like to thank Jianfeng Zhang for an instructive discussion.}.
    We extend \cite{ETZ_2014} in the following two ways:
   First, $\tau_0$ is more general than $\hb{\footnotesize H}$ so that our   result
   can be at least applied to identify an optimal triplet
   for robust Dynkin game. See also
   Example \ref{eg_rm} for $\tau_0$'s that are the first exit time   of
   $B$  from certain non-convex  domain.
   Second,  we impose a weaker stability under pasting assumption on the probability class
   than the {\it stability under finite pasting} used in \cite{ETZ_2014}.

 Since the seminal work \cite{Snell_1952},  the martingale approach became a primary tool
 in optimal stopping theory (see e.g. \cite{Neveu_1975}, \cite{El_Karoui_1981},
 Appendix D of \cite{Kara_Shr_MF}).  
 Like \cite{ETZ_2014}, we will take a martingale approach with respect to the nonlinear expectation
 $\ol{\sE}_0$. As   probabilities in $\cP$ are mutually singular, one can not define
 the conditional expectation of $\ol{\sE}_0$, and thus
 the Snell envelope of payoff process $\sY$, in  essential supremum sense.
 Instead, we use shifted processes and regular conditional probability distributions
 (see Subsection \ref{subsec:shift_proc} for details)
 to construct the Snell envelope $  \Xi $ of $\sY$ with respect to  pathwise-defined nonlinear expectations
  $  \ol{\sE}_t [\xi] (\o) \n :=  \n   \underset{\hP \in \cP_t }{\sup} \hE_\hP [ \xi^{t,\o}]  $,
  $(t,\o)  \n \in \n  [0,T]  \n \times \n  \O$.
 Here $\cP_t$ is a set of probabilities on the shifted canonical space
 $\O^t $ which includes all regular conditional probability distributions
 stemming from   $\cP$, see (P3).
 In demonstrating  the martingale property of $  \Xi $ with respect to the
 nonlinear expectations $\ol{\sE}  \n = \n  \{\ol{\sE}_t\}_{t \in [0,T]}$,
 we have  encountered  two  major technical difficulties: First,
 no dominating probability in $\cP$ means no  bounded convergence theorem
 for the  nonlinear expectations  $  \ol{\sE} $,
 then one can not follow the classical approach
 for optimal stopping in El Karoui \cite{El_Karoui_1981} to obtain the  $\ol{\sE}-$martingale property of $\Xi$.
  Second, the jump of  payoff process $\sY$   at the random maturity $\tau_0$
 and 
 the discontinuity of  each $\sY_t$ 
 over $\O$ (because of the discontinuity of $\tau_0$) bring technical subtleties  in deriving
 the dynamic programming principle of  $\Xi$, a necessity for
 the  $\ol{\sE}-$martingale property of $  \Xi $.

 To resolve the optimization problem \eqref{eq:RDOSRT}, we first consider the case
 $Y \= L \= U$, however, with a Lipschitz continuous stopping time $\wp$ as the random maturity.
 For the modified payoff process $ \wh{Y}_t \n := \n  Y_{\wp \land t} $, $t  \n \in \n  [0,T]$,
 we construct in Theorem \ref{thm_cst} an optimal pair $(\wh{\hP}, \wh{\nu} )  \n \in \n  \cP  \n \times \n  \cT$ of
 the corresponding  optimization problem
 \bea \label{eq:et431}
   \underset{(\hP, \ga) \in \cP \times \cT}{\sup} \,
      \hE_\hP  \Big[  \wh{Y}_\ga      \Big] = \hE_{\wh{\hP}}  \Big[  \wh{Y}_{\wh{\nu}}      \Big]
 \eea
 such that  $\wh{\nu}$ is the first time $ \wh{Y} $ meets
 its $ \ol{\sE}-$Snell envelope $Z$.
  Using the uniform continuity of $Y$ and
  the Lipschitz continuity of $\wp$, we first derive
  a   continuity estimate \eqref{eq:gd011} of  each $Z_t$ on $\O$, which  leads to
  a dynamic programming principle   \eqref{eq:ef011} of $Z$
  and thus a path continuity estimate \eqref{eq:es017} of process $Z$.
  In virtue of  \eqref{eq:ef011}, we show in Proposition \ref{prop_Z_martingale} that
  $Z$ is an $\ol{\sE}-$supermartingale  and that $Z$ is also an $\ol{\sE}-$submartingale up to
  each approximating stopping time $\nu_n$ of  $\wh{\nu}$, the latter of which shows that for some $\hP_n \n \in \n  \cP$
  \bea \label{eq:es014}
  Z_0 = \ol{\sE}_0 [Z_{\nu_n}] \le \hE_{\hP_n}[Z_{\nu_n}] + 2^{-n} .
  \eea
  Up to a subsequence, $\{\hP_n\}_{n \in \hN}$ has a limit $\wh{\hP} $
  in the weakly compact probability set $\cP$. Then as $n  \n \to \n  \infty$ in \eqref{eq:es014},
  we can deduce  $Z_0  \n = \n  \hE_{\wh{\hP} }[Z_{\wh{\nu} }] $ and thus
  \eqref{eq:et431} by leveraging the continuity estimates
  \eqref{eq:gd011}, \eqref{eq:es017}  of $Z$ as well as a similar argument
  to the one used in the proof of  \cite[Theorem 3.3]{ETZ_2014} that
  replaces  $ \nu_n $'s with a   sequence   of quasi-continuous random variables
  decreasing sequence to $\wh{\nu}$.

 To approximate the general payoff process $\sY$ in problem \eqref{eq:RDOSRT},
 we construct in Proposition \ref{prop_wp_n} an increasing sequence $\{\wp_n\}_{n \in \hN}$
  of Lipschitz continuous stopping times that   converges to $\tau_0$
  and satisfies
  \bea \label{eq:eu011}
   \wp_{n+1}  \n -  \n  \wp_n    \n  \le  \n  \frac{2T}{n+3} , \q n  \n \in \n  \hN .
   \eea
    This result
  together with its premises, Lemma \ref{lem_conti_st} and Lemma \ref{lem_sandwich2}, are
  among the main contributions of this paper. Given $n,k  \n \in \n  \hN$,
  connecting $L$ and $U$ near  $\wp_n$ with lines of slope $2^k$ yields a uniformly continuous process
  $ Y^{n,k}_t  \n :=  \n   L_t  \n + \n \big[ 1 \n \land  \n  (  2^k ( t  \n - \n  \wp_n )  \n - \n  1 )^+ \big]
  ( U_t  \n - \n   L_t )$, $t  \n \in \n  [0,T]$, see Lemma \ref{lem_Y_nkl}.
  Then one can apply Theorem \ref{thm_cst} to $ Y^{n,k} $ and  Lipschitz continuous stopping time
  $\wp^{n,k} \n := \n  (\wp_n \n + \n 2^{1-k})    \n  \land  \n    T $ to find a $\hP_{n,k}  \n \in \n  \cP$ such that
  the $\ol{\sE}-$Snell envelope $Z^{n,k}$ of
  process $ \wh{Y}^{n,k}_t \n := \n  Y^{n,k}_{\wp^{n,k} \land t} $, $t  \n \in \n  [0,T]$  satisfies
  \bea \label{eq:eu014}
   Z^{n,k}_0 =    \hE_{\hP_{n,k}} \Big[ Z^{n,k}_{ \nu_{n,k} \land \z} \Big] , \q \fa \z \in \cT ,
  \eea
 where $\nu_{n,k}$ is the first time $\wh{Y}^{n,k}$ meets $Z^{n,k}$.

 Since $\wh{Y}^{n,k} $ differs from process
 $  \sY^n_t \n := \n \lmt{k \to \infty} \wh{Y}^{n,k}_t  \n = \n    \b1_{\{ t \le \wp_n \}} L_t   \n + \n
 \b1_{\{ t > \wp_n     \}} U_{\wp_n}  $, $\fa t  \n \in \n  [0,T]$
  only over the stochastic interval $\[\wp_n,\wp_n \n + \n 2^{1-k}\]$
  (both processes are stopped after $\wp_n \n + \n 2^{1-k}$),
   the uniform  continuity  of   $L $ and $ U$ gives rise to
     an inequality \eqref{eq:eh137} on how
 $ \wh{Z}^{n,k} $ converges to the $\ol{\sE}-$Snell envelope $ \sZ^n $ of $\sY^n$
 in term of $ 2^{1-k} $.  Similarly, one can deduce from \eqref{eq:eu011} and the uniform  continuity  of   $L , U$
   an estimate \eqref{eq:eh137b}  on
 the distance between  $\sZ^n$ and $\sZ^{n+1}$, which further implies that
 for each $(t,\o) \n \in \n  [0,T]  \n \times \n  \O$,
  $\{\sZ^n_t (\o)\}_{n \in \hN} $ is a Cauchy sequence, and thus admits a limit $\sZ_t (\o)$, see \eqref{eq:et317}.
  We then show in Proposition \ref{prop_Z_ultim}   that $\sZ$ is an $\bF-$adapted continuous process
   that is above the $\ol{\sE}-$Snell envelope of
   the stopped payoff process $\sY^{\tau_0 }$ and
    stays at $U_{\tau_0}$ after the maturity $\tau_0$, so the first time $\ga_*$ when $\sZ$ meets $\sY$
   precedes $\tau_0$.

To prove our main result, Theorem \ref{thm_RDOSRT}, we let $    n  \n < \n  i   \n <  \n   \ell  \n < \n   m  $
so that the stopping time
$ \z_{i,\ell} \n : = \n   \inf\big\{t  \n \in \n  [0,T] \n : Z^{\ell,\ell}_t
  \n \le \n  L_t  \n + \n  1/i   \big\} $
satisfies
$ \z_{i,l}  \n \land \n  \wp_n   \n \le \n  \nu_{m,m}  \n \land \n  \wp_n $.
Applying  \eqref{eq:eh137}, \eqref{eq:eh137b} and \eqref{eq:eu014} with
$(n,k,\z) \n = \n \big( m,m,\z_{i,l} \land \wp_n \big)$ yields
 \bea \label{eq:eu016}
 \sZ_0    \n \le \n  Z^{m,m}_0 \n + \n    \ol{\e}_m   \n \le \n
 \hE_{\hP_{m,m}} \Big[  Z^{m,m}_{  \z_{i,\ell}  \land \wp_n } \Big] \n + \n    \ol{\e}_m
  \n \le \n  \hE_{\hP_{m,m}} \Big[ Z^{\ell,\ell}_{\z_{i,\ell}  \land \wp_n } \Big]
    \n + \n    \ol{\e}_m  \n + \n   \ol{\e}_\ell .
   \eea
   Let $\hP_*$ be the limit of $\{\hP_{m,m}\}_{m \in \hN}$ (up to a subsequence) in the weakly compact probability set $\cP$.
As $m  \n \to \n  \infty$ in \eqref{eq:eu016}, we can deduce
  $ \sZ_0    \n \le \n    \hE_{\hP_*} \Big[ \sZ^{\ell,\ell}_{\z_{i,\ell}  \land \wp_n } \Big]
    \n + \n   \ol{\e}_\ell
    \le   \hE_{\hP_*} \big[ \sZ_{\z_{i,\ell}  \land \wp_n } \big]
    \n + \n   \ol{\e}_\ell   $
 from \eqref{eq:eh137}, \eqref{eq:eh137b},
the continuity estimates \eqref{eq:gd011}, \eqref{eq:es017}   of $Z^{\ell,\ell}$ as well as a similar argument
  to the one used in the proof of  \cite[Theorem 3.3]{ETZ_2014} that approximates
    $ \z_{i,\ell} $ by a decreasing  sequence   of quasi-continuous random variables.
 Then sending $\ell,i,n$  to $\infty$ leads to
  \bea \label{eq:eu018}
 \sZ_0    \n \le \n   \hE_{\hP_*} \big[ \sZ_{\ga_*} \big]  \n = \n  \hE_{\hP_*} \big[ \sY_{\ga_* \land \tau_0} \big]
      \n \le \n     \underset{\hP \in \cP }{\sup} \,  \hE_\hP [\sY_{\ga_* \land \tau_0} ]
    \n \le \n   \underset{(\hP, \ga) \in \cP \times \cT}{\sup}  \,
 \hE_\hP [\sY_{\ga \land \tau_0} ] \n \le  \n  \sZ_0 ,
  \eea
   thus  \eqref{eq:RDOSRT} holds.

 Among our assumptions on the probability class $\{\cP_t\}_{t \in [0,T]}$,
  (P2) is a   continuity condition of the shifted canonical process $B^t$
  that is uniform at each  $\bF^t-$stopping time ($\bF^t$ denotes the natural filtration  of   $B^t$)
  and under each $\hP \n \in \n  \cP_t$.  This condition together with  the  uniform  continuity  of $L,U$
  implies the path continuity \eqref{eq:es017} of $\ol{\sE}-$envelope of any uniformly continuous process
  as well as the aforementioned estimates \eqref{eq:eh137}, \eqref{eq:eh137b}
  about the approximating Snell envelopes $Z^{n,k}$ and $\sZ^n$, all are crucial for the proof of Theorem
  \ref{thm_RDOSRT}.    Another important assumption we impose   on the probability class
  $\{ \cP_t\}_{t \in [0,T]  }$ is the ``weak stability under pasting" (P4),
  which is the key to the supersolution part of the dynamic programming principle \eqref{eq:ef011}
  for the $\ol{\sE}-$envelope of any uniformly continuous process. More precisely,
  (P4) allows us to assemble local  $\e-$optimal controls of the $\ol{\sE}-$envelope
  to form approximating  strategies.
  In Example \ref{eg_weak_formulation}, we show  that these two assumptions
  along with (P3)      are satisfied   by  controls in weak formulation
  i.e. $\cP$ contains all semimartingale measures  under   which $B$ has uniformly bounded drift and
  diffusion coefficients.

  \no    {\bf Relevant Literature.}
 The authors analyzed in \cite{OSNE1,OSNE2} an   optimal stopping problem under
 a non-linear  expectation $\underset{ i  \in \cI  }{\sup} \, \cE_i [\cd] $
 over a filtered probability space $\big(\wt{\O},\wt{\cF}, \wt{\hP}, \wt{\bF} \n = \n \big\{\wt{\cF}_t\big\}_{t \in
[0,T]}\big)$, where
 $  \{\cE_i[\cd |\cF_t]\}_{t \in [0,T]} $ is a   $\wt{\bF}-$consistent
(nonlinear) expectation under $ \wt{\hP} $ for each index $i \n \in \n \cI$.
 A notable example of   $\wt{\bF}-$consistent expectations are the   ``$g$-expectations"
 introduced by \cite{Peng-97}, which represent a fairly large
 class of  convex risk measures, thanks to \cite{CHMP,Pln}.
 If $\cE_i$'s are conditional expected values with   controls,   the
 optimal stopping problem under  $\underset{ i  \in \cI  }{\sup} \, \cE_i [\cd] $
 is exactly the classic  control  problem with discretionary stopping, whose
 general existence/characterization results  can be found in
 \cite{Dubins_Savage_1976, Krylov_CDP, El_Karoui_1981,Bensoussan_Lions_1982,Haussmann_Lepeltier_1990,
 Maitra_Sudderth_1996,Morimoto_2003,Ceci_Bassan_2004,Delbaen_2006,Kara_Zam_2006} among others.
 (For explicit solutions to   applications of such control  problems with discretionary stopping,
 e.g. target-tracking models and
 computation of the upper-hedging prices of American contingent claims under constraints,
 please refer to the literature in \cite{Kara_Zam_2006}.)
 See also \cite{Cadenillas_Sethi_1997, JLK_2004, Choi_Koo_2005}
 for the related optimal consumption-portfolio selection problem with discretionary stopping.
  When the nonlinear expectation becomes
 $\underset{ i  \in \cI  }{\inf} \, \cE_i [\cd] $, the optimal stopping problem considered in \cite{OSNE1,OSNE2}
 transforms to the robust optimal stopping under Knightian uncertainty or the closely related controller-stopper-game,
 which were also extensively studied over the past few decades:
    \cite{Kara_Zam_2005, Kara_Zam_2008, Follmer_Schied_2004, CDK-2006, Delbaen_2006, Riedel_2009,
    OS_CRM, OSNE1, OSNE2,  QRBSDE, riedel2012, Morlais_2008}
    and etc.

 All works cited in the last paragraph assumed that the  probability set $\cP$ is dominated by a single probability
 or that the controller is only allowed to affect the drift.
 When $\cP$ contains mutually singular probabilities  or the controller can influence not only the  drift
 but also the volatility,  there has been a little progress in research due to
 the technical subtleties caused by the mutual singularity of $\cP$
 such as the bounded/dominated convergence theorem generally fails  in this framework.
 Krylov \cite{Krylov_CDP} solved the control problem with discretionary stopping
 in an one-dimensional Markov model with uniformly non-degenerate diffusion,
 however, his approach that relies heavily  on the smoothness of the (deterministic) value function
 does not work in the general case.  In order to  extend  the notion of viscosity solutions
 to the fully nonlinear path-dependent PDEs, as developed in \cite{ETZ_2014_part1,ETZ_2014_part2},
 Ekren, Touzi and Zhang \cite{ETZ_2014} studied
 the optimal  stopping problem with the   random maturity $\hb{\footnotesize H}$  under
 the  nonlinear expectation $\ol{\sE}_0$. Our paper analyzed a similar problem, however, with
 allow for more general forms for $\tau_0$ as explained above.

 In spite of  following its  technical set-up, we adopt a quite different method than \cite{ETZ_2014}:
To estimate the difference between $t-$time Snell envelope values   along two paths $\o, \o' \n \in \n  \O$
satisfying $t  \n < \n  \hb{\footnotesize H}(\o)  \n \land \n   \hb{\footnotesize H}(\o') $, i.e.
  $ \D_t (\o,\o')  \n := \n  \underset{(\hP,\ga) \in \cP_t \times \cT^t}{\sup} \,
      \hE_\hP  \Big[  \ol{Y}^{t,\o}_\ga      \Big]  \n - \n  \underset{(\hP,\ga) \in \cP_t \times \cT^t}{\sup} \,
      \hE_\hP  \Big[  \ol{Y}^{t,\o'}_\ga      \Big]  $ with
      $\ol{Y}_s : = Y_{H \land s} $, $s  \n \in \n  [0,T]$,
 \cite{ETZ_2014}  focuses on all trajectories
 traveling along the straight line $\fl$ from $\o'(t)$   to $\o(t)$ over a short period $[t, t \n + \n \d] $.
 Using a ``stability of finite  pasting" assumption on the probability class $\{\cP_s\}_{s \in [0,T]}$
 \big(which  implies  (P4), see Remark \ref{rem_P2} (3)\big) and the assumption that
 $ \cP_t|_{[t,T-\d]} \subset \cP_{t+\d}$,  \cite{ETZ_2014}
  shifts distributions $\hP $ along these trajectories from time $t$ to time $t \n + \n \d$.
  As $\fl$  is still inside the convex open domain, the stopping time {\footnotesize H} can also be transferred along
 these trajectories with a delay of  $\d$. Then one can use the uniform continuity of $Y$ to estimate
 $ \big| \D_t (\o,\o') \big| $. On the other hand, as described above, we first solve the optimal stopping
  problem with Lipschitz continuous  random  maturity $\wp$   and then approximate the hitting time $\tau_0$
  by Lipschitz continuous stopping times.

As to the robust optimal stopping problem,
or the related controller-stopper-game, with respect to the set $\cP$ of mutually singular probabilities,
Nutz and Zhang \cite{NZ_2015} and Bayraktar and Yao \cite{ROSVU},
used different methods to obtain the existence of the game value and  its martingale property
under the nonlinear expectations
$  \ul{\sE}_t [\xi] (\o):=   \underset{\hP \in \cP_t }{\inf} \hE_\hP [ \xi^{t,\o}]  $, $(t,\o) \in [0,T] \times \O$
(see  the introduction of \cite{ROSVU} for its comparison with \cite{NZ_2015}).
 Such a robust optimal stopping problem   are also considered by, e.g.,
 \cite{Karatzas_Sudderth_2001} and \cite{Bayraktar_Huang_2013} for some particular cases,
    (see also   \cite{ROSVU} for a summary).

    Moreover,   Bayraktar and Yao \cite{RDG} analyzed a robust Dynkin game
 with respect to the set $\cP$ of mutually singular probabilities,
 they show that the Dynkin game has a value  and characterize its $\ul{\sE}-$martingale property.
 Applying the main result of the current paper, Theorem \ref{thm_RDOSRT},
 \cite{RDG} also reaches   an optimal triplet for the robust Dynkin game.

 Very recently, Ekren and Zhang \cite{EkrenZhang1016}  found that our results are useful
 for defining the viscosity solutions  of fully non-linear degenerate path dependent PDEs.

      The rest of the paper is organized as follows:
   Section~\ref{sec:preliminary}   introduces some notation and preliminary results
  such as the regular conditional probability distribution.
  In section \ref{sec:main}, we state   our main result
  on the   optimal stopping problem with random maturity $\tau_0$ under   nonlinear expectation  $\ol{\sE}_0$
  after we impose some assumptions on  the payoff process and the   classes
  $\{\cP_t\}_{t \in [0,T]}$ of mutually singular probabilities.
  In Section \ref{sec:DPP}, we first solve an auxiliary optimal stopping problem with
  uniformly continuous payoff process and Lipschitz continuous
  random maturity under the nonlinear expectation  $\ol{\sE}_0$  by exploring
  the properties of the corresponding $\ol{\sE}-$Snell envelope
  such as   dynamic programming principles it satisfies, the path regularity properties
  as well as the $\ol{\sE}-$martingale characterization.
  In Section \ref{sec:OSRM},
   we approximate the hitting time $\tau_0$  of the index process $\sX$ by Lipschitz continuous stopping times
  and approximate the general payoff process $ \sY $ with discontinuity at
   $\tau_0$ by uniformly continuous   processes. Then
  we show that the convergence of the Snell envelopes of the approximating uniformly continuous  processes
  and    derive the regularity of  their limit,  which is necessary to prove our main result.
  Section \ref{sec:proofs} contains proofs of our results while the demonstration of some auxiliary
  statements with starred labels in these proofs are deferred to the Appendix. We also include
  two technical lemmata in the appendix.

\section{Notation and Preliminaries} \label{sec:preliminary}


 Throughout this paper, we fix $d  \n  \in   \n  \hN$ and a time horizon $T \n \in \n  (0,\infty)  $.
 Let $t  \n \in \n  [0 , T ] $.

 We set $\O^t   \n := \n   \big\{\o    \n  \in   \n   \hC \big([t,T]; \hR^d \big)  \n  :   \o(t)   \n  =  \n   0 \big\}$
 as the   canonical space    over    period    $[t,T]$.
   Given   $ \o  \n \in \n   \O $,
    $\phi^\o_t (x)  \n : = \n  \sup \big\{|\o(r') \n - \n \o(r)| \n :
 r,r'  \n \in \n  [0,t] , \; 0  \n \le \n | r' \n - \n r | \n \le \n  x  \big\}$, $x \in [0,t]$ is clearly a
 modulus of continuity function satisfying  $ \lmtd{x \to 0+} \phi^\o_t (x)  \n = \n  0 $.
     For any $ s  \n \in \n  [ t , T ]  $,
     $ \|\o\|_{t,s}  \n := \n  \underset{r \in [t,s]}{\sup} |\o(r)|  $, $ \fa \o  \n \in \n  \O^t $
     defines  a semi-norm   on   $\O^t   $.
     In particular, $\|\cd\|_{t,T}$ is the {\it uniform} norm on $\O^t$, under which
     $\O^t$ is  a separable complete metric space.

     The canonical process  $ B^t $  of  $\O^t$
   is a  $d-$dimensional standard Brownian motion  under  the   Wiener measure $\hP^t_0$ of $\big(\O^t,  \cF^t_T\big)$.
    Let   $\bF^t \n  =   \n    \{ \cF^t_s     \}_{s \in [t,T]}$,
    with $\cF^t_s   \n  :=   \n  \si \big(B^t_r; r   \n  \in  \n   [t,s]\big)$,   be the  natural filtration of $ B^t $ and   let $\cT^t$   collect all $\bF^t-$stopping times.
    Also, let $\fP_t  $ collect all probabilities   on   $\big(\O^t,  \cF^t_T  \big)   $.
    For any $\hP \n \in \n \fP_t  $ and any sub-sigma-field  $\cG$  of $\cF^t_T$,  we denote by
     $L^1 (\cG , \hP) $    the space of all  real-valued,
$  \cG-$measurable random variables $\xi$ with $\|\xi\|_{L^1(\cG ,\hP)}
   :=   \hE_\hP  \big[ |\xi|  \big]  < \infty$.

    Given $s  \n \in \n  [t, T]$,
    we set $\cT^t_s  \n := \n  \{\tau  \n \in \n  \cT^t \n : \tau (\o)  \n \ge \n  s , \fa \o  \n \in \n  \O^t \}$
    and define   the  {\it truncation}  mapping $\Pi^t_s$  from
 $ \O^t $ to $  \O^s $   by
 $ \big(\Pi^t_s (\o)\big)(r) \n := \n \o (r) \n - \n  \o (s) $,
 $  \fa (r,\o)  \n \in \n  [s, T] \n \times \n  \O^t $.
  By Lemma A.1 of \cite{ROSVU},   $\tau  (\Pi^t_s) \= \tau \circ \Pi^t_s  \n  \in \n  \cT^t_s $,
  $\fa \tau  \n  \in \n  \cT^s$.
  For any $\d  \n > \n  0     $ and $ \o  \n \in \n  \O^t$,
  \bea   \label{eq:bb237}
 O^s_{\d} (\o)  := \big\{\o' \in \O^t:  \| \o'  -   \o  \|_{t,s} < \d \big\}
 \hb{ is an $  \cF^t_s-$measurable open set of }\O^t ,
    \eea
    and $\ol{O}^s_{\d} (\o) \n := \n  \big\{\o'  \n \in \n  \O^t:  \| \o'   \n - \n    \o  \|_{t,s}  \n \le \n  \d \big\}$
   is an $  \cF^t_s-$measurable closed set of $\O^t$
    \big(see e.g.   (2.1) of \cite{ROSVU}\big).
  In particular, we will simply denote $O^T_{\d} (\o)$ and $\ol{O}^T_{\d} (\o)$ by
 $O_{\d} (\o)$ and $\ol{O}_{\d} (\o)$ respectively.

      We will drop the superscript $t$ from the above notations if it is $0$.
  For example,   $  (\O,\cF)    \n =  \n  (\O^0,\cF^0) $.

  We say that $\xi$ is   a continuous random variable   on $\O$
  if for any $\o  \n \in \n  \O$ and $\e  \n > \n  0$,
  there exists a $\d  \n = \n  \d (\o,\e)  \n > \n  0$ such that $|\xi (\o') \n - \n \xi (\o)|  \n < \n  \e$
  for any $\o'  \n \in \n  O_\d (\o)$. Also, $\xi$ is called  a
  Lipschitz continuous random variable   on $\O$ if for some $\k  \n > \n  0$,
  $|\xi (\o') \n - \n \xi (\o)|  \n \le \n  \k \|\o' \n - \n \o\|_{0,T}$ holds for any $\o, \o'  \n \in \n  \O$.

  We say that a process $X$ is   bounded by some $C  > 0$
  if $ |X_t(\o)| \le C   $ for any $ (t,\o) \in [0,T] \times \O $. Also,
  a real-valued  process $X$ is called to be  uniformly   continuous
  on $  [0,T] \times \O$   with respect to some modulus of continuity function $\rho   $ if
       \bea  \label{eq:aa211}
        | X_{t_1}( \o_1) \n - \n   X_{t_2}( \o_2) |
   \n \le \n  \rho  \Big(  \bd_\infty \big((t_1, \o_1), (t_2, \o_2)\big)  \Big)  , \q
   \fa   (t_1, \o_1), (t_2, \o_2) \in [0, T] \times \O   ,
   \eea
 where $\bd_\infty \big((t_1, \o_1), (t_2, \o_2)\big)
 \n := \n  |t_1 \n - \n t_2|  \n + \n  \|\o_1 (\cd  \n \land \n  t_1)  \n - \n  \o_2(\cd  \n \land \n  t_2) \|_{0,T} $.
 For any $t  \n \in \n  [0,T]$, taking $t_1  \n = \n t_2  \n = \n t$ in \eqref{eq:aa211} shows that
  $ \big| X_t( \o_1) \n - \n   X_t( \o_2)\big|  \n \le \n  \rho \big( \|\o_1  \n - \n  \o_2\|_{0,t} \big) $,
  $  \o_1 ,  \o_2  \n \in \n  \O$, which implies the $\cF_t-$measurability of  $X_t$. So
  \bea \label{eq:et014}
  \hb{ $X$ is   indeed an $\bF-$adapted process with all continuous paths.}
  \eea

 Moreover, let $\fM$ denote all modulus of continuity functions $\rho$ such that
 for some $ C  \n > \n  0    $ and
 $ 0  \n < \n  p_1  \n \le \n    p_2  $,
 \bea \label{eq:gb017}
 \rho (x )  \n \le \n    C  (x^{p_1}  \n \vee \n  x^{p_2} )      ,
   \q \fa x  \n \in \n    [0, \infty) .
 \eea
 In this paper, we will frequently use  the convention $ \inf \es := \infty$ as well as  the inequalities
     \bea \label{eq:ec017}
     |x \land a - y \land a| \le |x-y| \q \hb{and} \q |x \vee a - y \vee a| \le |x-y| , \q \fa a, x, y \in \hR .
     \eea

 \subsection{Shifted Processes and Regular Conditional Probability Distributions}
 \label{subsec:shift_proc}

 In this subsection, we fix $ 0 \n \le \n  t   \n \le \n   s   \n \le \n   T$.
 The concatenation $\o  \n \otimes_s \n   \wt{\o}$ of  an  $\o  \n \in \n  \O^t$
 and an $ \wt{\o}  \n \in \n  \O^s$ at time $s$:
 \beas  
 \big(\o \otimes_s  \wt{\o}\big)(r)  :=   \o(r) \, \b1_{\{r \in [t,s)\}}
  + \big(\o(s) + \wt{\o}(r) \big) \, \b1_{\{r \in [s,T]\}} , \q \fa  r \in [t,T]
 \eeas
 defines another path in $\O^t$.
 Set $\o  \n \otimes_s \n  \es   \n = \n \es $ and $  \o  \n \otimes_s \n  \wt{A}  \n := \n
 \big\{ \o   \n \otimes_s \n  \wt{\o} \n : \wt{\o}  \n \in \n   \wt{A} \big\}$
 for  any non-empty  subset   $\wt{A} $ of $  \O^s $.


\begin{lemm}  \label{lem_element}
  If  $ A \in \cF^t_s$, then  $  \o \otimes_s \O^s   \subset A  $   for any $\o \in A  $.

  \end{lemm}

    For any    $ \cF^t_s-$measurable random variable $\eta$,
  since $   \{ \o'  \n \in \n  \O^t  \n : \eta (\o') \n
  = \n  \eta (\o) \} \n \in \n \cF^t_s$,
    Lemma \ref{lem_element} implies that
 \bea   \label{eq:bb421}
  \o  \n \otimes_s \n  \O^s \subset \{\o'  \n \in \n  \O^t  \n : \eta (\o') \n = \n  \eta (\o) \}
 \q \hb{i.e.,} \q
    \eta(  \o \otimes_s \wt{\o} )  \n = \n  \eta(\o), \q
  \fa \wt{\o}  \n \in \n  \O^s .
  \eea
  To wit, the value $ \eta(\o)  $  depends only on $\o|_{[t,s]}$.

   Let $\o  \n \in \n  \O^t$.  For  any   $A  \n \subset \n  \O^t$ we set $A^{s, \o}  \n := \n
   \{ \wt{\o}  \n \in \n  \O^s \n : \o  \n \otimes_s \n  \wt{\o}  \n \in \n  A  \} $
   as the  projection of $A$  on $\O^s $ along $\o$. In particular, $\es^{s,\o}  \n = \n  \es$.
 Given a random variable $\xi$  on $\O^t$,  define  the {\it shift} $\xi^{s,\o}$ of $\xi$ along $\o|_{[t,s]}$
 by  $ \xi^{s,\o}(\wt{\o})  \n := \n  \xi ( \o  \n \otimes_s \n   \wt{\o} ) $,  $  \fa  \wt{\o}  \n \in \n  \O^s $.
 Correspondingly, for a process $X  \n = \n  \{X_r\}_{r \in [t,T]}$ on $\O^t$,
     its {\it shifted} process $X^{s,\o}$ is
      \beas
       X^{s,\o}  (r, \wt{\o}) := (X_r)^{s,\o}(\wt{\o}) = X_r ( \o \otimes_s \wt{\o}) , 
      \q  \fa   ( r, \wt{\o} )  \in [s,T] \times \O^s   .
      \eeas

   Shifted random variables and shifted processes ``inherit" the measurability of   original ones:

 \begin{prop}  \label{prop_shift0}
   Let $ 0 \n \le \n  t   \n \le \n   s   \n \le \n   T$ and     $\o \in \O^t$.

 \no \(1\)  If a  real-valued random variable $\xi $ on $\O^t$ is   $\cF^t_r-$measurable for some $r \in [s,T]$,
   then   $\xi^{s,\o} $ is $  \cF^s_r-$measurable.

   \no \(2\)  For any  $\tau  \n \in \n  \cT^t $,
  if $\t(\o  \n \otimes_s \n  \O^s )  \n \subset \n  [r,T]$
  for some $r   \n \in \n  [s,T]$,  then $\t^{s,\o}  \n \in \n  \cT^s_r   $.

  \no \(3\)  If a real-valued   process $  \{X_r \}_{r \in [t, T]}$
  is $\bF^t-$adapted \(resp. $\bF^t-$progressively measurable\),
 then  $  X^{s,\o} $ is $\bF^s-$adapted \(resp. $\bF^s-$progressively measurable\).

 \end{prop}

  Let  $\hP \n \in \n  \fP_t $. In light of the
   {\it regular conditional probability distributions} (see e.g. \cite{Stroock_Varadhan}),
   we can follow Section 2.2 of \cite{ROSVU} to introduce
   a family of {\it shifted} probabilities $\{\hP^{s,\o}  \}_{\o \in \O^t } \n \subset \n \fP_s $,
   under which the corresponding shifted random variables  inherit the $\hP$ integrability of original ones:

   \begin{prop}  \label{prop_shift1}
  If   $\xi \n \in \n  L^1 \big( \cF^t_T,\hP \big)$ for some  $\hP  \n \in \n  \fP_t$,
  then it holds for $\hP -$a.s.     $\o  \n \in \n  \O^t$  that
   $\xi^{s,\o}  \n \in \n  L^1 \big( \cF^s_T ,  \hP^{s,\o}   \big) $ and
  \bea   \label{eq:f475}
 \hE_{\hP^{s,\o}}  \big[ \xi^{s,\o} \big]= \hE_\hP \big[\xi\big| \cF^t_s\big](\o) \in  \hR  .
    \eea

 \end{prop}

 This subsection was presented in   \cite{ROSVU} with more details and proofs.


 \section{Main Results}

 \label{sec:main}

  In this section,   after   imposing some assumptions on  the payoff process and the   classes
  $\{\cP_t\}_{t \in [0,T]}$ of mutually singular probabilities,
  we will present our main result, Theorem \ref{thm_RDOSRT},
  on the   optimal stopping problem    under the  nonlinear expectation
  $\ol{\sE}_0 [\cd] : = \underset{\hP \in \cP}{\sup} \, \hE_\hP  [\cd]$, whose
  random maturity is in form of the hitting time $\tau_0$ to level $0$ of some continuous index process $\sX$.
  More precisely, let $\sX$ be a  process  with $\sX_0  \n > \n  0$  such that all its paths are continuous and that
  for some modulus of continuity function $\rX$
 \bea \label{eq:es111}
 |\sX_t( \o) - \sX_t( \o')| \le \rX (\|\o-\o'\|_{0,t}) , \q \fa t \in [0,T], ~ \fa \o,\o' \in \O .
 \eea
  Clearly,  \eqref{eq:es111}  implies that the $\bF-$adaptedness of $\sX$.
  Then   $\tau_0  \n : = \n
  \inf\{t  \n \in \n  [0,T] \n : \sX_t  \n \le \n  0\} \n \land \n T  \n \in \n ( 0 , T ]    $
 is an $\bF-$stopping time and 
  \bea \label{def_tau_n}
  \tau_n \n : = \n
  \inf\big\{t  \n \in \n  [0,T] \n : \sX_t  \n \le \n
    \big(  \lceil \log_2(n \n + \n 2) \rceil  \n + \n \lfloor  \sX^{-1}_0 \n \rfloor \n - \n 1 \big)^{-1}  \big\}
  \n \land \n T  \n \in \n ( 0 , \tau_0 ]   , \q n \n \in \n \hN
  \eea
    is an increasing sequence of $\bF-$stopping times
    that converges to $\tau_0$.

   The following example shows that $\tau_0$ could   be the first exit time   of
   $B$  from some non-convex  domain.

   \begin{eg} \label{eg_rm}

1\) Let $d = 2$. Clearly,
$\sX_t = 1 + B^{(2)}_t + \big|B^{(1)}_t \big|$, $t \in [0,T]$ defines   a process with $\sX_0 = 1$ such that
 all its  paths are continuous and that for any $t \in [0,T]$ and $\o,\o' \in \O$,
 $ |\sX_t (\o)-\sX_t (\o') |   \le \big|B^{(1)}_t (\o)- B^{(1)}_t (\o')\big| + \big|B^{(2)}_t (\o)- B^{(2)}_t (\o')\big|
  \le 2 |B_t (\o)- B_t (\o')| \le 2 \|\o-\o'\|_{0,t} $.
 However,   $\tau_0 = \inf\{t \in [0,T]: \sX_t \le 0  \} \land T
 = \inf\{t \in [0,T]: B_t \notin  \U   \}  \land T $ is the first exit time   of
   $B$ from  $\U : = \{(x,y) \in \hR^2: y > -1 - |x|\}$, a non-convex subset of $\hR^2$.

 \no 2\)
Let $d =2$  and let $\G := \{(r \cos \th, r \sin \th): r \in [0,1], \,   \th \in [0, \frac32 \pi]\}$
be the $3/4$ unit disk in $\hR^2$ centered at the origin $(0,0)$.
Clearly, $\sX_t := 1/2- \hb{dist} (B_t,\G)$, $t \in [0,T]$
is a process with $\sX_0 = 1/2$ such that
 all its  paths are continuous and that for any $t \in [0,T]$ and $\o,\o' \in \O$,
 $ |\sX_t (\o)-\sX_t (\o') |  \le | \hb{dist}  (B_t (\o ),\G)-  \hb{dist}  (B_t (\o'),\G)|
 \le |B_t (\o)- B_t (\o')| \le \|\o-\o'\|_{0,t} $.
However, $\tau_0 = \inf\{t \in [0,T]: \sX_t \le 0  \} \land T  = \inf\{t \in [0,T]: B_t \notin \wt{\G}  \} \land T $
is the first exit time   of    $B$ from
 $\wt{\G} : = \{ (x,y) \in \hR^2: dist ((x,y),\G) < 1/2 \} $, another  non-convex subset of $\hR^2$.

 \if{0}
  \no 3\)
 Let $d =2$  and let $\G := \{(x, 0): x \in [0, 1]\} \cup \{(0,y  ): y \in [0, 1]\}$.
 Clearly, $\sX_t := 1/2- \hb{dist} (B_t,\G)$, $t \in [0,T]$
 is a process with $\sX_0 = 1/2$ such that
 all its  paths are continuous and that for any $t \in [0,T]$ and $\o,\o' \in \O$,
 $ |\sX_t (\o)-\sX_t (\o') |  \le | \hb{dist}  (B_t (\o ),\G)-  \hb{dist}  (B_t (\o'),\G)|
 \le |B_t (\o)- B_t (\o')| \le \|\o-\o'\|_{0,t} $.
 However, $\tau_0 = \inf\{t \in [0,T]: \sX_t \le 0  \} \land T  = \inf\{t \in [0,T]: B_t \notin \wt{\G}  \} \land T $
 is the first exit time   of    $B$ from
 $\wt{\G} : = \{ (x,y) \in \hR^2: dist ((x,y),\G) < 1/2 \} $, a  non-convex subset of $\hR^2$ again.
 \fi

   \end{eg}

    \subsection{Uniform Continuity of Payoff Processes}

   \no   {\bf Standing assumptions on payoff processes $(   L ,  U)$.}

   Let    $ L$ and $ U$ be two real-valued    processes bounded by some $M_0  \n > \n  0$ such that

  \no {\bf (A1)} both are   uniformly   continuous   on $  [0,T] \times \O$
   with respect to some   $\rho_0 \in \fM$
   such that $\rho_0$ satisfies \eqref{eq:gb017} with some $\fC>0$ and $0<\fp_1 \le \fp_2  $;

  \no {\bf (A2)}  $ L_t (\o) \le  U_t (\o)$,     $\fa (t,\o) \in [0,T) \times \O$
 and $ L_T (\o)  \n = \n   U_T (\o) $, $\fa \o \in \O$.

 We consider the following payoff process
   \bea \label{eq:et411}
   \sY_t : = \b1_{\{t < \tau_0 \}}  L_t + \b1_{\{t \ge  \tau_0 \}}  U_t , \q t \in [0,T] ,
   \eea
   Clearly, $\sY$ is an $\bF-$adapted process bounded by $M_0$
   whose paths are all continuous except a possible positive jump at $\tau_0$.

 \begin{eg} \label{eg_01}
 1\) \(American-type contingent claims for controllers\)
 Consider   an American-type contingent claim for an agent
 who is able to influence the probability model via certain controls \(e.g. an insider\):
 The claim pays the agent an endowment $ U_{\tau_{\overset{}{\fI}}}$ at
 the first time $\tau_{\overset{}{\fI}}$ when some financial index process $\fI$ rises to certain level $ \fra $
 \(Taking $\sX_t \n = \n  \fra  \n - \n  \fI_t $, $t  \n \in \n  [0,T]$ shows that
 $\tau_0  \n = \n  \tau_{\overset{}{\fI}}$\).
 If the agent chooses to exercise at an earlier time $\ga$ than $\tau_{\overset{}{\fI}}$, she will  receive $L_\ga$.
 Then the price of such an American-type contingent claim is
 $ \underset{(\hP, \ga) \in \cP \times \cT}{\sup}  \,  \hE_\hP [\sY_{\ga \land \tau_0} ] $.

\no  2\) \(robust Dynkin game\)
  \cite{RDG} analyzed a robust Dynkin game with respect to the set $\cP$ of mutually singular probabilities:
  Player 1 \(who  conservatively thinks that the {\it Nature} is not in favor of her\) will receive from Player 2
  a payoff  $ R(\tau,\ga)  \n := \n  \b1_{\{\tau \le \ga\}} \fL_{\tau} \n +  \n \b1_{\{\ga < \tau \}}  \fU_{\ga} $
  if they choose to exit the game at $\tau  \n \in \n  \cT$ and $\ga  \n \in \n  \cT$ respectively.
 The paper shows   that Player 1 has a value in the robust Dynkin game, i.e.
 $ V  \n = \n  \underset{\hP \in \cP  }{\inf} \,  \underset{\ga  \in \cT  }{\inf} \,
  \underset{\t  \in \cT  }{\sup} \,   \hE_\hP  \big[  R  (\t,\ga)    \big]
   \n = \n  \underset{\t  \in \cT  }{\sup} \,  \underset{\ga  \in \cT  }{\inf} \,
   \underset{\hP \in \cP  }{\inf} \,  \hE_\hP  \big[  R  (\t,\ga)    \big]   $ and
    identifies an optimal stopping time $\tau_*$ for Player 1,
    which is the first time Player 1's value process meets $\fL$ \(see Theorem 5.1 therein\). Then the
    robust Dynkin game reduces to the optimal stopping problem with random maturity $\tau_*$ under nonlinear
    expectation $\ol{\sE}_0 := \underset{ \hP  \in \cP  }{\sup}  \,  \hE_\hP [ ~ ] $, i.e.
    $   \underset{(\hP, \ga) \in \cP \times \cT}{\sup}  \,  \hE_\hP [\sY_{\ga \land \tau_*} ] $,
    where $L  \n := \n  - \fU $ and $U  \n := \n  - \fL$.

 \end{eg}

 \subsection{Weak Stability under Pasting}

 Let $\fS$ collect all pairs $(Y,\wp)$ such that

  \no  (i)  $Y$ is a real-valued process   bounded by $M_Y \n > \n 0$
 and uniformly    continuous  on $  [0,T]  \n \times \n  \O$   with respect to some $\rY  \n \in \n  \fM $;

  \no  (ii) $\wp  \n \in \n  \cT $ is a Lipschitz continuous stopping time on $\O$ with coefficient   $\k_\wp  \n > \n  0$:
 $|\wp(\o)  \n - \n  \wp(\o')|  \n \le \n  \k_\wp \|\o \n - \n \o'\|_{0,T}$,  $\fa \o,\o'  \n \in \n  \O  $.

    For any  $(Y,\wp) \in \fS$, we define
 \bea \label{def_whY}
   \wh{Y}_t : = Y_{\wp \land t}  , \q t \in [0,T] ,
 \eea
 which is clearly an $\bF-$adapted  process bounded by $M_Y  $ that has   all continuous paths.

     \no   {\bf Standing assumptions on probability class.}

    We   consider a family $\{\cP_t  \}_{ t \in [0,T]}$ of subsets   of $ \fP_t $, $t \in [0,T]$ such that

    \no \(\textbf{P1}\)    $ \cP := \cP_0        $ is a  weakly compact subset   of $\fP_0$.

   \no ({\bf P2})  For any $\rho \in \fM$,
       there exists another   $ \wh{\rho}  $ of $\fM$ such that
   \bea  \label{eq:ex015}
  \underset{(\hP,\z) \in \cP_t  \times \cT^t}{\sup}  \hE_\hP
 \bigg[  \rho  \Big(    \d +      \underset{ r \in [\z, (\z+\d) \land T]    }{\sup}
  \big|  B^{t}_r - B^{t}_\z  \big|  \Big) \bigg]  \le   \wh{\rho}  (\d) , \q \fa t \in [0,T) , ~
   \fa \d \in  ( 0 , \infty)   .
   \eea
 In particular, we require $\wh{\rho}_0$ to satisfy \eqref{eq:gb017}
 with some $\wh{\fC}>0$ and $1<\wh{\fp}_1 \le \wh{\fp}_2  $.

   \no \(\textbf{P3}\)  For any  $ 0 \le  t <s \le T     $,
 $ \o  \in \O$ and $\hP  \n \in \n  \cP_t$,
   there exists an   extension $(\O^t,\cF',\hP')$ of $(\O^t,\cF^t_T,\hP)$
   \big(i.e. $\cF^t_T \n \subset \n  \cF' $ and $ \hP' |_{\cF^t_T}  \n = \n  \hP$\big) and
   $\O' \in \cF'$ with $\hP' (\O') = 1$ such that
    $\hP^{s,   \wt{\o}} $ belongs to $  \cP_s  $ for any $\wt{\o} \in \O'$.

   \no ({\bf P4})  For any  $(Y,\wp)  \n \in \n  \fS$,
   there exists a modulus of continuity function  $\ol{\rho}_Y  $ such that
 the following statement holds
   for any  $ 0  \n \le \n   t  \n < \n s  \n \le \n  T     $, $ \o   \n \in \n  \O$
   and $\hP  \n \in \n  \cP_t$:
  Given $ \d \n \in \n  \hQ_+   $ and $\l \n  \in \n \hN$,
 let $\{\cA_j\}^\l_{j=0}$ be a $\cF^t_s-$partition of $\O^t$ such that for $j \n = \n 1,\cds \n , \l$,
  $\cA_j  \n \subset \n  O^s_{\d_j} (\wt{\o}_j)$ for some $\d_j  \n \in \n  \big((0,\d]  \n \cap \n  \hQ\big) \cup \{\d\}$ and
   $\wt{\o}_j  \n \in \n  \O^t $. Then for any
    $  \{ \hP_j  \}^\l_{j=1} \n   \subset   \n    \cP_s $,
 there is a $\wh{\hP} \= \wh{\hP} (Y,\wp)    \n \in \n  \cP_t $  such that

 \no (\,i) $\wh{\hP}  (A \cap \cA_0 ) \n =  \n  \hP   (A \cap \cA_0 ) $, $ \fa A \in \cF^t_T$;

   \no  (ii) For any $j \n = \n 1,\cds  \n ,\l $ and $A \in \cF^t_s$,
   $\wh{\hP}  (A \cap \cA_j ) =  \hP   (A \cap \cA_j ) $ and
  \bea    \label{eq:eb127}
    \hE_{\wh{\hP} }
  \Big[ \b1_{A \cap \cA_j} \wh{Y}^{t,\o}_{\ga(\Pi^t_s)}    \Big]
  \n  \ge \n  \hE_{ \hP  } \Big[ \b1_{\{\wt{\o} \in A \cap  \cA_j\}} \Big(
   \hE_{\hP_j}   \big[ \wh{Y}^{s, \o \otimes_t \wt{\o}}_\ga     \big]
    \n - \n  \ol{\rho}_Y (\d)  \Big) \Big] , \q \fa \ga \in \cT^s .
   \eea

 What follows is the main result of this paper on the solvability of the optimization problem  \eqref{eq:RDOSRT}.

 \begin{thm} \label{thm_RDOSRT}
   Assume \eqref{eq:es111},   \(A1\), \(A2\) and \(P1\)$-$\(P4\).
   Then the optimization problem  \eqref{eq:RDOSRT}
   admits an optimal pair  $(\hP_*, \ga_* ) \n \in \n  \cP  \n \times \n  \cT$, where
   the form of $\ga_*$ will be specified in Proposition \ref{prop_Z_ultim} \(4\).

 \end{thm}

  For any $\cF_T-$measurable random variable $\xi$ that is bounded by some $C \n > \n 0$, we define
 its nonlinear  expectations with respect to the  probability class $\{\cP_t\}_{t \in [0,T]  }$  by
    \beas
       \ol{\sE}_t [\xi] (\o):=   \underset{\hP \in \cP_t }{\sup} \hE_\hP [\xi^{t,\o}]  , \q
    \fa  (t,  \o) \in [0,T] \times \O    .
    \eeas
    Then \eqref{eq:RDOSRT} can be alternatively expressed as \eqref{eq:eu021},
    namely,  $\ga_*$ alone is  an optimal stopping time for the optimal stopping   under  nonlinear expectation
    $\ol{\sE}_0$.

 \begin{rem} \label{rem_P2}
 \(1\)  Clearly, $\vr (x) := x $, $\fa x \in [0,\infty)$ is a modulus of continuity function in $\fM$.
 Let     $ \wh{\vr} $ be its corresponding element  in
 $  \fM $  in \(P2\) and assume that $ \wh{\vr} $  satisfies \eqref{eq:gb017} for some $ C_{\n \vr}  \n > \n  0    $ and
 $ 0  \n < \n  q_1  \n \le \n    q_2  $.

  \no \(2\)  Based on \(P2\),    the expectation on the right-hand-side of \eqref{eq:eb127}
 is well-defined since    the mapping $ \wt{\o} \to \n
   \hE_{\wt{\hP}}   \big[ \wh{Y}^{s, \o \otimes_t \wt{\o}}_\ga     \big]$ is continuous  under norm $\|~\|_{t,T}$
   for any $\wt{\hP} \ins \fP_s$ and $\ga \ins \cT^s$.

  \no \(3\) Analogous to the assumption \(P2\) of \cite{ROSVU},
 the condition \(P4\)   can be regarded as a weak form of stability under pasting
 since it  is implied by   the  ``stability under finite pasting"
 \big(see e.g. \(4.18\) of \cite{STZ_2011b}\big):  for any  $0  \n \le \n  t  \n < \n  s  \n \le \n  T$,
 $ \o   \n \in \n  \O$, $\hP  \n \in \n  \cP_t$,
    $ \d  \n \in \n  \hQ_+   $ and $\l  \n \in \n  \hN$,
    let $\{\cA_j\}^\l_{j=0}$ be a $\cF^t_s-$partition of $\O^t$ such that for $j=1,\cds \n , \l$,
  $\cA_j \subset O^s_{\d_j} (\wt{\o}_j)$ for some $\d_j  \n \in \n  \big((0,\d]  \n \cap \n  \hQ\big) \cup \{\d\}$ and $\wt{\o}_j \in \O^t $.
Then for any
    $  \{ \hP_j  \}^\l_{j=1} \n   \subset   \n    \cP_s $,
 the probability defined by
 \bea   \label{eq:xxx131c}
    \wh{\hP} (A)   \n =  \n   \hP ( A \cap    \cA_0  \big)
     +   \sum^\l_{j=1} \hE_\hP   \n  \left[   \b1_{\{\wt{\o} \in \cA_j\}}  \hP_j \big( A^{s,\wt{\o}} \big) \right]
          , \q \fa  A \in \cF^t_T
   \eea
 is in $  \cP_t $.

\end{rem}

\begin{eg} \label{eg_weak_formulation}
 \(Controls of weak formulation\) Given $\ell > 0$, let $\{\cP^\ell_t\}_{t \in [0,T]}$
 be the family of semimartingale measures considered in \cite{ETZ_2014} such that
 $\cP^\ell_t$   collects all continuous semimartingale measures   on $(\O^t,\cF^t_T)$,
 whose drift and diffusion characteristics are bounded by $\ell$ and $\sqrt{2\ell}$ respectively.
 According to Lemma 2.3 therein,  $\{\cP^\ell_t\}_{t \in [0,T]}$ satisfies
 \(P1\), \(P3\) and stability under finite pasting \big(thus \(P4\) by Remark \ref{rem_P2} \(3\)\big).
 Also,   one can  deduce from the Burkholder-Davis-Gundy inequality that  $\{\cP^\ell_t\}_{t \in [0,T]}$ satisfies
 \(P2\), see Section \ref{sec:proofs} for details.
 \end{eg}

 \section{Optimal Stopping  with Lipschitz Continuous Random  Maturity  }
  \label{sec:DPP}

   To  solve the optimization problem \eqref{eq:RDOSRT}
    we first analyze in this section   an auxiliary optimal stopping problem with
  uniformly continuous payoff process and Lipschitz continuous
  random maturity under the nonlinear expectation $\ol{\sE}_0$.

    Let the     probability class $\{\cP_t\}_{t \in [0,T]  }$ satisfy (P2)$-$(P4).
  To solve \eqref{eq:RDOSRT}, we first consider the case
 $Y \= L \= U$   with random maturity $\wp$ for some   $(Y,\wp)  \n \in \n  \fS$.
    For any $(t,\o)  \n \in \n  [0,T]  \n \times \n  \O$,  define
       \beas
  Z_t (\o)  :=   \underset{(\hP,\ga) \in \cP_t \times \cT^t}{\sup} \,
      \hE_\hP  \Big[  \wh{Y}^{t,\o}_\ga      \Big]
    \eeas
  as a Snell envelope of the payoff process $\wh{Y}$
  with respect to the nonlinear expectations $ \ol{\sE} \n = \n \{\ol{\sE}_{\n t}\}_{t \in [0,T]} $
  given the historical path $\o|_{[0,t]}$. We will simply refer to $Z$ as the $ \ol{\sE}-$Snell envelope of $\wh{Y}$.
  Since the  $\bF-$adaptedness of $\wh{Y}$ and \eqref{eq:bb421} imply that
 $\wh{Y}^{t,\o}_t (\wt{\o})  \n = \n  \wh{Y}_t (\o \otimes_t \wt{\o})   \n = \n  \wh{Y}_t (\o) $,
 $\fa \wt{\o}  \n \in \n  \O^t$, one has
 \bea \label{eq:ef221}
  M_Y   \ge Z_t (\o) \ge    \underset{\hP \in \cP_t}{\sup}
  \hE_\hP \Big[\wh{Y}^{t,\o}_t \Big] = \wh{Y}_t (\o) \ge - M_Y     , \q
  \fa   (t,\o) \in [0,T] \times \O    .
 \eea

 Given $t \n \in \n  [0,T]$, we have the following      estimate on the continuity of random variable $Z_t$ at
each $ \o \in \O$, which is    not only in term of the distance from $\o$ under
$\|~\|_{0,t}$ but also in term of the path information of $\o$ up to time $t$.

     \begin{prop} \label{prop_Z_conti_in_o}
     Assume \(P2\).  Let $(Y,\wp)  \n \in \n  \fS$ and $(t,\o) \n \in \n  [0,T]  \n \times \n  \O $.
 It holds for any $\o'  \n \in \n  \O$
      \bea  \label{eq:gd011}
     |  Z_t (\o ) \n - \n  Z_t (\o') |
   \n \le  \n      \wrY \Big( (1 \n + \n \k_\wp) \|\o \n - \n \o'\|_{0,t}
  \n + \n  \underset{r \in [ t_1 , t_2 ]}{\sup} \big|    \o    (r  )   \n - \n    \o    ( t_1 ) \big|  \Big)
   \n  \le  \n      \wrY \Big( (1 \n + \n \k_\wp) \|\o \n - \n \o'\|_{0,t}
  \n + \n \phi^\o_t \big( \k_\wp   \|\o   \n  - \n  \o'  \|_{0,t} \big)   \Big)  , \q
   \eea
   where $t_1 : = \wp (\o) \land \wp (\o') \land t $ and $t_2 : = \big(\wp (\o) \vee \wp (\o')\big) \land t $.
   Consequently, $  Z_t(\o)$ is continuous in $\o$ under the norm $\|~\|_{0,t}$: i.e. for any $\e > 0$, there exists a
   $ \, \d = \d (t,\o) > 0$ such that
   \bea \label{eq:gd014}
   |  Z_t (\o' ) - Z_t (\o) | < \e , \q \fa \o' \in O^t_\d (\o) ,
   \eea
   and thus    $Z_t$ is   $\cF_t-$measurable.
    \end{prop}

 The resolution of the auxiliary  optimization problem \eqref{eq:et431}
 with the payoff process $Y$ and random maturity $\wp$
 relies on  the following  dynamic programming principle   for
 the $ \ol{\sE}-$Snell envelope $Z$ of $\wh{Y}$
 and a consequence of it,  a path continuity estimate   of process $Z$:

   \begin{prop} \label{prop_DPP_Z}
   Assume \(P2\)$-$\(P4\).
   Let $(Y,\wp) \in \fS$.   It holds for any $(t,\o) \in [0,T] \times \O$ and $\nu \in \cT^t$ that
 \bea    \label{eq:ef011}
 Z_t (\o) =   \underset{(\hP,\ga) \in \cP_t \times \cT^t}{\sup}   \,     \hE_\hP
\Big[ \b1_{\{  \ga < \nu \}}   \wh{Y}^{t,\o}_\ga  + \b1_{\{  \ga \ge \nu \}} Z^{t,\o}_\nu    \Big]  .
 \eea
 Consequently, $Z$ is   an $\bF-$adapted  process bounded by $M_Y  $ that has   all continuous paths.
 More precisely,    for    any   $ \o   \n \in \n   \O$ and
     $0  \n \le \n  t  \n \le \n  s  \n \le \n  T$,
  \bea   \label{eq:es017}
  \big| Z_t (\o) \n - \n  Z_s (\o) \big|    \le  2 C_{\n \vr} M_Y   \Big( (s \n - \n t)^{ \frac{q_1}{2}}  \n \vee  \n  (s \n - \n t)^{q_2-\frac{q_1}{2}} \Big)  \n + \n \wrY (s \n - \n t)
   \n + \n    \wrY  \big(   \d_{t,s} (\o)  \big)
    \n \vee \n  \wh{\wh{\rho}\,}_Y \big(   \d_{t,s} (\o)  \big)       ,
  \eea
  where  $\d_{t,s} (\o) \n := \n   (1 \n + \n \k_\wp)
 \Big(  (s \n - \n t)^{\frac{q_1}{2}}  \n + \n
  \underset{t\le r < r' \le s  }{\sup} \big| \o(r' ) \n - \n  \o(  r) \big|  \Big)   $
  and $ \wh{\wh{\rho}\,}_Y \ins \fM $ is the modulus of continuity function corresponding to $\wrY$ in \(P2\).
  See Remark \ref{rem_P2} \(1\) for the notations $C_{\n \vr}, q_1$ and $q_2$ here.
 \end{prop}

 In light of Proposition \ref{prop_DPP_Z},  the $\ol{\sE}-$Snell envelope $Z$ of $\wh{Y}$ has the following
 $\ol{\sE}-$martingale properties:

\begin{prop} \label{prop_Z_martingale}
 Assume \(P2\)$-$\(P4\).
 Let $ (Y,\wp) \in \fS $ and $n \in \hN$. Then $Z$ is an $\ol{\sE}-$supermartingale,
 and  $Z$ is an $\ol{\sE}-$submartingale over
$[0, \nu_n]$  in sense that for any $\z \in \cT$
\beas
    Z_{   \z \land t} (\o) \ge \ol{\sE}_{\n t} [Z_{  \z}] (\o) \q \hb{and} \q
      Z_{ \nu_n \land \z \land t} (\o) \le \ol{\sE}_{\n t} [Z_{ \nu_n \land \z}] (\o)
      , \q \fa (t,\o) \in [0,T] \times \O ,
\eeas
where $ \nu_n   \n : = \n  \inf \big\{ t  \n \in \n  [0,T] \n : Z_t  \n - \n  \wh{Y}_t
\n \le \n \frac{1}{n} \big\}  \n \in \n  \cT $.

\end{prop}

 Exploiting the $\ol{\sE}-$submartingale of $Z$ up to $\nu_n$ as well as the
 continuity estimates \eqref{eq:gd011}, \eqref{eq:es017} of $Z$, we can
 solve the optimization problem \eqref{eq:et431} by taking a similar argument
 to the one used in the proof of  \cite[Theorem 3.3]{ETZ_2014}.

 \begin{thm} \label{thm_cst}
  Assume  \(P1\)$-$\(P4\) and let   $(Y,\wp)  \n \in \n  \fS$. There exists
  a $\wh{\hP}  \n \in \n  \cP$ such that
  $  Z_0    \n = \n \ol{\sE}_{\n 0} \big[ Z_{\wh{\nu}  } \big]
   \n = \n  \hE_{\wh{\hP}} \big[ Z_{\wh{\nu}}   \big]
   \n = \n  \hE_{\wh{\hP}} \big[ \wh{Y}_{\wh{\nu}}   \big]   $,
  where $\wh{\nu} \n : = \n  \inf \big\{ t  \n \in \n  [0,T] \n : Z_t  \n = \n  \wh{Y}_t \big\}  \n \in \n  \cT$.
  To wit, $ \big(\wh{\nu},\wh{\hP} \big)$ solves
  the  optimization problem \eqref{eq:et431}  with the payoff process $\wh{Y}$.
  Moreover, it holds for any $\z \in \cT$ that   $ Z_0 = \ol{\sE}_{\n 0} \big[ Z_{  \wh{\nu} \land \z  } \big]
  = \hE_{\wh{\hP}} \big[ Z_{\wh{\nu} \land \z} \big] $.

 \if{0}
 The $\ol{\sE}-$martingale of $Z$ requires $\cP_t$ is weakly compact for any $t \in [0,T]$.
 \fi

 \end{thm}

 \section{Optimal Stopping  with Random  Maturity $\tau_0$  }

 \label{sec:OSRM}

  In this section,  we approximate the hitting time $\tau_0$ of the index process $\sX$
  by Lipschitz continuous stopping times  and approximate the general payoff process $ \sY $
  in \eqref{eq:et411} by uniformly continuous processes.
  We show   the convergence of the Snell envelopes of the approximating uniformly continuous  processes
  and    derive the regularity of  their limit,
  which is crucial for the proof of our main result, Theorem \ref{thm_RDOSRT}.

   To apply  Theorem  \ref{thm_RDOSRT}, we first  approximate $\tau_0$ by
   an increasing sequence $\{\wp_n\}_{n \in \hN}$ of Lipschitz continuous stopping times such that
   the increment  $\wp_{n+1}  \n - \n   \wp_n  $ uniformly decreases   to $0$ as $n  \n \to \n  \infty$:

 \begin{prop} \label{prop_wp_n}
 Assume \eqref{eq:es111}.
   There exist   an increasing sequence $\{\wp_n\}_{n \in \hN}$ in $\cT$
   and  an increasing sequence $\{\k_n\}_{n \in \hN}$ of positive numbers with
   $\lmtu{n \to \infty} \k_n \n = \n  \infty  $
   such that for any $n  \n \in \n  \hN$

  \no \(1\)  $ \tau_n (\o) \n \le  \n  \wp_n (\o) \n \le  \n  \tau_0 (\o) $ and
  $ 0 \le   \wp_{n+1} (\o) -  \wp_n (\o)  \le  \frac{2T}{n+3} $, $\fa \o \in \O$.
  In particular, if $\{t  \n \in \n  [0,T] \n : \sX_t  (\o')  \n \le \n  0\}$  is not empty for some $\o' \in \O$,
  then $\wp_n (\o') \n <  \n  \tau_0 (\o')$.

  \no \(2\)  Given   $\o_1, \o_2  \n \in \n  \O$,
 $ \big| \wp_n (\o_1)    \n  -  \n    \wp_n (\o_2) \big|
  \n \le \n    \k_n   \|\o_1    \n  -  \n    \o_2 \|_{0,t_0} 
  $
 holds  for any  $t_0  \n \in \n   \big\{ t  \n \in \n  [a_n, T) \n : t  \n \ge \n  a_n
   \n + \n   \k_n  \| \o_1   \n - \n   \o_2 \|_{0,t} \big\} \cup \{T\} $,
 where $ a_n  \n := \n  \wp_n (\o_1)  \n \land \n  \wp_n (\o_2) $.

 \end{prop}

 Let $n, k \in \hN $ and let  $\wp_n  $ be the $\bF-$stopping time stated in Proposition \ref{prop_wp_n}.
 We use lines of slope $2^k$  to connect $L$ and $U$ near  $\wp_n$ as follows:  For any $t \in [0,T]$,
 \bea
 \hspace{-3mm}  Y^{n,k}_t  & \tn \dn :=  & \tn \dn   L_t + \big[ 1 \land  (  2^k ( t - \wp_n ) - 1 )^+ \big]
  ( U_t -  L_t ) \label{eq:eh110} \\
    & \tn  \dn   =  & \tn  \dn    \b1_{\{  t \le \wp_n   +   2^{-k}   \}}  L_t
 \n + \n  \b1_{\{    \wp_n +   2^{-k}     < t < \wp_n +   2^{1-k}     \}}
   \Big\{  \big[ 1 \n - \n   2^k  ( t \n - \n \wp_n  \n - \n    2^{-k}   )  \big]    L_t  \n + \n
   2^k ( t \n - \n \wp_n  \n - \n    2^{-k} )     U_t   \Big\}
   \n  + \n  \b1_{\{ t \ge \wp_n +   2^{1-k}   \}}  U_t ,  \qq  \label{eq:eh111}
 \eea
 where the set $\{    \wp_n  (\o)  \n + \n    2^{-k}       \n < \n  t  \n < \n  \wp_n  (\o)  \n + \n    2^{1-k}  \}$
 \big(resp.  $\{ t  \n \ge \n  \wp_n  (\o)  \n + \n    2^{1-k}   \}$\big) may be empty if
 $\wp_n (\o)   \n + \n    2^{-k}     \n \ge \n  T$
 \big(resp.  $ \wp_n (\o)   \n + \n    2^{1-k}  \n > \n  T$\big) for some $\o  \n \in \n  \O$.

   Clearly, the process $Y^{n,k}$ is also bounded by $M_0$, and it is uniformly continuous on $  [0,T]  \n \times \n  \O$
   with respect to some $\rho_{n,k} \in \fM$:

  \begin{lemm} \label{lem_Y_nkl}
 Assume   \eqref{eq:es111} and \(A1\). For any  $n ,   k \n \in \n  \hN$,
   $  Y^{n,k} $ is uniformly continuous on $  [0,T]  \n \times \n  \O$
 with respect to the modulus of continuity function $\rho_{n,k} (x)
  \n := \n  6 \rho_0 (2x)   \n + \n  2^{1+k}  M_0 (1 \n + \n \k_n)  x
   \n \le \n  C_{n,k} \big( x^{\fp_1 \land 1}  \n \vee \n  x^{\fp_2 \vee 1} \big) $,
     $\fa x  \n \in \n  [ 0 , \infty ) $,
  where $C_{n,k}  \n : = \n  6  \n \cd \n  2^{\fp_2} \fC \n + \n  2^{1+k}  M_0 (1 \n + \n \k_n) $
  and $\{\k_n\}_{n \in \hN}$ is the   increasing sequence  of positive numbers
    in Proposition \ref{prop_wp_n}.

   \end{lemm}

 Applying Proposition \ref{prop_wp_n} (2) with $t_0 \n = \n T$ shows that $\wp_n$
 is a Lipschitz continuous stopping time on $\O$ with coefficient $\k_n$, so is
 $\wp^{n,k} \n := \n  (\wp_n \n + \n 2^{1-k})    \n  \land  \n    T   $
 by \eqref{eq:ec017}. Then  we define
  \bea
 \wh{Y}^{n,k}_t & \tn \dn := & \tn  \dn  Y^{n,k}_{\wp^{n,k} \land t}
 \n = \n   Y^{n,k}_{  (\wp_n+2^{1-k}) \land t }
   \n   =  \n   \b1_{\{  t \le \wp_n   +   2^{-k}   \}} L_t
 \n + \n  \b1_{\{    \wp_n +   2^{-k}     < t < \wp_n +   2^{1-k}     \}}
   \Big\{  \big[ 1 \n - \n   2^k  ( t \n - \n \wp_n  \n - \n    2^{-k}   )  \big]   L_t
    \n +  \n    2^k ( t \n - \n \wp_n  \n - \n    2^{-k} )    U_t   \Big\}  \nonumber \\
 & \tn  & \tn  +    \b1_{\{ t \ge \wp_n +   2^{1-k}   \}} U_{(\wp_n+2^{1-k}) \land T} ,
   \q \fa t  \n \in \n  [0,T]  ,       \label{eq:et031}
 \eea
 and its $\ol{\sE}-$Snell envelope:
 \beas
  Z^{n,k}_t (\o) : = \underset{(\hP,\ga) \in \cP_t \times \cT^t}{\sup}
  \hE_\hP \Big[ \big(\wh{Y}^{n,k} \big)^{t,\o}_\ga \Big] , \q  \fa  (t,\o) \in [0,T] \times \O .
 \eeas
 As $L $ and $  U$ are bounded by $  M_0 $, so are $Y^{n,k}$ and $Z^{n,k}$ by \eqref{eq:ef221}.
 In light of Lemma \ref{lem_Y_nkl}, we can apply the results in Section \ref{sec:DPP} to each
 $Z^{n,k}$, $n,k \n \in \n \hN$.

   Given  $  t  \n \in \n  [0,T]$,
    \eqref{eq:eh110} shows that  $\lmtu{k \to \infty} Y^{n,k}_t  \n = \n    \b1_{\{ t \le \wp_n \}} L_t   \n + \n
 \b1_{\{ t > \wp_n     \}} U_t  $. Since
 \beas
  \lmt{k \to \infty} \big(\wh{Y}^{n,k}_t - Y^{n,k}_t \big)
 = \lmt{k \to \infty}   \b1_{\{ t \ge \wp_n +   2^{1-k}   \}} (U_{(\wp_n+2^{1-k}) \land T} - U_t)
 =  \b1_{\{ t > \wp_n     \}} (U_{ \wp_n   } - U_t)
 \eeas
  by the continuity of $U$, we see that
 $  \sY^n_t \n := \n  \lmt{k \to \infty} \wh{Y}^{n,k}_t  \n = \n    \b1_{\{ t \le \wp_n \}} L_t   \n + \n
 \b1_{\{ t > \wp_n     \}} U_{\wp_n}  $, $\fa t  \n \in \n  [0,T]$,
 which is   an $\bF-$adapted process with all c\`agl\`ad paths.
 For any $(t,\o)  \n \in \n  [0,T]  \n \times \n  \O$,   Proposition \ref{prop_shift0} (3) shows that
 $(\sY^n)^{t,\o}$ is an $\bF^t-$adapted process with all c\`agl\`ad paths
 and thus an $\bF^t-$progressively measurable process.
 Then we can consider the following $\ol{\sE}-$Snell envelope of $\sY^n$:
 \beas
 \sZ^n_t(\o)  : = \underset{(\hP,\ga) \in \cP_t \times \cT^t}{\sup} \,
 \hE_\hP \big[ (\sY^n)^{t,\o}_{\ga} \big] , \q \fa (t,\o) \in [0,T] \times \O .
 \eeas
 Again,   $\sY^n$ and $\sZ^n$ are   bounded by $  M_0$.

     The next two inequalities show  how $ Z^{n,k} $ converges to   $ \sZ^n $
 in term of $ 2^{1-k} $ and how $\sZ^n$ differs from $\sZ^{n+1}$,   both inequalities also depend on
 the historical path of process $U$.

 \begin{prop} \label{lem_Z_approx}
 Assume \eqref{eq:es111}, \(A1\), \(A2\), \(P2\) and let $n,k \in \hN$.
 It holds for any  $(t,\o) \in [0,T] \times \O$ that
   \bea \label{eq:eh137}
   - 2 \wh{\rho}_0 \big( 2^{1-k}   \big) & \tn  \le & \tn
     Z^{n,k}_t(\o) -  \sZ^n_t(\o) - U \big( ( \wp_n (\o)  \n + \n  2^{1-k} ) \land t  ,   \o    \big)
      +     U  (\wp_n (\o) \land t , \o   )
  \le    \wh{\rho}_0 \big( 2^{1-k}   \big)   \qq           \\
  \hb{and} \q -  2  \wh{\rho}_0 \big( \hb{$\frac{2T}{n+3}$} \big) & \tn  \le & \tn
     \sZ^{n+1}_t(\o) \n - \n  \sZ^n_t(\o)
  \n - \n  U \big(  \wp_{n+1} (\o)    \n \land \n  t , \o \big)
  \n + \n U \big(  \wp_n (\o)    \n \land \n  t , \o \big)
        \le          \wh{\rho}_0 \big( \hb{$\frac{2T}{n+3}$} \big) .  \label{eq:eh137b}
 \eea

 \end{prop}

 As $\wh{\rho}_0$ satisfies \eqref{eq:gb017} with some $\wh{\fC}>0$ and $1<\wh{\fp}_1 \le \wh{\fp}_2  $
 by (P2), we see from \eqref{eq:eh137b} that
  for each $(t,\o) \n \in \n  [0,T]  \n \times \n  \O$,
  $\{\sZ^n_t (\o)\}_{n \in \hN} $ is a Cauchy sequence, and thus admits a limit $\sZ_t (\o)$.
 The following results shows that $\sZ$ is an $\bF-$adapted continuous  process
 above the Snell envelope of the stopped payoff process $\sY^{\tau_0 }$ and that the first time
 $\sZ$ meets $\sY$ is exactly the optimal stopping time expected in Theorem \ref{thm_RDOSRT}.

 \begin{prop} \label{prop_Z_ultim}
 Assume \eqref{eq:es111}, \(A1\), \(A2\) and \(P2\)$-$\(P4\).

  \no \(1\) For any $n \n \in  \n  \hN$,
 $\sZ^n$ is an $\bF-$adapted process bounded by $M_0$ that has all continuous paths.

  \no \(2\)  For any  $(t,\o) \in [0,T] \times \O$, the limit $\sZ_t (\o) : = \lmt{n \to \infty} \sZ^n_t (\o) $ exists
 and satisfies
 \bea \label{eq:et317}
 -  2 \e_n  \le   \sZ_t(\o) \n - \n  \sZ^n_t(\o)
  \n - \n  U \big(  \tau_0 (\o)    \n \land \n  t , \o \big)
  \n + \n U \big(  \wp_n (\o)    \n \land \n  t , \o \big)   \le   \e_n , \q \fa n \in \hN ,
 \eea
 where $\e_n : = \sum^\infty_{i=n} \wh{\rho}_0 \big( \hb{$\frac{2T}{i+3}$} \big) $
 decreases to $0$ as   $n \to \infty$.

  \no \(3\)    $\sZ $ is   an $\bF-$adapted process bounded by $M_0$ that has all continuous paths.
 Set $ \wh{\sY}_t   \n  : = \n  \sY_{\tau_0 \land t} $, $t \ins [0,T]$.
   It holds for any $\o \n \in \n \O$ that
   \bea \label{eq:et335}
   \wh{\sY}_t (\o)   \n \le  \n
    \underset{(\hP,\ga) \in \cP_t \times \cT^t}{\sup} \,
 \hE_\hP \Big[ \wh{\sY}\,^{t,\o}_{\ga} \Big]  \n \le \n  \sZ_t (\o) , ~   \fa t  \n \in \n  [0,T] ~ \hb{and} ~
 \sZ \big(t,\o\big)  \n = \n  U \big(\tau_0(\o),\o\big) , ~   \fa t  \n \in \n  [\tau_0(\o),T].
   \eea

  \no \(4\) $\ga_*  \n : = \n  \inf \big\{t  \n \in \n  [0, T] \n : \sZ_t    \n = \n  \wh{\sY}_t \, \big\}
   \n = \n  \inf\{t  \n \in \n  [0,\tau_0): \sZ_t  \n = \n  L_t \}   \n \land \n  \tau_0$  is an $\bF-$stopping time.

 \end{prop}

  \section{Proofs}
  \label{sec:proofs}

    \subsection{Proofs of Results in Section \ref{sec:main}}

     \label{subsect:proof_ROSVU}

  \no {\bf Proof of Remark \ref{rem_P2}:}
  {\bf 2)} Let $(Y,\wp)  \n \in \n  \fS$,  $\wt{\hP} \ins \fP_s$  and $\ga \ins \cT^s$.
  Given  $ \wt{\o}    \n \in \n  \O^t  $,
  Lemma \ref{lem_Y_diff} shows that
 \beas
 \hE_{\wt{\hP}} \Big[ \big| \wh{Y}^{s,\o \otimes_t \wt{\o}}_\ga \n - \n  \wh{Y}^{s,\o \otimes_t \wt{\o}'}_\ga \big| \Big]
  &  \tn \le &  \tn    \wrY \Big( (1 \n + \n \k_\wp) \|\o  \n \otimes_t \n  \wt{\o}
  \n - \n \o  \n \otimes_t \n  \wt{\o}'\|_{0,s}
  \n + \n \phi^{\o \otimes_t \wt{\o}}_s
  \big( \k_\wp   \|\o  \n \otimes_t \n  \wt{\o}   - \o  \n \otimes_t \n  \wt{\o}'  \|_{0,s} \big)   \Big) \\
  &  \tn \le &  \tn   \wrY \Big( (1 \n + \n \k_\wp) \|  \wt{\o} \n - \n   \wt{\o}'\|_{t,T}
  \n + \n \phi^{\o \otimes_t \wt{\o}}_s  \big( \k_\wp   \|  \wt{\o}   \n  -  \n   \wt{\o}'  \|_{t,T} \big)   \Big)  ,
  \q \fa \wt{\o} '    \n \in \n  \O^t .
 \eeas
 Hence,
 the mapping $\wt{\o}  \to \hE_{\wt{\hP}} \Big[ \wh{Y}^{s,\o \otimes_t \wt{\o}}_\ga  \Big] $
 is continuous   under norm $\|~\|_{t,T}$ and thus $\cF^t_T-$measurable.

  \no {\bf 3)} Similar to the proof of Remark 3.3 (2) in \cite{ROSVU}, one can   show that
 the probability $\wh{\hP}$ defined in \eqref{eq:xxx131c} satisfies (P4) (i) and the first part of (P4) (ii): i.e.
 $ \wh{\hP}  (A  \n \cap \n  \cA_0 ) \n =  \n  \hP   (A  \n \cap \n  \cA_0 )  $, $ \fa   A  \n \in \n  \cF^t_T $, and $
 \wh{\hP}  (A  \n \cap \n  \cA_j )  \n = \n   \hP  (A  \n \cap \n  \cA_j )$, $ \fa  j \n = \n 1,\cds  \n ,\l $,
   $ \fa    A  \n \in  \n  \cF^t_s $.
 \if{0}
 Given $  A \n \in \n  \cF^t_T$, for any $j  \n = \n  1,  \cds \n , \l $ and $\wt{\o}  \n \in \n  \cA_j $,
    since $\cA_j  \n \in \n  \cF^t_s$,
    Lemma \ref{lem_element}  shows that $( \cA_j  )^{s , \wt{\o}} = \O^s $,
    which implies that  $(A \cap \cA_0 )^{s , \wt{\o}} = \es$. Then one can easily calculate  that
     $\wh{\hP}  (A \cap \cA_0 ) =  \hP   (A \cap \cA_0 ) $.

       Next, let $j = 1,  \cds \n , \l $ and  $A \in \cF^t_s$.
  We see from Lemma \ref{lem_element}  again that
  \bea  \label{eq:xxx614}
 \hb{ if $\wt{\o} \in A \cap \cA_j$ (resp. $\notin A \cap \cA_j$),
 then $( A \cap \cA_j )^{s , \wt{\o}}  = \O^s  $ (resp. $= \es $). }
 \eea
    It follows that
\beas
\wh{\hP}  (A \cap \cA_j ) =  \sum^\l_{j'=1} \hE_\hP   \n  \left[   \b1_{\{\wt{\o} \in \cA_{j'}\}}
  \hP_{j'}  \big( (A \cap \cA_j )^{s,\wt{\o}}  \big)    \right]
 = \sum^\l_{j'=1} \hE_\hP   \n  \left[ \b1_{\{\wt{\o} \in A \cap \cA_j \}}  \b1_{\{\wt{\o} \in \cA_{j'} \}}
  \hP_{j'}  \big(\O^s\big)  \right]   = \hP(A \cap \cA_j) .
\eeas
 \fi

 To see   $\wh{\hP}$  satisfying \eqref{eq:eb127} for some $(Y,\wp)  \n \in \n  \fS$,
 we fix  $j \n = \n  1,  \cds \n , \l $,   $A  \n \in \n  \cF^t_s$,  and $\ga  \n \in \n  \cT^s  $.
 By   Lemma \ref{lem_element},  $( A \cap \cA_j )^{s , \wt{\o}}  = \O^s  $ (resp. $= \es $),
 when   $\wt{\o} \in A \cap \cA_j$ (resp. $\notin A \cap \cA_j$).
 Then we can deduce   that
   \beas
     \hE_{\wh{\hP} }  \Big[  \b1_{ A \cap \cA_j}  \wh{Y}^{t,\o}_{\ga  (\Pi^t_s)}  \Big]
   & \tn   =     & \tn   \sum^\l_{j' = 1}    \hE_{ \hP  }  \bigg[
  \b1_{ \{ \wt{\o} \in  \cA_{j'}  \} } \hE_{ \hP_{j'}}
  \bigg[ \Big( \b1_{   A \cap   \cA_j}  \wh{Y}^{t,\o}_{\ga  (\Pi^t_s)} \Big)^{s,\wt{\o}} \bigg]   \bigg]
 \n = \n \sum^\l_{j' = 1}    \hE_{ \hP  }  \bigg[
  \b1_{ \{ \wt{\o} \in  A \cap \cA_j  \} } \b1_{ \{ \wt{\o} \in  \cA_{j'}  \} } \hE_{ \hP_{j'}}
  \bigg[ \Big(   \wh{Y}^{t,\o}_{\ga  (\Pi^t_s)} \Big)^{s,\wt{\o}} \bigg]   \bigg] \\
   & \tn   =  & \tn  \hE_{ \hP  }  \bigg[  \b1_{ \{ \wt{\o} \in  A \cap  \cA_j  \} }
 \hE_{\hP_j} \Big[   \wh{Y}^{s, \o   \otimes_t    \wt{\o}}_\ga   \Big]  \bigg]   ,
     \eeas
 where we used the fact that  for any $\wh{\o}  \n \in \n  \O^s$,
 $ \Big(   \wh{Y}^{t,\o}_{\ga  (\Pi^t_s)} \Big)^{s,\wt{\o}} (\wh{\o})
      \n = \n   \Big(  \wh{Y}^{t,\o}_{\ga  (\Pi^t_s)}   \Big)  (\wt{\o} \n \otimes_s \n  \wh{\o}   )
      \n = \n  \wh{Y} \big( \ga \big( \Pi^t_s (\wt{\o}  \n \otimes_s \n  \wh{\o} )   \big) ,
     \o  \n \otimes_t \n  ( \wt{\o}  \n \otimes_s \n  \wh{\o}) \big)
      \n = \n  \wh{Y} \big( \ga (  \wh{\o} )  , (\o  \n \otimes_t \n  \wt{\o} )  \n \otimes_s \n  \wh{\o} \big)
      \n = \n   \wh{Y}^{s, \o    \otimes_t    \wt{\o} }_\ga   (\wh{\o}) $.  \qed

  \no {\bf Proof of Example \ref{eg_weak_formulation}:}
 Let $\rho \n \in \n  \fM$ satisfies   \eqref{eq:gb017} with some $ C \n > \n 0$ and $0 \n < \n  p_1  \n \le \n   p_2  $.
 \if{0}
 Given  $t  \n \in \n  [0,T]$ and $\hP  \n \in \n  \cP^\ell_t $,
  suppose that $B^t $ can be decomposed as the sum of two $\bF^t-$adapted processes $ A^\hP  $ and $ M^\hP  $ with all continuous paths   such that $ A^\hP  $ is
   a $\hP-$finite variation process
   with $\big|\frac{dA^\hP_t}{dt} \big|  \n \le \n  \ell$, \dtp ~
    and that $ M^\hP  $ is a $\hP-$martingale $ M^\hP  $ with
    $ trace   \big(\frac{d \lan M^\hP \ran_t}{dt} \big)
     \n \le \n  2 \ell$, \dtp

  Set $ \ol{\O}^t : = \O^t \times  \O^t \times  \O^t $.
  For any $ \ol{\o} = ( \o^{(1)},  \o^{(2)},  \o^{(3)}) \in \ol{\O}^t $, we
   define $\Pi^{(i)}_t (\ol{\o}) : = \o^{(i)} (t)$, $t \in [0,T]$, $i=1,2,3$. Then the Coordinator process on $\ol{\O}^t$
   is in form of
 \beas
 \ol{B}_t (\ol{\o}) = \big( \Pi^{(1)}_t (\ol{\o}) , \Pi^{(2)}_t (\ol{\o}) , \Pi^{(3)}_t (\ol{\o})  \big)
 = \big( \o^{(1)} (t) ,  \o^{(2)} (t) ,  \o^{(3)} (t)    \big), \q (t, \ol{\o}) \in [0,T] \times \ol{\O} .
 \eeas
 Define $\Phi : \O^t \to \ol{\O}^t $ by
   $\Phi  (\o) : = \big(\o, A^\hP (\o), M^\hP (\o) \big) $, $\fa \o \in \O^t$. Based on $\hP$,
 we can induce a probability $\ol{\hP}$ on the product measurable space  $\big(\ol{\O}^t, \ol{\cF}^t_T  \big) :=
\big(  \ol{\O}^t ,  \cF^t_T \otimes  \cF^t_T \otimes  \cF^t_T
    \big)$ via $\Phi$:
    \beas
    \ol{\hP} (A) : = \hP \big(\Phi^{-1}(A)\big), \q \fa A \in \ol{\cF}^t_T .
    \eeas

  \no (i)  $\Pi^{(1)}  = \Pi^{(2)}  + \Pi^{(3)}$, $\ol{\hP}-$a.s.
  \beas
\q && \hspace{-1.5cm} \ol{\hP} \big\{ \ol{\o} \in \ol{\O}^t: \Pi^{(1)}_s (\ol{\o})
 = \Pi^{(2)}_s (\ol{\o}) + \Pi^{(3)}_s (\ol{\o}) , ~ \fa s \in [t,T] \big\}
 = \hP \big\{ \o \in \O^t : \Pi^{(1)}_s (\Phi(\o) )
 = \Pi^{(2)}_s (\Phi(\o)) + \Pi^{(3)}_s (\Phi(\o)) , ~ \fa s \in [t,T] \big\} \\
&&   = \hP \big\{ \o \in \O^t : B^t_s (\o)
 = A^\hP_s (\o) +  M^\hP_s (\o)  , ~ \fa s \in [t,T] \big\} = 1 .
 \eeas

  \no (ii) $\Pi^{(2)}$ is a finite variation process under $\ol{\hP}$.
  \beas
\q && \hspace{-1.5cm} \ol{\hP} \big\{ \ol{\o} \in \ol{\O}^t: \Pi^{(2)}_\cd  (\ol{\o}) \hb{ is a path of finite variation}
  \big\}
 = \hP \big\{ \o \in \O^t : \Pi^{(2)}_\cd  (\Phi(\o)) \hb{ is a path of finite variation} \big\} \\
&&   = \hP \big\{ \o \in \O^t : A^\hP_\cd (\o)  \hb{ is a path of finite variation} \big\} = 1 .
 \eeas

  \no (iii) $\Pi^{(3)}$ is a martingale with respect to   $\big(\ol{\bF}^t , \ol{\hP}\big)$,
 where $ \ol{\bF}^t = \big\{\ol{\cF}^t_s : = \cF^t_s \otimes  \cF^t_s \otimes  \cF^t_s   \big\}_{s \in [t,T]} $:
 Let $t \le s < r \le T$. For any $D_i \in \cF^t_s$, $i=1,2,3$,
 \beas
 \int_{\{ \ol{\o} \in D_1 \times D_2 \times  D_3 \}} \Pi^{(3)}_r (\ol{\o}) \ol{\hP} (d \ol{\o})
 & \tn \dn =& \tn \dn   \int_{\{ \o \in \O^t: \Phi (\o)   \in D_1 \times D_2 \times  D_3 \}}
 \Pi^{(3)}_r (\Phi (\o))  \hP  (d \o)
 =  \int_{(B^t)^{-1} (D_1) \cap (A^\hP)^{-1} (D_2)  \cap  (M^\hP)^{-1} (D_3)  } M^\hP_r (\o)  \hP  (d \o) \\
 & \tn \dn =& \tn \dn  \int_{(B^t)^{-1} (D_1)  \cap  (A^\hP)^{-1} (D_2)  \cap  (M^\hP)^{-1} (D_3)  } M^\hP_s (\o)  \hP  (d \o)
 = \int_{\{ \ol{\o} \in D_1 \times D_2 \times  D_3 \}} \Pi^{(3)}_s (\ol{\o}) \ol{\hP} (d \ol{\o}) ,
 \eeas
 where we used in the third equality the fact that   $(B^t)^{-1} (D_1), (A^\hP)^{-1} (D_2), (M^\hP)^{-1} (D_3) $ are all in $ \cF^t_s  $. So $\Th : \{D_1 \times D_2 \times  D_3: D_i \in \cF^t_s, i=1,2,3 \} \subset
 \L:= \{ \ol{D} \in  \ol{\cF}^t_T:
  \int_{\ol{D}} \Pi^{(3)}_r  d \ol{\hP} = \int_{\ol{D}} \Pi^{(3)}_s  d \ol{\hP}   \}$. Clearly,
  $\Th$ is closed under intersection and $\L$ is a Dynkin system. Dynkin System Theorem then shows that
  \beas
 \int_{\ol{D}} \Pi^{(3)}_r  d \ol{\hP} = \int_{\ol{D}} \Pi^{(3)}_s  d \ol{\hP}, \q \fa \ol{D} \in
 \ol{\cF}^t_s : = \cF^t_s \otimes  \cF^t_s \otimes  \cF^t_s ,
  \eeas
 so $\hE_{\ol{\hP}} \big[ \Pi^{(3)}_r  \big| \ol{\cF}^t_s  \big] = \Pi^{(3)}_s $, $\ol{\hP}-$a.s.

  \no (iv)
  \beas
\q && \hspace{-1.5cm}
 (ds \otimes d \ol{\hP}) \Big\{ (s,\ol{\o}) \in [t,T] \times \ol{\O}^t:  \Big|\frac{d \Pi^{(2)}_s}{ds} (\ol{\o}) \Big|  \n > \n  \ell \Big\}
 = (ds \otimes d  \hP ) \Big\{ (s, \o ) \in [t,T] \times \O^t : \Big|\frac{d \Pi^{(2)}_s}{ds} (\Phi (\o)) \Big|  \n > \n  \ell \Big\}  \\
 &&   = (ds \otimes d  \hP ) \Big\{ (s, \o ) \in [t,T] \times \O^t : \big|\frac{dA^\hP_s}{ds} (\o)
    \big|  \n > \n  \ell \Big\} = 0 .
 \eeas

  \no (v)
    Given $i,j \in \{1,\cds,d\}$, one can deduce that
 \beas
   0& \dn  \tn =& \dn  \tn \ol{\hP} \n - \n \lmt{n \to \infty} \,
  \underset{ s \in [t,T] }{\sup} \bigg| \lan  \Pi^{3,i} , \Pi^{3,j} \ran_s
   \n - \n  \sum^{\lfloor 2^n \n  s\rfloor}_{k=1 }  
   \Big(\Pi^{3,i}_{ \frac{k}{2^n}  }  \n - \n  \Pi^{3,i}_{ \frac{k-1}{2^n}} \Big)
    \Big(\Pi^{3,j}_{ \frac{k}{2^n}  } \n  - \n  \Pi^{3,j}_{ \frac{k-1}{2^n}}   \Big)  \bigg|    \nonumber   \\
 & \dn  \tn = & \dn  \tn \hP  - \lmt{n \to \infty} \,
   \underset{ s \in [t,T] }{\sup} \bigg| \lan \Pi^{3,i}, \Pi^{3,j} \ran_s  ( \Phi   )
    \n - \n   \sum^{\lfloor 2^n \n  s\rfloor}_{k=1 }  
   \Big(M^{\hP,i}_{ \frac{k}{2^n}  }  \n -  \n  M^{\hP,i}_{ \frac{k-1}{2^n}} \Big)\Big(M^{\hP,j}_{ \frac{k}{2^n}  }
   \n -   \n   M^{\hP,j}_{ \frac{k-1}{2^n}}   \Big)  \bigg|  \,  .   
 \eeas
 Thus $
   \hP \big\{ \o \in \O^t :  \lan \Pi^{(3)}  \ran_s  ( \Phi (\o)  )
 = \lan M^\hP    \ran_s (\o)  , ~ \fa s \in [t,T]\} =1$ and it follows that
 \beas
&& \hspace{-1.5cm} (ds \otimes d \ol{\hP}) \Big\{ (s,\ol{\o}) \in [t,T] \times \ol{\O}^t:
  trace   \Big(\frac{d \lan \Pi^{(3)} \ran_s}{ds} (\ol{\o}) \Big)   \n > \n  2 \ell \Big\}
  =  (ds \otimes d \hP) \Big\{ (s,\o ) \in [t,T] \times \O^t:
  trace   \Big(\frac{d \lan \Pi^{(3)} \ran_s}{ds} (\Phi (\o)) \Big)   \n > \n  2 \ell \Big\} \\
 && =  (ds \otimes d \hP) \Big\{ (s,\o ) \in [t,T] \times \O^t:
  trace   \Big(\frac{d \lan M^\hP \ran_s}{ds}   (\o)  \Big)   \n > \n  2 \ell \Big\} = 0 .
 \eeas

 For any $A \in \cF^t_T$, one has
\beas
\ol{\hP} \circ \big( \Pi^{(1)} \big)^{-1} (A) =  \ol{\hP} \big\{ \ol{\o} \in \ol{\O}^t: \Pi^{(1)}  (\ol{\o})
 \in A \big\} =   \hP  \big\{  \o  \in  \O^t: \Pi^{(1)}  \big(\Phi( \o )\big)
 \in A \big\}  =   \hP  \big\{  \o  \in  \O^t:   \o   \in A \big\} = \hP (A) .
\eeas

 \fi
  Fix  $t  \n \in \n  [0,T)$ and $\d  \n \in \n  (0,\infty)$. We  consider   an enlarged canonical space
  $\ol{\O}^t  \n : = \n  \O^t  \n \times \n   \O^t  \n \times \n   \O^t   $ with canonical processes
 \beas
 \ol{B}_t (\ol{\o}) = \big( X_t (\ol{\o}) , A_t (\ol{\o}) , M_t (\ol{\o})  \big)
 = \big( x (t) ,  a (t) ,  m (t)    \big), \q \fa \ol{\o} = (x,a,m) \in \ol{\O}^t , ~
 \fa t \in [0,T]  .
 \eeas
 Given  $\hP  \n \in \n  \cP^\ell_t$, there exists
  an extension $\ol{\hP}$ of $\hP$ on $\ol{\O}^t$
 such that

  \no (i) 
  $ \ol{\hP} \big\{\ol{\o} \ins \ol{\O}^t : X (\ol{\o}) \ins A \big\} \= \hP (A)$ for any $A \ins \cF^t_T$;

  \no (ii)   $X  \n = \n  K  \n + \n  M$, $\ol{\hP}-$a.s., in which $K$ is an absolutely continuous process   with
 $\big|\frac{d K_t}{dt} \big|  \n \le \n  \ell$, $\ol{\hP}-$a.s., and $M$ is a $\ol{\hP} -$martingale    with
 $    \hb{trace}  \big(\frac{d \lan M \ran_t}{dt} \big)   \n \le \n  2 \ell $, $\ol{\hP}-$a.s.

  Let $\z \n \in \n  \cT^t$ and set  $ \eta  \n  : =  \n  \underset{ r \in [\z (X), \, (\z (X)+\d) \land T]    }{\sup}
  \big|  M_r  \n - \n  M_{\z (X)}  \big|  \n = \n  \underset{ r \in [t,   T]    }{\sup}
  \big|  M_{(\z (X)+\d) \land r}  \n - \n  M_{\z (X) \land r}  \big| $.
  Given   $p > 0$,   since
 \bea  \label{eq:ex001}
 \big( 1 \land n^{p-1} \big) \sum_{i=1}^n a_i^p \le
\bigg(\sum_{i=1}^n a_i \bigg)^p \le \big( 1 \vee n^{p-1} \big) \sum_{i=1}^n a_i^p  ,
\q \fa n \in \hN, ~ \fa \{a_i\}^n_{i=1} \subset [0,\infty) ,
\eea
  one can deduce from   the Burkholder-Davis-Gundy inequality  that
  \beas
  \hspace{-0.5cm}
   \hE_{\ol{\hP}} \big[  \eta^p \big] & \tn \dn  \le  & \tn \dn    \big( 1  \n \vee \n  d^{\frac{p}{2}-1} \big)  \sum_{i=1}^d
  \hE_{\ol{\hP}} \bigg[  \underset{ r \in [t,   T]    }{\sup}
  \big|  M^i_{(\z (X)+\d) \land r}  \n - \n  M^i_{\z (X) \land r}  \big|^p \bigg]
       \ls        c_p \big( 1  \n \vee \n  d^{\frac{p}{2}-1} \big)  \sum_{i=1}^d
  \hE_{\ol{\hP}} \bigg[ \Big( \int_t^T \n \b1_{\{ \z (X) \le r \le \z (X)+\d \}} d \lan M^i,M^i \ran_r \Big)^{\frac{p}{2}}  \bigg] \\
  & \tn  \dn  \le  & \tn  \dn  c_p \frac{1 \ve  d^{\frac{p}{2}-1}}{1 \ld d^{\frac{p}{2}-1}} \hE_{\ol{\hP}} \bigg[ \Big( \int_t^T \n \b1_{\{ \z (X) \le r \le \z (X)+\d \}} \hb{trace}  \Big(\frac{d \lan M \ran_r}{dr} \Big) dr   \Big)^{\frac{p}{2}}  \bigg]
   \le  c_p \frac{1 \ve  d^{\frac{p}{2}-1}}{1 \ld d^{\frac{p}{2}-1}} (2\ell \d )^{\frac{p}{2}}   ,
  \eeas
  where $ c_p $ is a constant depending on $p$.   Then we see from (i), (ii) and \eqref{eq:ex001} that
 \beas
  &&  \hspace{-1cm}  \hE_\hP \bigg[  \rho  \Big(    \d +      \underset{ r \in [\z, (\z+\d) \land T]    }{\sup}  \big|  B^{t}_r - B^{t}_\z  \big|  \Big) \bigg]
  = \hE_{\ol{\hP}} \bigg[  \rho  \Big(    \d +     \underset{ r \in [\z (X), (\z (X)+\d) \land T]    }{\sup}
  \big|  X_r - X_{\z (X)}  \big|  \Big) \bigg]
  \le C \sum^2_{i=1} \hE_{\ol{\hP}} \Big[    \big(     (1 \n +  \n    \ell ) \d +  \eta  \big)^{p_i} \Big] \\
  && \le   C \sum^2_{i=1} (1 \ve  2^{p_i-1})
  \big( (1 \n +  \n    \ell )^{p_i} \d^{p_i} \+  \hE_{\ol{\hP}}  [   \eta^{p_i}  ] \big)
  \ls   C (1 \ve  2^{p_2-1}) \sum^2_{i=1}
  \Big( (1 \n +  \n    \ell )^{p_i} \d^{p_i} \+  \d^{p_i/2}
  \+ \d^{-1/2} \hE_{\ol{\hP}}  \big[  \b1_{\{\eta \ge \sqrt{\d}\}}   \eta^{1+p_i }  \big] \Big) \\
   &&   \ls \frac14 \wh{C}     \sum^2_{i=1}
   (       \d^{p_i} \+   \d^{p_i/2}  ) \ls    \wh{C} \big( \d^{p_1/2} \ve  \d^{p_2} \big)
  \eeas
  for some   constant $\wh{C}$  depending on $C$, $d$,  $\ell$, $p_1$ , $p_2$  and $c_{p_2}$.
  Hence, \eqref{eq:ex015} holds for
  $\wh{\rho} (\d) : = \wh{C} \big( \d^{p_1/2} \vee \d^{p_2} \big) $.    \qed


 \subsection{Proofs of   Results in   Section \ref{sec:DPP}}

  \no {\bf Proof of Proposition \ref{prop_Z_conti_in_o}:}
 Fix $(Y,\wp)  \n \in \n  \fS$ and $(t,\o) \n \in \n  [0,T]  \n \times \n  \O $.
 Let  $\o' \n \in \n  \O$.
  We   set $t_1  \n : = \n  \wp (\o)  \n \land \n  \wp (\o')  \n \land \n  t $,
 $t_2  \n : = \n  (\wp (\o)  \n \vee \n  \wp (\o'))  \n \land \n  t $.
 Given $(\hP, \ga) \n \in \n \cP_t \n \times \n  \cT^t$, we see from
 Lemma \ref{lem_Y_diff} that
  \beas
  \hE_\hP \big[ \wh{Y}^{t,\o}_\ga    \big] \n - \n Z_t(\o')  & \tn \dn \le  & \tn  \dn
   \hE_\hP \big[ \wh{Y}^{t,\o}_\ga   \n - \n   \wh{Y}^{t,\o'}_\ga \big]
   \n  \le  \n    \wrY \Big( (1 \n + \n \k_\wp) \|\o \n - \n \o'\|_{0,t}
  \n + \n  \underset{r \in [ t_1 , t_2 ]}{\sup} \big|    \o    (r  )   \n - \n    \o    ( t_1 ) \big|  \Big) , \\
 \hb{and} \q   \hE_\hP \big[ \wh{Y}^{t,\o'}_\ga    \big] \n - \n Z_t(\o)   & \tn  \dn  \le  & \tn  \dn
   \hE_\hP \big[ \wh{Y}^{t,\o'}_\ga   \n - \n   \wh{Y}^{t,\o}_\ga \big]
   \n  \le  \n    \wrY \Big( (1 \n + \n \k_\wp) \|\o \n - \n \o'\|_{0,t}
  \n + \n  \underset{r \in [ t_1 , t_2 ]}{\sup} \big|    \o    (r  )   \n - \n    \o    ( t_1 ) \big|  \Big) .
 \eeas
 Taking supremum over $(\hP, \ga) \n \in \n \cP_t \times \cT^t$ on the left-hand-sides
 of both inequalities  leads to \eqref{eq:gd011}.

  For any $\e > 0$, there exists a $\l > 0$ such that $\wrY (x) < \e  $, $\fa x \in [0,\l)$.
 One can also find a $\wt{\l} (t,\o) > 0  $ such that $\phi^\o_t (y) < \l/2$, $\fa y \in \big[0, \wt{\l} (t,\o)\big) $.
 Now, taking $\d (t,\o) : = \frac{\l}{2(1+\k_\wp)} \land \frac{\wt{\l} (t,\o)}{\k_\wp}$,
 we will obtain \eqref{eq:gd014}. \qed

 \no {\bf Proof of Proposition \ref{prop_DPP_Z}:} Fix $(Y,\wp)  \n \in \n  \fS$.

  \no {\bf 1)} {\it We first show \eqref{eq:ef011}   for stopping time $\nu$ taking finitely many values. }

  Fix $ (t,  \o)  \n \in \n [0,T] \n \times \n \O $ and let
     $\nu \ins \cT^t$  take values in some finite  subset  $\{  t_1 < \cds < t_m  \}$ of $[t,T]$.
 We simply   denote
 \bea  \label{eq:wtY_wtZ}
 \cY_r := \wh{Y}^{t,\o}_r \q \hb{and} \q \cZ_r := Z^{t,\o}_r  , \q \fa r \in [t,T].
 \eea
     Proposition \ref{prop_shift0} (3) and \eqref{def_whY} show that
 $\cY$ is an $\bF^t-$adapted bounded process with all continuous paths.

    \no {\bf 1a)} {\it In the first step, we show
     \bea \label{eq:ef114}
 Z_t (\o) 
  \le  \underset{(\hP,\ga) \in \cP_t \times \cT^t}{\sup}   \,
  \hE_\hP \big[ \b1_{\{\ga  < \nu \}} \cY_{\ga   }
  \n + \n  \b1_{\{\ga  \ge \nu\}} \cZ_\nu  \big]
 \eea
 for the $\bF^t-$stopping time $\nu$ taking finitely many values.}

     Let $(\hP, \ga) \n \in \n \cP_t    \n \times \n \cT^t  $ and let $i = 1, \cds , m$.
 In light of  \eqref{eq:f475}, there exists  a $\hP-$null set  $ \cN_i  $  such that
      \bea \label{eq:ef111}
     \hE_\hP \big[\cY_{\ga  \vee t_i} \big|\cF^t_{t_i}\big]  (\wt{\o})
    =  \hE_{\hP^{t_i,\wt{\o}}}  \Big[  ( \cY_{\ga  \vee t_i}  )^{t_i,\wt{\o}}  \Big]
    =  \hE_{\hP^{t_i,\wt{\o}}}  \Big[   \wh{Y}^{t_i,\o \otimes_t \wt{\o}}_{(\ga   \vee   t_i)^{t_i,\wt{\o}}}  \Big] , \q
    \fa \wt{\o} \in  \cN^c_i ,
  \eea
  where we used the fact that for any $\wt{\o} \in \O^t$ and $\wh{\o} \in \O^{t_i}$
  \beas
  (\cY_{\ga  \vee t_i}  )^{t_i,\wt{\o}} (\wh{\o})
  \n =  \n   \cY_{\ga  \vee t_i} ( \wt{\o} \otimes_{t_i} \n  \wh{\o})
   \n = \n  \wh{Y} \big( (\ga \n \vee \n t_i)   (   \wt{\o} \otimes_{t_i} \n  \wh{\o})      ,
    \o  \n \otimes_t \n  ( \wt{\o} \otimes_{t_i} \n  \wh{\o} ) \big)
    \n  =  \n  \wh{Y} \big( (\ga \n \vee \n t_i)^{t_i, \wt{\o}} (   \wh{\o}    )  ,
   ( \o  \n \otimes_t  \n  \wt{\o} )  \n \otimes_{t_i} \n  \wh{\o}  \big)
   \n = \n     \wh{Y}^{t_i, \o \otimes_t  \wt{\o} }_{(\ga   \vee   t_i)^{t_i, \wt{\o}}} ( \wh{\o} ) .
  \eeas

    By (P3),   there exist
  an extension $(\O^t,\cF^{(i)},\hP^{(i)})$ of $(\O^t,\cF^t_T,\hP)$
  and $\O^{(i)} \n \in \n  \cF^{(i)}$ with $\hP^{(i)}(\O^{(i)})  \n = \n  1$
  such that for any $\wt{\o}  \n \in \n  \O^{(i)}$, $\hP^{t_i,   \wt{\o}}  \n \in \n  \cP_{t_i}  $.
  Given $\wt{\o} \n \in  \n \O^{(i)} \n \cap \n   \cN^c_i $,
  since $(\ga   \vee   t_i)^{t_i, \wt{\o}}  \n \in \n  \cT^{t_i}$
  by Proposition \ref{prop_shift0} (2),   we see from \eqref{eq:ef111} that
     \beas
     \hE_\hP \big[\cY_{\ga  \vee t_i} \big|\cF^t_{t_i}\big]  (\wt{\o})
    =  \hE_{\hP^{t_i,\wt{\o}}}  \Big[   \wh{Y}^{t_i,\o \otimes_t \wt{\o}}_{(\ga   \vee   t_i)^{t_i, \wt{\o}}}  \Big]
    \le Z (t_i,\o \otimes_t \wt{\o}) 
     = \cZ_{t_i} (\wt{\o}) .
  \eeas
  So $\O^{(i)} \n \cap \n   \cN^c_i \subset A_i  \n : = \n
  \{  \hE_\hP \big[\cY_{\ga  \vee t_i} \big|\cF^t_{t_i}\big]  \n   \le \n  \cZ_{t_i}   \}  $.
   The $\bF-$adaptedness of $Z$ by Proposition \ref{prop_Z_conti_in_o} as well as
  Proposition \ref{prop_shift0} (3) imply that $\cZ_{t_i}$ is $\cF^t_{t_i}-$measurable and thus
  $A_i \n \in \n \cF^t_{t_i}$.
  It follows that $\hP (A_i)  \n = \n  \hP^{(i)} (A_i)  \n \ge \n  \hP^{(i)} \big( \O^{(i)}
  \n \cap \n   \cN^c_i \big)  \n = \n  1 $. Namely,
 \bea  \label{eq:ef121}
  \hE_\hP \big[\cY_{\ga  \vee t_i} \big|\cF^t_{t_i}\big] \le  \cZ_{t_i}   , \q  \pas
  \eea

   Setting $A_i : = \{\nu = t_i\} \in \cF^t_{t_i}$, as
      \bea \label{eq:ef131}
   \b1_{\{\ga < t_i\}} \cY_{\ga} =  \b1_{\{\ga < t_i\}} \cY_{\ga \land t_i}   \in \cF^t_{\ga \land  t_i}
   \subset  \cF^t_{  t_i}  ,
   \eea
    we can deduce from \eqref{eq:ef121} that
  \beas
  \q  \hE_\hP[ \b1_{A_i}\cY_{\ga} ]
  & \tn =  & \tn   \hE_\hP \big[ \hE_\hP[ \b1_{A_i} \b1_{\{\ga < t_i\}} \cY_{\ga}
  \n + \n \b1_{A_i} \b1_{\{\ga \ge t_i\}} \cY_{\ga  \vee t_i} |\cF^t_{t_i} ]  \big]
 \n = \n    \hE_\hP \big[ \b1_{A_i} \b1_{\{\ga < t_i\}} \cY_{\ga}
  \n + \n  \b1_{A_i}\b1_{\{\ga \ge t_i\}} \hE_\hP[  \cY_{\ga  \vee t_i} |\cF^t_{t_i} ]  \big] \\
   & \tn  \le  & \tn   \hE_\hP \big[ \b1_{A_i}\b1_{\{\ga < t_i\}} \cY_{\ga }
  \n + \n  \b1_{A_i} \b1_{\{\ga \ge t_i\}} \cZ_{t_i}  \big]
  = \hE_\hP \big[ \b1_{A_i}\b1_{\{\ga < \nu \}} \cY_{\ga }
  \n + \n  \b1_{A_i} \b1_{\{\ga \ge \nu\}} \cZ_\nu  \big] ,
 \eeas
 and similarly that
   $    \hE_\hP[ \b1_{A_i}\cY_{\ga} ]
   \n =   \n     \hE_\hP \big[ \b1_{A_i} \b1_{\{\ga \le t_i\}} \cY_{\ga}
  \n + \n  \b1_{A_i}\b1_{\{\ga > t_i\}} \hE_\hP[  \cY_{\ga  \vee t_i} |\cF^t_{t_i} ]  \big]
    \n  \le   \n   \hE_\hP \big[ \b1_{A_i}\b1_{\{\ga \le \nu \}} \cY_{\ga  }
  \n + \n  \b1_{A_i} \b1_{\{\ga > \nu\}} \cZ_\nu  \big] $.
 Summing them  up over $i \in \{1, \cds, m\}$ yields that
 \bea \label{eq:ef211}
   \hE_\hP[ \cY_{\ga}  ]   \n \le \n  \hE_\hP \big[ \b1_{\{\ga  < \nu \}} \cY_{\ga   }
  \n + \n  \b1_{\{\ga  \ge \nu\}} \cZ_\nu  \big]
  \q \hb{and} \q
  \hE_\hP[ \cY_{\ga}  ]   \n \le \n  \hE_\hP \big[ \b1_{\{\ga  \le \nu \}} \cY_{\ga   }
  \n + \n  \b1_{\{\ga  > \nu\}} \cZ_\nu  \big]  .
 \eea
 Taking supremum of the former over $(\hP, \ga) \n \in \n \cP_t    \n \times \n \cT^t  $ leads to \eqref{eq:ef114}.

   \no {\bf 1b)} {\it To demonstrate the inverse inequality of \eqref{eq:ef114},
   we shall paste the local approximating $\hP-$maximizers of $Z^{t,\o}_{t_i}$'s
  according to \(P4\) and then make some estimations. }

  Fix  $(\hP, \ga) \n \in \n \cP_t  \n \times \n  \cT^t$,   $\e  \n > \n 0$
  and let $\d  \n \in \n \hQ_+$ satisfy $\ol{\rho}_Y (\d )  \n < \n  \e / 4$.
  For any $\wt{\o}  \n \in \n  \O^t $, let $\d(\wt{\o})  \n \in \n  \big((0,\d]  \n \cap \n  \hQ\big) \n \cup \n  \{\d\} $
  such that
   \bea \label{eq:ge011}
   \wrY \Big( (1 \n + \n \k_\wp) \d(\wt{\o})
      \n + \n  \phi^{\o \otimes_t \wt{\o} }_T  \big( \k_\wp \, \d(\wt{\o})   \big) \Big) < \e / 4 .
   \eea
   Since the canonical space $\O^t$ is separable and thus Lindel\"of,
   there exists    a sequence $\{\wt{\o}_j\}_{j \in \hN}$ of $\O^t$ such that
   $ \dis \underset{j \in \hN}{\cup} O_{ \d_j }(\wt{\o}_j) \n = \n  \O^t $
   with $ \d_j  :=  \d (\wt{\o}_j)   $.

   Let $i \n = \n  1, \cds  \n  , m$ and $j  \n \in \n  \hN$.
   By \eqref{eq:bb237},  $   \cA^i_j  \n := \n
  \{ \nu  \n = \n  t_i \}    \cap    \Big(  O^{t_i}_{\d_j }(\wt{\o}_j)
  \backslash \underset{j' < j}{\cup}  O^{t_i}_{\d_{j'} } \n \big(\wt{\o}_{j'} \big) \Big)
   \n \in \n  \cF^t_{t_i} $.  We can find a pair    $ (\hP^i_j, \ga^i_j       )  \n  \in \n
  \cP_{t_i}   \n \times \n \cT^{t_i} $ such that
      \bea   \label{eq:bb417}
   Z_{t_i} (\o \otimes_t \wt{\o}_j )
   \le    \hE_{  \hP^i_j  }  \Big[  \wh{Y}^{t_i,\o \otimes_t \wt{\o}_j}_{\ga^i_j}  \Big] + \e/4 .
   \eea
 Given $\wt{\o} \in
  O^{\raisebox{1.5pt}{\scriptsize $t_i$}}_{\d_j } (  \wt{\o}_j )   $,
 applying Lemma \ref{lem_Y_diff} with $(t,\o,\o',\hP,\ga) = \big(t_i,\o \otimes_t \wt{\o}_j,\o \otimes_t \wt{\o},
 \hP^i_j,\ga^i_j \big)$, we see from   \eqref{eq:ge011}   that
 \beas
 \hspace{-2mm}  \hE_{  \hP^i_j }  \Big[  \wh{Y}^{t_i, \o \otimes_t \wt{\o}_j}_{\ga^i_j}
   \n - \n   \wh{Y}^{t_i, \o \otimes_t \wt{\o}}_{\ga^i_j} \Big]
   & \tn \dn \le  & \tn  \dn
     \wrY   \Big( (1 \n + \n \k_\wp) \|   \o  \n \otimes_t \n   \wt{\o}_j
     \n -  \n  \o  \n \otimes_t  \n  \wt{\o} \|_{0,t_i}
      \n + \n  \phi^{\o \otimes_t \wt{\o}_j}_{t_i}
      \big( \k_\wp \|   \o  \n \otimes_t \n   \wt{\o}_j
       \n -  \n  \o  \n \otimes_t \n   \wt{\o} \|_{0,t_i}  \big) \Big) \\
     & \tn  \dn  =  & \tn  \dn    \wrY   \Big( (1 \n + \n \k_\wp) \|    \wt{\o}_j   \n -  \n   \wt{\o} \|_{t,t_i}
     \dn  + \n  \phi^{\o \otimes_t \wt{\o}_j}_{t_i}
      \big( \k_\wp \|   \wt{\o}_j   \n -  \n   \wt{\o} \|_{t,t_i}  \big) \Big)
       \n \le  \n       \wrY \Big( (1 \n + \n \k_\wp) \d_j
      \n + \n  \phi^{\o \otimes_t \wt{\o}_j}_T  \big( \k_\wp \, \d_j  \big) \Big)  \n < \n  \e / 4   .
         \eeas
 Then applying \eqref{eq:gd011} with $(t,\o,\o' ) = \big(t_i,\o \otimes_t \wt{\o}_j,\o \otimes_t \wt{\o}  \big)$,
   one can deduce from   \eqref{eq:bb417} and \eqref{eq:ge011} again that
        \bea
        \cZ_{t_i}(   \wt{\o}) & \tn  \dn   =& \tn  \dn     Z_{t_i}( \o  \n \otimes_t \n   \wt{\o})
        \le  Z_{t_i} \big(\o  \n \otimes_t \n  \wt{\o}_j \big)
    \n + \n    \wrY   \Big( (1 \n + \n \k_\wp) \|   \o  \n \otimes_t \n   \wt{\o}_j
     \n -  \n  \o  \n \otimes_t  \n  \wt{\o} \|_{0,t_i}
      \n + \n  \phi^{\o \otimes_t \wt{\o}_j}_{t_i}
      \big( \k_\wp \|   \o  \n \otimes_t \n   \wt{\o}_j
       \n -  \n  \o  \n \otimes_t \n   \wt{\o} \|_{0,t_i}  \big) \Big) \nonumber \\
     & \tn  \dn     =  & \tn  \dn     Z_{t_i}  (\o  \n \otimes_t \n  \wt{\o}_j  )  \n + \n
        \wrY   \Big( (1 \n + \n \k_\wp) \|    \wt{\o}_j   \n -  \n   \wt{\o} \|_{t,t_i}
      \dn + \n  \phi^{\o \otimes_t \wt{\o}_j}_{t_i}
      \big( \k_\wp \|   \wt{\o}_j   \n -  \n   \wt{\o} \|_{t,t_i}  \big) \Big)
       \n \le  \n    Z_{t_i}  (\o  \n \otimes_t \n  \wt{\o}_j  )  \n + \n     \wrY \Big( (1 \n + \n \k_\wp) \d_j
      \n + \n  \phi^{\o \otimes_t \wt{\o}_j}_T  \big( \k_\wp \, \d_j  \big) \Big)    \nonumber \\
       & \tn  \dn  <   & \tn  \dn    \hE_{  \hP^i_j  }  \Big[  \wh{Y}^{t_i,\o \otimes_t \wt{\o}_j}_{\ga^i_j} \Big]
      \n + \n  \e/2 \n <  \n     \hE_{  \hP^i_j }
   \Big[  \wh{Y}^{t_i, \o \otimes_t \wt{\o}}_{\ga^i_j} \Big]    \n + \n  \frac34 \e  .
        \label{eq:aa103}
    \eea

      Now, fix     $\l \in \hN  $. Setting $\hP^\l_{m+1} := \hP $,
      we recursively pick up  $\hP^\l_i  $, $  i = m , \cds \n , 1 $
    from $  \cP_t$ such that
    (P4) holds for $\Big(s, \wh{\hP} , \hP , \big\{(\cA_j, \d_j, \wt{\o}_j, \hP_j) \big\}^\l_{j=1}   \Big)
    = \Big(t_i, \hP^\l_i, \hP^\l_{i+1},   \big\{(\cA^i_j, \d_j, \wt{\o}_j, \hP^i_j) \big\}^\l_{j=1}   \Big)$
    and  $ \cA_0 \= \cA^i_0 \df \Big( \underset{j=1}{\overset{\l}{\cup}} \cA^i_j \Big)^c    \in \cF^t_{t_i}   $. Then
    \bea
    \hspace{-3mm}
     \hE_{ \hP^\l_i} [\xi]  \n = \n  \hE_{\hP^\l_{i+1}} [\xi]  ,  \;
          \fa \xi \n \in \n  L^1  \n ( \cF^t_{t_i}, \hP^\l_i )
           \n \cap \n  L^1  \n  \big( \cF^t_{t_i}, \hP^\l_{i+1} \big)
          \hb{ and }       \hE_{ \hP^\l_i } [\b1_{\cA^i_0}\xi] \n = \n  \hE_{\hP^\l_{i+1}} [\b1_{\cA^i_0}\xi]  ,
    \; \fa \xi  \n \in \n  L^1  \n  ( \cF^t_T, \hP^\l_i )
     \n \cap \n  L^1  \n  \big( \cF^t_T, \hP^\l_{i+1} \big)  . \q  \label{eq:ff024b}
         \eea

    For any $i \n = \n 1,\cds \n ,m$, as Lemma A.1 of \cite{ROSVU} shows that
       $\ga^i_j  (\Pi^t_{t_i})  \n \in \n  \cT^t_{t_i}$, stitching $\ga$ with $  \ga^i_j  (\Pi^t_{t_i})   $'s
       forms a new $\bF^t-$stopping time
    \bhe
    \bea \label{eq:a107}
     \wh{\ga}_\l  \n : = \n  \b1_{\{ \ga  < \nu \}}  \ga
    \n + \n  \b1_{\{  \ga    \ge   \nu \}}
         \Big(  \b1_{   \underset{i=1}{\overset{m}{\cap}}  \cA^i_0   } \ga
    + \sum^{m}_{i=1} \sum^\l_{j=1}   \b1_{   \cA^i_j   }  \ga^i_j  (\Pi^t_{t_i}) \Big) .
    \eea
    \ehe
   We see from    \eqref{eq:ff024b} that
  \bea   \label{eq:ff137}
  \hE_{\,\hP^\l_1} \Big[   \b1_{ \underset{i=1}{\overset{m}{\cap}}  \cA^i_0  } \cY_{\wh{\ga}_\l}    \Big]
  = \hE_{\,\hP^\l_2} \Big[ \b1_{ \underset{i=1}{\overset{m}{\cap}}  \cA^i_0   }  \cY_{\wh{\ga}_\l}    \Big]
  = \cds
  = \hE_{\,\hP^\l_m} \Big[ \b1_{ \underset{i=1}{\overset{m}{\cap}}  \cA^i_0   }  \cY_{\wh{\ga}_\l}    \Big]
  = \hE_{\,\hP^\l_{m+1}} \Big[ \b1_{ \underset{i=1}{\overset{m}{\cap}}  \cA^i_0   }  \cY_{\wh{\ga}_\l}    \Big]
  = \hE_\hP  \Big[    \b1_{ \underset{i=1}{\overset{m}{\cap}}  \cA^i_0   }  \cY_\ga       \Big] .
  \eea
  On the other hand,  for any  $ (i,j) \n  \in \n  \{1,\cds \n , m\}  \n \times \n  \{1,\cds \n , \l \} $,
 as $  \cA^i_j \n \subset \n   \cA^{i'}_0$ for
  $i'  \n \in \n  \{1,\cds \n , m \}\backslash \{i\}  $, we  can  deduce from
  \eqref{eq:ef131}, \eqref{eq:ff024b}, \eqref{eq:eb127} and \eqref{eq:aa103}    that
  \beas
  && \hspace{-1cm} \hE_{\,\hP^\l_1 }  \big[  \b1_{   \cA^i_j}   \cY_{\wh{\ga}_\l}  \big]
    \n  =   \n   \hE_{\,\hP^\l_2 } \Big[ \b1_{   \cA^i_j}   \cY_{\wh{\ga}_\l}   \Big] =
        \cds =  \hE_{ \hP^\l_{i-1} }  \big[ \b1_{   \cA^i_j}   \cY_{\wh{\ga}_\l}  \big]
         =  \hE_{ \hP^\l_i }  \big[ \b1_{   \cA^i_j}   \cY_{\wh{\ga}_\l}  \big]
  \n  =   \n \hE_{ \hP^\l_i }  \Big[   \b1_{    \{ \ga < \nu \} \cap  \cA^i_j}   \cY_\ga    +
    \b1_{      \{ \ga \ge  \nu \} \cap   \cA^i_j}   \cY_{\ga^i_j  (\Pi^t_{t_i})}       \Big] \nonumber \\
     & &     =   \n \hE_{ \hP^\l_i }  \Big[   \b1_{    \{ \ga < t_i \} \cap  \cA^i_j}   \cY_{\ga }
 +  \b1_{      \{  \ga \ge  t_i \} \cap   \cA^i_j}   \cY_{\ga^i_j  (\Pi^t_{t_i})}       \Big]
    \n  \ge  \n  \hE_{ \hP^\l_{i+1} } \n  \bigg[  \b1_{    \{ \ga < t_i \} \cap  \cA^i_j}   \cY_{\ga }
 \n + \n \b1_{  \{\ga (\wt{\o}) \ge t_i\} \cap \{ \wt{\o} \in  \cA^i_j\}}
  \Big( \,  \hE_{\hP^i_j}   \Big[ \wh{Y}^{t_i, \o \otimes_t  \wt{\o} }_{\ga^i_j}   \Big]
        \n  -  \n  \ol{\rho}_Y (\d)  \Big) \bigg]  \nonumber \\
 &&  \ge  \n  \hE_{ \hP^\l_{i+1} }  \n  \big[   \b1_{    \{ \ga < t_i \} \cap  \cA^i_j}   \cY_\ga
    +     \b1_{ \{  \ga   \ge t_i \} \cap   \cA^i_j  }
   ( \cZ_{t_i}   -  \e   )  \big]
      = \cds = \hE_{ \hP^\l_{m+1} } \big[  \b1_{    \{ \ga < \nu \} \cap  \cA^i_j}   \cY_\ga
 + \b1_{\{ \ga     \ge  \nu\} \cap   \cA^i_j  }
          \big( \cZ_\nu    -  \e \big) \big]  \\
 &&   = \hE_\hP   \big[  \b1_{    \{ \ga < \nu \} \cap  \cA^i_j}   \cY_\ga
 + \b1_{\{ \ga     \ge  \nu\} \cap   \cA^i_j  }
          \big( \cZ_\nu    -  \e \big) \big]  .
   \eeas
   Taking summation   over $ (i,j) \in \{1,\cds \n, m\} \times \{1,\cds \n , \l \}$
   and then combining with \eqref{eq:ff137}   yield that
         \bea
        Z_t (\o)   & \tn \ge& \tn \hE_{\,\hP^\l_1 }  \big[   \cY_{\wh{\ga}_\l}  \big] \ge
         \hE_{\hP}  \big[ \b1_{\{\ga < \nu\}} \cY_\ga
      \n + \n     \b1_{\{\ga \ge \nu\}}  \cZ_\nu   \big]
        +   \hE_{\hP}  \Big[ \b1_{\{  \ga  \ge  \nu \}}
          \b1_{ \underset{i=1}{\overset{m}{\cap}}  \cA^i_0   } ( \cY_\ga - \cZ_\nu )    \Big] -\e     \nonumber \\
       & \tn  \ge  & \tn   \hE_\hP  \Big[ \b1_{\{\ga < \nu\}} \cY_{\ga}
      \n + \n     \b1_{\{\ga \ge \nu\}}  \cZ_\nu  \Big]
      \n - \n  2 M_Y  \hP \Big(\underset{i=1}{\overset{m}{\cap}}  \cA^i_0 \Big)   \n - \n  \e  ,
     \label{eq:ef135}
        \eea
     where $\underset{i=1}{\overset{m}{\cap}}  \cA^i_0 \n = \n  \Big( \underset{i=1}{\overset{m}{\cup}}
     \underset{j=1}{\overset{\l}{\cup}}  \cA^i_j \Big)^c$.
     Since $ \underset{j \in \hN}{\cup}       O^{t_i}_{\d_j }(\wt{\o}_j)
     \supset \underset{j \in \hN}{\cup}       O_{\d_j }(\wt{\o}_j)   = \O^t $ for each $i \n \in \n \{ 1, \cds , m \}$,
     we see that
   $\underset{i=1}{\overset{m}{\cup}} \underset{j \in \hN}{\cup} \cA^i_j  \n = \n
     \underset{i=1}{\overset{m}{\cup}} \bigg[ \{\nu = t_i\} \cap
     \Big( \underset{j \in \hN}{\cup}       O^{t_i}_{\d_j }(\wt{\o}_j) \Big) \bigg]
      \n = \n  \underset{i=1}{\overset{m}{\cup}} \{\nu = t_i\}  \n = \n  \O^t $,
         letting $\l  \n \to \n  \infty$  and then letting $\e  \n \to \n  0$ in \eqref{eq:ef135} yield that
  \bea \label{eq:ef233}
     Z_t (\o)   \n \ge \n   \hE_\hP  \Big[ \b1_{\{\ga < \nu\}} \cY_\ga
      \n + \n     \b1_{\{\ga \ge \nu\}}  \cZ_\nu  \Big]     .
  \eea
 Taking supremum over $ (\hP, \ga) \n \in \n \cP_t    \n \times \n \cT^t  $ and combining with \eqref{eq:ef114}
 prove \eqref{eq:ef011}   for stopping times $\nu$ taking finitely many values.

      \no {\bf 2)} {\it Next, let us   show \eqref{eq:es017} and thus the continuity of process $Z$.}

 Fix   $ \o   \n \in \n   \O$ and
     $0  \n \le \n  t  \n \le \n  s  \n \le \n  T$.
     If $t  \n = \n  s$, then \eqref{eq:es017} trivially holds. So we assume   $ t  \n < \n  s $.

   \no {\bf 2a)} {\it Let us start by proving   an auxiliary inequality:}
  \bea \label{eq:ge317}
 \hE_\hP \Big[   \big| Z_s^{t,\o}   \n   -  \n   Z_s(\o)    \big| \Big] \le
  2 C_{\n \vr} M_Y   \Big( (s \n - \n t)^{ \frac{q_1}{2}}  \n \vee  \n  (s \n - \n t)^{q_2-\frac{q_1}{2}} \Big)
   \n + \n    \wrY  \big(   \d_{t,s} (\o)  \big)
    \n \vee \n  \wh{\wh{\rho}\,}_Y \big(   \d_{t,s} (\o)  \big)   : = \wh{\phi}_{t,s} (\o) .
 \eea

  For any $\wt{\o} \in \O^t$,  applying \eqref{eq:gd011} with
 $(t,\o,\o') = (s, \o \otimes_t \wt{\o}, \o  )$  yields that
   \bea \label{eq:ge021}
       \big|  Z (s, \o \otimes_t \wt{\o} )    \n    -  \n  Z_s(\o) \big|
      \n  \le \n   \wrY  \Big( (1 \n + \n \k_\wp) \| \o \otimes_t \wt{\o}   \n - \n \o  \|_{0,s}
  \n + \n  \underset{r \in [ s_1 (\wt{\o}) , s_2 (\wt{\o}) ]}{\sup} \big|  ( \o  \n \otimes_t \n  \wt{\o} )  (r  )
    \n - \n  ( \o  \n \otimes_t \n  \wt{\o} )   ( s_1 (\wt{\o}) ) \big|  \Big) , \qq
 \eea
 where $s_1 (\wt{\o})  \n : = \n  \wp (\o  \n \otimes_t \n  \wt{\o})  \n \land \n  \wp (\o)  \n \land \n  s   $
 and $s_2 (\wt{\o})  \n : = \n  \big(\wp (\o  \n \otimes_t \n  \wt{\o})  \n \vee \n  \wp (\o)\big)  \n \land \n  s    $.

   Let $ \hP \in \cP_t $  and set $ \dis  A \n  : = \n  \Big\{ \, \underset{r \in [t,s]}{\sup} |B^t_r|
  \n \le \n   (s-t)^{\frac{q_1}{2}} \Big\}$. As $B^t_t = 0 $,
  one can deduce from \eqref{eq:ef221} and \eqref{eq:ex015} that
 \bea
 \hE_\hP \big[ \b1_{A^c}  | Z^{t,\o}_s   \n - \n   Z_s(\o) |  \big] & \tn \dn \le & \tn  \dn  2 M_Y \hP (A^c)
  \n \le  \n  2 M_Y (s \n - \n t)^{-\frac{q_1}{2}}
   \hE_\hP \bigg[ \underset{r \in [t,s]}{\sup} |B^t_r \n - \n B^t_t| \bigg]
   \n \le  \n  2 M_Y (s \n - \n t)^{-\frac{q_1}{2}}
   \hE_\hP \Big[ \vr \Big( (s \n - \n t)
    \n + \n  \underset{r \in [t,s]}{\sup} |B^t_r \n - \n B^t_t| \Big) \Big] \nonumber \\
    & \tn  \dn  \le   & \tn  \dn   2 M_Y (s \n - \n t)^{-\frac{q_1}{2}} \wh{\vr} (s \n - \n t)
   \n \le  \n 2 C_{\n \vr} M_Y   \Big( (s \n - \n t)^{ \frac{q_1}{2}}  \n \vee  \n  (s \n - \n t)^{q_2-\frac{q_1}{2}} \Big) .
   \label{eq:ge311}
 \eea
 As to $ \hE_\hP \big[  \b1_A | Z^{t,\o}_s   \n - \n   Z_s(\o) |  \big] $,
 we shall estimate it by two cases on values of $\wp(\o)$:

  \no (i) When $\wp(\o) \n \le \n  t$, let   $\wt{\o}  \n \in \n A $.
 Applying Lemma \ref{lem_shift_stopping_time} with $(t,s,\tau) \n = \n (0,t,\wp)$ yields that
 $\wp (\o    \otimes_t    \wt{\o})  \n = \n  \wp (\o)$, thus $s_1 (\wt{\o}) \n = \n s_2 (\wt{\o})  \n = \n  \wp (\o)  \n \land \n  s
  \n = \n  \wp (\o) $. Since
  \bea \label{eq:ge023}
 \| \o \otimes_t \wt{\o}  \n - \n  \o  \|_{0,s}
  \n = \n    \underset{r \in [t, s]}{\sup} \big|  \wt{\o}(r) \n  + \n \o(t)  \n  - \n  \o(r)  \big|
  \n \le \n     \underset{r \in [t, s]}{\sup} \big| \wt{\o}(r)    \big|    \n + \n
     \underset{r \in [t, s]}{\sup} \big|  \o(r)  \n - \n  \o(t)     \big|
   \n \le \n   (s \n - \n t)^{\frac{q_1}{2}}
    \n + \n    \underset{r \in [t, s]}{\sup} \big|  \o(r)  \n - \n  \o(t)     \big|  ,
 \eea
 we can deduce from \eqref{eq:ge021}   that
   \bea \label{eq:ge314}
   \hE_\hP \Big[ \b1_A  \big| Z_s^{t,\o}   \n   -  \n   Z_s(\o)    \big| \Big]
   \le   \wrY  \big(   \d_{t,s} (\o)  \big) .
   \eea

  \no (ii) When $\wp(\o) \n > \n  t$,
 applying Lemma \ref{lem_shift_stopping_time} again shows that
 $\wp(\o  \n \otimes_t \n \O^t  )  \n \subset \n  (t,T] $
 and that 
 $\z  \n : = \n  \wp^{t,\o}    \n \land \n  \wp (\o)  \n \land \n  s  \n \in \n  \cT^t $.
 Let   $\wt{\o}  \n \in \n  A$.
 Since  $\z (\wt{\o})   \n  = \n  \wp (\o  \n \otimes_t \n  \wt{\o})  \n \land \n  \wp (\o)  \n \land \n  s
  \n = \n  s_1 (\wt{\o})  \n > \n  t $, we have
    \beas
   \underset{r \in [ s_1 (\wt{\o}) , s_2 (\wt{\o}) ]}{\sup} \big|  ( \o \otimes_t \wt{\o} )  (r  )
    \n - \n  ( \o \otimes_t \wt{\o} )  ( s_1 (\wt{\o}) ) \big|
    =    \underset{r \in [ s_1 (\wt{\o}) , s_2 (\wt{\o}) ]}{\sup} \big|  \wt{\o}    (r  )   \n - \n  \wt{\o}    ( s_1 (\wt{\o}) ) \big|
    =  \underset{r \in [ \z (\wt{\o})  , s_2 (\wt{\o})]}{\sup} \big| B^t_r ( \wt{\o}   )
     \n - \n B^t_{\z   } ( \wt{\o}  ) \big|  .
   \eeas
 By \eqref{eq:ec017} and \eqref{eq:ge023},
 $ s_2 (\wt{\o})  \n - \n \z (\wt{\o}) \n = \n s_2 (\wt{\o})  \n - \n s_1 (\wt{\o})
  \n \le  \n   \wp (\o    \otimes_t    \wt{\o})    \ve     \wp (\o)
   \n - \n  \wp (\o    \otimes_t    \wt{\o})    \ld    \wp (\o)
    \n = \n  |   \wp (\o    \otimes_t    \wt{\o})  \n - \n   \wp (\o)  |
  \n \le \n  \k_\wp \|   \o    \otimes_t    \wt{\o}  \n - \\   \o \|_{0,s}   \n < \n    \d_{t,s} (\o) $.
 So $  \underset{r \in [ s_1 (\wt{\o}) , s_2 (\wt{\o}) ]}{\sup} \big|  ( \o  \n \otimes_t \n  \wt{\o} )  (r  )
    \n - \n  ( \o  \n \otimes_t \n  \wt{\o} )  ( s_1 (\wt{\o}) ) \big|
    \n  \le  \n   \underset{r \in [ \z(\wt{\o}) , (\z(\wt{\o})+\d_{t,s} (\o)) \land T ]}{\sup} \big| B^t_r ( \wt{\o} )
    \n - \n B^t_\z ( \wt{\o}  )   \big| $.
 Then \eqref{eq:ge021}, \eqref{eq:ge023} and \eqref{eq:ex015} imply  that
   $   \hE_\hP \Big[  \b1_A  \big| Z_s^{t,\o}   \n   -  \n   Z_s(\o)    \big| \Big]
    \ls   \hE_\hP \bigg[  \b1_A  \wrY  \Big(  \d_{t,s} (\o)
   \+  \underset{r \in [ \z  , (\z +\d_{t,s} (\o)) \land T ]}{\sup} \big| B^t_r
    \n - \n B^t_\z     \big|   \Big) \bigg]
  \ls   \wh{\wh{\rho}\,}_Y \big(   \d_{t,s} (\o)  \big) $,
 which together with \eqref{eq:ge311} and \eqref{eq:ge314} leads to \eqref{eq:ge317}.

   \no {\bf 2b)}  {\it Now, we shall use \eqref{eq:ef233},
   \eqref{eq:ge317}, \eqref{eq:ef211} as well as  \eqref{eq:ex015} to derive \eqref{eq:es017}.}

    For any  $\hP \n \in \n \cP_t$, applying \eqref{eq:ef233} with $ \nu = s $
 and $\ga = s  $, we see from \eqref{eq:ge317}  that
 \bea
  Z_t(\o)  \n - \n   Z_s(\o)
  \n  \ge    \n
  \hE_\hP \big[    Z^{t,\o}_s   \n - \n   Z_s(\o)   \big] \ge - \wh{\phi}_{t,s} (\o) .   \label{eq:bb439}
 \eea
  As to the inverse inequality, let us fix $\e \n > \n  0$.  There exists a pair
  $(\hP, \ga) 
  \n \in \n \cP_t    \n \times \n \cT^t   $ such that
  $  Z_t(\o)  \n \le \n     \hE_\hP  \big[    \wh{Y}^{t,\o}_\ga        \big] \n + \n  \e  $. Applying
  the first inequality of \eqref{eq:ef211} with $\nu  \n = \n  s$ yields that
 \bea \label{eq:ge321}
   Z_t(\o)  \n \le \n     \hE_\hP  \big[    \wh{Y}^{t,\o}_\ga        \big] \n + \n  \e
 \le  \hE_\hP  \big[  \b1_{\{\ga < s\}}  \wh{Y}^{t,\o}_\ga  +  \b1_{\{\ga \ge s\}} Z^{t,\o}_s     \big]
 \n + \n  \e    .
 \eea

 For any $ \wt{\o} \ins \O^t $,
 $ \wh{Y}^{t,\o}_\ga  (\wt{\o}) \n - \n    \wh{Y}^{t,\o}_s  (\wt{\o})
  \n = \n     \wh{Y}  \big( \ga  (\wt{\o}),   \o  \n \otimes_t \n  \wt{\o} \big)
  \n - \n   \wh{Y}  \big( s,   \o  \n \otimes_t \n  \wt{\o} \big)
   \n = \n      Y  \big( \fs_1 (\wt{\o}) ,   \o  \n \otimes_t \n  \wt{\o} \big)
  \n - \n   Y  \big( \fs_2 (\wt{\o}) ,   \o  \n \otimes_t \n  \wt{\o} \big)   $,
 where $ \fs_1 (\wt{\o})  \n : = \n  \ga  (\wt{\o})  \n \land \n  \wp(\o  \n \otimes_t \n  \wt{\o}) $
 and $ \fs_2 (\wt{\o})  \n : = \n  s  \n \land \n  \wp(\o  \n \otimes_t \n  \wt{\o}) $.
 Let us  show by two cases that
 \bea \label{eq:a121}
   \hE_\hP \Big[ \b1_{\{\ga < s\}} \big|   \wh{Y}^{t,\o}_\ga  (\wt{\o})
 \n - \n    \wh{Y}^{t,\o}_s  (\wt{\o})  \big| \Big]   \n \le \n  \wrY (s \n - \n t)  .
 \eea
  If $\wp(\o) \n \le \n  t$, for any $\wt{\o} \ins \O^t$,
  since Lemma \ref{lem_shift_stopping_time} shows that
 $\wp (\o  \oti_t    \O^t)  \n = \n  \wp (\o)$,
 we see that $\fs_1 (\wt{\o}) \n = \n \fs_2 (\wt{\o})  
  \n = \n  \wp (\o) $ and thus that
    $\big|   \wh{Y}^{t,\o}_\ga  (\wt{\o}) \n - \n    \wh{Y}^{t,\o}_s  (\wt{\o})  \big|  \n = \n  0$.
 Otherwise,  if $\wp(\o) \n > \n  t$,
 applying Lemma \ref{lem_shift_stopping_time} again gives that
   $\wh{\ga}  \n : = \n \ga  \n \land \n  \wp^{t,\o}  \n \in \n  \cT^t $.
  For any $\wt{\o}  \n \in \n  \{\ga  \n < \n  s\} $,
  since $\wh{\ga} (\wt{\o})   \n  = \n  \ga (\wt{\o})  \n \land \n \wp (\o  \n \otimes_t \n  \wt{\o})
  \n = \n  \fs_1 (\wt{\o})  \n \ge \n  t $
  and since $ 
  \fs_2 (\wt{\o})  \n - \n  \fs_1 (\wt{\o}) 
  \n \le \n  s  \n - \n  t $,
     \eqref{eq:aa211} implies that
 \beas
  \big|   \wh{Y}^{t,\o}_\ga  (\wt{\o}) \n - \n    \wh{Y}^{t,\o}_s  (\wt{\o})  \big|
 & \tn  \dn  \le & \tn  \dn  \rY \Big( \big(\fs_2 (\wt{\o}) \n - \n \fs_1 (\wt{\o})\big)
  \n + \n  \underset{r \in [0,T]}{\sup} \Big| (\o \n \otimes_t \n \wt{\o}) \big(r  \n \land \n  \fs_1 (\wt{\o})\big)
   \n - \n  (\o \n \otimes_t \n \wt{\o}) \big(r  \n \land \n  \fs_2 (\wt{\o}) \big) \Big|   \Big)  \nonumber  \\
  & \tn  \dn  \le  & \tn  \dn  \rY \Big( (s  \n - \n t)
   \n + \n  \underset{r \in [ \wh{\ga} (\wt{\o})  , \fs_2 (\wt{\o})]}{\sup} \big|  \wt{\o}  (r  )
   \n - \n    \wt{\o}  \big(  \wh{\ga} (\wt{\o}) \big) \big| \Big)
       \n  \le  \n  \rY \Big( (s  \n - \n t)
       \n + \n  \underset{r \in [ \wh{\ga} (\wt{\o})  , (\wh{\ga} (\wt{\o})+s-t)   \land    T]}{\sup}
       \big| B^t_r ( \wt{\o}   )   \n - \n B^t_{ \wh{\ga}   } ( \wt{\o}  ) \big| \Big)    .
 \eeas
 Then \eqref{eq:a121} follows from \eqref{eq:ex015}.
 Plugging \eqref{eq:a121} into  \eqref{eq:ge321}, we can deduce from  \eqref{eq:ef221} and \eqref{eq:ge317} that
  \beas
  \hspace{-3mm}  Z_t(\o) \n - \n  Z_s (\o)
   \n \le \n  \hE_\hP  \big[  \b1_{\{\ga < s\}}  \wh{Y}^{t,\o}_s
    \n + \n   \b1_{\{\ga \ge s\}} Z^{t,\o}_s   \dn - \n  Z_s (\o)   \big]
 \n + \n \wrY (s \n - \n t)  \n + \n  \e
   \n \le \n  \hE_\hP  \big[    Z^{t,\o}_s   \n - \n  Z_s (\o)   \big]
 \n + \n \wrY (s \n - \n t) \n + \n  \e
  \n \le  \n  \wh{\phi}_{t,s} (\o)  \n + \n \wrY (s \n - \n t)  \n + \n  \e  .
 \eeas
 Letting $\e \to 0$ and combining with \eqref{eq:bb439} yield that
 $ 
  \big| Z_t (\o) \n - \n  Z_s (\o) \big|    \le  \wh{\phi}_{t,s} (\o)  \n + \n \wrY (s \n - \n t)      $,
  i.e. \eqref{eq:es017}.

   As $ \lmtd{t \nearrow s} \d_{t,s} (\o) = \lmtd{s \searrow t} \d_{t,s} (\o) = 0$, we see that
    $ \lmtd{t \nearrow s} \wh{\phi}_{t,s} (\o) = \lmtd{s \searrow t} \wh{\phi}_{t,s} (\o) = 0$,
    which together with \eqref{eq:ef221} and \eqref{eq:gd014} shows that
    $Z$ is   an $\bF-$adapted  process bounded by $M_Y  $ and with   all continuous paths.

       \no {\bf 3)} {\it Finally,  we   show \eqref{eq:ef011}   for general stopping time $\nu$.}

  Fix $ (t,  \o)  \n \in \n [0,T] \n \times \n \O $, $\nu \in \cT^t$ and $(\hP, \ga) \n \in \n \cP_t    \n \times \n \cT^t $.
  We still take the simple notation \eqref{eq:wtY_wtZ}.
  For any   $k  \n \in \n  \hN$, let us   set    $t^k_i    \n  :=  \n   t  \n \vee \n  (i 2^{-k} T)     $,
     $i  \n   =    \n    0 , \cds   \n   , 2^k  $ and
   define
   \bea \label{eq:ef231}
   \nu_k := \b1_{\{\nu=t\}} t + \sum^{2^k}_{i=1} \b1_{\{ t^k_{i-1}  <  \nu
   \le   t^k_i  \}} t^k_i       \in \cT^t   .
   \eea
 Applying the second inequality of \eqref{eq:ef211} with $\nu = \nu_k$ yields that
 $\hE_\hP[ \cY_\ga ]   \n \le \n  \hE_\hP \big[ \b1_{\{\ga \le \nu_k \}} \cY_{\ga  }
  \n + \n  \b1_{\{\ga > \nu_k \}} \cZ_{\nu_k}   \big] $.
  Since $\lmtd{k \to \infty} \nu_k = \nu $ and
  since
  \bea \label{eq:uxu011}
  \hb{ the function $x \to \b1_{\{x \ge a\}}$ is right-continuous for any $a \in \hR$, }
   \eea
  letting $k \to \infty$  we can deduce from the continuity of   $Z$ by part 2),
  the bounded convergence theorem   and \eqref{eq:ef221} that
 \bea \label{eq:ef235}
   \hE_\hP[ \cY_\ga ]   \n \le \n  \hE_\hP \big[ \b1_{\{\ga \le \nu  \}} \cY_{\ga  }
  \n + \n  \b1_{\{\ga > \nu  \}} \cZ_\nu   \big] \le \hE_\hP \big[ \b1_{\{\ga < \nu  \}} \cY_{\ga  }
  \n + \n  \b1_{\{\ga \ge \nu  \}} \cZ_\nu   \big] .
 \eea

 Next, let $n, k \ins \hN$ with $n \< k$. We define $\ga_n \df   \b1_{\{\ga=t\}} t \n + \n
  \sum^{2^n}_{i=1} \b1_{\{ t^n_{i-1}  <  \ga
   \le   t^n_i  \}} t^n_i      \ins  \cT^t  $ and still consider $\nu_k$ defined  in \eqref{eq:ef231}.
   Applying \eqref{eq:ef233} with $(\hP,\ga,\nu) = (\hP,\ga_n,\nu_k) $ gives that
  \bea \label{eq:cc727}
     Z_t (\o)   \n \ge \n   \hE_\hP  \big[ \b1_{\{\ga_n < \nu_k\}} \cY_{\ga_n}  \n + \n  \b1_{\{\ga_n \ge \nu_k\}}
     \cZ_{\nu_k}   \big]     .
  \eea
  Clearly, $ \{\ga_n < \nu \} \subset \{\ga_n < \nu_k\} $. To see the reverse inclusion, we let
  $\o \in \{\ga_n < \nu_k\}$. There exist $i \in \{0,\cds \n, 2^n\}$ and $j \in \{1,\cds \n, 2^k\}$ such that
  $ t^n_i = \ga_n (\o) < \nu_k (\o) = t^k_j $. Since $\{t^n_\ell\}^{2^n}_{\ell=0}
  \subset \{t^k_\ell\}^{2^k}_{\ell=0}  $, one has
  $  \ga_n (\o) = t^n_i \le t^k_{j-1} < \nu (\o) $. Thus $ \{\ga_n < \nu \} = \{\ga_n < \nu_k\} $ and
  \eqref{eq:cc727} becomes
  $  Z_t (\o)   \n \ge \n   \hE_\hP  \big[ \b1_{\{\ga_n < \nu \}} \cY_{\ga_n}  \n + \n  \b1_{\{\ga_n \ge \nu \}}
     \cZ_{\nu_k}   \big]     $.
  As $k \to \infty $,   the continuity of   $Z$ by part 2), \eqref{eq:ef221}
   and   the bounded convergence theorem imply   that
  \beas 
  \hE_\hP   \big[ \b1_{\{\ga_n < \nu  \}} \cY_{\ga_n}
  \n +  \n     \b1_{\{\ga_n \ge \nu \}} \cZ_\nu  \big]
  =  \lmt{k \to \infty}  \hE_\hP   \big[ \b1_{\{\ga_n < \nu  \}} \cY_{\ga_n}
  \n +  \n     \b1_{\{\ga_n \ge \nu \}} \cZ_{\nu_k}  \big] \le   Z_t (\o) .
  \eeas
  Since $\lmtd{n \to \infty} \ga_n \= \ga $,
  letting $n \to \infty $, we can deduce from \eqref{eq:uxu011},  the   continuity of $Y$,   \eqref{eq:ef221},
  the bounded convergence theorem as well as \eqref{eq:ef235}  that
  \beas
  \hE_\hP[ \cY_\ga ]   \n \le \n
  \hE_\hP   \big[ \b1_{\{\ga <  \nu  \}} \cY_\ga
  \n +  \n     \b1_{\{\ga \ge  \nu \}} \cZ_\nu   \big]
  \n = \n   \lmt{n \to \infty}   \hE_\hP   \big[ \b1_{\{\ga_n < \nu  \}} \cY_{\ga_n}
  \n +  \n     \b1_{\{\ga_n \ge \nu \}} \cZ_\nu  \big] \le  Z_t (\o)  .
  \eeas
  Taking supremum over $ (\hP, \ga) \n \in \n \cP_t    \n \times \n \cT^t  $ proves \eqref{eq:ef011}. \qed

 \no {\bf Proof of Proposition \ref{prop_Z_martingale}:} Fix $ (Y,\wp) \in \fS $ and $n \in \hN$.
Since  both $\wh{Y}$ and $Z$ are $\bF-$adapted  processes with  all continuous paths by Proposition \ref{prop_DPP_Z}
 and since
 \bea \label{eq:es011}
 Z_T (\o)  \n = \n    \underset{\hP \in \cP_T }{\sup} \,
      \hE_\hP  \Big[  \wh{Y}^{T,\o}_T      \Big]  \n = \n  \wh{Y}_T (\o)  , \q \fa \o \in \O ,
 \eea
 we see that
 \bea \label{eq:et345}
  \nu_n  \n : = \n
 \inf\big\{t  \n \in \n  [0,T] \n : Z_t  \n - \n  \wh{Y}_t  \n \le \n  1/n    \big\}
 \eea
 is an $\bF-$stopping time.
 Let us also fix   $(t,\o)  \n \in \n  [0,T]  \n \times \n  \O$.

 \no {\bf 1)} {\it Given $\z \n \in \n  \cT$,
 let us first show that }
 \bea  \label{eq:xax_011}
   Z_{   \z \land t} (\o) \ge \ol{\sE}_{\n t} [Z_{  \z}] (\o)  .
 \eea

 If $ \wh{t}  \n := \n     \z (\o)  \n \le \n  t $,
 applying Lemma \ref{lem_shift_stopping_time} with
 $(t,s,\tau)  \n = \n  (0,t,   \z)$ shows that
 $    \z  (\o  \n \otimes_t \n  \O^t)  \n \equiv \n \wh{t} $.
 Since $Z_{\wh{t}} \in \cF_{\wh{t}} \subset \cF_t$ by Proposition \ref{prop_DPP_Z},
 using \eqref{eq:bb421} with $ (t,s,\eta) = \big(0,t,Z_{\wh{t}}\big)  $  shows that
 \bea \label{eq:gk111}
  ( Z_\z  )^{t,\o} (\wt{\o}) = Z_\z (\o \otimes_t \wt{\o})
 = Z \big(\wh{t}, \o \otimes_t \wt{\o} \big) = Z \big(\wh{t}, \o  \big)
= Z\big(   \z (\o) \land t , \o \big) , \q \fa \wt{\o}  \n \in \n  \O^t .
 \eea
It follows that
$\ol{\sE}_{\n t} [Z_\z] (\o)
=  \underset{\hP \in \cP_t}{\sup}   \,     \hE_\hP \big[ (Z_\z )^{t,\o} \big]
= Z_{  \z \land t} (\o) $.

On the other hand, if $     \z (\o)  \n > \n  t $, 
as $\z^{t,\o} \ins  \cT^t$ by Lemma \ref{lem_shift_stopping_time},
applying \eqref{eq:ef011} with $\ga \= \nu \= \z^{t,\o} $ yields that
\beas
 Z_{\z \land t} (\o) =
 Z_t (\o) =   \underset{(\hP,\ga) \in \cP_t \times \cT^t}{\sup}   \,     \hE_\hP
\Big[ \b1_{\{  \ga < \z^{t,\o} \}}   \wh{Y}^{t,\o}_\ga  + \b1_{\{  \ga \ge \z^{t,\o} \}} Z^{t,\o}_{\z^{t,\o}}    \Big]
\ge   \underset{\hP \in \cP_t}{\sup}   \,     \hE_\hP  \big[   (Z_\z)^{t,\o}    \big] =
\ol{\sE}_{\n t} [Z_\z] (\o) .
\eeas

 \no {\bf 2)} {\it Let $\z \n \in \n  \cT$. We next   show that $  Z_{ \nu_n \land \z \land t} (\o)
 \n \le \n  \ol{\sE}_{\n t} [Z_{ \nu_n \land \z}] (\o) $.}

 If $  \nu_n (\o)   \land   \z (\o)  \n \le \n  t $,
 using similar arguments that lead  to \eqref{eq:gk111} yields that
 \if{0}
 applying Lemma \ref{lem_shift_stopping_time} with
 $(t,s,\tau)  \n = \n  (0,t,\nu_n    \n \land \n  \z)$ shows that
 $ (\nu_n    \n \land \n  \z) (\o  \n \otimes_t \n  \O^t)  \n \equiv \n \wh{t} $.
 In light of Proposition \ref{prop_DPP_Z}, using \eqref{eq:bb421} with $ (t,s) = (0,t)  $
 and $ \eta = Z_{\wh{t}} \in \cF_{\wh{t}} \subset \cF_t $ shows that
 \beas
  ( Z_{\nu_n   \land \z}  )^{t,\o} (\wt{\o}) = Z_{\nu_n   \land \z} (\o \otimes_t \wt{\o})
  \n = \n  Z(\wh{t}, \o \otimes_t \wt{\o}) = Z(\wh{t}, \o  )
  \n = \n  Z\big( \nu_n (\o) \land \z (\o) \land t , \o \big) , \q \fa \wt{\o}  \n \in \n  \O^t .
 \eeas
  \fi
 $( Z_{\nu_n    \land   \z}  )^{t,\o} (\wt{\o})
 \n = \n  Z\big( \nu_n (\o)   \land   \z (\o)   \land   t , \o \big) $, $ \fa \wt{\o}  \n \in \n  \O^t $
and thus
$\ol{\sE}_{\n t} [Z_{ \nu_n \land \z}] (\o)
  \n = \n  Z_{ \nu_n \land \z \land t} (\o) $.

On the other hand, suppose that $   \nu_n (\o)     \land       \z (\o)  \n > \n  t $.
We see from Lemma \ref{lem_shift_stopping_time} again  that
$\z_n  \n := \n  (\nu_n   \n  \land  \n   \z)^{t,\o}  \n \in \n  \cT^t$. Let $\e  \n > \n  0$.
Applying \eqref{eq:ef011} with $  \nu  \n = \n  \z_n $, one can find a pair
$(\hP_\e, \ga_\e ) \= (\hP^n_\e, \ga^n_\e ) 
 \n \in \n  \cP_t    \n \times \n \cT^t $  such that
\bea \label{eq:gk114}
 Z_t (\o) \n =  \n   \underset{(\hP,\ga) \in \cP_t \times \cT^t}{\sup}   \,     \hE_\hP
\Big[ \b1_{\{  \ga < \z_n \}}   \wh{Y}^{t,\o}_\ga   \n + \n  \b1_{\{  \ga \ge \z_n \}} Z^{t,\o}_{\z_n}    \Big]
 \n \le  \n    \hE_{\hP_\e}
\Big[ \b1_{\{ \ga_\e < \z_n \}}   \wh{Y}^{t,\o}_{\ga_\e}
 \n + \n  \b1_{\{ \ga_\e \ge \z_n \}} Z^{t,\o}_{\z_n}    \Big]   \n + \n  \e .
\eea
For any $\wt{\o} \in \{ \ga_\e < \z_n \} $, since $ \ga_\e (\wt{\o}) < \z_n (\wt{\o}) =
(\nu_n   \n  \land  \n   \z) (\o \otimes_t \wt{\o}) \le \nu_n (\o \otimes_t \wt{\o}) $,
 the definition of $\z_n$ shows that $\frac{1}{n} < Z \big(\ga_\e (\wt{\o}), \o \otimes_t \wt{\o} \big)
- \wh{Y} \big(\ga_\e (\wt{\o}), \o \otimes_t \wt{\o} \big)
= Z^{t,\o}_{\ga_\e} (\wt{\o}) - \wh{Y}^{t,\o}_{\ga_\e} (\wt{\o}) $.
It follows from   \eqref{eq:gk114}   that
\bea \label{eq:gk115}
 Z_t (\o)  \n \le  \n    \hE_{\hP_\e}
\Big[ \b1_{\{ \ga_\e < \z_n \}}   \wh{Y}^{t,\o}_{\ga_\e}
 \n + \n  \b1_{\{ \ga_\e \ge \z_n \}} Z^{t,\o}_{\z_n}    \Big]   \n + \n  \e
  \n \le \n   \hE_{\hP_\e}
\Big[    Z^{t,\o}_{ \ga_\e \land \z_n}   \n - \n  \frac{1}{n}\b1_{\{ \ga_\e < \z_n \}}  \Big]  \n + \n  \e  .
\eea
 Since   $ \ga_\e (\Pi^0_t)  \n  \in \n  \cT_t  $ by Lemma A.1 of \cite{ROSVU},
 applying \eqref{eq:xax_011} with $\z \= \ga_\e (\Pi^0_t)  \land \nu_n \land \z $ yields that
\bea \label{eq:gk117}
 Z_t (\o) = Z_{\ga_\e (\Pi^0_t) \land \nu_n \land \z \land t} (\o)
 \ge \ol{\sE}_{\n t} [Z_{\ga_\e (\Pi^0_t)  \land \nu_n \land \z}] (\o)
 \ge \hE_{\hP_\e} \Big[ \big( Z_{\ga_\e (\Pi^0_t)  \land \nu_n \land \z} \big)^{t,\o} \Big]
 = \hE_{\hP_\e}  \Big[    Z^{t,\o}_{ \ga_\e \land \z_n  } \Big]   ,
\eea
where we used the fact that for any $\wt{\o} \in \O^t $
\beas
\big( Z_{\ga_\e (\Pi^0_t)  \land \nu_n \land \z} \big)^{t,\o} (\wt{\o})
 \n = \n  Z \big( \ga_\e \big(\Pi^0_t (\o  \n \otimes_t \n  \wt{\o})\big)
  \n \land \n  \nu_n (\o  \n \otimes_t \n  \wt{\o})
 \n \land \n  \z (\o  \n \otimes_t \n  \wt{\o}) , \o  \n \otimes_t \n  \wt{\o} \big)
 \n = \n  Z \big( \ga_\e ( \wt{\o})    \n \land \n  \z_n ( \wt{\o}) , \o  \n \otimes_t \n  \wt{\o} \big)
 \n = \n  Z^{t,\o}_{ \ga_\e \land \z_n  } ( \wt{\o}) .
\eeas

Putting \eqref{eq:gk115} and  \eqref{eq:gk117} together shows  that
 $\hP_\e \{\ga_\e \n < \n  \z_n\}  \n \le \n  n\e$. Then we can deduce from   \eqref{eq:gk114}   and
 \eqref{eq:ef221}  that
\beas
Z_{  \nu_n \land \z \land t} (\o) & \tn = &  \tn
 Z_t (\o) \n \le  \n       \hE_{\hP_\e}
\Big[ \b1_{\{ \ga_\e < \z_n \}}  ( \wh{Y}^{t,\o}_{\ga_\e}
 \n - \n   Z^{t,\o}_{\z_n} )   \n + \n    Z^{t,\o}_{\z_n}    \Big] \n + \n \e
 \n \le \n  2 M_Y \hP_\e \{\ga_\e < \z_n\}  \n + \n  \hE_{\hP_\e}[Z^{t,\o}_{\z_n}] \n + \n \e \\
 & \tn  \le & \tn  \hE_{\hP_\e} \big[(Z_{\nu_n \land \z})^{t,\o} \big]  \n + \n  (1 \n + \n 2 n M_Y )\e
 \le \ol{\sE}_{\n t} [Z_{\nu_n \land \z}] (\o) \n + \n  (1 \n + \n 2 n M_Y )   \e .
\eeas
Letting $\e \to 0$ yields that $Z_{  \nu_n \land \z \land t} (\o) \le \ol{\sE}_{\n t} [Z_{\nu_n \land \z}] (\o) $. \qed

  \no {\bf Proof of Theorem \ref{thm_cst}:}
 Fix $(Y,\wp) \n  \in \n  \fS$.
 Since  both $\wh{Y}$ and $Z$ are $\bF-$adapted  processes
 with  all continuous paths by Proposition \ref{prop_DPP_Z},
 we see from 
 \eqref{eq:ef221} and \eqref{eq:es011} that
 $\wh{\nu} \n : = \n  \inf \big\{ t  \n \in \n  [0,T] \n : Z_t  \n = \n  \wh{Y}_t \big\}$
 is an $\bF-$stopping time. For any $n \ins \hN$,
 let $\nu_n$ be the $\bF-$stopping time defined in \eqref{eq:et345}.
 Since  $Z$ is an $\ol{\sE}-$martingale over  $[0, \nu_n]$ by Proposition \ref{prop_Z_martingale},
 one can find  a $ \hP_n \ins  \cP $ satisfying \eqref{eq:es014}.
 By (P1), $\{\hP_n\}^\infty_{n = 2}$ 
 has a weakly convergent subsequence $\{\hP_{m_j}\}_{j \in \hN}$
 with limit $\wh{\hP} \n \in \n  \cP$.

 \if{0}
  If $Z_0  \n = \n  \wh{Y}_0  \n = \n  Y_0   $, then $\wh{\nu}  \n = \n  0$ and it holds for any $\hP \in \cP$ that
 $   \hE_\hP \big[ Z_0   \big]   \n = \n  Z_0 $. So the conclusion   holds.

 Next, let us assume $Z_0  \n > \n  \wh{Y}_0  \n = \n  Y_0 $ and set $   n_0  \n : = \n   1  \n + \n
   \big\lfloor (Z_0  \n - \n  \wh{Y}_0)^{-1}  \big\rfloor \n > \n (Z_0  \n - \n  \wh{Y}_0)^{-1} $.
 For any integer $n \ge n_0$,  as $Z_0 - \wh{Y}_0 > \frac{1}{n_0} \ge \frac{1}{n }$,
   we see from  \ref{lem_stopping_time} again that  $ \nu_n  \n : = \n
 \inf\big\{t  \n \in \n  [0,T] \n : Z_t  \n - \n  \wh{Y}_t  \n \le \n   \frac{1}{n}  \big\} \n > \n 0 $
 is also an $\bF-$stopping time.
 \fi

 When   $m_j  \n \ge \n  n  $,  (so $ \nu_n  \n \le \n   \nu_{m_j}    $),
 applying Lemma \ref{lem_P_supermg}
 with $(\hP, \tau,\ga)  \n = \n  (\hP_{m_j},\nu_n, \nu_{m_j} )$, we see from \eqref{eq:es014} that
 \bea   \label{eq:es021}
  Z_0    \le   \hE_{\hP_{m_j}} \big[ Z_{\nu_{m_j}} \big] + 2^{-m_j}
 \le \hE_{\hP_{m_j}} \big[ Z_{\nu_n} \big] + 2^{-m_j}    .
 \eea

  \no {\bf 1)} {\it Before sending  $j$ to $\infty$   in order to approximate the distribution $ \wh{\hP} $
  in \eqref{eq:es021},  we need to approach
     $\{\nu_n\}_{n \in \hN}$ by a sequence
  $\big\{\wh{\th}_n \big\}_{n \in \hN}$ of   Lipschitz  continuous random variables. }

   Fix integer $n  \n \ge \n  2$.  
   There exists a $\l_n  \n > \n  0$ such that
   $\rY (x)  \n \vee \n  \wrY (x)  \n \le \n  \frac{1}{2 n (n+1)}$, $\fa x  \n \in \n  [0, \l_n]$.
   Let    $\o  \n \in \n  \O$,
  set $ \dis \d_n ( \o )  \n : = \n  \frac{\l_n}{2(1 \+ \k_\wp)}  \n \land \n  \frac{(\phi^\o_T)^{-1}(\l_n/2)}{\k_\wp}$
  with  $(\phi^\o_T)^{-1} (x)
 \n : = \n  \inf\{y  \n > \n  0 \n : \phi^\o_T (y)  \n = \n  x \} $, $\fa x  \n > \n  0 $,
 and let   $\o'  \n \in \n  \ol{O}_{\d_n ( \o )} (\o)$.
 Given $t  \n \in \n  [0,T]$,   set $s \n := \n   \wp(\o)  \n \land \n t $
 and $s' \n := \n   \wp(\o') \n \land \n t $.
    By \eqref{eq:ec017}, $ |s \n - \n s' |  \n \le \n  \big| \wp(\o)  \n - \n  \wp(\o') \big|
   \n \le \n   \k_\wp   \|\o    \n - \n  \o'  \|_{0,T}   $.
 Then   \eqref{eq:aa211} implies that
 \bea
 && \hspace{-1cm}  \big|\wh{Y}  ( t,\o ) \n - \n  \wh{Y}  ( t,\o') \big| \n  = \n
 | Y  ( s,\o ) - Y ( s',\o') |
  \le    \rY \Big( |s-s'  |
 + \underset{r \in [0,T]}{\sup} \big|  \o   (r \land s)
  -  \o'   (r \land s' ) \big|   \Big)  \nonumber \\
  &  & \le      \rY \Big(  \k_\wp   \|\o \n  - \n  \o'  \|_{0,T}
  \n + \n  \underset{r \in [0,T]}{\sup}  |   \o   (r  \n \land \n  s  )
  \n - \n   \o  (r  \n \land \n  s'  )  |
  \n + \n  \underset{r \in [0,T]}{\sup} \big|  \o (r  \n \land \n  s' )
  \n - \n  \o' (r  \n \land \n  s' ) \big|   \Big) \nonumber \\
  &  & \le
  \rY \Big(  (1 \n + \n \k_\wp)   \|\o \n  - \n  \o'  \|_{0,T}
  \n + \n   \phi^\o_T \big( |s' \n - \n s| \big)   \Big)
   \n \le \n  \rY \Big(  (1 \n + \n \k_\wp)   \|\o \n  - \n  \o'  \|_{0,T}
  \n + \n   \phi^\o_T \big( \k_\wp    \|\o \n  - \n  \o'  \|_{0,T} \big)   \Big)
   \n \le \n  \frac{1}{2 n (n+1)}  .  \qq   \label{eq:gk011}
 \eea

Taking $t= \nu_n (\o) $, we see from \eqref{eq:gd011} that
 \beas
&& \hspace{-1cm} \big|  (Z \n - \n \wh{Y}) ( \nu_n (\o),\o)  \n - \n   (Z \n - \n \wh{Y}) ( \nu_n (\o),\o') \big|
  \n \le  \n  \big|   Z ( \nu_n (\o),\o)  \n - \n    Z ( \nu_n (\o),\o') \big|
  \n + \n   \big| \wh{Y}  ( \nu_n (\o),\o)  \n - \n \wh{Y}  ( \nu_n (\o),\o') \big|  \\
&&  \n \le \n   \wrY \Big( (1 \n + \n \k_\wp) \|\o \n - \n \o'\|_{0,\nu_n (\o)}
  \n + \n \phi^\o_{\nu_n (\o)} \big( \k_\wp   \|\o   \n  - \n  \o'  \|_{0,\nu_n (\o)} \big)  \Big)
  \n + \n \frac{1}{2 n (n \n + \n 1)}  \\
&&  \n \le \n    \wrY \Big( (1 \n + \n \k_\wp) \|\o \n - \n \o'\|_{0,T}
  \n + \n \phi^\o_T \big( \k_\wp   \|\o   \n  - \n  \o'  \|_{0,T} \big)  \Big)
  \n + \n \frac{1}{2 n (n \n + \n 1)}
   \n \le  \n  \frac{1}{  n (n \n + \n 1)}  \n < \n  \frac{1}{ (n \n - \n 1) n }  .
\eeas
 As   the continuity of $Z  \n - \n \wh{Y}$    shows that
 \bea \label{eq:es041}
 \big( Z  \n - \n \wh{Y} \big) ( \nu_n (\o),\o)  \n \le \n  \frac{1}{n } ,
 \eea
 it follows that
 $ (Z  \n - \n \wh{Y}) ( \nu_n (\o),\o')
  \n \le \n  
  \frac{1}{n }  \n + \n  \frac{1}{ (n-1) n  }  \n = \n  \frac{1}{n-1 }   $,
 so $\nu_{n-1} (\o')  \n \le \n  \nu_n (\o)$.
 Analogously,
 taking $t= \nu_{n+1} (\o') $ in \eqref{eq:gk011} yields that
  \beas
&& \hspace{-0.8cm} \big|  (Z \n - \n \wh{Y}) ( \nu_{n+1} (\o'),\o)  \n - \n   (Z \n - \n \wh{Y}) ( \nu_{n+1} (\o'),\o') \big|
  \n \le  \n  \big|   Z ( \nu_{n+1} (\o'),\o)  \n - \n    Z ( \nu_{n+1} (\o'),\o') \big|
  \n + \n   \big| \wh{Y}  ( \nu_{n+1} (\o'),\o)  \n - \n \wh{Y}  ( \nu_{n+1} (\o'),\o') \big|  \\
&&  \n \le \n   \wrY \Big( (1 \n + \n \k_\wp) \|\o \n - \n \o'\|_{0,T}
  \n + \n \phi^\o_T \big( \k_\wp   \|\o   \n  - \n  \o'  \|_{0,T} \big)  \Big)
  \n + \n \frac{1}{2 n (n \n + \n 1)}
  \n \le  \n  \frac{1}{  n (n \n + \n 1)}    ,
\eeas
and that  $ (Z  \n - \n \wh{Y}) ( \nu_{n+1} (\o'),\o)
  \n < \n   (Z  \n - \n \wh{Y}) ( \nu_{n+1} (\o'),\o')  \n + \n  \frac{1}{ n(n+1)  }
  \n \le \n  \frac{1}{n+1 }  \n + \n  \frac{1}{ n(n+1) }  \n = \n  \frac{1}{n }  $,
 which shows that  $\nu_n (\o)  \n \le \n  \nu_{n+1} (\o') $.

   Now, we can   apply Lemma \ref{lem_sandwich}  with
 $ (\O_0,\ul{\th}, \th, \ol{\th},I, \d (\o),\e )
  \n = \n  (\O, \nu_{n-1},   \nu_n,   \nu_{n+1}, [0,T], \d_n (\o), 2^{-n} )$
 to find an   open subset $\wh{\O}_n$  of $\O$   and
  a Lipschitz  continuous random variable   $\wh{\th}_n \n : \O \to [0,T]$ such that
 \bea \label{eq:eh211}
  \underset{\hP \in \cP}{\sup} \,  \hP \big( \wh{\O}^c_n \big) \le 2^{-n} , \q
  \nu_{n-1} - 2^{-n}   \n < \n  \wh{\th}_n  \n < \n  \nu_{n+1} + 2^{-n}
 \; \hb{ on } \; \wh{\O}_n .
 \eea

  \no {\bf 2)} {\it Next, let us estimate the expected difference
   $\hE_{\hP_{m_j}} \big[ \big| Z_{\wh{\th}_n} \-  Z_{\nu_n} \big| \big]$.}

 Given $\o \n \in \n  \wh{\O}_{n-1}   \n   \cap   \n   \wh{\O}_{n+1} $,
 as $ \wh{\th}_{n-1}  \n - \n  2^{-n+1}  \n < \n  \nu_n  \n < \n  \wh{\th}_{n+1}  \n + \n  2^{-n-1} $,
   $t  \n := \n  \wh{\th}_n (\o)    \n  \land  \n    \nu_n (\o) $
 and $s  \n := \n  \wh{\th}_n (\o)    \n  \vee  \n    \nu_n (\o)$  satisfy
 \bea
 s \n - \n t & \tn \dn =  & \tn \dn  \big| \nu_n (\o) \n - \n \wh{\th}_n (\o) \big|
  \n <  \n  (\wh{\th}_{n-1}  \n - \n  \wh{\th}_n  \n - \n  2^{-n+1})^-
  \n \vee \n  (\wh{\th}_{n+1}  \n - \n  \wh{\th}_n  \n + \n  2^{-n-1})^+
  \n \le \n |\wh{\th}_{n-1}  \n - \n  \wh{\th}_n  \n - \n  2^{-n+1}|
  \n \vee \n  |\wh{\th}_{n+1}  \n - \n  \wh{\th}_n  \n + \n  2^{-n-1}|  \nonumber  \\
   & \tn \dn \le  & \tn \dn
   |\wh{\th}_{n-1} (\o)  \n - \n  \wh{\th}_n (\o) |  \n + \n  |\wh{\th}_{n+1} (\o)   \n - \n  \wh{\th}_n (\o) |
  \n  + \n  2^{-n+1}  \n := \n  \wh{\d}_n (\o) . \label{eq:uxu017}
 \eea
 Set  $\phi_n (\o) \n := \n   (1 \n + \n \k_\wp)
 \Big(  \big( \wh{\d}_n (\o) \big)^{\frac{q_1}{2}}  \n + \n \phi^\o_T \big( \wh{\d}_n (\o) \big)  \Big)   $.
  Then \eqref{eq:es017} shows that
 \beas
 \big| Z_{\wh{\th}_n} (\o) \n  -  \dn  Z_{\nu_n}  (\o) \big|
  \n \= \n  \big| Z  (t,\o)  \n - \dn   Z   (s,\o) \big|
   \ls 2 C_{\n \vr} M_Y   \Big( \n \big(\wh{\d}_n (\o)\big)^{ \frac{q_1}{2}}
   \n \ve     \big(\wh{\d}_n (\o)\big)^{q_2-\frac{q_1}{2}} \n \Big)
   \dn + \n \wrY \big(  \wh{\d}_n (\o)  \big)
   \dn + \n    \wrY  \big(   \phi_n (\o)  \big)
    \n \vee \n  \wh{\wh{\rho}\,}_Y \big(   \phi_n (\o)  \big)   \df \xi_n (\o)     .
  \eeas

  Let $j \ins \hN$ with $m_j  \n \ge \n  n  $.
 We see from  \eqref{eq:es021}, \eqref{eq:ef221} and \eqref{eq:eh211} that
 \bea
   Z_0 \- 2^{-m_j}  & \tn  \ls  & \tn  \hE_{\hP_{m_j}} \big[ Z_{\wh{\th}_n} \big] \+
 \hE_{\hP_{m_j}} \big[ \big| Z_{\wh{\th}_n} \-  Z_{\nu_n} \big| \big]  \ls
  \hE_{\hP_{m_j}} \big[  Z_{\wh{\th}_n}  \+  \b1_{\wh{\O}_{n-1}   \cap   \wh{\O}_{n+1}}  ( \xi_n \ld 2 M_Y ) \big]
      \n + \n  2 M_Y \hP_{m_j} \big(\wh{\O}^c_{n-1}   \cup   \wh{\O}^c_{n+1}\big) \nonumber  \\
      & \tn   \ls  & \tn   \hE_{\hP_{m_j}} \big[  Z_{\wh{\th}_n}  \+   ( \xi_n \ld 2 M_Y ) \big]
        \n + \n  5 M_Y 2^{-n}  . \q  \label{eq:es035}
 \eea

 The random variables  $\wh{\th}_{n-1}, \wh{\th}_n, \wh{\th}_{n+1}$ are Lipschitz  continuous   on $\O$, so is $\wh{\d}_n $.
 Then one can deduce that
 \bhe
 \bea \label{eq:a114}
  \hb{ $\o \n \to \n  \phi^\o_T \big(\wh{\d}_n (\o) \big)$ is  a continuous random variable on $\O$,  }
 \eea
 \ehe
 which together
 with the Lipschitz continuity of $ \wh{\d}_n $ implies that $\phi_n$
 and thus $\xi_n$ are   continuous random variables on $\O$. Moreover, the
 Lipschitz continuity of random variable $\wh{\th}_n$ and the continuity of process $Z$
 implies that
 \bhe
 \bea \label{eq:a117}
 \hb{$Z_{\wh{\th}_n}$ is also a   continuous random variable on $\O$.}
 \eea
 \ehe
 Letting   $j \to \infty$ in \eqref{eq:es035},
 we see from the continuity of random variables $Z_{\wh{\th}_n}$  and $\xi_n$ that
 \bea   \label{eq:es043}
   Z_0   \le  \hE_{\wh{\hP}} \big[ Z_{\wh{\th}_n} \n + \n ( \xi_n \land 2 M_Y ) \big]  \n + \n  5 M_Y 2^{-n}  ,
   \q \fa n \ge 2 .
 \eea

 \no {\bf 3)} {\it Finally, we   use the convergence of $ \wh{\th}_n $ to $\wh{\nu} $
 and    the   continuity of $Z$ to derive the $\ol{\sE}-$martingality of $Z$ over $[0,\wh{\nu}]$.}

 Set  $\wh{\nu}' \df \lmtu{n \to \infty} \nu_n   \ls \wh{\nu} $.
   The continuity of $Z  \n - \n \wh{Y}$,   \eqref{eq:es041} and \eqref{eq:ef221} imply that
 $    Z_{\wh{\nu}'}  \n - \n \wh{Y}_{\wh{\nu}'}   \n = \n  0 $, thus
 $ \wh{\nu} \= \wh{\nu}' \= \lmtu{n \to \infty} \nu_n    $. Then we can deduce from \eqref{eq:eh211} that
 $ \lmt{n \to \infty} \wh{\th}_n (\o) \= \wh{\nu} (\o) $, $ \fa
 \o \ins \underset{n = 3}{\overset{\infty}{\cup}} \, \underset{k \ge n}{\cap} \wh{\O}_k $.
 As $ \dis  \sum^\infty_{n = 3} \wh{\hP} \big( \wh{\O}^c_n \big)
 \le \sum^\infty_{n = 3}  \underset{\hP \in \cP}{\sup} \,  \hP \big( \wh{\O}^c_n \big) \le \frac14   $,
 the Borel-Cantelli Lemma implies that
 $ \wh{\hP} \Big( \underset{n = 3}{\overset{\infty}{\cup}} \, \underset{k \ge n}{\cap} \wh{\O}_k \Big) = 1  $.
 So
 \bea    \label{eq:et359}
  \lmt{n \to \infty} \wh{\th}_n   = \wh{\nu}    , \q   \wh{\hP}- a.s.
 \eea
 It follows that
 $ \lmt{n \to \infty} \wh{\d}_n   = 0 $, $\wh{\hP}-$a.s.
 and thus $ \lmt{n \to \infty} \xi_n   = 0 $, $\wh{\hP}-$a.s.
 Eventually, letting $n \to \infty$ in \eqref{eq:es043}, we can deduce from the continuity of process $Z$, $\wh{Y}$ and
    the bounded dominated convergence theorem    that
 \beas
   Z_0   \le  \hE_{\wh{\hP}} \big[ Z_{\wh{\nu}} \big] 
   \le 
   \ol{\sE}_{\n 0} \big[ Z_{\wh{\nu}} \big] =    \ol{\sE}_{\n 0} \big[ \wh{Y}_{\wh{\nu}} \big]
   \le \underset{(\hP, \ga) \in \cP \times \cT}{\sup}  \hE_\hP \big[ \wh{Y}_\ga \big]  = Z_0 .
 \eeas
 Hence, $ Z_0 \n = \n  \ol{\sE}_{\n 0} \big[ Z_{\wh{\nu}} \big]
   \n = \n  \hE_{\wh{\hP}} \big[ Z_{\wh{\nu}}   \big]
   \n = \n  \hE_{\wh{\hP}} \big[ \wh{Y}_{\wh{\nu}}   \big]  $.

      Next, let $\z \in \cT$. For any $\hP \in \cP$, we see from Lemma
   \ref{lem_P_supermg} that
   \bea \label{eq:es573}
   Z_0 = \hE_\hP [ Z_0 ] \ge \hE_\hP \big[ Z_{\wh{\nu} \land \z} \big] \ge \hE_\hP \big[ Z_{\wh{\nu}} \big]  .
   \eea
   Taking supremum over $\hP \in \cP$ yields that
   $Z_0 \ge \ol{\sE}_{\n 0} \big[ Z_{\wh{\nu} \land \z} \big] \ge \ol{\sE}_{\n 0} \big[ Z_{\wh{\nu}} \big] = Z_0 $.
   In particular, taking $\hP = \wh{\hP}$ in \eqref{eq:es573} shows that
   $  Z_0 \ge \hE_{\wh{\hP}} \big[ Z_{\wh{\nu} \land \z} \big] \ge \hE_{\wh{\hP}} \big[ Z_{\wh{\nu}} \big] = Z_0 $.
   \qed

 \subsection{Proofs of   results in Section  \ref{sec:OSRM}}

  \no {\bf Proof of Proposition \ref{prop_wp_n}:}
  Set $   \fn_0  \n : = \n    1 \n + \n   \lfloor  \sX^{-1}_0  \rfloor  \n > \n  \sX^{-1}_0   $. Given
  $k \in \hN \cup \{0\}$, since   $\sX$ is an $\bF-$adapted process with all continuous paths
  and since $\sX_0  \n > \n \frac{1}{   \fn_0   }  \n \ge \n  \frac{1}{k + \fn_0   } $,
 we see 
    that $ \wh{\tau}_k  \n : = \n
 \inf\big\{t  \n \in \n  [0,T] \n : \sX_t  \n \le \n   \frac{1}{k+\fn_0}  \big\} \n \land \n T    $
 is an $\bF-$stopping time satisfying $0  \n < \n  \wh{\tau}_k (\o)  \n \le \n  \tau_0 (\o) $, $\fa \o \in \O$.
 In particular, if $\{t  \n \in \n  [0,T] \n : \sX_t  (\o')  \n \le \n  0\}$  is not empty for some $\o' \in \O$,
 then $\wh{\tau}_k (\o') \n <  \n  \tau_0 (\o')$.
 Let $ \{ \d_k   \}_{k \in \hN} $  be   a  sequence decreasing to $0$  such that
 $   \rX (\d_k)  \n \le \n   \frac{1}{  (k+\fn_0) (k+\fn_0+1)}$, $\fa k \in \hN$.

   \no {\bf a)} {\it First, we construct an auxiliary increasing sequence $\{\vth_\ell\}_{\ell \in \hN}$ of Lipschitz continuous stopping times.  }

     Fix $k \n \in \n  \hN$.
   For $i \n = \n k \n - \n 1$, $k$, let $\o, \o'  \n \in \n  \O$ with
   $\|\o' \n - \n \o\|_{0,\wh{\tau}_{i+1} (\o)}\n \le \n  \d_k$,
  \eqref{eq:es111} shows that
 \beas
 \big|  \sX ( \wh{\tau}_{i+1} (\o),\o')  \n - \n   \sX ( \wh{\tau}_{i+1} (\o),\o) \big|
  \n \le \n  \rX \big( \|\o' \n - \n \o\|_{0,\wh{\tau}_{i+1} (\o)} \big)
  \n \le \n  \rX (\d_k) \le \frac{1}{ (k \n + \n \fn_0) (k \n + \n \fn_0 \n + \n 1) } .
\eeas
If $ \sX_t (\o) \n > \n  \frac{1}{i+\fn_0+1}$ for all $t  \n \in \n  [0,T]$, then
$ \wh{\tau}_{i+1} (\o)  \n = \n  T  \n \ge \n  \wh{\tau}_i (\o') $. On the other hand, if
the set $\big\{t  \n \in \n  [0,T] \n : \sX_t (\o)  \n \le \n   \frac{1}{i+\fn_0+1}  \big\}$is not empty,
   the continuity of $\sX$ imply that
 $ \sX ( \wh{\tau}_{i+1} (\o),\o)  = \frac{1}{i+\fn_0+1} $,
 it follows that  $ \sX ( \wh{\tau}_{i+1} (\o),\o')
  \n \le \n 
  \frac{1}{i+\fn_0+1}  \n + \n  \frac{1}{ (i+\fn_0) (i+\fn_0+1) }  \n = \n  \frac{1}{i+\fn_0 }  $,
 so one still has $\wh{\tau}_i (\o')  \n \le \n  \wh{\tau}_{i+1} (\o)$.
 Then we can   apply Lemma \ref{lem_sandwich2}  with
 $ (\th_1, \th_2, \th_3, \d,\k)  \n = \n  \big(\wh{\tau}_{k-1},   \wh{\tau}_k,   \wh{\tau}_{k+1}, \d_k,
 2T/\d_k \big)$ to find   a $\wh{\wp}_k  \n \in \n  \cT$ such that
 \bea \label{eq:eh211b}
  \wh{\tau}_{k-1} (\o)   \n \le \n  \wh{\wp}_k (\o)   \n \le \n  \wh{\tau}_{k+1} (\o)   \n \le \n \tau_0 (\o) ,
 \q \fa \o \in \O ,
 \eea
 (the last inequality is strict  if the set
 $\{t  \n \in \n  [0,T] \n : \sX_t  (\o)  \n \le \n  0\}$  is not empty) and that given  $\o_1, \o_2 \in \O$,
 \bea \label{eq:ec014}
 \big| \wh{\wp}_k (\o_1) - \wh{\wp}_k (\o_2) \big| \le 2T \d^{-1}_k   \|\o_1 - \o_2 \|_{0,t_0}
 \eea
 holds for any    $t_0  \dn \in \dn   \big[   \wh{b}_k,T  \big]
        \cup         \big\{ t  \n \in \n  [\wh{a}_k, \wh{b}_k) \n :
  t  \n \ge \n  \wh{a}_k   \n + \n  2T \d^{-1}_k  \| \o_1   \n -   \n    \o_2 \|_{0,t} \n \big\} $,
 where $ \wh{a}_k  \n := \n  \wh{\wp}_k (\o_1)    \land    \wh{\wp}_k (\o_2) $
 and  $ \wh{b}_k  \n := \n  \wh{\wp}_k (\o_1)    \vee    \wh{\wp}_k (\o_2) $.

  Let $\ell \n \in \n  \hN$.
  We   define and $\bF-$stopping time $\vth_\ell \n  := \n    \underset{k=1,\cds , \ell}{\max} \,  \wh{\wp}_k   $.
  Let $\o_1, \o_2  \n \in \n  \hN$ and set $\fra_\ell  \n : = \n  \vth_\ell (\o_1)  \n \land \n  \vth_\ell (\o_2)$,
  $\fb_\ell  \n : = \n  \vth_\ell (\o_1)  \n \vee \n  \vth_\ell (\o_2)$.
  To see that
      \bea \label{eq:ga011}
        \big| \vth_\ell (\o_1) \n - \n \vth_\ell (\o_2) \big|
    \n \le \n  2  T \d^{-1}_\ell  \|\o_1  \n - \n  \o_2\|_{0,t_0}
    \eea
    holds for any  $t_0  \n \in \n  [\fb_\ell,T ]   \cup   \big\{ t  \n \in \n  [\fra_\ell, \fb_\ell) \n :
  t  \n \ge \n  \fra_\ell   \n + \n  2T \d^{-1}_\ell  \| \o_1   \n -     \o_2 \|_{0,t} \big\} $,
  we first let   $t_0 \in  [\fb_\ell, T  ]$. For any   $k  \n = \n 1, \cds  \n , \ell$,
  since $ \wh{b}_k 
  \n \le \n    \vth_\ell (\o_1)  \n \vee \n  \vth_\ell (\o_2)  \n = \n \fb_\ell
  \n \le \n  t_0   $,
   applying \eqref{eq:ec014}   yields that
   \bea \label{eq:ej171}
    \big|\wh{\wp}_k (\o_1) \n - \n \wh{\wp}_k (\o_2) \big|
    \n \le \n  2  T \d^{-1}_k  \|\o_1  \n - \n  \o_2\|_{0,t_0}
    \n \le \n  2  T \d^{-1}_\ell  \|\o_1  \n - \n  \o_2\|_{0,t_0} .
    \eea
    It follows that $ \wh{\wp}_k (\o_1)  \n \le \n  \wh{\wp}_k (\o_2) \n + \n
     2  T \d^{-1}_\ell  \|\o_1  \n - \n  \o_2\|_{0,t_0}
     \n \le \n  \vth_\ell (\o_2)  \n + \n  2  T \d^{-1}_\ell  \|\o_1  \n - \n  \o_2\|_{0,t_0}  $. Taking
    maximum over $k  \n = \n 1, \cds  \n , \ell$ shows that
    $ \vth_\ell (\o_1)    \n  \le \n  \vth_\ell (\o_2)  \n + \n  2  T \d^{-1}_\ell  \|\o_1  \n - \n  \o_2\|_{0,t_0} $.
    Then exchanging the roles of $\o_1$ and $\o_2$ yields \eqref{eq:ga011}.

  We  next  suppose that the set   $    \big\{ t  \n \in \n  [\fra_\ell, \fb_\ell) \n :
  t  \n \ge \n  \fra_\ell   \n + \n  2T \d^{-1}_\ell  \| \o_1   \n -     \o_2 \|_{0,t} \big\} $
  is not empty and contains $t_0$.
  Given $k  \n = \n 1, \cds  \n , \ell$, since $ t_0 \ins [\fra_\ell, \fb_\ell)
  \sb \big[ \wh{a}_k , T \big]  $ and since
   \beas
   \wh{\wp}_k (\o_1)  \n \land \n  \wh{\wp}_k (\o_2)   \n + \n  2T\d^{-1}_k \|\o_1  \n - \n  \o_2\|_{0,t_0}
    \n  \le  \n   \vth_\ell (\o_1)  \n \land \n   \vth_\ell (\o_2)  \n + \n  2T\d^{-1}_\ell \|\o_1  \n - \n  \o_2\|_{0,t_0}
      \n \le \n  t_0 ,
   \eeas
   applying \eqref{eq:ec014}   yields \eqref{eq:ej171} and thus leads to \eqref{eq:ga011} again.

   Now, fix $n \n \in \n  \hN$. We set $\ell  \n := \n  \lceil \log_2 (n \n + \n 2 ) \rceil  \n \ge \n 2 $,
 $ j  \n : = \n  n \n + \n 2    \n - \n  2^{\ell-1}$ and
 define $ \wp_n  \n : = \n  (\vth_{\ell-1}  \n + \n  j 2^{1-\ell} T)  \n \land \n  \vth_\ell  \n \in \n  \cT  $.

  \no {\bf b)} {\it In this step, we show that $\wp_n$'s is   the increasing sequence of
  Lipschitz continuous stopping times in quest such that the increment $\wp_{n+1} \- \wp_n$ is bounded by $\frac{2T}{n+3}$.}

 Since  $\ell \n - \n 1  \n < \n  \log_2(n \n + \n 2)  \n \le \n \ell $, 
 we see that $1  \n \le \n  j  \n \le \n  2^{\ell-1}$.
 If $j \n < \n  2^{\ell-1}$,  as $n \n + \n 2  \n = \n  2^{\ell-1} \n + \n j  \n \le \n  2^\ell \n - \n 1$,
 one has $ \ell  \n = \n  \lceil \log_2 (n \n + \n 2 ) \rceil  \n \le \n
  \lceil \log_2(n \n + \n 3) \rceil \n \le \n \ell $, so $  \lceil \log_2(n \n + \n 3) \rceil \n = \n \ell $.
 Then \eqref{eq:ec017} implies that
  \beas
 \q  0 \n \le \n  \wp_{n+1} (\o)  \n - \n  \wp_n  (\o)  \n = \n
   \big(\vth_{\ell-1}  (\o)  \n + \n  (j \n + \n 1) 2^{1-\ell} T \big)  \n \land \n  \vth_\ell  (\o)
 \n - \n  \big(\vth_{\ell-1}  (\o)  \n + \n  j 2^{1-\ell} T\big)  \n \land \n  \vth_\ell  (\o)
  \n \le \n  2^{1-\ell} T  \n \le \n  \frac{2T}{n+3}   , \q \fa \o  \n \in \n  \O .
 \eeas
 On the other hand, if $j \n = \n 2^{\ell-1} $, i.e. $n \n + \n 2  \n = \n 2^\ell$, then
 $ \wp_n  \n = \n  (\vth_{\ell-1}  \n +  \n   T)  \n \land \n  \vth_\ell  \n = \n  \vth_\ell   $
 and $ \lceil \log_2(n \n + \n 3) \rceil \n = \n
 \lceil \log_2(2^\ell \n + \n 1) \rceil \n = \n \ell \n + \n  1 $. Applying \eqref{eq:ec017} again yields that
 \beas
 0  \n \le \n  \wp_{n+1} (\o)   \n - \n  \wp_n (\o)
  \n  = \n  \big(\vth_\ell  (\o)  \n  +  \n   2^{-\ell} T \big)  \n \land \n  \vth_{\ell+1}  (\o)
   \n - \n   \vth_\ell  (\o)   \n \land \n  \vth_{\ell+1}  (\o)
   \n \le \n  2^{- \ell } T  \n = \n  \frac{T}{n+2}  \n < \n  \frac{2T}{n+3}   , \q \fa \o  \n \in \n  \O .
 \eeas

  Since $   \wh{\tau}_{\ell-2}
  \n = \n  \inf\big\{t  \n \in \n  [0,T] \n : \sX_t  \n \le \n   \frac{1}{\ell-2 +\fn_0}  \big\} \n \land \n T
  \n = \n  \inf\big\{t  \n \in \n  [0,T] \n : \sX_t
  \n \le \n (\lceil \log_2(n \n + \n 2) \rceil   \n + \n  \lfloor  \sX^{-1}_0  \rfloor
   \n - \n 1  )^{-1}  \big\}  \n \land \n T  \= \tau_n $ by \eqref{def_tau_n}, we can deduce from
   \eqref{eq:eh211b} that
 \beas
 \q  \tau_n (\o) \n = \n  \wh{\tau}_{\ell-2}  (\o)   \n \le \n  \wh{\wp}_{\ell-1} (\o)
   \n \le \n  \vth_{\ell-1} (\o)
     \n \le \n  (\vth_{\ell-1}  (\o)  \n + \n  j 2^{1-\ell} T)  \n \land \n  \vth_\ell (\o)
      \n = \n \wp_n  (\o)  \n \le \n  \vth_\ell (\o)
      \n = \n   \underset{i=1,\cds , \ell}{\max} \, \wh{\wp}_k  (\o)  \n \le \n  \tau_0  (\o)    ,
      ~ \fa \o  \n \in \n  \O ,
  \eeas
  where  the last inequality is strict  if  the set
 $\{t  \n \in \n  [0,T] \n : \sX_t  (\o)  \n \le \n  0\}$  is not empty.

\no {\bf c)} {\it It remains to show the Lipschitz continuity of $\wp_n$. }

 Set $\k_n \n := \n  2T \d_\ell^{-1}  \n = \n  2T \big( \d_{\lceil \log_2(n  +  2) \rceil} \big)^{-1} $,
 which is   increasing in $n$ and converges to $\infty$.
 Let $\o_1, \o_2  \n \in \n  \hN$ and set $a_n  \n : = \n  \wp_n (\o_1)  \n \land \n  \wp_n (\o_2)$.
  We assume  without loss of generality  that
  $ a_n  \n = \n \wp_n (\o_1)  \n \le \n  \wp_n (\o_2) $ and discuss by two cases:   

  \no i)
  When $\wp_n (\o_1)  \n = \n  \vth_{\ell-1} (\o_1)  \n + \n  j2^{1-\ell} T $,
  one has
  \bea \label{eq:es147}
    \wp_n (\o_2)  \n - \n  \wp_n (\o_1)  \n = \n  \wp_n (\o_2)  \n - \n  \vth_{\ell-1} (\o_1)  \n - \n  j 2^{1-\ell} T
   \n \le \n   \vth_{\ell-1} (\o_2) \n - \n  \vth_{\ell-1} (\o_1) .
   \eea
  Applying \eqref{eq:ga011}   with $t_0 \n = \n T$ shows that
  $  \wp_n (\o_2)  \n - \n  \wp_n (\o_1)  \n \le \n  2T \d^{-1}_{\ell-1}    \| \o_1   \n - \n   \o_2 \|_{0,T}
     \n \le \n  \k_n  \| \o_1   \n - \n   \o_2 \|_{0,T}  $.
  On the other hand,  suppose that the set  $ \big\{ t  \n \in \n  [a_n, T) \n : t  \n \ge \n  a_n
   \n + \n   \k_n  \| \o_1   \n - \n   \o_2 \|_{0,t} \big\}     $
     is not empty and contains $t_0$.
   since   $ \vth_{\ell-1} (\o_1)  \n = \n  \wp_n (\o_1)  \n - \n  j2^{1-\ell} T
    \n \le \n  \wp_n (\o_2)  \n - \n  j2^{1-\ell} T
    \n \le \n  \vth_{\ell-1} (\o_2)  $,
   we see that $\fra_{\ell-1} = \vth_{\ell-1} (\o_1)$ and  can deduce that
   \beas
   \hspace{-5mm}
   t_0    \n \ge \n  a_n  \n + \n   \k_n  \| \o_1   \n - \n   \o_2 \|_{0,t_0}
    \n =  \n   \vth_{\ell-1} (\o_1)  \n + \n  j2^{1-\ell} T   \n + \n   \k_n  \| \o_1   \n - \n   \o_2 \|_{0,t_0}
    \n > \n  \vth_{\ell-1} (\o_1)    \n + \n   2T \d^{-1}_{\ell-1}  \| \o_1   \n - \n   \o_2 \|_{0,t_0}
    \n = \n  \fra_{\ell-1}  \n + \n  2T \d^{-1}_{\ell-1}  \| \o_1   \n - \n   \o_2 \|_{0,t_0}     .
    \eeas
  Then \eqref{eq:es147} and  \eqref{eq:ga011}    imply  that
  $  \wp_n (\o_2)  \n - \n  \wp_n (\o_1)      \n \le \n  2T \d^{-1}_{\ell-1}    \| \o_1   \n - \n   \o_2 \|_{0,t_0}
     \n \le \n  \k_n  \| \o_1   \n - \n   \o_2 \|_{0,t_0}  $.

  \no ii) When $\wp_n (\o_1) \n  = \n  \vth_\ell (\o_1)  $,
 applying \eqref{eq:ga011}   with $t_0 \n = \n T$ shows that
  $\wp_n (\o_2)  \n - \n  \wp_n (\o_1)   \le  \vth_\ell (\o_2)  \n - \n  \vth_\ell (\o_1)
  \le 2T \d^{-1}_\ell    \| \o_1   \n - \n   \o_2 \|_{0,T} = \k_n  \| \o_1   \n - \n   \o_2 \|_{0,T}$.
  Next,  suppose that the set  $ \big\{ t  \n \in \n  [a_n, T) \n : t  \n \ge \n  a_n
   \n + \n   \k_n  \| \o_1   \n - \n   \o_2 \|_{0,t} \big\}     $
     is not empty and contains $t_0$.
  Since
  $ \vth_\ell (\o_1)  \n = \n  \wp_n (\o_1)  \n \le \n  \wp_n (\o_2)   \n \le \n  \vth_\ell (\o_2)  $.
  we see that $ \fra_\ell \n = \n \vth_\ell (\o_1) $ and can deduce that
    $t_0    \n \ge \n  a_n  \n + \n   \k_n  \| \o_1   \n - \n   \o_2 \|_{0,t_0}
    \n = \n   \vth_\ell (\o_1)  \n + \n   \k_n  \| \o_1   \n - \n   \o_2 \|_{0,t_0}
     \n = \n  \fra_\ell  \n + \n  2T \d^{-1}_\ell  \| \o_1   \n - \n   \o_2 \|_{0,t_0}    $.
 Applying \eqref{eq:ga011} again  yields that
  $  \wp_n (\o_2)  \n - \n  \wp_n (\o_1) \le  \vth_\ell (\o_2)  \n - \n  \vth_\ell (\o_1)
  \le 2T \d^{-1}_\ell    \| \o_1   \n - \n   \o_2 \|_{0,t_0} = \k_n  \| \o_1   \n - \n   \o_2 \|_{0,t_0}  $. \qed

     \no {\bf Proof of Lemma \ref{lem_Y_nkl}:} Fix $n ,   k \n \in \n  \hN$.
    We define
     $  H_t    \n := \n   1 \n \land  \n  (  2^k ( t  \n - \n  \wp_n )  \n - \n  1 )^+  $ and
     $ \D_t   \n : = \n   U_t     \n - \n   L_t      $, $ t   \n \in \n  [0,T]   $.
     Let $(t_1,\o_1), (t_2,\o_2)  \n \in \n  [0,T]  \n \times \n  \O$.
          We set $ \bd^{1,2}  \n : = \n  \bd_\infty \big((t_1, \o_1), (t_2, \o_2)\big) $ and assume
     without loss of generality  that $t_1  \n \le \n  t_2$.

        Since \eqref{eq:ec017} shows that $ \big| H_{t_1} (\o_2) \n - \n   H_{t_2} (\o_2) \big|
      \n \le \n    \big| \big( 2^k ( t_1 \n - \n \wp_n (\o_2) ) \n - \n 1 )^+  \n - \n
      \big( 2^k ( t_2 \n - \n \wp_n (\o_2) ) \n - \n 1 \big)^+ \big| \le 2^k |t_1 \n - \n t_2|  $,
     \eqref{eq:aa211} implies that
   \bea
   \big|Y^{n,k}_{t_1} (\o_2) \n - \n  Y^{n,k}_{t_2} (\o_2) \big|
   & \tn \dn \le &  \tn  \dn  \big| L_{t_1} (\o_2) -  L_{t_2} (\o_2)\big|
   + \big| H_{t_1} (\o_2) -  H_{t_2} (\o_2) \big|
   | \D_{t_1} (\o_2)   |
   + H_{t_2} (\o_2)  | \D_{t_1} (\o_2) - \D_{t_2} (\o_2)  | \nonumber  \\
    & \tn \dn  \le & \tn \dn   \rho_0 \big(  \bd_\infty \big((t_1, \o_2), (t_2, \o_2)\big)    \big)
    \n + \n  2^{1+k}   M_0 |t_1 \n - \n t_2|
    \n + \n  2 \rho_0 \big(  \bd_\infty \big((t_1, \o_2), (t_2, \o_2)\big)    \big) . \qq
    \label{eq:ec114}
   \eea
   Since
     \beas
     \underset{r \in [t_1,t_2]}{\sup} \big|\o_2 (r) \n - \n  \o_2 (t_1) \big|
     & \tn \le  & \tn   |\o_1 (t_1)  \n - \n  \o_2 (t_1)|
       \n + \n  \underset{r \in [t_1,t_2]}{\sup} |\o_2 (r)  \n - \n  \o_1 (t_1) |
     \ls   2  \Big(    \|\o_1 \n - \n \o_2\|_{0,t_1}   \n  \vee \n  \underset{r \in [t_1,t_2]}{\sup}
      |\o_1 (t_1)  \n - \n  \o_2 (r)|   \Big) \nonumber \\
       & \tn  =   & \tn   2  \|\o_1 (\cd  \n \land \n  t_1)  \n - \n  \o_2(\cd  \n \land \n  t_2) \|_{0,T} ,
     \eeas
     one can deduce that
     $
 \bd_\infty \big((t_1, \o_2), (t_2, \o_2)\big) \n =  \n  |t_1 \n - \n t_2| \n + \n  \underset{r \in [t_1,t_2]}{\sup}
    \big|\o_2 (r) \n - \n  \o_2 (t_1) \big|  \n \le  \n  2 \big( |t_1 \n - \n t_2| \n + \n
    \|\o_1 (\cd \land t_1)  \n - \n  \o_2(\cd \land t_2) \|_{0,T} \big)  \n = \n   2 \bd^{1,2} $.
 Then it follows from \eqref{eq:ec114} that
    \bea
   \big|Y^{n,k}_{t_1} (\o_2) \n - \n  Y^{n,k}_{t_2} (\o_2) \big|
   \n \le \n  3  \rho_0 \big(  2 \bd^{1,2}   \big)
    \n + \n  2^{1+k}   M_0 \bd^{1,2}  .
    \label{eq:es271}
   \eea

 Since \eqref{eq:ec017}, Proposition \ref{prop_wp_n} (2) imply that
\bhe
\bea \label{eq:a119}
 \big|H_{t_1} (\o_1) \n - \n  H_{t_1} (\o_2) \big|    \n \le \n  2^k \k_n   \|\o_1  \n - \n  \o_2\|_{0,t_1} ,
\eea
\ehe
and since $ \|\o_1  \n - \n  \o_2\|_{0,t_1}
     \n \le \n     \|\o_1 (\cd  \n \land \n  t_1)  \n - \n  \o_2(\cd  \n \land \n  t_2) \|_{0,T}
     \n \le \n     \bd^{1,2}$,   we can further deduce that
   \beas
   \big|Y^{n,k}_{t_1} (\o_1) \n - \n  Y^{n,k}_{t_1} (\o_2) \big|
   & \tn \le &  \tn  \big| L_{t_1} (\o_1) -  L_{t_1} (\o_2)\big|
   + \big| H_{t_1} (\o_1) -  H_{t_1} (\o_2) \big|
   | \D_{t_1} (\o_1)   |
   + H_{t_1} (\o_2)  | \D_{t_1} (\o_1) - \D_{t_1} (\o_2)  | \\
    & \tn \le & \tn  \rho_0 \big( \|\o_1  \n - \n  \o_2\|_{0,t_1} \big)
    \n + \n  2^{1+k}   M_0 \k_n  \|\o_1  \n - \n  \o_2\|_{0,t_1}
    \n + \n  2 \rho_0 \big( \|\o_1  \n - \n  \o_2\|_{0,t_1} \big)
    \n  \le  \n   3 \rho_0 (2\bd^{1,2})  \n + \n  2^{1+k}  M_0   \k_n   \bd^{1,2}     ,
   \eeas
   which together with \eqref{eq:es271}   leads to that
   $   \big|Y^{n,k}_{t_1} (\o_1) \n - \n  Y^{n,k}_{t_2} (\o_2) \big|
   \le 6 \rho_0 (2\bd^{1,2}) + 2^{1+k}  M_0 (1 \n + \n \k_n)  \bd^{1,2}
   = \rho_{n,k} (\bd^{1,2}) $.   \qed

      \no {\bf Proof of \eqref{eq:eh137}:}
     Fix  $ (t,\o) \n \in \n  [0,T]  \n \times \n  \O  $.
 We will simply denote $ 2^{1-k}    $ by   $\d   $ and denote the term
 $  U \big( ( \wp_n (\o)  \n + \n  \d ) \land t  ,   \o    \big)
    \n - \n   U  (\wp_n (\o) \land t , \o   )  $ by $\D^U $.
  Let    $(\hP, \ga,\nu)  \n \in \n  \cP_t    \n \times \n \cT^t \ti \cT^t $ and define
  \beas
  \q   J_{\ga,\nu} (\wt{\o}) \n := \n  \b1_{\{\ga (\wt{\o}) > \wp_n (\o \otimes_t \wt{\o})\}}
    \Big( U \big( (\wp_n (\o  \n \otimes_t \n  \wt{\o})  \n + \n  \d  )
   \n \land  \n  (  \nu (\wt{\o}) \n \vee \n \wp_n (\o  \n \otimes_t \n  \wt{\o}) )  ,
     \o  \n \otimes_t \n  \wt{\o} \big)  \n - \n
  U \big(   \wp_n (\o  \n \otimes_t \n  \wt{\o}) ,
  \o  \n \otimes_t \n  \wt{\o} \big) \Big)  , \q   \fa \wt{\o} \in \O^t .
  \eeas

  \no  {\bf 1)} {\it We first  show by three cases that }
 \bea \label{eq:et011}
   \hE_\hP \big[ | J_{\ga,\nu} - \D^U | \big]   \le     \wh{\rho}_0 (\d)   .
 \eea

   \no (i) When $  \wp_n (\o) \n < \n  t  \n - \n  \d$, applying Lemma \ref{lem_shift_stopping_time}
  with $(t,s,\tau)  \n = \n  (0,t,\wp_n)$ yields that
  $t_n  \n : = \n  \wp_n (\o)  \n = \n  \wp_n ( \o  \n \otimes_t \n  \wt{\o} )  $,
  $\fa \wt{\o}  \n \in \n  \O^t$.
  Since  $U$ is an $\bF-$adapted process by (A1) and \eqref{eq:et014},
  one has $U_{t_n}  \n \in \n  \cF_{t_n}    \n \subset \n  \cF_t $ and
   $U_{t_n+\d}  \n \in \n  \cF_{t_n+\d}  \n \subset  \n   \cF_t $.
  Let $\wt{\o} \n \in \n  \O^t$.  Using \eqref{eq:bb421} with $(t,s,\eta)  \n = \n  (0,t,U_{t_n})$ and
  $(t,s,\eta)  \n = \n  (0,t,U_{t_n+\d})$     respectively shows that
  $  U  (t_n   ,  \o  \n \otimes_t \n  \wt{\o}  )   \n = \n  U  (t_n ,    \o  )    $ and
 $ U \big(t_n  \n + \n \d ,  \o  \n \otimes_t \n  \wt{\o} \big)
 \n  = \n  U \big(t_n  \n + \n \d ,    \o \big) $.
 As $   t_n  \n + \n  \d  \n < \n t  \n \le \n  \ga (\wt{\o}) \n \land \n \nu (\wt{\o})   $,
 one has
 \beas
 J_{\ga,\nu} (\wt{\o})
 & \tn \dn  = &  \tn \dn  \b1_{\{\ga (\wt{\o}) > t_n\}} \big( U \big( (t_n  \n + \n  \d  )
   \n \land  \n   ( \nu (\wt{\o}) \n \vee \n t_n )   ,   \o  \n \otimes_t \n  \wt{\o} \big)  \n - \n
  U  (   t_n   , \o  \n \otimes_t \n  \wt{\o}  ) \big)
  =   U \big(  t_n  \n + \n  \d   ,   \o    \big)  \n - \n   U  (t_n , \o   )
  = \D^U  .
 \eeas

   \no (ii) When $  t  \n - \n  \d  \n \le  \n   \wp_n (\o)  \n < \n  t$,
  we still have    $ t_n  \dn   = \n  \wp_n (\o)  \n = \n  \wp_n ( \o  \n \otimes_t \n  \wt{\o} ) $
  and $ U  (t_n   ,  \o  \n \otimes_t \n  \wt{\o}  )  \n  = \n  U  (t_n ,    \o  )$, $\fa \wt{\o} \in \O^t$.
  Set  $\nu_n   \n : = \n  (t_n  \n + \n  \d)   \n  \land  \n    \nu  \n \in \n  \cT^t$.
  For any $ \wt{\o} \n \in \n \O^t$, we see   from
  $ t_n \dn < \n t    \n \le \n  \ga (\wt{\o}) \n \land \n   \nu (\wt{\o})   $  that
  \beas
     J_{\ga,\nu} (\wt{\o}) \n - \n \D^U   \n =  \n    \b1_{\{\ga (\wt{\o}) > t_n \}}
    \Big(   U  \n \big( (t_n  \dn + \n  \d  )
   \n \land  \n   ( \nu (\wt{\o}) \n \vee \n t_n )  ,   \o  \n \otimes_t \n  \wt{\o} \big)  \n - \n
  U  \n    (    t_n   , \o  \n \otimes_t \n  \wt{\o}  )   \Big)  \n - \n     U  \n  (t , \o   )
   \n + \n   U   (t_n , \o   )
   \dn =  \n      U   \big(   \nu_n  (\wt{\o})   ,   \o  \n \otimes_t \n  \wt{\o} \big)  \n - \n
  U  \n  (t , \o   ) .
  \eeas
 Since    $ 
  t \n \le \n \nu_n (\wt{\o}) \n \le \n (t_n  \n + \n  \d) \n \land \n T
  \n \le \n (t  \n + \n  \d) \n \land \n T $, one can further deduce from \eqref{eq:aa211}   that
  \beas
  \big|  J_{\ga,\nu} (\wt{\o}) \n - \n \D^U  \big|
   & \tn  \dn    \le  & \tn  \dn        \rho_0 \Big( (\nu_n (\wt{\o})  \n - \n  t )
    \n + \n  \underset{r \in [0,T]}{\sup} \big|
   (\o  \n  \otimes_t \wt{\o})  \big( r  \n \land \n  \nu_n (\wt{\o}) \big)
    \n - \n   \o   (r  \n \land \n  t ) \big| \Big)
    \n \le \n  \rho_0 \Big( \d  \n + \n  \underset{r \in [t , (t  +  \d ) \land T]}{\sup}
    |   \wt{\o}  (r  )   | \Big) \\
     & \tn  \dn  =    & \tn  \dn
     \rho_0 \Big( \d  \n + \n  \underset{r \in [t, (t    +    \d) \land T]}{\sup}
   \big|    B^t_r (\wt{\o})    \n - \n  B^t_t (\wt{\o}) \big|  \Big)   .
   \eeas
   Taking expectation $\hE_\hP [~]$, we see from   \eqref{eq:ex015}   that
   $  \hE_\hP \big[ | J_{\ga,\nu} - \D^U | \big]    \le  \wh{\rho}_0 (\d)   $.

   \no (iii) When $  \wp_n (\o)  \n \ge \n  t$, we see that
  $\D^U \n = \n  U(t,\o)   \-     U(t,\o)  \n = \n  0 $.
  As Lemma  \ref{lem_shift_stopping_time} shows that
  $ \wp^{t,\o}_n \n \in \n \cT^t $,
 $ \z_n  \n : = \n  (\wp^{t,\o}_n  \n + \n  \d)      \land      (\nu      \vee      \wp^{t,\o}_n   )  $
 is also an  $\bF^t-$stopping time.
 Given   $\wt{\o}  \n \in \n  \O^t$, we set
   $s^1_n  \n := \n  \wp^{t,\o}_n ( \wt{\o})  \n \le \n  \z_n (\wt{\o})  \n : = \n  s^2_n   $.
 Since  $ s^2_n      \n \le \n
   \wp^{t,\o}_n ( \wt{\o})  \n  +  \n  \d   \n = \n  s^1_n   \n  +  \n  \d  $,
   applying   \eqref{eq:aa211} again yields that
 \beas
  \big| J_{\ga,\nu} (\wt{\o}) \n - \n \D^U \big| & \tn \dn = & \tn \dn  \big| J_{\ga,\nu} (\wt{\o})   \big| =
     \b1_{\{\ga (\wt{\o}) > \wp^{t,\o}_n (  \wt{\o})\}} \Big| U \big( (\wp^{t,\o}_n (   \wt{\o})  \n + \n  \d  )
   \n \land  \n   ( \nu (\wt{\o}) \n \vee \n  \wp^{t,\o}_n (   \wt{\o}) )   ,
     \o  \n \otimes_t \n  \wt{\o} \big)  \n - \n
  U \big(   \wp^{t,\o}_n (   \wt{\o}) , \o  \n \otimes_t \n  \wt{\o} \big)  \Big| \\
    & \tn \dn  \le   & \tn \dn     \big| U \big( s^2_n ,   \o  \n \otimes_t \n  \wt{\o} \big)  \n - \n
  U \big( s^1_n , \o  \n \otimes_t \n  \wt{\o} \big) \big|
    \n  \le  \n   \rho_0 \Big( (s^2_n  \n - \n  s^1_n)  \n + \n  \underset{r \in [0,T]}{\sup} \big|
  (\o  \n \otimes_t \n  \wt{\o}) (r  \n \land \n  s^2_n)
   \n - \n  (\o  \n \otimes_t \n  \wt{\o}) (r  \n \land \n  s^1_n) \big| \Big) \\
  & \tn \dn  =  & \tn \dn  \rho_0 \Big( (s^2_n  \n - \n  s^1_n)  \n + \n  \underset{r \in [s^1_n,s^2_n]}{\sup} \big|
    \wt{\o}  (r  )  \n - \n    \wt{\o}  (  s^1_n) \big| \Big)
      \n \le  \n   \rho_0 \Big( \d
       \n + \n  \underset{r \in [\wp^{t,\o}_n (\wt{\o}) , (\wp^{t,\o}_n (\wt{\o}) + \d ) \land T
    ]}{\sup} \big|   B^t_r ( \wt{\o}   )  \n - \n   B^t \big(\wp^{t,\o}_n (\wt{\o}) , \wt{\o} \big) \big| \Big) .
 \eeas
 Taking expectation $\hE_\hP [~]$ and using \eqref{eq:ex015} yield that
 $  \hE_\hP \big[ | J_{\ga,\nu} - \D^U | \big]  \le
   \wh{\rho}_0 (\d)  $.
 Hence, we proved \eqref{eq:et011}.

 \no {\bf 2)} {\it Next, we use \eqref{eq:et011} to verify \eqref{eq:eh137}.}

  \no {\bf 2a)}  For any $(t',\o')  \n \in \n  [0,T]  \n \times \n  \O$,
  since \eqref{eq:eh111} and (A2) imply     that
  $  L (t',\o') \n \le \n   Y^{n,k} (t',\o')    \n \le \n  U (t',\o')    $,
  \bea \label{eq:et021}
  \wh{Y}^{n,k} (t',\o')
  & \tn \dn =  & \tn  \dn   \b1_{\{t' \le \wp_n (\o')  \}} L (t',\o') \n + \n  \b1_{\{t' > \wp_n (\o')\}} Y^{n,k}
   \big(  \wp^{n,k} (\o')  \n \land \n  t', \o' \big) \nonumber \\
  & \tn \dn  \le   & \tn \dn   \b1_{\{t' \le \wp_n (\o')\}} L (t',\o') \n + \n  \b1_{\{t' > \wp_n (\o')\}}
   U \big((\wp_n (\o')  \n + \n  \d)  \n \land \n    t' , \o' \big) .
 \eea
Given $(\hP, \ga)  \n \in \n  \cP_t    \n \times \n \cT^t $ and $ \wt{\o}  \n \in \n  \O^t$, taking $(t',\o') \n = \n
\big( \ga (\wt{\o}), \o  \n \otimes_t \n  \wt{\o}  \big) $ in    \eqref{eq:et021} yields that
 \beas
 && \hspace{-0.8cm} \big(\wh{Y}^{n,k}\big)^{t,\o}_\ga  (\wt{\o}) \n -  \n    (\sY^n)^{t,\o}_\ga  (\wt{\o})
 \n = \n  \wh{Y}^{n,k} \big( \ga (\wt{\o}), \o  \n \otimes_t \n  \wt{\o}  \big)
  \n - \n  \sY^n \big( \ga (\wt{\o}), \o  \n \otimes_t \n  \wt{\o}  \big) \nonumber \\
 & &    \le  \n    \b1_{\{\ga (\wt{\o}) > \wp_n (\o \otimes_t \wt{\o})\}}
  \Big( U \big( (\wp_n (\o  \n \otimes_t \n  \wt{\o})  \n + \n  \d  )
   \n \land  \n   ( \ga (\wt{\o}) \n \vee \n  \wp_n (\o  \n \otimes_t \n  \wt{\o}) )   ,   \o  \n \otimes_t \n  \wt{\o} \big)  \n - \n
  U \big(  \wp_n (\o  \n \otimes_t \n  \wt{\o}) , \o  \n \otimes_t \n  \wt{\o} \big) \Big)
  \n = \n J_{\ga,\ga} (\wt{\o}) .
 \eeas
  It then follows from \eqref{eq:et011} that
  $  \hE_\hP \Big[ \big(\wh{Y}^{n,k}\big)^{t,\o}_\ga  \Big]
  \n \le \n  \hE_\hP \big[   (\sY^n)^{t,\o}_\ga  \n + \n J_{\ga,\ga} \big]
     \n \le \n  \sZ^n_t (\o)  \n + \n \D^U   \n + \n     \wh{\rho}_0 (\d)$.
  Taking supremum over $(\hP, \ga)  \n \in \n  \cP_t    \n \times \n \cT^t $
 on the left-hand-side  leads to that
 $ Z^{n,k}_t (\o)  \n \le \n  \sZ^n_t (\o)  \n + \n \D^U   \n + \n     \wh{\rho}_0 (\d)   $.

  \no {\bf 2b)} To show the left-hand-side of  \eqref{eq:eh137},
  we let $(\hP, \ga)  \n \in \n  \cP_t    \n \times \n \cT^t $ and set
     $\wt{\ga}  \n : = \n  \big( \ga  \n + \n  \d \big)   \n \land \n  T  \n \in \n  \cT^t $.
     Also, let $(t',\o') \n \in \n [0,T] \n \times \n \O$, one has
     \bea \label{eq:et041}
      \sY^n (t',\o') \n \le  \n    \b1_{\{ t' \le \wp_n (\o') \}} U (t',\o')   \n + \n
      \b1_{\{ t' > \wp_n (\o')    \}} U \big( \wp_n (\o'), \o' \big) \n = \n
      U  \big( \wp_n (\o')  \n  \land  \n  t' , \o' \big) .
     \eea
 If    $t' \n \le \n  T  \n - \n  \d $, since
   \beas
 && \hspace{-1cm} \wh{Y}^{n,k} ( t' \n + \n  \d , \o')  \n = \n    \b1_{\{  t'\le \wp_n (\o')   - 2^{-k} \}}
    L(t'+  \d, \o') +    \b1_{\{ t'\ge \wp_n (\o')      \}}
           U \big( (\wp_n (\o')  +  \d)  \n \land \n  T , \o' \big) \\
    &    & \q   +    \b1_{\{    \wp_n (\o')     - 2^{-k}   < t'< \wp_n (\o')        \}}
    \Big\{ \big[ 1 \n - \n   2^k  ( t'\n + \n 2^{-k} \n - \n  \wp_n (\o')        )  \big]
   L(t' \n + \n   \d, \o')  \n + \n  2^k ( t'\n + \n 2^{-k} \n - \n  \wp_n (\o') )
   U(t' \n + \n   \d, \o')    \Big\}  \\
         & & \ge \n \b1_{\{  t'\le \wp_n (\o')    \}}   L(t' \n + \n   \d, \o')
        \n  +  \n   \b1_{\{ t'> \wp_n (\o')      \}}
        U \big( (\wp_n (\o')  \n  +  \n  \d)  \n \land \n  T , \o' \big)    ,
 \eeas
 we can obtain that
 \bea \label{eq:eh131}
 \sY^n(t' , \o') \n - \n  \wh{Y}^{n,k}(t' \n + \n  \d , \o')
 & \tn \dn  \le  & \tn  \dn     \b1_{\{  t'\le \wp_n (\o')    \}}
 \big( L(t' , \o') \n - \n   L(t' \n + \n  \d , \o')    \big) \nonumber \\
   & \tn  \dn  & \tn  \dn   +    \b1_{\{ t'> \wp_n (\o')    \}}  \big( U(\wp_n (\o') , \o')
    \n - \n U \big( (\wp_n (\o')  \n  +  \n  \d)  \n \land \n  T , \o' \big)  \big) .
 \eea
 Also, \eqref{eq:et031} and (A2) imply   that
 \bea \label{eq:et043}
   \wh{Y}^{n,k}   (T,\o')  \n = \n \b1_{\{   \wp_n (\o')  > T -      \d   \}} U (T,\o')   \n + \n
  \b1_{\{   \wp_n (\o')  \le T -   \d   \}} U \big( ( \wp_n (\o')  \n + \n \d )  \n \land \n  T , \o' \big)
    \n = \n  U \big( ( \wp_n (\o')  \n + \n \d )  \n \land \n  T , \o' \big)    . ~ \;
  \eea

     Let   $ \wt{\o}  \n \in \n  \{ \ga    \n > \n  T \n - \n \d \}  $, so $\wt{\ga} (\wt{\o})  \n = \n  T $.
     Taking $(t',\o') \n = \n \big( \ga (\wt{\o}), \o  \n \otimes_t \n  \wt{\o}  \big) $ in
     \eqref{eq:et041}, \eqref{eq:et043} and using \eqref{eq:aa211} yield that
 \bea
 &&  \hspace{-1.2cm}  (\sY^n)^{t,\o}_\ga  (  \wt{\o}  )
 \n - \dn  \big(\wh{Y}^{n,k}\big)^{t,\o}_{\wt{\ga}} (  \wt{\o}  )
  \n = \n       \sY^n  \big(\ga (\wt{\o}) ,  \o  \n \otimes_t \n   \wt{\o} \big)
  \n - \n   \wh{Y}^{n,k}   ( T ,  \o  \n \otimes_t \n   \wt{\o}  )
    \n  \le  \n  U \n  \big(\wp_n (\o \n \otimes_t  \n  \wt{\o})  \n \land \n  \ga (\wt{\o}) ,
  \o  \n \otimes_t  \n  \wt{\o} \big)
  \n - \n  U  \n  \big( (\wp_n(\o  \n \otimes_t \n   \wt{\o}) \n + \n \d)  \n \land \n  T ,
   \o  \n \otimes_t  \n  \wt{\o} \big)  \nonumber  \\
 &&  \hspace{-0.7cm}   = \n \b1_{\{ \ga (\wt{\o}) \le \wp_n  \n  (\o   \otimes_t     \wt{\o})    \}}
 \Big( \n  U  \n   \big(\ga ( \wt{\o})   ,
  \o  \n \otimes_t  \n  \wt{\o} \big)
  \dn - \n  U \n  \big(  T ,   \o  \n \otimes_t  \n  \wt{\o} \big) \n  \Big)
  \dn  + \n  \b1_{\{  \ga (\wt{\o}) > \wp_n  \n  (\o   \otimes_t     \wt{\o}) \}}
 \Big( \n  U  \n   \big(\wp_n (\o \n \otimes_t  \n  \wt{\o})   ,
  \o  \n \otimes_t  \n  \wt{\o} \big)
  \dn - \n  U \n  \big( (\wp_n(\o  \n \otimes_t \n   \wt{\o}) \n + \n \d)  \n \land \n  T ,
   \o  \n \otimes_t  \n  \wt{\o} \big) \n  \Big)  \nonumber   \\
 &&  \hspace{-0.7cm}   \le     \rho_0 \Big( \big(T \n - \n  \ga (\wt{\o})\big)  \n + \n  \underset{r \in [\ga (\wt{\o}),T]}{\sup}
 \big|   \wt{\o}  ( r  )    \n - \n    \wt{\o}  (   \ga (\wt{\o})  )  \big| \Big)
   \n - \n  J_{\ga,T} (\wt{\o})
   \n \le  \n      \rho_0 \Big( \d
    \n + \n  \underset{r \in [\ga (\wt{\o}),(\ga (\wt{\o})+\d) \land T]}{\sup}
 \big|  B^t_r  (   \wt{\o}   )      \n - \n   B^t_\ga  (   \wt{\o}  )  \big| \Big)
   \n - \n  J_{\ga,T} (\wt{\o})   .   \label{eq:eh133}
 \eea
 On the other hand, let  $\wt{\o}  \n \in \n  \{ \ga    \n \le \n  T \n - \n \d \}$.
 applying \eqref{eq:eh131}
 with $(t',\o') \n =  \n  \big(\ga  (\wt{\o}) , \o  \n \otimes_t  \n  \wt{\o} \big)  $
 and using \eqref{eq:aa211} yield  that
 \beas
  &&  \hspace{-1cm}  (\sY^n)^{t,\o}_\ga  (  \wt{\o}  )
 - \big(\wh{Y}^{n,k}\big)^{t,\o}_{\wt{\ga}} (  \wt{\o}  )
 =    \sY^n  \big(\ga  (\wt{\o}) , \o \otimes_t  \wt{\o} \big)
 -  \wh{Y}^{n,k}  \big(\ga  (\wt{\o}) \n + \n  \d ,  \o \otimes_t   \wt{\o} \big) \\
 && \hspace{-0.3cm} \le    \b1_{\{  \ga(\wt{\o})\le \wp_n (\o \otimes_t \wt{\o})    \}}
 \Big( L(\ga(\wt{\o}) , \o  \n \otimes_t \n  \wt{\o})
 \n - \n   L(\ga(\wt{\o}) \n + \n  \d , \o  \n \otimes_t \n  \wt{\o})    \Big) \\
 && \hspace{-0.3cm}   \q   +    \b1_{\{ \ga(\wt{\o})> \wp_n (\o \otimes_t \wt{\o})    \}}
     \Big( U(\wp_n (\o  \n \otimes_t \n  \wt{\o}) , \o  \n \otimes_t \n  \wt{\o})
    \n - \n U \big( (\wp_n (\o  \n \otimes_t \n  \wt{\o})  \n  +  \n  \d)  \n \land \n  T ,
    \o  \n \otimes_t \n  \wt{\o} \big)  \Big) \\
 && \hspace{-0.3cm}  \le \n  \rho_0 \Big( \d  \n + \n
  \underset{r \in [\ga  (\wt{\o}), (\ga  (\wt{\o})   +    \d) \land T ]}{\sup}
 \big|    \wt{\o}  (r   )   \n - \n     \wt{\o}  ( \ga  (\wt{\o}) )  \big| \Big)  \n - \n  J_{\ga,T} (\wt{\o})
    \n = \n  \rho_0 \Big( \d
     \n + \n  \underset{r \in [\ga  (\wt{\o}), (\ga  (\wt{\o})   +    \d) \land T ]}{\sup}
 \big| B^t_r (   \wt{\o})   \n   - \n  B^t_\ga ( \wt{\o})   \big| \Big)  \n - \n  J_{\ga,T} (\wt{\o})  .
 \eeas
 Combining this with \eqref{eq:eh133},
 we see from \eqref{eq:et011} and \eqref{eq:ex015}  that
 \beas
  \hE_\hP \big[  ( \sY^n)^{t,\o}_{\ga}  \big]
    \n \le \n  \hE_\hP \Big[  \big(\wh{Y}^{n,k}\big)^{t,\o}_{\wt{\ga}} \n - \n J_{\ga,T}  \Big]
     \n + \n   \wh{\rho}_0(\d)
  \n \le \n   Z^{n,k}_t(\o) \n - \n \D^U   \n + \n    2 \wh{\rho}_0 (\d)   .
 \eeas
   Then taking supremum over $(\hP, \ga)  \n \in \n  \cP_t    \n \times \n \cT^t $
 on the left-hand-side  leads to that
 $    \sZ^n_t(\o)   \n \le \n   Z^{n,k}_t(\o)  \n - \n \D^U    \n + \n   2 \wh{\rho}_0 (\d)   $. \qed

 \no {\bf Proof of \eqref{eq:eh137b}:}
  Fix  $ (t,\o) \n \in \n  [0,T]  \n \times \n  \O  $.
  We will simply  denote $ \frac{2T}{n+3} $ by $ \d $ and
  denote     the term $ U \big(  \wp_{n+1} (\o)   \n \land \n  t , \o \big)
  \n - \n U \big(  \wp_n (\o)  \n \land \n  t , \o \big) $ by $\wt{\D}^U$.
  Let  $(\hP, \ga, \nu)  \n \in \n  \cP_t    \n \times \n \cT^t \ti \cT^t $ and   define
  \beas
  \q   J_{\ga,\nu} (\wt{\o}) \n := \n  \b1_{\{\ga (\wt{\o}) > \wp_n (\o \otimes_t \wt{\o})\}}
    \Big( U \big(  \wp_{n+1} (\o  \n \otimes_t \n  \wt{\o})
   \n \land  \n  (  \nu (\wt{\o}) \n \vee \n \wp_n (\o  \n \otimes_t \n  \wt{\o}) )  ,
     \o  \n \otimes_t \n  \wt{\o} \big)  \n - \n
  U \big(   \wp_n (\o  \n \otimes_t \n  \wt{\o}) ,
  \o  \n \otimes_t \n  \wt{\o} \big) \Big)  , \q   \fa \wt{\o} \in \O^t .
  \eeas

   In light of Proposition \ref{prop_wp_n} (1),
  one can  deduce \eqref{eq:et011} again
  by three cases:  $  \wp_{n+1} (\o) \n < \n  t   $,
  $  \wp_n (\o) \n < \n  t   \n \le  \n   \wp_{n+1} (\o)  $ and $  \wp_n (\o) \n  \ge \n  t  $.

  \if{0}

   \no (i) When $  \wp_{n+1} (\o) \n < \n  t   $,
  applying Lemma \ref{lem_shift_stopping_time}
  with $(t,s,\tau)  \n = \n  (0,t,\wp_n)  $   and $(t,s,\tau)  \n = \n  (0,t,\wp_{n+1})  $ respectively yields that
  $t_n  \n : = \n  \wp_n (\o)   \n = \n  \wp_n ( \o  \n \otimes_t \n  \wt{\o} )  $
  and $t_{n+1}  \n : = \n  \wp_{n+1} (\o)    \n = \n  \wp_{n+1} ( \o  \n \otimes_t \n  \wt{\o} )   $,
  $\fa \wt{\o}  \n \in \n  \O^t$.
  Since  $U$ is an $\bF-$adapted process by (A1) and \eqref{eq:et014},
  one has   $U_{t_n}  \n \in \n  \cF_{t_n}    \n \subset \n  \cF_t $
  and  $U_{t_{n+1}}  \n \in \n  \cF_{t_{n+1}}  \n \subset  \n   \cF_t $.
  Let $\wt{\o} \n \in \n  \O^t$.  Using \eqref{eq:bb421} with
    $(t,s,\eta)  \n = \n  (0,t,U_{t_n})$ and $(t,s,\eta)  \n = \n  (0,t,U_{t_{n+1}})$   respectively shows that
 $  U  (t_n   ,  \o  \n \otimes_t \n  \wt{\o}  )   \n = \n  U  (t_n ,    \o  )    $ and
  $ U \big(t_{n+1} ,  \o  \n \otimes_t \n  \wt{\o} \big)  \n  = \n  U \big(t_{n+1} ,    \o \big) $.
 As  $ t_n \n \le \n   t_{n+1} \n < \n t  \n \le \n  \ga (\wt{\o}) \n \land \n   \nu (\wt{\o})  $,
 we can deduce from \eqref{eq:aa211} that
 \beas
 J_{\ga,\nu} (\wt{\o})
 & \tn \dn  = &  \tn \dn  \b1_{\{\ga (\wt{\o}) > t_n \}}  \Big(  U \big(  t_{n+1}
   \n \land  \n   ( \nu (\wt{\o}) \n \vee \n t_n )   ,   \o  \n \otimes_t \n  \wt{\o} \big)  \n - \n
  U \big( t_n , \o  \n \otimes_t \n  \wt{\o} \big) \Big)
  =   U \big(  t_{n+1}    ,   \o    \big)  \n - \n   U  (t_n  , \o   ) = \wt{\D}^U    .
 \eeas

   \no (ii) When $  \wp_n (\o) \n < \n  t   \n \le  \n   \wp_{n+1} (\o)  $,
  we still have    $ t_n  \n   = \n  \wp_n (\o)
   \n = \n  \wp_n ( \o  \n \otimes_t \n  \wt{\o} )  $
  and $ U  (t_n   ,  \o  \n \otimes_t \n  \wt{\o}  )  \n  = \n  U  (t_n ,    \o  )$, $\fa \wt{\o} \in \O^t$.
  Lemma \ref{lem_shift_stopping_time} implies that
    $\ol{\nu}_n   \n : = \n   \wp^{t,\o}_{n+1}     \n  \land  \n    \nu  $ is an $\bF^t-$stopping time.
  For any $ \wt{\o} \n \in \n \O^t$,
  we see from  $ t_n \n < \n t    \n \le \n  \ga (\wt{\o}) \n \land \n   \nu (\wt{\o})   $ that
  \beas
  \hspace{-3mm}
  J_{\ga,\nu} (\wt{\o}) \n - \n \wt{\D}^U  \n  = \n  \b1_{\{\ga (\wt{\o}) > t_n \}}
      \Big( \n  U \n  \big(  \wp_{n+1} (\o  \n \otimes_t \n  \wt{\o})
   \n \land  \n   ( \nu (\wt{\o}) \n \vee \n t_n )  ,   \o  \n \otimes_t \n  \wt{\o} \big)  \n - \n
  U \n  \big(    t_n   , \o  \n \otimes_t \n  \wt{\o} \big) \n  \Big)  \dn - \n  U(t,\o)  \n + \n  U (t_n,\o)
   \n =  \n     U \n  \big(   \ol{\nu}_n  (\wt{\o})   ,   \o  \n \otimes_t \n  \wt{\o} \big)  \n - \n
  U  \n  (t  , \o   ) .
  \eeas
  As $ \wp_{n+1} (\o \n \otimes_t \n  \wt{\o}) \n \le \n \wp_n (\o \n \otimes_t \n  \wt{\o}) \n + \n \d
   \n = \n  t_n  \n + \n \d  $ by Proposition \ref{prop_wp_n} (1), one can deduce that
   $ 
    t \n \le \n \ol{\nu}_n (\wt{\o}) \n \le \n   \wp^{t,\o}_{n+1} (\wt{\o})
  \n \le \n (t_n  \n + \n  \d) \n \land \n T
  \n \le \n (t   \n + \n  \d) \n \land \n T   $. It then follows from \eqref{eq:aa211}  that
    \beas
    \big|  J_{\ga,\nu} (\wt{\o}) \n - \n \wt{\D}^U  \big|
    & \tn    \le  & \tn        \rho_0 \Big( (\ol{\nu}_n (\wt{\o})  \n - \n  t  )
    \n + \n  \underset{r \in [0,T]}{\sup} \big|
   (\o  \n  \otimes_t \wt{\o})  \big( r  \n \land \n  \ol{\nu}_n (\wt{\o}) \big)
    \n - \n   \o   (r  \n \land \n  t ) \big| \Big)
   \n  \le   \n       \rho_0 \Big( \d  \n + \n  \underset{r \in [t , (t  +  \d ) \land T]}{\sup}
    |    \wt{\o}   (r  )    | \Big)  \\
    & \tn    =  & \tn   \rho_0 \Big( \d  \n + \n  \underset{r \in [t, (t    +    \d) \land T]}{\sup}
   \big|    B^t_r (\wt{\o})    \n - \n  B^t_t (\wt{\o}) \big|  \Big)   .
   \eeas
   Taking expectation $\hE_\hP [~]$, we see from   \eqref{eq:ex015}   that
   $  \hE_\hP \big[ |J_{\ga,\nu}-\wt{\D}^U| \big]    \le
     \wh{\rho}_0 (\d)   $.

   \no (iii) When $  \wp_n (\o) \n  \ge \n  t  $, we see that
  $\wt{\D}^U  \n = \n  U(t,\o )  \n - \n  U(t,\o )   \n = \n  0  $.
  As Lemma  \ref{lem_shift_stopping_time} shows that
  $ \wp^{t,\o}_n , \wp^{t,\o}_{n+1} \n \in \n \cT^t $,
 $ \z_n  \n : = \n   \wp^{t,\o}_{n+1}   \n  \land  \n    ( \nu     \n  \vee  \n    \wp^{t,\o}_n )  $
 is also an  $\bF^t-$stopping time.
 Given   $\wt{\o}  \n \in \n  \O^t$, we set
   $s^1_n  \n := \n  \wp^{t,\o}_n ( \wt{\o})  \n \le \n  \z_n (\wt{\o})  \n : = \n  s^2_n   $.
 Since Proposition \ref{prop_wp_n} (1) implies that  $ s^2_n  \n \le  \n  \wp^{t,\o}_{n+1} ( \wt{\o})  \n  \le  \n
 \wp^{t,\o}_n ( \wt{\o})  \n + \n  \d\n = \n  s^1_n  \n + \n  \d$,
   applying   \eqref{eq:aa211} again yields that
 \beas
   \big| J_{\ga,\nu} (\wt{\o}) \n - \n \wt{\D}^U \big| & \tn  =  & \tn
   \big| J_{\ga,\nu} (\wt{\o}) \big| \n  = \n
    \b1_{\{\ga (\wt{\o}) > \wp^{t,\o}_n (  \wt{\o})\}} \Big| U \big( \wp^{t,\o}_{n+1} (   \wt{\o})
   \n \land  \n   ( \nu (\wt{\o})  \n \vee \n  \wp^{t,\o}_n (   \wt{\o}) )   ,
     \o  \n \otimes_t \n  \wt{\o} \big)  \n - \n
  U \big( \wp^{t,\o}_n (   \wt{\o}) , \o  \n \otimes_t \n  \wt{\o} \big) \Big| \\
     & \tn  \le   & \tn     \big| U \big( s^2_n ,   \o  \n \otimes_t \n  \wt{\o} \big)  \n - \n
  U \big( s^1_n , \o  \n \otimes_t \n  \wt{\o} \big) \big|
    \n  \le  \n   \rho_0 \Big( (s^2_n  \n - \n  s^1_n)  \n + \n  \underset{r \in [0,T]}{\sup} \big|
  (\o  \n \otimes_t \n  \wt{\o}) (r  \n \land \n  s^2_n)
   \n - \n  (\o  \n \otimes_t \n  \wt{\o}) (r  \n \land \n  s^1_n) \big| \Big) \\
  & \tn   =  & \tn  \rho_0 \Big( (s^2_n  \n - \n  s^1_n)  \n + \n  \underset{r \in [s^1_n,s^2_n]}{\sup} \big|
    \wt{\o}  (r  )  \n - \n    \wt{\o}  (  s^1_n) \big| \Big)
      \n \le  \n   \rho_0 \Big(\d
       \n + \n  \underset{r \in [\wp^{t,\o}_n  (\wt{\o}) , (\wp^{t,\o}_n  (\wt{\o}) + \d) \land T
    ]}{\sup} \big|   B^t_r ( \wt{\o}   )  \n - \n   B^t \big(\wp^{t,\o}_n  (\wt{\o}) , \wt{\o} \big) \big| \Big) .
 \eeas
 Taking expectation $\hE_\hP [~]$ and using \eqref{eq:ex015} yield that
 $  \hE_\hP \big[ |J_{\ga,\nu}-\wt{\D}^U| \big]  \le   \wh{\rho}_0 (\d) $.
 Hence, we proved \eqref{eq:et011} again.

 \fi

 \no {\bf 1)} {\it Let us show the right-hand-side of \eqref{eq:eh137b} first.}

    For any $(t',\o')  \n \in \n  [0,T]  \n \times \n  \O$,
 since an analogy to \eqref{eq:et041} shows that
  $  \sY^{n+1} (t',\o')  \le    U \big( \wp_{n+1} (\o')    \n \land \n    t' , \o' \big) $, we have
  \bea \label{eq:et121}
  \sY^{n+1} (t',\o')
  & \tn \dn =  & \tn  \dn   \b1_{\{t' \le \wp_n (\o')  \}} L (t',\o') \n + \n  \b1_{\{t' > \wp_n (\o')\}} \sY^{n+1}
   \big(  t', \o' \big) \nonumber \\
  & \tn \dn  \le   & \tn \dn   \b1_{\{t' \le \wp_n (\o')\}} L (t',\o') \n + \n  \b1_{\{t' > \wp_n (\o')\}}
   U \big( \wp_{n+1} (\o')   \n \land \n    t' , \o' \big) .
 \eea
Given $(\hP, \ga)  \n \in \n  \cP_t    \n \times \n \cT^t $ and $ \wt{\o}  \n \in \n  \O^t$, taking $(t',\o') \n = \n
\big( \ga (\wt{\o}), \o  \n \otimes_t \n  \wt{\o}  \big) $ in    \eqref{eq:et121} yields that
 \beas
 && \hspace{-0.8cm} \big(\sY^{n+1}\big)^{t,\o}_\ga  (\wt{\o}) \n -  \n    (\sY^n)^{t,\o}_\ga  (\wt{\o})
 \n = \n  \sY^{n+1} \big( \ga (\wt{\o}), \o  \n \otimes_t \n  \wt{\o}  \big)
  \n - \n  \sY^n \big( \ga (\wt{\o}), \o  \n \otimes_t \n  \wt{\o}  \big) \nonumber \\
 & &    \le  \n    \b1_{\{\ga (\wt{\o}) > \wp_n (\o \otimes_t \wt{\o})\}}
  \Big( U \big(  \wp_{n+1} (\o  \n \otimes_t \n  \wt{\o})
   \n \land  \n   ( \ga (\wt{\o})
   \n  \vee \n  \wp_n (\o  \n \otimes_t \n  \wt{\o}) )   ,   \o  \n \otimes_t \n  \wt{\o} \big)  \n - \n
  U \big(  \wp_n (\o  \n \otimes_t \n  \wt{\o}) , \o  \n \otimes_t \n  \wt{\o} \big) \Big)
  \n = \n J_{\ga,\ga} (\wt{\o}) .
 \eeas
  It then follows from \eqref{eq:et011} that
  $  \hE_\hP \big[ \big(\sY^{n+1}\big)^{t,\o}_\ga  \big]
  \n \le \n  \hE_\hP \big[   (\sY^n)^{t,\o}_\ga  \n + \n J_{\ga,\ga} \big]
     \n \le \n  \sZ^n_t (\o)  \n + \n \wt{\D}^U   \n + \n     \wh{\rho}_0 (\d)$.
  Taking supremum over $(\hP, \ga)  \n \in \n  \cP_t    \n \times \n \cT^t $
 on the left-hand-side  leads to that
 $ \sZ^{n+1}_t (\o)  \n \le \n  \sZ^n_t (\o)  \n + \n \wt{\D}^U   \n + \n     \wh{\rho}_0 (\d)   $.

 \no {\bf 2)} To show the left hand side of \eqref{eq:eh137b}, we
    let $(\hP, \ga) \n  \in  \n   \cP_t    \n \times \n \cT^t $ and set
     $\wt{\ga}   \n  : =  \n   \big( \ga   \n  +  \n   \d \big)    \n   \land   \n   T   \n  \in  \n   \cT^t $.
     We also let $(t',\o')  \n  \in \n   [0,T]  \n  \times  \n  \O$.
 If    $t'  \n  \le  \n   T    \n - \n    \d $,
 since $ \wp_{n+1} (\o')   \n  \le  \n   \wp_n  (\o')   \n  +  \n   \d  $ by
  Proposition \ref{prop_wp_n} (1), one can deduce that
  \beas
 \sY^{n+1} ( t' \dn + \n  \d , \o')
 & \tn  =   & \tn      \b1_{\{  t' + \d \le \wp_{n+1} (\o')   \}}   L(t' \dn + \n   \d, \o')
  \n  +  \n    \b1_{\{ t'+\d> \wp_{n+1} (\o')   \}} U \big( \wp_{n+1} (\o') , \o' \big)   \\
   & \tn  \ge  & \tn     \b1_{\{  t' + \d \le \wp_{n+1} (\o')   \}}   L(t' \dn + \n   \d, \o')
  \n + \n  \b1_{\{ \wp_{n+1} (\o') - \d  < t'   \le  \wp_n  (\o')    \}}
  L \big( \wp_{n+1} (\o')  , \o' \big)
    \n  + \n  \b1_{\{ t'  > \wp_n  (\o')   \}} U \big( \wp_{n+1} (\o')  , \o' \big)    ,
\eeas
 and thus that
 \bea
    \hspace{-5mm} \sY^n  (t',\o')  \n - \n   \sY^{n+1} (t' \dn + \n \d , \o')
    & \tn \dn  \le  & \tn \dn   \b1_{\{  t' + \d \le \wp_{n+1} (\o')   \}}
    \big( L (t',\o')   \n - \dn    L(t' \dn + \n   \d, \o')  \big)
     \n + \n  \b1_{\{ \wp_{n+1} (\o') - \d  < t'   \le  \wp_n  (\o')    \}}
      \Big( \n  L (t',\o')   \n  -  \dn   L \big( \wp_{n+1} (\o')  \n \vee \n   t'  , \o' \big) \n  \Big)  \nonumber  \\
  & \tn \dn  & \tn \dn
   \n +   \b1_{\{ t'  > \wp_n  (\o')   \}}
      \Big( U \n \big( \wp_n (\o')  , \o' \big)
      \dn - \n U \n  \big( \wp_{n+1} (\o') , \o' \big) \Big)   . \qq  \label{eq:eh131b}
 \eea
 Also,   (A2) implies   that
 \bea
   \sY^{n+1}   (T,\o') & \tn =  & \tn  \b1_{\{ T =  \wp_{n+1} (\o')    \}} L (T,\o')   \n + \n
  \b1_{\{ T >  \wp_{n+1} (\o')  \}} U \big( \wp_{n+1} (\o'), \o' \big) \nonumber  \\
    & \tn  =  & \tn  \b1_{\{ T =  \wp_{n+1} (\o')    \}} U (T,\o')   \n + \n
  \b1_{\{ T >  \wp_{n+1} (\o')  \}} U \big( \wp_{n+1} (\o'), \o' \big)
    \n = \n  U \big( \wp_{n+1} (\o')   , \o' \big)    . ~ \; \label{eq:et043b}
  \eea

     Let   $ \wt{\o}  \n \in \n  \{ \ga    \n > \n  T \n - \n \d \}  $, so $\wt{\ga} (\wt{\o})  \n = \n  T $.
     Taking $(t',\o') \n = \n \big( \ga (\wt{\o}), \o  \n \otimes_t \n  \wt{\o}  \big) $ in
     \eqref{eq:et041} and  \eqref{eq:et043b}    yields that
 \bea
   (\sY^n)^{t,\o}_\ga  (  \wt{\o}  )
 \n - \dn  \big(\sY^{n+1}\big)^{t,\o}_{\wt{\ga}} (  \wt{\o}  )
 & \tn \dn =  & \tn \dn        \sY^n  \big(\ga (\wt{\o}) ,  \o  \n \otimes_t \n   \wt{\o} \big)
  \n - \n   \sY^{n+1}   ( T ,  \o  \n \otimes_t \n   \wt{\o}  )
    \n  \le  \n  U \n  \big(\wp_n (\o \n \otimes_t  \n  \wt{\o})  \n \land \n  \ga (\wt{\o}) ,
  \o  \n \otimes_t  \n  \wt{\o} \big)
  \n - \n  U  \n  \big(  \wp_{n+1} (\o  \n \otimes_t \n   \wt{\o})  ,
   \o  \n \otimes_t  \n  \wt{\o} \big)  \nonumber  \\
  & \tn \dn    = & \tn \dn    \b1_{\{ \ga (\wt{\o}) \le \wp_n  \n  (\o   \otimes_t     \wt{\o})    \}}
 \Big(    U  \n   \big(\ga ( \wt{\o})   ,
  \o  \n \otimes_t  \n  \wt{\o} \big)
  \dn - \n  U \n  \big(   \wp_{n+1} (\o \n \otimes_t  \n  \wt{\o})
  \vee \ga ( \wt{\o}) ,   \o  \n \otimes_t  \n  \wt{\o} \big)    \Big)  \nonumber  \\
  & \tn \dn  & \tn \dn   +    \b1_{\{  \ga (\wt{\o}) > \wp_n  \n  (\o   \otimes_t     \wt{\o}) \}}
 \Big(    U  \n   \big(\wp_n (\o \n \otimes_t  \n  \wt{\o})   ,
  \o  \n \otimes_t  \n  \wt{\o} \big)
  \dn - \n  U \n  \big(  \wp_{n+1} (\o  \n \otimes_t \n   \wt{\o})   ,
   \o  \n \otimes_t  \n  \wt{\o} \big)    \Big)  \nonumber   \\
  & \tn \dn    \le   & \tn \dn     \rho_0 \Big( \d
    \n + \n  \underset{r \in [\ga (\wt{\o}),(\ga (\wt{\o})+\d) \land T]}{\sup}
 \big|  B^t_r  (   \wt{\o}   )      \n - \n   B^t_\ga  (   \wt{\o}  )  \big| \Big)
   \n - \n  J_{\ga,T} (\wt{\o})   ,    \label{eq:eh134}
 \eea
 where we obtained from  \eqref{eq:aa211} that
 \beas
  && \hspace{-1.5cm}   U    \big(\ga ( \wt{\o})   ,
  \o  \n \otimes_t  \n  \wt{\o} \big)
  \dn - \n  U   \big(   \wp_{n+1} (\o \n \otimes_t  \n  \wt{\o})
  \vee \ga ( \wt{\o}) ,   \o  \n \otimes_t  \n  \wt{\o} \big) \nonumber \\
  & \tn \le   & \tn   \rho_0 \Big( \big(  \wp_{n+1} (\o \n \otimes_t  \n  \wt{\o})  \vee \ga ( \wt{\o})
  \n - \n  \ga (\wt{\o})\big)  \n + \n  \underset{r \in [0,T]}{\sup}
 \big|   (\o  \n \otimes_t  \n  \wt{\o})
  \big( r  \land ( \wp_{n+1} (\o \n \otimes_t  \n  \wt{\o})  \vee \ga ( \wt{\o})) \big)
  \n - \n   (\o  \n \otimes_t  \n  \wt{\o})  (  r \land \ga (\wt{\o})  )  \big| \Big)  \nonumber  \\
   & \tn \le   & \tn  \rho_0 \Big( \big(  T \n - \n  \ga (\wt{\o})\big)
 \n + \n  \underset{r \in [\ga (\wt{\o}),T]}{\sup}
 \big|     \wt{\o}   ( r   )
  \n - \n     \wt{\o}   (    \ga (\wt{\o})  )  \big| \Big)
  \le \rho_0 \Big( \d  \n + \n  \underset{r \in [\ga (\wt{\o}),(\ga (\wt{\o})+\d) \land T]}{\sup}
 \big|    B^t_r  (   \wt{\o}   )      \n - \n   B^t_\ga  (   \wt{\o}  )  \big| \Big) . 
 \eeas

 On the other hand, let  $\wt{\o}  \n \in \n  \{ \ga    \n \le \n  T \n - \n \d \}$.
 applying \eqref{eq:eh131b}
 with $(t',\o') \n =  \n  \big(\ga  (\wt{\o}) , \o  \n \otimes_t  \n  \wt{\o} \big)  $
   yields   that
 \bea
  &&  \hspace{-1cm}  (\sY^n)^{t,\o}_\ga  (  \wt{\o}  )
 - \big(\sY^{n+1}\big)^{t,\o}_{\wt{\ga}} (  \wt{\o}  )
 =    \sY^n  \big(\ga  (\wt{\o}) , \o \otimes_t  \wt{\o} \big)
 -  \sY^{n+1}  \big(\ga  (\wt{\o}) \n + \n  \d ,  \o \otimes_t   \wt{\o} \big)  \nonumber  \\
 &&   \le    \b1_{\{  \ga(\wt{\o}) + \d \le \wp_{n+1} (\o \otimes_t \wt{\o})    \}}
 \Big( L(\ga(\wt{\o}) , \o  \n \otimes_t \n  \wt{\o})
 \n - \n   L(\ga(\wt{\o}) \n + \n  \d , \o  \n \otimes_t \n  \wt{\o})    \Big) \nonumber  \\
 && \q   +    \b1_{\{ \wp_{n+1} (\o \otimes_t \wt{\o})  - \d <
 \ga(\wt{\o}) \le \wp_n (\o \otimes_t \wt{\o})    \}}
 \Big( L(\ga(\wt{\o}) , \o  \n \otimes_t \n  \wt{\o})
 \n - \n   L \big( \wp_{n+1} (\o \otimes_t \wt{\o})  \n \vee \n   \ga  (\wt{\o}) ,
 \o  \n \otimes_t \n  \wt{\o} \big)    \Big)  \nonumber \\
 &&     \q   +    \b1_{\{ \ga(\wt{\o})> \wp_n (\o \otimes_t \wt{\o})    \}}
     \Big( U(\wp_n (\o  \n \otimes_t \n  \wt{\o}) , \o  \n \otimes_t \n  \wt{\o})
    \n - \n U \big(  \wp_{n+1} (\o  \n \otimes_t \n  \wt{\o})     ,
    \o  \n \otimes_t \n  \wt{\o} \big)  \Big)  \nonumber \\
 &&    \le \n  \rho_0 \Big( \d  \n + \n
  \underset{r \in [\ga  (\wt{\o}), (\ga  (\wt{\o})   +    \d) \land T ]}{\sup}
 \big|   B^t_r  (   \wt{\o}   )      \n - \n   B^t_\ga  (   \wt{\o}  )  \big| \Big)  \n - \n  J_{\ga,T} (\wt{\o}) ,
 \label{eq:a123}
 \eea
  where we derived from  \eqref{eq:aa211} that if $   \wp_{n+1} (\o \oti \wt{\o})     \<
 \ga (\wt{\o}) \+ \d      $,
 \beas
  && \hspace{-1.5cm}   L    \big(\ga ( \wt{\o})   ,
  \o  \n \otimes_t  \n  \wt{\o} \big)
  \dn - \n  L   \big(   \wp_{n+1} (\o \n \otimes_t  \n  \wt{\o})
  \vee \ga ( \wt{\o}) ,   \o  \n \otimes_t  \n  \wt{\o} \big) \nonumber \\
  & \tn \le   & \tn   \rho_0 \Big( \big(  \wp_{n+1} (\o \n \otimes_t  \n  \wt{\o})  \vee \ga ( \wt{\o})
  \n - \n  \ga (\wt{\o})\big)  \n + \n  \underset{r \in [0,T]}{\sup}
 \big|   (\o  \n \otimes_t  \n  \wt{\o})
  \big( r  \land ( \wp_{n+1} (\o \n \otimes_t  \n  \wt{\o})  \vee \ga ( \wt{\o})) \big)
  \n - \n   (\o  \n \otimes_t  \n  \wt{\o})  (  r \land \ga (\wt{\o})  )  \big| \Big)  \nonumber  \\
   & \tn \le   & \tn  \rho_0 \Big( \d
 \n + \n  \underset{r \in [\ga (\wt{\o}),(\ga (\wt{\o})+\d) \land T]}{\sup} \big|     \wt{\o}   ( r   )
  \n - \n     \wt{\o}   (    \ga (\wt{\o})  )  \big| \Big)
  \le \rho_0 \Big( \d  \n + \n  \underset{r \in [\ga (\wt{\o}),(\ga (\wt{\o})+\d) \land T]}{\sup}
 \big|    B^t_r  (   \wt{\o}   )      \n - \n   B^t_\ga  (   \wt{\o}  )  \big| \Big) .
 \eeas
 Combining \eqref{eq:eh134} with  \eqref{eq:a123},
 we see from \eqref{eq:et011} and  \eqref{eq:ex015}  that
 \beas
  \hE_\hP \big[  ( \sY^n)^{t,\o}_{\ga}  \big]
    \n \le \n  \hE_\hP \Big[  \big(\sY^{n+1}\big)^{t,\o}_{\wt{\ga}} \n - \n J_{\ga,T}  \Big]
     \n + \n   \wh{\rho}_0(\d)
  \n \le \n   \sZ^{n+1}_t(\o) \n - \n \wt{\D}^U   \n + \n    2 \wh{\rho}_0 (\d)   .
 \eeas
   Then taking supremum over $(\hP, \ga)  \n \in \n  \cP_t    \n \times \n \cT^t $
 on the left-hand-side  leads to that
 $    \sZ^n_t(\o)   \n \le \n   \sZ^{n+1}_t(\o)  \n - \n \wt{\D}^U    \n + \n   2 \wh{\rho}_0 (\d)   $. \qed

  \no {\bf Proof of Proposition \ref{prop_Z_ultim}:}
 {\bf 1)} Let $n \n \in \n  \hN$. Lemma \ref{lem_Y_nkl} and Proposition \ref{prop_DPP_Z} show that
  $Z^{n,k}$, $k  \n \in \n  \hN$ are $\bF-$adapted processes with all continuous paths.
 For any $(t,\o)  \n \in \n  [0,T] \n \times \n \O$, as
 $k \n \to \n  \infty$ in \eqref{eq:eh137},   the continuity of  $U$ implies that
 \bea \label{eq:et311}
    \lmt{k \to \infty} Z^{n,k}_t (\o) = \sZ^n_t (\o)  .
 \eea
     Then the $\bF-$adaptedness of
 $\{Z^{n,k}\}_{k \in \hN}$ shows that process $\sZ^n$ is also   $\bF-$adapted.

 Given $(s,\o)  \n \in \n  [0,T] \n \times \n \O$,
 letting $t  \n \to \n  s$ in \eqref{eq:eh137},
 we can deduce from the continuity of processes $U$,  $\{Z^{n,k}\}_{k \in \hN}$   that
 \beas
 && \hspace{-1.5cm}  Z^{n,k}_s (\o) - \wh{\rho}_0(2^{1-k}) - U \big( ( \wp_n (\o)  \n + \n  2^{1-k} ) \land s  ,   \o    \big)
      +     U  (\wp_n (\o) \land s , \o   )   \le \linf{t \to s}  \,  \sZ^n_s (\o)  \le \lsup{t \to s}  \,  \sZ^n_s (\o)  \\
 &&  \le    Z^{n,k}_s (\o) + 2 \wh{\rho}_0(2^{1-k})
   - U \big( ( \wp_n (\o)  \n + \n  2^{1-k} ) \land s  ,   \o    \big)   +     U  (\wp_n (\o) \land s , \o   ) ,
    \q  \fa k \in \hN    .
 \eeas
 As $k \n \to \n  \infty$, \eqref{eq:et311} and the continuity of $U$ imply that
 $ \lmt{t \to s} \, \sZ^n_t (\o) = \lmt{k \to \infty} Z^{n,k}_s (\o) = \sZ^n_s (\o) $.
 Hence, the process $\sZ^n$ has all continuous paths.

  \no {\bf 2)}
 Fix $(t,\o)  \n \in \n  [0,T] \n \times \n \O$.
 For any integers $n \n < \n m $, adding \eqref{eq:eh137b} up
 from $i \n = \n n$ to $i  \n = \n  m \n - \n 1$ shows that
 \bea \label{eq:et314}
 -  2  \sum^{m-1}_{i=n} \wh{\rho}_0 \big( \hb{$\frac{2T}{i+3}$} \big) & \tn  \le & \tn
     \sZ^{m}_t(\o) \n - \n  \sZ^{n}_t(\o)
  \n - \n  U \big(  \wp_{m} (\o)    \n \land \n  t , \o \big)
  \n + \n U \big(  \wp_{n} (\o)    \n \land \n  t , \o \big)
        \le          \sum^{m-1}_{i=n}  \wh{\rho}_0 \big( \hb{$\frac{2T}{i+3}$} \big) .
 \eea
 Since  $\wh{\fp}_1 \n > \n  1$ by (P2),   \eqref{eq:gb017} gives that
 $ \sum^\infty_{i= 0}   \wh{\rho}_0 \big( \hb{$\frac{2T}{i+3}$} \big)  \n \le \n
 \sum^{n_0-1}_{i= 0 }   \wh{\rho}_0 \big( \hb{$\frac{2T}{i+3}$} \big) \+
 \wh{\fC} \sum^\infty_{i= n_0}  \n  \big(\hb{$\frac{2T}{i+3}$} \big)^{\wh{\fp}_1}  \dn < \n  \infty$,
 where   $ n_0  \n := \n 1 \n + \n \lfloor (2T \n - \n 3)^+ \rfloor $. Then
 we see from the continuity of $U$ and \eqref{eq:et314} that $\big\{\sZ^n_t(\o)\big\}_{n \in \hN}$
 is a Cauchy sequence of $\hR$. Let $\sZ_t (\o)$ be the limit of $\big\{\sZ^n_t(\o)\big\}_{n \in \hN}$, i.e.
 $\sZ_t (\o)  \n := \n  \lmt{n \to \infty} \sZ^n_t(\o) $.
 As  $ \lmtu{m \to \infty} \tau_m (\o)  \n = \n  \tau_0 (\o) $,   Proposition \ref{prop_wp_n} (1)
 shows that $ \lmtu{m \to \infty} \wp_m (\o)  \n = \n  \tau_0 (\o) $.
 Letting $m \n \to \n \infty$ in \eqref{eq:et314} and using the continuity of $U$ yield \eqref{eq:et317}.

     \no {\bf 3a)} {\it Let us now show the first inequality of \eqref{eq:et335}. }

     Clearly, the $\bF-$adaptedness of
 $\{\sZ^n\}_{n \in \hN}$ implies that of    $\sZ $ and the boundedness of
 $\{\sZ^n\}_{n \in \hN}$ by $M_0$ implies that of    $\sZ $.
  \if{0}
 Moreover,  for any $(s,\o)  \n \in \n  [0,T] \n \times \n \O$,
 as $t \to s$ in \eqref{eq:et317},
   the continuity of processes $\{\sZ^n\}_{n \in \hN}$ and $U$ implies that
   for any $n \in \hN$
 \beas
    \sZ^n_s (\o) \- 2 \e_n \+     U  (\tau_0 (\o) \ld s , \o   )
 \- U \big(   \wp_n (\o)   \ld s  ,   \o    \big)
        \ls \linf{t \to s}  \,  \sZ_s (\o)  \ls \lsup{t \to s}  \,  \sZ_s (\o)
    \ls    \sZ^n_s (\o) \+   \e_n \+     U  (\tau_0 (\o) \ld s , \o   )
   \- U \big(   \wp_n (\o)   \ld s  ,   \o    \big)      .
 \eeas
 Since $ \lmt{n \to \infty} \wp_n (\o) = \tau_0 (\o) $,
 letting  $n \n \to \n  \infty$, we see from   the continuity of $U$   that
 $ \lmt{t \to s} \, \sZ_t (\o) = \lmt{n \to \infty} \sZ^n_s (\o) = \sZ_s (\o) $.
 Therefore, the process $\sZ$ has all continuous paths.
  \fi
 Similar to the argument used in part 1), letting  $t \to s$ in \eqref{eq:et317},
 we can deduce from   the continuity of processes  $\{\sZ^n\}_{n \in \hN}$,  $U$
 and $\lmtu{n \to \infty} \wp_n    \n = \n  \tau_0$  that
   the process $\sZ$ has all continuous paths.

 Let $(t,\o)  \n \in \n  [0,T] \n \times \n \O$.
   Given $\e  \n > \n  0$, there exists $( \hP_\e , \ga_\e    )  \n \in \n  \cP_t    \n \times \n \cT^t  $ such that
   $ \underset{(\hP, \ga) \in \cP_t \times \cT^t}{\sup}  \,
  \hE_\hP \Big[ \wh{\sY}\,^{t,\o}_{\ga} \Big]
  \le \hE_{\hP_\e} \Big[ \wh{\sY}\,^{t,\o}_{\ga_\e} \Big]  \n + \n  \e $.
   Since $ \lmtu{n \to \infty} \tau_n    \n = \n  \tau_0 $, one can deduce from the continuity of $U$ that
  \bhe
  \bea \label{eq:et415}
  \lmt{n \to \infty} \sY^n_{t'} (\o') \n = \n  \wh{\sY}_{t'} (\o')  , \q \fa (t',\o')  \n \in \n  [0,T] \n \times \n \O .
  \eea
  \ehe
  It follows that
  $  \lmt{n \to \infty} (\sY^n)^{t,\o}_{\ga_\e} (\wt{\o})
   \= \lmt{n \to \infty} \sY^n \big( \ga_\e  (\wt{\o}) , \o \otimes_t \wt{\o} \big)
   \= \wh{\sY} \big( \ga_\e  (\wt{\o}) , \o \otimes_t \wt{\o} \big) \= \wh{\sY}\,^{t,\o}_{\ga_\e} (\wt{\o}) $,
   $\fa \wt{\o} \n \in \n \O^t$.
 As $\sY^n$'s are all bounded by $M_0$,
   applying the bounded convergence theorem yields that
 \beas
 \underset{(\hP,\ga)  \in \cP_t \times \cT^t}{\sup}   \,
 \hE_\hP \Big[\wh{\sY}\,^{t,\o}_{\ga} \Big] \n \le \n  \hE_{\hP_\e} \Big[\wh{\sY}\,^{t,\o}_{\ga_\e} \Big]  \n + \n  \e
  \n = \n  \lmt{n \to \infty} \hE_{\hP_\e} \big[ (\sY^n)^{t,\o}_{\ga_\e} \big]  \n + \n  \e
  \n \le \n 
    \lmt{n \to \infty} \sZ^n_t (\o)  \n + \n  \e
   \n = \n  \sZ_t (\o)  \n + \n  \e .
 \eeas
 Then letting $\e \n \to \n  0$ leads to that $ \sZ_t (\o)  \n \ge \n  \underset{(\hP, \ga) \in \cP_t \times \cT^t}{\sup}  \,
 \hE_\hP \Big[\wh{\sY}\,^{t,\o}_{\ga} \Big]
  \n \ge \n  \underset{\hP \in \cP_t }{\sup}  \, \hE_\hP \Big[\wh{\sY}\,^{t,\o}_t \Big]
  \n = \n  \wh{\sY}_t (\o) $, where we used the $\bF-$adaptedness of $\wh{\sY}$
   and \eqref{eq:bb421} in the last equality.

  \no {\bf 3b)} {\it Let $(t,\o)  \n \in \n  [0,T] \n \times \n \O$. We verify the third equality
  of \eqref{eq:et335} by two cases. }

   If $\tau_0 (\o) \= T$, \eqref{eq:et043b}
 and the continuity of $U$ imply that
 \beas
 \sZ_T (\o) \= \lmt{n \to \infty} \sZ^n_T (\o)
 \= \lmt{n \to \infty} \, \underset{\hP \in \cP_t}{\sup} \hE_\hP \big[ (\sY^n)^{T,\o}_T \big]
 \= \lmt{n \to \infty} \sY^n_T (\o) \= \lmt{n \to \infty} U \big(\wp_n (\o), \o\big)
 \= U \big(\tau_0 (\o), \o\big) .
 \eeas

 Suppose next that $ \tau_0 (\o) \< T $. By the definition of $\tau_0 (\o)$,   the set
 $\{t  \n \in \n  [0,T] \n : \sX_t  (\o )  \n \le \n  0\}$  is not empty.
 So Proposition \ref{prop_wp_n} shows that    $\wp_n (\o) \n <  \n  \tau_0 (\o)$.

  Let $t \ins [\tau_0 (\o),T]$ and  $n \n \in \n \hN$.
 As $t_n  \n := \n  \wp_n (\o)  \n < \n   \tau_0 (\o)  \n \le \n  t $,
 Lemma \ref{lem_shift_stopping_time} implies that
 $ \wp_n (\o \n \otimes_t \n \O^t) \n = \n   \wp_n (\o )  \n = \n t_n $.
 Let $\ga  \n \in \n  \cT^t$.
 Since   $U$ is an $\bF-$adapted process by (A1) and \eqref{eq:et014},
 one has $U_{t_n}  \n \in \n  \cF_{t_n}    \n \subset \n  \cF_t $. Given   $\wt{\o}  \n \in \n  \O^t$,
 using \eqref{eq:bb421} with $(t,s,\eta)  \n = \n  (0,t,U_{t_n})$  shows that
  $  U  (t_n   ,  \o  \n \otimes_t \n  \wt{\o}  )   \n = \n  U  (t_n ,    \o  )    $.
  Then we can deduce from   $\ga (\wt{\o})  \n \ge \n  t  \n > \n  t_n  \n = \n \wp_n (\o \n \otimes_t \n \wt{\o}) $ that
 \beas
 (\sY^n)^{t,\o}_\ga (\wt{\o})
   =   \sY^n \big( \ga  (\wt{\o}) , \o \oti_t \wt{\o} \big)
   = U \big(\wp_n (\o \n \otimes_t \n \wt{\o}), \o \n \otimes_t \n \wt{\o} \big)
   = U (t_n, \o \n \otimes_t \n \wt{\o}) = U (t_n, \o) ,
 \eeas
 which leads to that
 $   \sZ^n (t,\o)  \=  \underset{(\hP, \ga) \in \cP_t \times \cT^t}{\sup} \,
    \hE_\hP  \big[  (\sY^n)^{ t , \o  }_\ga     \big]
  \= U (t_n, \o) \= U \big( \wp_n (\o) , \o \big) $.
 Letting $n \to \infty$, we obtain from the continuity of $U$ that
 $ \sZ \big( t , \o  \big) = U \big( \tau_0 (\o) , \o  \big)  $.

  \no {\bf 4)}  By \eqref{eq:et411} and the continuity of $\sZ$ obtained in part 3a),
  $D_t \n := \n  \sZ_t  \n - \n \wh{\sY}_t  \n \ge \n  0$, $t \in [0,T]$
is an $\bF-$adapted process whose paths are all continuous except a possible negative jump at $\tau_0$.
In particular, each path of $D$ is lower-semicontinuous and right-continuous.
It follows that  $\ga_*$ is an $\bF-$stopping time (see   Lemma A.13 in the ArXiv version of \cite{ROSVU} for a proof).

\if{0}
 Let  $ t  \n \in \n  [0,T) $.   We claim that for any $\o \n \in \n  \O $
 \bea \label{eq:et331}
 \hb{ if   $ D_s (\o)  \n > \n   0 ,~ \fa s  \n \in \n  [0,t]$, then
 $ \underset{s \in [0,t]}{\inf} D_s (\o) \n > \n   0  $. }
 \eea
 Assume not, i.e. there exists a $\o' \n \in \n \O$ such that
  $ D_s (\o')  \n > \n   0 $, $ \fa s  \n \in \n  [0,t]$ and
 $ \underset{s \in [0,t]}{\inf} D_s (\o') \n \le \n   0  $. Then one can find   a
 sequence $\{s_n  \n = \n  s_n (t,\o')\}_{n \in \hN}$ of $[0,t]$ such that
 $ \lmtd{n \to \infty} D_{s_n} (\o')    \n  = \n  \underset{s \in [0,t]}{\inf} D_s (\o')  $.
 Clearly, $\{s_n  \}_{n \in \hN}$ has a convergent subsequence
 $\{s_{n_i}  \}_{i \in \hN}$ with limit $s_*  \n \in \n  [0,t]$.
 We can deduce from the lower-semicontinuity of $D$   that
 $0 \n < \n D_{s_*} ( \o' )  \n \le \n  \linf{s \to s_*} D_s(\o')
  \n \le \n  \lmtd{i \to \infty} D_{s_{n_i}} (\o' )  \n = \n  \underset{s \in [0,t]}{\inf} D_s (\o')  \n \le \n  0 $.
 An contradiction appears. So \eqref{eq:et331} holds and it follows that
         \bea \label{eq:et333}
     \{ \ga_* \n > \n  t\} & \tn =& \tn  \{\o  \n \in \n  \O \n :
      D_s (\o)  \n > \n   0,~ \fa s  \n \in \n  [0,t] \}
      \n = \n  \underset{n    \in    \hN}{\cup} \{\o  \n \in \n  \O \n :
      D_s (\o)  \n \ge \n    1/n,~ \fa s  \n \in \n  [0,t] \} .
      \eea
 For any $n \n \in \n  \hN$, the right-continuity of $D$ implies  that
 $\{\o  \n \in \n  \O \n :
      D_s (\o)  \n \ge \n   1/n, \; \fa s  \n \in \n  [0,t] \}
       \n = \n  \{\o  \n \in \n  \O \n :
     D_s (\o)  \n \ge \n   1/n, \; \fa s  \n \in \n  \hQ_{0,t} \} $,
     where $\hQ_{0,t}  \n := \n  \big(  [0,t]  \n \cap \n  \hQ \big)  \n \cup \n  \{t \} $.
     Putting these equalities back into \eqref{eq:et333} yields that
             \beas
     \{ \ga_* \n > \n  t \} & \tn =& \tn  \underset{n    \in    \hN}{\cup} \{\o  \n \in \n  \O \n :
     D_s (\o)  \n \ge \n   1/n,~ \fa s  \n \in \n  \hQ_{0,t} \}
      \n = \n  \underset{n    \in    \hN}{\cup} \, \underset{s \in \hQ_{0,t}}{\cap}
      \{\o  \n \in \n  \O \n :
      D_s (\o)  \n \ge  \n   1/n  \}  \n \in \n  \cF_t .
     \eeas
     Hence,    $\ga_*$ is an $\bF-$stopping time.

\fi

     As $\sZ_t \n = \n  U_{\tau_0}
        \n = \n  \wh{\sY}_t   $, $ \fa t  \n \in \n  [\tau_0 , T]$ by \eqref{eq:et335}, one can deduce that
     $\ga_*  \n = \n  \ga_*  \ld \tau_0 \=
     \inf \big\{t  \n \in \n  [0,\tau_0) \n : \sZ_t  \n = \n  \wh{\sY}_t \, \big\} \ld \tau_0
      \n = \n   \inf \big\{t  \n \in \n  [0,\tau_0) \n : \sZ_t  \n = \n   \sY_t  \,  \big\} \ld \tau_0
      \n = \n   \inf\{t  \n \in \n  [0,\tau_0) \n : \sZ_t  \n = \n  L_t \}   \n \land \n  \tau_0 $.  \qed

 \subsection{Proof of Theorem \ref{thm_RDOSRT}}

    For any $m \n \in \n  \hN$, applying Theorem \ref{thm_cst} with $(Y,\wp) \n = \n (Y^{m,m},\wp^{m,m})$ shows that
    there exists a $\hP_m  \n \in \n  \cP$ such that
 \bea \label{eq:et343}
 Z^{m,m}_0 =    \hE_{\hP_m} \Big[ Z^{m,m}_{ \nu_m \land \z} \Big] , \q \fa \z \in \cT ,
 \eea
 where $\nu_m \n : = \n  \inf\big\{ t  \n \in \n  [0,T] \n : Z^{m,m}_t  \n = \n  \wh{Y}^{m,m}_t \big\}  \n \in \n  \cT$.
 By (P1), $\{\hP_m\}_{m \in \hN}$ has a weakly convergent sequence $\{\hP_{m_j}\}_{i \in \hN}$ with limit $\hP_*$.

 \no {\bf 1)} {\it First,  we    use \eqref{eq:eh137}, \eqref{eq:et317} and similar arguments to those proving
 Theorem \ref{thm_cst} to show that
  \bea \label{eq:et389}
  \sZ_0   \n  \le  \n
  \hE_{\hP_*} \Big[ \lsup{n \to \infty} \, \lsup{i \to \infty} \,
  \lsup{\ell \to \infty} \sZ_{\z_{i,\ell} \land \wp_n} \Big]  ,
  \eea
  where $ \z_{i,\ell} \df  \inf\big\{t  \n \in \n  [0,T] \n : Z^{\ell,\ell}_t
  \n \le \n  L_t  \n + \n  1/i   \big\}    \land    T $. This part is relatively lengthy, we will split it into
  several steps.}

 \no {\bf 1a)} {\it We start with an auxiliary inequality:
 for any $n, k \n \in \n  \hN$ with  $k \ge n$ and $\o  \n \in \n  \O$,}
 \bea \label{eq:et341}
 \big|  Z^{k,k}_t (\o)  \n - \n  \sZ_t(\o) \big| \le    \ol{\e}_k
 \df 2 \wh{\rho}_0 \big( 2^{1-k}   \big) \n + \n 2 \sum^\infty_{i = k} \wh{\rho}_0 (\hb{$\frac{2T}{i+3}$}) ,
 \q \fa t \in [0,\wp_n (\o)] .
\eea

  Let $n, k \n \in \n  \hN$ with  $k \ge n$ and   let $\o  \n \in \n  \O$.
 For any $t  \n \in \n  [0,T]$, we see from
  \eqref{eq:eh137} and \eqref{eq:et317}    that
 $ -2  \wh{\rho}_0 \big( 2^{1-k}   \big)  \n \le \n  Z^{k,k}_t (\o)  \n - \n  \sZ^k_t (\o)
 \n - \n U \big( ( \wp_k (\o)  \n + \n  2^{1-k} )  \n \land \n  t  ,   \o    \big)
 \n + \n U \big(   \wp_k (\o)    \n  \land \n  t  ,   \o    \big)
  \n \le \n
  \wh{\rho}_0 \big( 2^{1-k}   \big)$ and that
 $  -  \n   \sum^\infty_{i = k} \wh{\rho}_0 (\hb{$\frac{2T}{i+3}$})
   \n \le \n  \sZ^k_t(\o)  \n - \n \sZ_t(\o) \n - \n U \big(   \wp_k (\o)    \n  \land \n  t  ,   \o    \big)
   \n + \n U \big(   \tau_0 (\o)    \n  \land \n  t  ,   \o    \big)
     \n \le \n  2 \sum^\infty_{i = k} \wh{\rho}_0 (\hb{$\frac{2T}{i+3}$})  $.
 Adding them together yields that
\bea   \label{eq:et340}
 -  \ol{\e}_k  \n  \le  \n   Z^{k,k}_t (\o)  \n - \n  \sZ_t(\o)
    \n - \n U \big( ( \wp_k (\o)  \n + \n  2^{1-k} )  \n \land \n  t  ,   \o    \big)
      \n  +  \n     U \big(  \tau_0 (\o)    \n \land \n  t , \o \big) \le  \ol{\e}_k , \q \fa t \in [0,T] .
\eea
In particular, for any $t  \n \in \n  [0,\wp_n (\o)]$,
 since $t  \n \le \n  \wp_n (\o)  \n \le \n  \wp_k (\o)  \n \le \n  \tau_0 (\o) $
 by Proposition \ref{prop_wp_n} (1), one has
 $  U \big( ( \wp_k (\o)  \n + \n  2^{1-k} )  \n \land \n  t  ,   \o    \big)
 \n = \n  U \big(  \tau_0 (\o)    \n \land \n  t , \o \big)   \n = \n  U(t,\o) $.
 Then \eqref{eq:et341} directly follows from   \eqref{eq:et340}.

 Now, fix integers    $ 1 \ls n  \< i \<  \ell  \< \a         $ such that
 $    \ol{\e}_\ell \ls \frac{1}{2i} $ and fix
 $j \ins \hN$   such that $  m_j \gs \a     $.
 Since Lemma \ref{lem_Y_nkl},
 Proposition \ref{prop_DPP_Z}, (A1) and \eqref{eq:et014} show that $Z^{\ell,\ell} \n - \n L$
 is an $\bF-$adapted process with all continuous paths, 
    \bea \label{eq:et361}
   \z^\a_{i,\ell} : =  \inf\big\{t \in [0,T]: Z^{\ell,\ell}_t \le L_t + 1/i + 1/\a \big\} \land T ~
   \hb{  defines an $\bF-$stopping time. }
    \eea
 Similar to $\nu_n$ in \eqref{eq:et345},
 $  \wh{\z}^\a_{i,\ell} \n : =  \n  \inf\big\{ t  \n \in \n  [0,T] \n :
 Z^{\ell,\ell}_t  \n \le \n  \wh{Y}^{\ell,\ell}_t  \n + \n  1/i  \n + \n  1/\a \big\}  $
 is also an $\bF-$stopping time  satisfying
 \bhe
 \bea \label{eq:a125}
 \wh{\z}^\a_{i,\ell} \n \land \n  \wp_n \=  \z^\a_{i,\ell} \n \land \n  \wp_n   \n \le \n  \nu_{m_j}  \n \land \n  \wp_n .
 \eea
 \ehe
 Then  applying \eqref{eq:et341} with $(k,t) \n = \n (m_j,0)$, $(k,t)
 \n = \n \big(m_j , \wh{\z}^\a_{i,\ell}   \land  \wp_n \big) $
 and $(k,t) \n = \n \big(\ell , \wh{\z}^\a_{i,\ell}  \land \wp_n \big) $ respectively
 as well as applying \eqref{eq:et343}
 with $(m,\z) \n = \n \big( m_j, \wh{\z}^\a_{i,\ell} \n \land \n  \wp_n \big)$, we obtain
 \bea
 \sZ_0 \n - \n     \ol{\e}_{m_j}   \n \le \n  Z^{m_j,m_j}_0   \n = \n
   \hE_{\hP_{m_j}} \Big[ Z^{m_j,m_j}_{\nu_{m_j} \land \wh{\z}^\a_{i,\ell} \land \wp_n } \Big]
    \n = \n  \hE_{\hP_{m_j}} \Big[  Z^{m_j,m_j}_{ \wh{\z}^\a_{i,\ell}  \land \wp_n } \Big]
  \n \le \n  \hE_{\hP_{m_j}} \Big[ \sZ_{\wh{\z}^\a_{i,\ell}  \land \wp_n } \Big]
    \n + \n    \ol{\e}_{m_j}  \n \le \n
 \hE_{\hP_{m_j}} \Big[ Z^{\ell,\ell}_{\wh{\z}^\a_{i,\ell}  \land \wp_n } \Big]
    \n + \n    \ol{\e}_{m_j}  \n + \n   \ol{\e}_\ell .  \q  \label{eq:et349}
   \eea

  \no {\bf 1b)}  {\it Before sending  $j$ to $\infty$   in order to approximate the distribution $  \hP_*  $
  in \eqref{eq:es021},  we need to approach
     $\big\{ \wh{\z}^\a_{i,\ell} \big\}_{\a \in \hN}$ by a sequence
  $\big\{  \th^\a_{i,\ell}  \big\}_{\a  \in \hN}$ of   Lipschitz  continuous random variables
  and estimate the expected difference $ \hE_{\hP_{m_j}} \Big[ \Big| Z^{\ell,\ell}_{\wh{\z}^\a_{i,\ell}  \land \wp_n }
  \n - \n   Z^{\ell,\ell}_{\th^\a_{i,\ell}  \land \wp_n } \Big| \Big] $. }

    Recall from Lemma \ref{lem_Y_nkl} and the remark following it that
 $ Y^{\ell,\ell}$ is   uniformly continuous on $  [0,T]  \n \times \n  \O$
 with respect to the modulus of continuity function $\rho_{\ell,\ell}$ and that
  $\wp^{\ell,\ell}  $  is a Lipschitz continuous stopping time on $\O$ with coefficient $\k_\ell$.
 Replacing $(Z,\wh{Y}, \nu_n)$   by
 $\big(Z^{\ell,\ell},\wh{Y}^{\ell,\ell}, \wh{\z}^\a_{i,\ell} \big)$
 in  the arguments   
   that lead to \eqref{eq:eh211}, we can
 \if{0}
  Fix   $\o  \n \in \n  \O$.
   There exists a $\l^\ell_\a  \n > \n  0$ such that
   $\rho_{\overset{}{\ell,\ell}} (x)  \n \vee \n  \wh{\rho}_{\overset{}{\ell,\ell}} (x)  \n \le \n  \frac{1}{2 \a (\a \n + \n 1)}$, $\fa x  \n \in \n  [0, \l^\ell_\a]$.
 Set $ \dis \d^\ell_\a ( \o )  \n : = \n  \frac{\l^\ell_\a}{2(1 \n + \n \k_\ell)}  \n \land \n  \frac{(\phi^\o_T)^{-1}(\l^\ell_\a/2)}{\k_\ell}$ with $(\phi^\o_T)^{-1} (x)
 \n : = \n  \inf\{y  \n > \n  0 \n : \phi^\o_T (y)  \n = \n  x \} $, $\fa x  \n > \n  0 $,
 and let   $\o'  \n \in \n  \ol{O}_{\d^\ell_\a ( \o )} (\o)$.
 Given $t  \n \in \n  [0,T]$,   set $s \n := \n  t  \n \land \n  \wp^{\ell,\ell}(\o)$
 and $s' \n := \n  t  \n \land \n  \wp^{\ell,\ell}(\o')$.
    By \eqref{eq:ec017}, $ |s \n - \n s' |  \n \le \n  \big| \wp^{\ell,\ell}(\o)  \n - \n  \wp^{\ell,\ell}(\o') \big|
   \n \le \n   \k_\ell   \|\o    \n - \n  \o'  \|_{0,T}   $.
 Then   \eqref{eq:aa211} implies that
 \bea
 && \hspace{-1cm}  \big|\wh{Y}^{\ell,\ell}  ( t,\o ) \n - \n  \wh{Y}^{\ell,\ell}  ( t,\o') \big| \n  = \n
 | Y^{\ell,\ell} ( s,\o )  \n - \n  Y^{\ell,\ell} ( s',\o') |
   \n \le \n     \rho_{\overset{}{\ell,\ell}} \Big( |s \n - \n s'  |
  \n + \n  \underset{r \in [0,T]}{\sup} \big|  \o   (r  \n \land \n  s)
   \n - \n   \o'   (r  \n \land \n  s' ) \big|   \Big)  \nonumber \\
  &  & \le      \rho_{\overset{}{\ell,\ell}} \Big(  \k_\ell   \|\o \n  - \n  \o'  \|_{0,T}
  \n + \n  \underset{r \in [0,T]}{\sup}  |   \o   (r  \n \land \n  s  )
  \n - \n   \o  (r  \n \land \n  s'  )  |
  \n + \n  \underset{r \in [0,T]}{\sup} \big|  \o (r  \n \land \n  s' )
  \n - \n  \o' (r  \n \land \n  s' ) \big|   \Big) \nonumber \\
  &  & \le
  \rho_{\overset{}{\ell,\ell}} \Big(  (1 \n + \n \k_\ell)   \|\o \n  - \n  \o'  \|_{0,T}
  \n + \n   \phi^\o_T \big( |s' \n - \n s| \big)   \Big)
   \n \le \n  \rho_{\overset{}{\ell,\ell}} \Big(  (1 \n + \n \k_\ell)   \|\o \n  - \n  \o'  \|_{0,T}
  \n + \n   \phi^\o_T \big( \k_\ell    \|\o \n  - \n  \o'  \|_{0,T} \big)   \Big)
   \n \le \n  \frac{1}{2 \a (\a \n + \n 1)}  .  \qq   \label{eq:gk011}
 \eea
Taking $t= \wh{\z}^\a_{i,\ell} (\o) $, we see from \eqref{eq:gd011} that
 \beas
&& \hspace{-1cm} \big|  (Z^{\ell,\ell} \n - \n \wh{Y}^{\ell,\ell}) ( \wh{\z}^\a_{i,\ell} (\o),\o)  \n - \n   (Z^{\ell,\ell} \n - \n \wh{Y}^{\ell,\ell}) ( \wh{\z}^\a_{i,\ell} (\o),\o') \big|
  \n \le  \n  \big|   Z^{\ell,\ell} ( \wh{\z}^\a_{i,\ell} (\o),\o)  \n - \n    Z^{\ell,\ell} ( \wh{\z}^\a_{i,\ell} (\o),\o') \big|
     \+      \big| \wh{Y}^{\ell,\ell}  ( \wh{\z}^\a_{i,\ell} (\o),\o)
   \n - \n \wh{Y}^{\ell,\ell}  ( \wh{\z}^\a_{i,\ell} (\o),\o') \big|  \\
&&  \n \le \n   \wh{\rho}_{\overset{}{\ell,\ell}} \Big( (1 \n + \n \k_\ell) \|\o \n - \n \o'\|_{0,\wh{\z}^\a_{i,\ell} (\o)}
  \n + \n \phi^\o_{\wh{\z}^\a_{i,\ell} (\o)} \big( \k_\ell   \|\o   \n  - \n  \o'  \|_{0,\wh{\z}^\a_{i,\ell} (\o)} \big)  \Big)
  \n + \n \frac{1}{2 \a (\a \n + \n 1)}  \\
&&  \n \le \n    \wh{\rho}_{\overset{}{\ell,\ell}} \Big( (1 \n + \n \k_\ell) \|\o \n - \n \o'\|_{0,T}
  \n + \n \phi^\o_T \big( \k_\ell   \|\o   \n  - \n  \o'  \|_{0,T} \big)  \Big)
  \n + \n \frac{1}{2 \a (\a \n + \n 1)}
   \n \le  \n  \frac{1}{  \a (\a \n + \n 1)} \ls \frac{1}{(\a \- 1)\a} .
\eeas
 As   the continuity of $Z^{\ell,\ell}  \n - \n \wh{Y}^{\ell,\ell}$    shows that
 \bea \label{eq:es041}
 (Z^{\ell,\ell}  \n - \n \wh{Y}^{\ell,\ell}) \big( \wh{\z}^\a_{i,\ell} (\o),\o \big)  \n \le \n
 \frac{1}{i} + \frac{1}{\a } ,
 \eea
 it follows that
 $ (Z^{\ell,\ell}  \n - \n \wh{Y}^{\ell,\ell}) ( \wh{\z}^\a_{i,\ell} (\o),\o')
  \n \le \n
  \frac{1}{i} \n + \n \frac{1}{\a }  \n + \n  \frac{1}{ (\a \- 1)\a  }
   \n = \n   \frac{1}{i} \n + \n  \frac{1}{\a \- 1  }   $,
 so $\wh{\z}^{\a-1}_{i,\ell} (\o')  \n \le \n  \wh{\z}^\a_{i,\ell} (\o)$. Analogously,
 taking $t= \wh{\z}^{\a+1}_{i,\ell} (\o') $ in \eqref{eq:gk011} yields that
  \beas
&& \hspace{-1.2cm} \big|  (Z^{\ell,\ell} \n - \n \wh{Y}^{\ell,\ell}) ( \wh{\z}^{\a+1}_{i,\ell} (\o'),\o)  \n - \n   (Z^{\ell,\ell} \n - \n \wh{Y}^{\ell,\ell}) ( \wh{\z}^{\a+1}_{i,\ell} (\o'),\o') \big| \\
&&   \n \le  \n  \big|   Z^{\ell,\ell} ( \wh{\z}^{\a+1}_{i,\ell} (\o'),\o)  \n - \n    Z^{\ell,\ell} ( \wh{\z}^{\a+1}_{i,\ell} (\o'),\o') \big|
  \n + \n   \big| \wh{Y}^{\ell,\ell}  ( \wh{\z}^{\a+1}_{i,\ell} (\o'),\o)  \n - \n \wh{Y}^{\ell,\ell}  ( \wh{\z}^{\a+1}_{i,\ell} (\o'),\o') \big|  \\
&&  \n \le \n   \wh{\rho}_{\overset{}{\ell,\ell}} \Big( (1 \n + \n \k_\ell) \|\o \n - \n \o'\|_{0,T}
  \n + \n \phi^\o_T \big( \k_\ell   \|\o   \n  - \n  \o'  \|_{0,T} \big)  \Big)
  \n + \n \frac{1}{2 \a (\a \n + \n 1)}
  \n \le  \n  \frac{1}{  \a (\a \n + \n 1)}   ,
\eeas
and that  $ (Z^{\ell,\ell}  \n - \n \wh{Y}^{\ell,\ell}) ( \wh{\z}^{\a+1}_{i,\ell} (\o'),\o)
  \n < \n   (Z^{\ell,\ell}  \n - \n \wh{Y}^{\ell,\ell}) ( \wh{\z}^{\a+1}_{i,\ell} (\o'),\o')
   \n + \n  \frac{1}{ \a (\a+1)   }
  \n \le \n   \frac{1}{i} \n + \n  \frac{1}{\a+1 }  \n + \n  \frac{1}{ \a(\a+1)   }
   \n = \n   \frac{1}{i} \n + \n  \frac{1}{ \a  }  $,
 which shows that  $\wh{\z}^\a_{i,\ell} (\o)  \n \le \n  \wh{\z}^{\a+1}_{i,\ell} (\o') $.

   Now, we can   apply Lemma \ref{lem_sandwich}  with
 $ (\O_0,\ul{\th}, \th, \ol{\th},I, \d (\o),\e )
  \n = \n  (\O, \wh{\z}^{\a-1}_{i,\ell},   \wh{\z}^\a_{i,\ell},   \wh{\z}^{\a+1}_{i,\ell}, [0,T], \d^\ell_\a (\o), 2^{-\a} )$
 to
 \fi
 find  an   open subset $\O^\a_{i,\ell}$  of $\O$   and
  a Lipschitz  continuous random variable   $\th^\a_{i,\ell} \n : \O \to [0,T]$ such that
 \bea \label{eq:et353}
  \underset{\hP \in \cP}{\sup} \,  \hP \big( (\O^\a_{i,\ell})^c \big) \le 2^{-\a} , \q
  \wh{\z}^{\a-1}_{i,\ell} - 2^{-\a}   \n < \n  \th^\a_{i,\ell}  \n < \n  \wh{\z}^{\a+1}_{i,\ell} + 2^{-\a}
 \; \hb{ on } \; \O^\a_{i,\ell} .
 \eea

    Given $\o \n \in \n  \wh{\O}^{\a-1}_{i,\ell}   \n   \cap   \n   \wh{\O}^{\a+1}_{i,\ell} $,
 since $ \th^{\a-1}_{i,\ell}  \n - \n  2^{-\a+1}  \n < \n  \wh{\z}^\a_{i,\ell}
  \n < \n  \th^{\a+1}_{i,\ell}  \n + \n  2^{-\a-1} $, \eqref{eq:ec017}
  and an analogy to \eqref{eq:uxu017} imply that
   $t  \n := \n  \th^\a_{i,\ell} (\o)    \n  \land  \n    \wh{\z}^\a_{i,\ell} (\o)  \n  \land  \n  \wp_n  $
 and $s  \n := \n \Big( \th^\a_{i,\ell} (\o)    \n  \vee  \n    \wh{\z}^\a_{i,\ell} (\o) \Big)
  \n  \land  \n  \wp_n $  satisfy
 \beas
 s \n - \n t  & \tn \dn =  & \tn \dn \big|  \wh{\z}^\a_{i,\ell} (\o) \n  \land  \n  \wp_n (\o)
  \n - \n \th^\a_{i,\ell}  (\o)  \n  \land  \n  \wp_n  (\o)  \big|
  \n\le \n \big|  \wh{\z}^\a_{i,\ell}   (\o)   \n - \n \th^\a_{i,\ell}  (\o)   \big|
   \n  <   \n   |\th^{\a-1}_{i,\ell} (\o)  \n - \n  \th^\a_{i,\ell} (\o) |
    \n + \n  |\th^{\a+1}_{i,\ell} (\o)   \n - \n  \th^\a_{i,\ell} (\o) |
  \n  + \n  2^{-\a+1}  \n := \n  \d^\a_{i,\ell} (\o) .
 \eeas
 Set   $\phi^\a_{i,\ell} (\o) \n := \n   (1 \n + \n \k_\ell)
 \Big(  \big( \d^\a_{i,\ell} (\o) \big)^{\frac{q_1}{2}}  \n + \n \phi^\o_T \big( \d^\a_{i,\ell} (\o) \big)  \Big)   $. An application of \eqref{eq:es017} to $Z \n = \n Z^{\ell,\ell} $ shows that
 \beas
 && \hspace{-1.5cm} \big| Z^{\ell,\ell}_{\th^\a_{i,\ell} \land \wp_n} (\o)
 \n  -  \n  Z^{\ell,\ell}_{\wh{\z}^\a_{i,\ell} \land \wp_n}  (\o) \big|
  \n = \n  \big| Z^{\ell,\ell}  (t,\o)  \n - \n   Z^{\ell,\ell}   (s,\o) \big| \\
 & &  \le  \n  2 C_{\n \vr} M_0   \Big( \big(\d^\a_{i,\ell} (\o)\big)^{ \frac{q_1}{2}}
   \n \vee  \n  \big(\d^\a_{i,\ell} (\o)\big)^{q_2-\frac{q_1}{2}} \Big)
   \n + \n \wh{\rho}_{\ell,\ell} \big(  \d^\a_{i,\ell} (\o)  \big)
   \n + \n    \wh{\rho}_{\ell,\ell}  \big(   \phi^\a_{i,\ell} (\o)  \big)
    \n \vee \n  \wh{\wh{\rho}\,}_{\ell,\ell} \big(   \phi^\a_{i,\ell} (\o)  \big)   : = \xi^\a_{i,\ell} (\o) .
  \eeas

  As $Z^{\ell,\ell}$ is bounded by $M_0$,
   \eqref{eq:et349} and \eqref{eq:et353} imply that
 \bea
 \sZ_0 \n - \n   2 \ol{\e}_{m_j}  \n - \n   \ol{\e}_\ell
  & \tn \le  & \tn  \hE_{\hP_{m_j}} \Big[ Z^{\ell,\ell}_{\th^\a_{i,\ell}  \land \wp_n } \Big]
      \n +     \hE_{\hP_{m_j}} \Big[ \Big| Z^{\ell,\ell}_{\wh{\z}^\a_{i,\ell}  \land \wp_n } \n - \n
     Z^{\ell,\ell}_{\th^\a_{i,\ell}  \land \wp_n } \Big| \Big] \nonumber \\
 & \tn \le  & \tn   \hE_{\hP_{m_j}} \Big[ Z^{\ell,\ell}_{\th^\a_{i,\ell} \land \wp_n} \Big]
       \n +    \hE_{\hP_{m_j}} \Big[ \b1_{\wh{\O}^{\a-1}_{i,\ell}   \cap   \wh{\O}^{\a+1}_{i,\ell}}
   ( \xi^\a_{i,\ell} \land 2 M_0 ) \Big]
      \n + \n  2 M_0 \hP_{m_j} \Big( \big( \wh{\O}^{\a-1}_{i,\ell} \big)^c
       \cup  \big( \wh{\O}^{\a+1}_{i,\ell} \big)^c \Big) \nonumber \\
      & \tn  \le  & \tn   \hE_{\hP_{m_j}} \Big[ Z^{\ell,\ell}_{\th^\a_{i,\ell} \land \wp_n}
      \n + \n ( \xi^\a_{i,\ell} \land 2 M_0 ) \Big]   \n + \n  5 M_0 2^{-\a}  .   \label{eq:et355}
 \eea
 The random variables $\th^{\a-1}_{i,\ell}, \th^\a_{i,\ell}, \th^{\a+1}_{i,\ell}$ are  Lipschitz continuous  on $\O$,
 so is $\d^\a_{i,\ell} $. Similar   to \eqref{eq:a114},
 one  can show that $\o \to \phi^\o_T \big(\d^\a_{i,\ell} (\o) \big)$ is also a continuous random variable on $\O$,
  \if{0}
 Let $ \k^\a_{i,\ell} $ be the Lipschitz coefficient of $ \d^\a_{i,\ell} $.
 We claim that $\o \n \to \n  \phi^\o_T \big(\d^\a_{i,\ell} (\o) \big)$
 is   a continuous random variable on $\O$:
 Let  $\o  \n \in \n  \O$, $\e  \n > \n  0$ and set  $ \l^\a_{i,\ell}
  \n = \n  \l^\a_{i,\ell} (\o,\e)
   \n := \n   \frac{ \e}{3}  \n \land \n  \frac{(\phi^\o_T)^{-1} (\e/3)}{\k^\a_{i,\ell}} $.
 Given $\o'  \n \in \n  O_{\l^\a_{i,\ell}} (\o)$, we first show that
 $\phi^\o_T  \big( \d^\a_{i,\ell} (\o) \big)  \n \le \n
  \phi^{\o'}_T \n \big(   \d^\a_{i,\ell} (\o') \big)  \n + \n  \e$.
 Let $0  \n \le \n  r  \n \le \n  r'  \n \le \n  T$ with $r' \n - \n r  \n \le \n  \d^\a_{i,\ell} (\o)$.
 If $\d^\a_{i,\ell} (\o)  \n \le \n  \d^\a_{i,\ell} (\o')$, then
 \bea  \label{eq:es031}
 |\o(r') \n - \n \o(r)|  \n \le \n |\o(r') \n - \n \o'(r')| \n + \n |\o'(r') \n - \n \o'(r)|
  \n + \n    |\o'(r ) \n - \n \o(r )|
  \n \le \n  \phi^{\o'}_T \n \big(   \d^\a_{i,\ell} (\o) \big)  \n + \n  2 \|\o'  \n - \n \o \|_{0,T}
  \n < \n  \phi^{\o'}_T \n \big(   \d^\a_{i,\ell} (\o') \big)  \n + \n  \frac23 \e .
 \eea
 Otherwise if $\d^\a_{i,\ell} (\o') \n < \n  \d^\a_{i,\ell} (\o)$,
 we set $s'  \n : = \n  r'  \n \land \n  \big(r \n + \n \d^\a_{i,\ell} (\o') \big)$.
 Since \eqref{eq:ec017} shows that $r' \n - \n s'  \n = \n  r'  \n \land \n  (r \n + \n \d^\a_{i,\ell} (\o))
  \n - \n  r'  \n \land \n  (r \n + \n \d^\a_{i,\ell} (\o')) \\
  \n \le \n  \d^\a_{i,\ell} (\o)  \n - \n  \d^\a_{i,\ell} (\o')
   \n \le \n \k^\a_{i,\ell} \|\o  \n - \n \o' \|_{0,T} $
 and that $s' \n - \n r  \n = \n  r'  \n \land \n  (r \n + \n \d^\a_{i,\ell} (\o'))
  \n - \n  r'  \n \land \n  r  \n \le \n  \d^\a_{i,\ell} (\o') $, we can deduce that
  \beas
 |\o(r') \n - \n \o(r)| & \tn \le & \tn  |\o (r') \n - \n \o (s')|  \n + \n  |\o (s') \n - \n \o(r )|
  \n \le \n  \phi^\o_T  \big( \k^\a_{i,\ell} \|\o  \n - \n \o' \|_{0,T} \big)
  \n + \n  |\o(s') \n - \n \o'(s')|
  \n + \n  |\o'(s') \n - \n \o'(r)|  \n + \n  |\o'(r ) \n - \n \o(r )| \\
  & \tn \le & \tn  \phi^\o_T  \big( \k^\a_{i,\ell} \|\o  \n - \n \o' \|_{0,T} \big)
  \n + \n  \phi^{\o'}_T \n \big(   \d^\a_{i,\ell} (\o') \big)  \n + \n  2 \|\o'  \n - \n \o \|_{0,T}
  \n < \n  \phi^{\o'}_T \n \big(   \d^\a_{i,\ell} (\o') \big)  \n + \n  \e .
 \eeas
 Combining it with \eqref{eq:es031}   and taking supremum over the pair $(r,r')$ yields that
 $ \phi^\o_T  \big( \d^\a_{i,\ell} (\o) \big) \le \phi^{\o'}_T \n \big(   \d^\a_{i,\ell} (\o') \big) + \e $.

 On the other hand,   let $0 \ls \wt{r} \ls \wt{r}' \le T$ with $\wt{r}' \- \wt{r} \ls \d^\a_{i,\ell} (\o')$.
 If $\d^\a_{i,\ell} (\o') \ls \d^\a_{i,\ell} (\o)$, an analogy to \eqref{eq:es031} shows that
   \bea \label{eq:es033}
 |\o'(\wt{r}')-\o'(\wt{r})| 
 < \phi^\o_T \big(   \d^\a_{i,\ell} (\o) \big) + \frac23 \e .
 \eea
 Otherwise if $\d^\a_{i,\ell} (\o) < \d^\a_{i,\ell} (\o')$, one can deduce that
 \beas
  |\o'(\wt{r}') \n - \n \o'(\wt{r})| & \tn \dn \le & \tn \dn
   |\o'(\wt{r}') \n - \n \o(\wt{r}')| \n + \n  |\o(\wt{r}') \n - \n \o(\wt{r})|
   \n + \n   |\o(\wt{r} ) \n - \n \o'(\wt{r} )|
  \n \le \n  \phi^\o_T  \big( \d^\a_{i,\ell} (\o')   \big)   \n + \n  2 \|\o  \n - \n \o' \|_{0,T}  \\
  & \tn \dn  \le & \tn \dn   \phi^\o_T  \big( \d^\a_{i,\ell} (\o')  \n - \n  \d^\a_{i,\ell} (\o) \big)
  \n + \n  \phi^\o_T \big(   \d^\a_{i,\ell} (\o) \big)  \n + \n  2 \|\o  \n - \n \o' \|_{0,T}
  \n < \n  \phi^\o_T  \big( \k^\a_{i,\ell} \|\o  \n - \n \o' \|_{0,T}  \big)
  \n + \n  \phi^\o_T \big(   \d^\a_{i,\ell} (\o) \big)
  \n + \n  \frac23 \e
  \n \le \n  \phi^\o_T \big(   \d^\a_{i,\ell} (\o) \big)  \n + \n  \e .
 \eeas
 Combining it with \eqref{eq:es033}   and taking supremum over the pair $\big(\wt{r},\wt{r}'\big)$ yields that
 $ \phi^{\o'}_T \n \big( \d^\a_{i,\ell} (\o') \big) \le \phi^\o_T  \big(   \d^\a_{i,\ell} (\o) \big) + \e $.
 Hence $\o \to \phi^\o_T \big(\d^\a_{i,\ell} (\o) \big)$ is   a continuous random variable on $\O$,
 \fi
  which together
 with the Lipschitz  continuity of $ \d^\a_{i,\ell} $ implies that $\phi^\a_{i,\ell}$
 and thus $\xi^\a_{i,\ell}$ are also continuous random variables on $\O$.
 Analogous to \eqref{eq:a117}, we can deduce from the Lipschitz continuity of random variable
 $\th^\a_{i,\ell} \land \wp_n$ and the continuity of process  $Z$ that
 $Z^\ell_{\th^\a_{i,\ell} \land \wp_n}$ is a   continuous random variable on $\O$.

 \if{0}
  Next, we show that $Z^\ell_{\th^\a_{i,\ell} \land \wp_n}$ is a   continuous random variable:
  Let $\o, \o' \n \in \n \O$ and set
 $t  \n : = \n  \th^\a_{i,\ell} (\o )  \n \land \n  \wp_n (\o )$,
 $s  \n : = \n  \th^\a_{i,\ell} (\o' )  \n \land \n  \wp_n (\o' )$.
 We see from \eqref{eq:gd011} and \eqref{eq:es017} that
       \beas
     \big|  Z^{\ell,\ell}_s (\o  ) \n - \n  Z^{\ell,\ell}_s (\o') \big|
  & \tn \dn \le   & \tn \dn       \wh{\rho}_{\ell,\ell}  \Big( (1 \n + \n \k_\ell) \|\o  \n - \n \o' \|_{0,T}
  \n + \n \phi^\o_T \big( \k_\ell   \|\o   \n  - \n  \o'  \|_{0,T} \big)   \Big)  , \q \hb{and} \\
  \big| Z^{\ell,\ell}_t (\o) \n - \n  Z^{\ell,\ell}_s (\o) \big|
    & \tn \dn  = & \tn \dn  \big| Z^{\ell,\ell}_{t \land s} (\o) \n - \n  Z^{\ell,\ell}_{t \vee s} (\o) \big|
    \n  \le  \n  2 C_{\n \vr} M_0   \Big( |s \n - \n t|^{ \frac{q_1}{2}}
     \n \vee  \n  |s \n - \n t|^{q_2-\frac{q_1}{2}} \Big)  \n + \n \wh{\rho}_{\ell,\ell} \big(|s \n - \n t|\big)
   \n + \n    \wh{\rho}_{\ell,\ell}  \big(   \d'_{t,s} (\o)  \big)
    \n \vee \n  \wh{\wh{\rho}\,}_{\ell,\ell} \big(   \d'_{t,s} (\o)  \big)       ,
  \eeas
  where  $\d'_{t,s} (\o) \n := \n   (1 \n + \n \k_\ell)
 \Big(  |s \n - \n t|^{\frac{q_1}{2}}  \n + \n \phi^\o_T \big(|s \n - \n t|\big) \Big)       $.
 Adding them up, one can deduce from   the Lipschitz  continuity of random variable $\th^\a_{i,\ell} \land \wp_n$   that
  $Z^{\ell,\ell}_{\th^\a_{i,\ell} \land \wp_n}$ is a  continuous   random variable on $\O$.
  \fi
 As Proposition \ref{prop_Z_ultim} (2) shows that
 \bea \label{eq:et393}
 \lmtd{m \to \infty} \ol{\e}_m = 0  ,
 \eea
  letting   $j \to \infty$ in \eqref{eq:et355}, we see from the continuity of
 random variables $Z^{\ell,\ell}_{\th^\a_{i,\ell} \land \wp_n}$ and $\xi^\a_{i,\ell}$ that
 \bea   \label{eq:et357}
   \sZ_0   \le  \hE_{ \hP_* } \Big[ Z^{\ell,\ell}_{\th^\a_{i,\ell} \land \wp_n}
   \n + \n ( \xi^\a_{i,\ell} \land 2 M_0 ) \Big] \n + \n \ol{\e}_\ell \n + \n  5 M_0 2^{-\a}  .
 \eea

 \no {\bf 1c)} {\it Next, we will use the convergence of $\th^\a_{i,\ell}$ to $\wh{\z}_{i,\ell}$,
 the continuity of $Z^{\ell,\ell}$ as well as  \eqref{eq:et341}   to derive \eqref{eq:et389}. }

    Since the continuity of $Z^{\ell,\ell}  \n - \n \wh{Y}^{\ell,\ell}$ implies that
  \bhe
  \bea \label{eq:a127}
  \lmtu{\a \to \infty} \wh{\z}^\a_{i,\ell} = \wh{\z}_{i,\ell}
 : =  \inf \big\{t \in [0,T]: Z^{\ell,\ell}_t \le \wh{Y}^{\ell,\ell}_t \+ 1/i   \big\} \n \in \n \cT ,
  \eea
  \ehe
  using an analogy   to \eqref{eq:et359}  we can deduce from \eqref{eq:et353} and the Borel-Cantelli Lemma that
 \if{0}
 Then we can deduce from \eqref{eq:et353} that
 $ \lmt{\a \to \infty} \th^\a_{i,\ell} (\o) = \wh{\z}_{i,\ell} (\o) $, $ \fa
 \o \in \underset{\a = \ell + 1}{\overset{\infty}{\cup}} \, \underset{k \ge \a}{\cap}  \O^k_{i,\ell} $.
 As $ \dis  \sum^\infty_{\a = \ell + 1} \hP_* \Big( \big( \O^\a_{i,\ell} \big)^c \Big)
 \le  \sum^\infty_{\a = \ell + 1}  \underset{\hP \in \cP}{\sup} \,
  \hP \Big( \big( \O^\a_{i,\ell} \big)^c \Big) \le 2^{-\ell}  $,
 the Borel-Cantelli Lemma implies that
 $ \hP_* \left( \underset{\a = \ell + 1}{\overset{\infty}{\cup}} \,
 \underset{k \ge \a}{\cap}  \O^k_{i,\ell} \right) = 1  $. So
 \fi
 $ \lmt{\a \to \infty} \th^\a_{i,\ell}   = \wh{\z}_{i,\ell}   $, $\hP_*-$a.s. It follows that
 $ \lmt{\a \to \infty} \d^\a_{i,\ell}   = 0 $, $\hP_*-$a.s.
 and thus $ \lmt{\a \to \infty} \xi^\a_{i,\ell}   = 0 $, $\hP_*-$a.s.
 As Proposition \ref{prop_DPP_Z} shows that
  $Z^{\ell,\ell}$ is an $\bF-$adapted process bounded by $M_0$ that has all continuous paths,
  letting $\a \to \infty$ in \eqref{eq:et357} we see from   the bounded dominated convergence theorem that
  \bea \label{eq:et365}
  \sZ_0   \le  \hE_{\hP_*} \Big[ Z^{\ell,\ell}_{\wh{\z}_{i,\ell} \land \wp_n} \Big] + \ol{\e}_\ell .
  \eea

 Similar to $\z^\a_{i,\ell}$ in \eqref{eq:et361},
 $  \z_{i,\ell}   $
 is an $\bF-$stopping time satisfying $ \z_{i,\ell} \n \land \n  \wp_n \= \wh{\z}_{i,\ell} \n \land \n  \wp_n $.
 \if{0}
   And an analogy to \eqref{eq:et363} shows that
  \beas
  \z_{i,\ell} \n \land \n  \wp_n 
    \n = \n   \inf\big\{t  \n \in \n  [0,\wp_n) \n : Z^{\ell,\ell}_t  \dn \le \n  L_t  \n + \n  1/i    \big\}
   \n \land \n  \wp_n
   \n = \n  \inf\big\{t  \n \in \n  [0,\wp_n) \n : Z^{\ell,\ell}_t
    \dn \le \n  \wh{Y}^{\ell,\ell}_t  \n + \n  1/i  \big\}
   \n \land \n  \wp_n 
   \n = \n \wh{\z}_{i,\ell} \n \land \n  \wp_n .
  \eeas
  \fi
  Applying \eqref{eq:et341} with $(k, t) \n = \n \big(\ell , \z_{i,\ell} \land \wp_n \big) $
  and using \eqref{eq:et365} yield  that
 $  \sZ_0   \n  \le  \n  \hE_{\hP_*} \Big[ Z^{\ell,\ell}_{\z_{i,\ell} \land \wp_n} \Big]  \n + \n  \ol{\e}_\ell
     \n \le  \n  \hE_{\hP_*} \Big[ \sZ_{\z_{i,\ell} \land \wp_n} \Big]  \n + \n  2 \ol{\e}_\ell $.
 Since Proposition \ref{prop_Z_ultim} (3) shows that $\sZ$ is bounded by $M_0$,
 letting $\ell \n \to \n  \infty$, using the Fatou's Lemma and \eqref{eq:et393} yield that
 $ \sZ_0   \n  \le  \n \lsup{\ell \to \infty} \hE_{\hP_*} \Big[  \sZ_{\z_{i,\ell} \land \wp_n} \Big]
 \n  \le  \n  \hE_{\hP_*} \Big[ \lsup{\ell \to \infty} \sZ_{\z_{i,\ell} \land \wp_n} \Big] $.
 Similarly, letting $ i  \n \to \n  \infty$ and then letting $ n  \n \to \n  \infty$,
 we derive   \eqref{eq:et389} from Fatou's Lemma again.

 \no {\bf 2)} {\it   In the second part, we show that for any $i \ins \hN$}
 \bea \label{eq:et375}
  \ga_i    \n \le \n  \linf{\ell \to \infty} \z_{i,\ell}
  \n \le \n  \lsup{\ell \to \infty} \z_{i,\ell}   \n \le \n  \ga_{2i}  \; \;   \hb{  holds on $\O$, }
  \eea
    where $  \ga_i  \n : = \n   \inf\big\{t  \n \in \n  [0,T] \n : \sZ_t
    \n \le \n  L_t  \n + \n  1/i   \big\}  \n \land \n  T $.

 Fix $i \ins \hN$.
 Since  Proposition \ref{prop_Z_ultim} (3), (A1) and \eqref{eq:et014} show that $\sZ  \n - \n L$
 is an $\bF-$adapted process with all continuous paths, 
 $\ga_i$ is an $\bF-$stopping time that satisfies
 \bhe
 \bea \label{eq:et367}
  \ga_i \n  = \n  \lmtu{h \to \infty} \ga^h_i ,
 \eea
 \ehe
 where $  \ga^h_i  \n : = \n   \inf \big\{t  \n \in \n  [0,T] \n :
   \sZ_t  \n \le \n  L_t  \n + \n  1/i  \n + \n  1/h  \big\}  \n \land \n  T \ins \cT $.

 Fix $\o \ins  \O$ and   define $\phi^\o_U (x)  \n : = \n  \sup \big\{| U_{r'} (\o) \n - \n U_r (\o)| \n :
 r,r'  \n \in \n  [0,T] , ~ 0  \n \le \n | r' \n - \n r | \n \le \n  x  \big\}$, $\fa x  \n \in \n  [0,T]$.
 For any $\ell \ins \hN$,
  since \eqref{eq:ec017} implies that
   $ \big| U \big( ( \wp_\ell (\o)  \n + \n  2^{1-\ell} )  \n \land \n   \z_{i,\ell} (\o)   ,   \o    \big)
       \n - \n      U \big(  \tau_0 (\o)    \n \land \n  \z_{i,\ell} (\o)   , \o \big) \big|
       \n \le \n  \phi^\o_U \big(\big|  ( \wp_\ell (\o)  \n + \n  2^{1-\ell} )  \n \land \n   \z_{i,\ell} (\o)  \n - \n
      \tau_0 (\o)    \n \land \n  \z_{i,\ell} (\o)  \big|\big)
       \n \le \n  \phi^\o_U \big(\big|  ( \wp_\ell (\o)  \n + \n  2^{1-\ell} )  \n \land \n  T    \n - \n
      \tau_0 (\o)     \big|\big)  $,
  applying \eqref{eq:et340} with $(k,t) \n = \n  \big(\ell, \z_{i,\ell} (\o) \big) $ implies that
 \bea \label{eq:et371}
 \big|  Z^{\ell,\ell} \big(\z_{i,\ell} (\o)   , \o \big)  \n - \n   \sZ \big(\z_{i,\ell} (\o)   , \o \big) \big|
  \n \le  \n    \ol{\e}_\ell  \n  + \n  \phi^\o_U \big(\big|  ( \wp_\ell (\o)  \n + \n  2^{1-\ell} )  \n \land \n  T
    \n - \n   \tau_0 (\o)     \big|\big) , \q \fa \ell \ins \hN .
 \eea
As $ \lmtu{\ell \to \infty} \wp_\ell (\o)  \n = \n  \tau_0 (\o) $ by Proposition \ref{prop_wp_n} (1),
    the uniform continuity of the path $U_\cd (\o)$ implies that
    \bea \label{eq:et373}
     \lmt{\ell \to \infty} \phi^\o_U \big(\big|  ( \wp_\ell (\o)  \n + \n  2^{1-\ell} )  \n \land \n  T
    \n - \n   \tau_0 (\o)     \big|\big) \n = \n 0  .
    \eea

 To see the first inequality of \eqref{eq:et375},   we assume without loss of generality that
 $  \linf{\ell \to \infty} \z_{i,\ell} (\o)  \n < \n  T $.
 There exists a subsequence
 $\{\ell_\l \n = \n  \ell_\l (i,\o) \}_{\l \in \hN}$ of $\hN$ such that
 $   \lmt{\l \to \infty} \z_{i,\ell_\l} (\o)  \n = \n
 \linf{\ell \to \infty} \z_{i,\ell} (\o)  \n < \n  T $.

    Let $h \n \in \n  \hN$.
   Since  $\lmtd{\ell \to \infty}  \ol{\e}_\ell = 0 $ and because of  \eqref{eq:et373},
   there exists a  $\wh{\l}_h  \n = \n  \wh{\l}_h (i,\o)  \n \in \n  \hN$ such that
 for any integer $ \l  \n \ge \n  \wh{\l}_h $, one has $\z_{i,\ell_\l} (\o) \n < \n T$ and  $ \ol{\e}_{\ell_\l}
  \n  + \n  \phi^\o_U \big(\big|  ( \wp_{\ell_\l} (\o)  \n + \n  2^{1-\ell_\l} )  \n \land \n  T
  \n - \n   \tau_0 (\o)     \big|\big)  \n \le \n  1/h $.
 Given $\l  \n \in \n  \hN$ with $ \l  \n \ge \n  \wh{\l}_h $, as $\z_{i,\ell_\l} (\o) \n < \n T$,
 the set $\big\{t  \n \in \n  [0,T] \n : Z^{\ell_\l,\ell_\l}_t (\o)
  \n \le \n  L_t (\o)  \n + \n  1/i   \big\}$ is not empty. So   the continuity of the path
  $ Z^{\ell_\l,\ell_\l}_\cd (\o) \n - \n L_\cd (\o) $ implies
  that $ Z^{\ell_\l,\ell_\l}  \big( \z_{i,\ell_\l} (\o) , \o\big)
  \n \le \n  L   \big( \z_{i,\ell_\l} (\o) , \o\big)  \n + \n  1/i $.
  Applying \eqref{eq:et371} with $\ell  \n = \n  \ell_\l$ yields that
 \beas
 \sZ  \big(\z_{i,\ell_\l} (\o) , \o \big)  \n \le \n  Z^{\ell_\l,\ell_\l} \big(\z_{i,\ell_\l} (\o) , \o \big)
  \n + \n  \ol{\e}_{\ell_\l}
  \n  + \n  \phi^\o_U \big(\big|  ( \wp_{\ell_\l} (\o)  \n + \n  2^{1-\ell_\l} )  \n \land \n  T
    \n - \n   \tau_0 (\o)     \big|\big)
  \n \le \n  L \big(\z_{i,\ell_\l} (\o) , \o \big)  \n + \n  1/i    \n + \n  1/h ,
 \eeas
 which shows that $ \ga^h_i (\o) \ls \z_{i,\ell_\l} (\o)    $.
 As $\l  \n \to \n  \infty$, we obtain
 $\ga^h_i (\o)  \n \le \n  \lmt{\l \to \infty} \z_{i,\ell_\l} (\o)  \n = \n   \linf{\ell \to \infty} \z_{i,\ell} (\o)   $.
 Then letting $h \to \infty$ and using \eqref{eq:et367} yield that
 $ \ga_i (\o) \n  = \n  \lmtu{h \to \infty} \ga^h_i (\o) \le  \linf{\ell \to \infty} \z_{i,\ell} (\o) $.

 As to the third inequality of \eqref{eq:et375}, we assume without loss of generality that $ \ga_{2i} (\o) \n < \n T $,
 or equivalently,  the set $\big\{t  \n \in \n  [0,T] \n : \sZ_t (\o)
  \n \le \n  L_t (\o)  \n + \n  \frac{1}{2i}   \big\}$ is not empty.
  Then one can deduce from  the continuity of the path
  $ \sZ_\cd (\o) \n - \n L_\cd (\o) $
  that
  \bea \label{eq:et377}
   \sZ  \big( \ga_{2i} (\o) , \o\big)
  \n \le \n  L   \big( \ga_{2i} (\o) , \o\big)  \n + \n  \frac{1}{2i}   .
  \eea
 Applying \eqref{eq:et340} with $(k,t) \n = \n  \big(\ell, \ga_{2i} (\o) \big) $
 and using a similar argument to the one that leads  to \eqref{eq:et371} yield that
 \bea  \label{eq:et379}
   \big|  Z^{\ell,\ell} \big(\ga_{2i} (\o)   , \o \big)  \n - \n   \sZ \big(\ga_{2i} (\o)   , \o \big) \big|
  \n \le  \n    \ol{\e}_\ell  \n  + \n  \phi^\o_U \big(\big|  ( \wp_\ell (\o)  \n + \n  2^{1-\ell} )  \n \land \n  T
    \n - \n   \tau_0 (\o)     \big|\big)  .
    \eea
 For any   $\ell \n \in \n  \hN$ such that
$ \ol{\e}_\ell   \n  + \n  \phi^\o_U \big(\big|  ( \wp_\ell (\o)  \n + \n  2^{1-\ell} )  \n \land \n  T
  \n - \n   \tau_0 (\o)     \big|\big)  \n \le \n  \frac{1}{2i} $,  \eqref{eq:et377} and \eqref{eq:et379}
  imply that
 $ Z^{\ell,\ell}   \big(\ga_{2i} (\o) , \o \big) \ls \sZ \big(\ga_{2i} (\o) , \o \big) \+  \frac{1}{2i}
 \ls  L \big(\ga_{2i} (\o) , \o \big) \+ 1/i $,
 which shows that $\z_{i,\ell} (\o) \n  \le \n  \ga_{2i} (\o)$.
 As $\ell \to \infty$, we obtain $ \lsup{\ell \to \infty} \z_{i,\ell} (\o) \n  \le \n  \ga_{2i} (\o) $.

  \no {\bf 3)} {\it Finally,  we show that
 $ \lsup{n \to \infty} \,     \lsup{i \to \infty} \,   \lsup{\ell \to \infty}
  \sZ (\z_{i,\ell} (\o) \land \wp_n (\o) , \o )
   \n  = \n  \sZ (\ga_* (\o)   ,\o)$, $\fa \o  \n \in \n  \O$.
   The conclusion thus follows. }

 Let $ 1 \ls n  \< i   $ and $\o \ins \O$.
 We set $t_\ell \n = \n  t_\ell (n,i,\o )  \n := \n  \big( \z_{i,\ell}    \n  \land \n  \wp_n \big) (\o) $,
 $\fa \ell \> i$.
 Let $ \big\{t_{\wt{\ell}} \big\}_{\wt{\ell} \in \hN}$ be the subsequence of
 $\{t_\ell\}^\infty_{\ell = i+1}$   such that
 $ 
  \lsup{\ell \to \infty}  \sZ (t_\ell , \o )
  \n = \n  \lmt{\wt{\ell} \to \infty}  \sZ (t_{\wt{\ell}} , \o ) $.
 The sequence $\big\{t_{\wt{\ell}} \big\}_{\wt{\ell} \in \hN}$ in turn has a
 convergent subsequence $\big\{ t_{\wt{\ell}'} \big\}_{\wt{\ell}' \in \hN}$ with limit $\ft  \n \in \n  [0,\wp_n (\o)]$.
 The continuity of   path $ \sZ_\cd (\o)$ shows that
 $ \sZ (\ft, \o )  \n = \n  \lmt{\wt{\ell}' \to \infty}  \sZ (t_{\wt{\ell}'} , \o )
   \n = \n  \lsup{\ell \to \infty}  \sZ (t_\ell , \o )$.
 Also, \eqref{eq:et375} implies that
 $ \big( \ga_i  \n \land \n  \wp_n \big) (\o)
  \n \le \n    \linf{\ell \to \infty} \big( \z_{i,\ell}   \n \land \n  \wp_n \big) (\o)
  \n = \n  \linf{\ell \to \infty} t_\ell   \n \le \n
 \ft  \n = \n  \lmt{\wt{\ell}' \to \infty} t_{\wt{\ell}'}  \n \le \n  \lsup{\ell \to \infty} t_\ell
  \n = \n   \lsup{\ell \to \infty} \big( \z_{i,\ell}   \n \land \n  \wp_n \big) (\o)
   \n \le \n  \big( \ga_{2i}    \n \land \n  \wp_n \big) (\o)$. Hence
\bea \label{eq:et385}
    \underset{t \in \cJ_{n,i} (\o)}{\inf}   \sZ(t,\o) \le  \sZ (\ft, \o )  =
  \lsup{\ell \to \infty} \sZ \big( \z_{i,\ell} (\o)   \land \wp_n (\o) , \o \big)
  \le \underset{t \in \cJ_{n,i} (\o)}{\sup}   \sZ(t,\o)   ,
\eea
where $\cJ_{n,i} (\o) \n := \n  \big[ (\ga_i    \n \land \n  \wp_n) (\o), (\ga_{2i}    \n \land \n  \wp_n ) (\o) \big]$.

 An analogy to \eqref{eq:et367} shows that
 $  \ga_\sharp (\o) \n : = \n   \inf\big\{t  \n \in \n  [0,T] \n : \sZ_t (\o) \n \le \n  L_t (\o) \big\}  \n \land \n  T
  \n  = \n  \lmtu{i \to \infty} \ga_i (\o)      $.
 Since $\wh{\sY}_t (\o) \n = \n \sY_t (\o) \n = \n  L_t (\o) $ over the interval $ \big[0, \tau_0 (\o)\big)
  \n \supset \n  \big[0, \wp_n (\o)\big) $ by Proposition \ref{prop_wp_n} (1), we can deduce from \eqref{eq:et335} that
  \bea
  \hspace{-3mm} \lmtu{i \to \infty}   (\ga_i    \n \land \n  \wp_n) (\o)
   & \dn \dn =  & \dn \dn  (\ga_\sharp     \n \land \n  \wp_n) (\o)
    \n = \n  \inf\{t  \n \in \n  [0,\wp_n  (\o)) \n : \sZ_t  (\o)  \n \le \n  L_t  (\o)  \}  \n \land \n  \wp_n  (\o)
   \n = \n  \inf \big\{t  \n \in \n  [0,\wp_n  (\o)) \n :
   \sZ_t  (\o)  \n \le \n  \wh{\sY}_t  (\o)  \big\}  \n \land \n  \wp_n  (\o) \nonumber \\
    & \dn \dn  =  & \dn \dn
    \inf \big\{t  \n \in \n  [0, \wp_n  (\o)) \n : \sZ_t  (\o)  \n = \n  \wh{\sY}_t  (\o)  \big\}  \n \land \n  \wp_n  (\o)
   \n = \n ( \ga_*    \n \land \n   \wp_n ) (\o)  .   \label{eq:et381}
  \eea
 It follows from the continuity of   path $\sZ_\cd (\o)$   that
 \bhe
 \bea \label{eq:et383}
  \lmt{i \to \infty} \, \underset{t \in \cJ_{n,i} (\o)}{\inf}   \sZ(t,\o)
  \n = \n   \lmt{i \to \infty}  \,  \underset{t \in \cJ_{n,i} (\o)}{\sup}   \sZ(t,\o)
  \n = \n   \sZ \big(  \ga_* (\o)  \n \land \n   \wp_n (\o) , \o \big)  .
  \eea
  \ehe

  Then letting $i  \n \to \n  \infty$ in \eqref{eq:et385}  yields that
  \bea \label{eq:et387}
     \lsup{i \to \infty} \,   \lsup{\ell \to \infty} \sZ (\z_{i,\ell} (\o)  \n \land \n  \wp_n (\o) , \o )
    \n = \n  \lmt{i \to \infty} \,   \lsup{\ell \to \infty} \sZ (\z_{i,\ell} (\o)  \n \land \n  \wp_n (\o) , \o )
    \n = \n  \sZ (\ga_* (\o)  \n \land \n   \wp_n (\o) ,\o)  .
    \eea
 Since Proposition \ref{prop_wp_n} (1) and Proposition \ref{prop_Z_ultim} (4) imply that
  $ \lmt{n \to \infty} (\ga_*    \n \land \n   \wp_n ) (\o)
   = (\ga_*    \n \land \n   \tau_0 ) (\o)   = \ga_*  (\o)  $,
   letting $n \to \infty$ in \eqref{eq:et387}, we see from the continuity of   path $\sZ_\cd (\o)$ again that
 \beas
  \lsup{n \to \infty} \,     \lsup{i \to \infty} \,   \lsup{\ell \to \infty}
  \sZ (\z_{i,\ell} (\o) \land \wp_n (\o) , \o )
   = \lmt{n \to \infty} \, \lsup{i \to \infty} \,   \lsup{\ell \to \infty}
   \sZ (\z_{i,\ell} (\o) \land \wp_n (\o) , \o )
   = \sZ (\ga_* (\o)   ,\o) , \q \fa \o \n \in \n  \O.
\eeas
  Putting this back into \eqref{eq:et389} and using Proposition \ref{prop_Z_ultim} (3) yield that
 $  \underset{(\hP, \ga) \in \cP \times \cT}{\sup}  \,
 \hE_\hP [\sY_{\ga \land \tau_0} ] \n = \n  \underset{(\hP, \ga) \in \cP \times \cT}{\sup}  \,
 \hE_\hP \big[ \wh{\sY}_{\ga} \big] \n \le  \n  \sZ_0
    \n \le  \n  \hE_{\hP_*} \big[ \sZ_{\ga_*} \big]$.
   Since the continuity of $\sZ$ and the right-continuity of $\wh{\sY}$ imply that
   $ \sZ_{\ga_*} (\o)  \n = \n   \wh{\sY}_{\ga_*} (\o) = \sY_{\ga_* \land \tau_0} (\o) $, $\fa \o  \n \in \n  \O$,
   one can further deduce \eqref{eq:eu018} and thus \eqref{eq:RDOSRT}. \qed

\appendix
\renewcommand{\thesection}{A}
\refstepcounter{section}
\makeatletter
\renewcommand{\theequation}{\thesection.\@arabic\c@equation}
\makeatother

\section{Appendix}

\subsection{Technical Lemmata}

 \if{0}

 \begin{lemm} \label{lem_stopping_time}
 Let $t \in [0,T]$, $\d \in \hR$ and let $X$ be an $\bF^t-$adapted process.
  If all paths of $X$ are continuous, then
 $\tau_\d  \n :=  \n  \inf\big\{s  \n \in \n  [t,T] \n : X_s  \n \le  \n   \d    \big\}
 \land T $ is an $\bF^t-$stopping time.

 \end{lemm}

  \no {\bf Proof:} Let $t \n \in \n  [0,T]$, $\d  \n \in \n  \hR$ and let $X$ be an $\bF^t-$adapted process.
 If all paths of $X$ are continuous,
 for any $ s  \n \in \n  [t,T) $   we can deduce that
             \beas
     \{ \tau_\d \n > \n  s\} & \tn =& \tn  \{\o  \n \in \n  \O^t \n :
      X_r (\o)  \n > \n   \d,~ \fa r  \n \in \n  [t,s] \}
      \n = \n  \underset{n    \in    \hN}{\cup} \{\o  \n \in \n  \O^t \n :
      X_r (\o)  \n \ge \n   \d  \n + \n  1/n,~ \fa r  \n \in \n  [t,s] \} \\
     & \tn =& \tn  \underset{n    \in    \hN}{\cup} \{\o  \n \in \n  \O^t \n :
     X_r (\o)  \n \ge \n   \d  \n + \n  1/n,~ \fa r  \n \in \n  \hQ_{t,s} \}
      \n = \n  \underset{n    \in    \hN}{\cup} \, \underset{r \in \hQ_{t,s}}{\cap}
      \{\o  \n \in \n  \O^t \n :
      X_r (\o)  \n \ge  \n  \d  \n + \n  1/n  \}  \n \in \n  \cF^t_s ,
     \eeas
     where $\hQ_{t,s} := \big(  [t,s] \cap \hQ \big) \cup \{t, s\} $. So
   $\tau_\d$ is an $\bF^t-$stopping time.     \qed

 \fi

 \begin{lemm} \label{lem_shift_stopping_time}
 Given $t  \n \in \n  [0,T]$, let $\tau  \n \in \n  \cT^t$ and $(s,\o)  \n \in \n  [t,T]  \n \times \n  \O^t$.
 If $\tau (\o)  \n \le \n  s$, then
 $ \tau (\o  \n \otimes_s \n  \O^s )  \n \equiv \n  \tau (\o) $; if $\tau (\o)  \n \ge \n  s$
 \(resp.  $ \n > \n  s$\), then
 $ \tau (\o  \n \otimes_s \n  \wt{\o})  \n \ge \n  s $ \(resp.  $ \n > \n  s$\), $\fa \wt{\o}  \n \in \n  \O^s$
 and thus $\tau^{s,\o}  \n \in \n  \cT^s$
 by Proposition \ref{prop_shift0} \(2\).
 \end{lemm}

  \no {\bf Proof: } Let  $t  \n \in \n  [0,T]$,   $\tau  \n \in \n  \cT^t$ and $(s,\o)  \n \in \n  [t,T]  \n \times \n  \O^t$. When $\wh{s}  \n : = \n  \tau (\o)  \n \le \n  s$,
 since $\o  \n \in \n  A  \n : = \n  \{  \tau   \n = \n  \wh{s}\}  \n \in \n  \cF^t_{\wh{s}}
  \n \subset \n  \cF^t_s$, Lemma \ref{lem_element} shows that
  $ \o  \n \otimes_s \n  \O^s  \n \subset \n  A $, i.e. $\tau (\o  \n \otimes_s \n  \O^s )  \n \equiv \n
  \wh{s}  \n   = \n  \tau (\o) $.

  On the other hand, when $\tau (\o)  \n \ge \n  s$ \(resp.  $ \n > \n  s$\), as $ \o    \n \in \n
  A'  \n : = \n  \{  \wp    \n  \ge  \n  s \} $  (resp.  $ \{  \wp    \n  >  \n  s \} \)
   \n \in \n     \cF^t_s $, applying Lemma \ref{lem_element} again yields that
   $ \o  \n \otimes_s \n  \O^s   \n \in \n  A'$.
   So $  \tau (\o  \n \otimes_s \n  \wt{\o})  \n  \ge  \n  s $ \(resp.  $ \n > \n  s$\), $\fa \wt{\o}  \n \in \n  \O^s$. \qed

        \begin{lemm} \label{lem_Y_diff}
     Assume \(P2\).  Let $(Y,\wp)  \n \in \n  \fS$ and $(t,\o) \n \in \n  [0,T]  \n \times \n  \O $.
 It holds for any $\o'  \n \in \n  \O$, $  \hP   \n \in \n     \cP_t $
  and  $ \ga    \n \in \n  \cT^t    $    that
  \beas
   \hE_\hP \Big[ \big| \wh{Y}^{t,\o}_\ga   \n - \n   \wh{Y}^{t,\o'}_\ga \big| \Big]
   \n  \le  \n     \wrY \Big( (1 \n + \n \k_\wp) \|\o \n - \n \o'\|_{0,t}
  \n + \n  \underset{r \in [ t_1 , t_2 ]}{\sup} \big|    \o    (r  )   \n - \n    \o    ( t_1 ) \big|  \Big)
   \n  \le  \n     \wrY \Big( (1 \n + \n \k_\wp) \|\o \n - \n \o'\|_{0,t}
  \n + \n \phi^\o_t \big( \k_\wp   \|\o   \n  - \n  \o'  \|_{0,t} \big)   \Big) ,
 \eeas
 where $t_1  \n : = \n  \wp (\o)  \n \land \n  \wp (\o')  \n \land \n  t $
 and $t_2  \n : = \n  \big(\wp (\o)  \n \vee \n  \wp (\o')\big)  \n \land \n  t $.

 \end{lemm}

    \no {\bf Proof: 1)}
 Fix     $\o' \n \in \n  \O$.
 We   set $t_1  \n : = \n  \wp (\o)    \land    \wp (\o')    \land    t $,
 $t_2  \n : = \n  \big(\wp (\o)    \vee    \wp (\o')\big)    \land    t $
 and $ \d \n : = \n   (1 \+ \k_\wp) \|\o   -   \o'\|_{0,t}
  \n + \n  \underset{r \in [ t_1 , t_2 ]}{\sup} \big|    \o    (r  )   \n - \n    \o    ( t_1 ) \big|   $.
 Fix also $  \hP   \n \in \n     \cP_t $ and  $ \ga    \n \in \n  \cT^t    $.
 Let $\wt{\o}   \n \in \n  \O^t $. One has
 \beas
 \Big| \wh{Y}^{t,\o} \big(\ga (\wt{\o}), \wt{\o} \big) \n - \n  \wh{Y}^{t,\o'} \big(\ga (\wt{\o}), \wt{\o} \big)  \Big|
    \n = \n     \Big|  \wh{Y} \big(\ga (\wt{\o}), \o  \n \otimes_t \n  \wt{\o} \big)
  \n - \n  \wh{Y} \big(\ga (\wt{\o}), \o'  \n \otimes_t \n  \wt{\o} \big)  \Big|
    \n  =  \n     \Big|   Y \big(s_1 (\wt{\o}), \o  \n \otimes_t \n  \wt{\o} \big)
  \n - \n  Y \big(s_2 (\wt{\o}) , \o'  \n \otimes_t \n  \wt{\o} \big)  \Big|  ,
 \eeas
 where   $  s_1(\wt{\o}) := \ga (\wt{\o}) \ld   \wp(\o \otimes_t \wt{\o}) \ld \wp(\o' \otimes_t \wt{\o}) $ and
   $  s_2 (\wt{\o}) := \ga (\wt{\o}) \ld \big( \wp(\o \otimes_t \wt{\o}) \ve \wp(\o' \otimes_t \wt{\o}) \big)   $.
   Since \eqref{eq:ec017} implies that
  \bea \label{eq:ga311}
    s_2 (\wt{\o}) -s_1 (\wt{\o})   \le \big| \wp(\o \otimes_t \wt{\o}) - \wp(\o' \n  \otimes_t \wt{\o}) \big|
  \le \k_\wp   \|\o \otimes_t \wt{\o} - \o' \n \otimes_t \wt{\o}\|_{0,T}
  = \k_\wp   \|\o   - \o'  \|_{0,t} < \d ,
  \eea
   one can deduce from \eqref{eq:aa211} that
 \bea
  && \hspace{-1.2cm} \Big| \wh{Y}^{t,\o} \big(\ga (\wt{\o}), \wt{\o} \big)
 - \wh{Y}^{t,\o'} \big(\ga (\wt{\o}), \wt{\o} \big) \Big|
     \ls   \rY \Big( \big( s_2 (\wt{\o}) \-s_1(\wt{\o}) \big)
 + \underset{r \in [0,T]}{\sup} \Big| (\o \otimes_t \wt{\o}) \big(r \ld s_1 (\wt{\o})\big)
  - (\o' \oti_t \wt{\o}) \big(r \ld s_2 (\wt{\o}) \big) \Big|   \Big)  \nonumber \\
  &&  \hspace{-0.6cm}   \le \n    \rY \n \Big(  \k_\wp   \|\o \n  - \n  \o'  \|_{0,t}
  \n + \n \cI (\wt{\o})
  \n + \n  \underset{r \in [0,T]}{\sup} \Big|  (\o  \n \otimes_t \n  \wt{\o}) \big(r  \n \land \n  s_2 (\wt{\o}) \big)
  \n - \n  (\o'  \n \otimes_t \n  \wt{\o}) \big(r  \n \land \n  s_2 (\wt{\o}) \big) \Big|   \Big)
     \ls    \rY \n  \Big( ( 1 \+ \k_\wp )  \|\o \n  - \n  \o'  \|_{0,t}
  \n + \n \cI (\wt{\o})   \Big)   , \qq  \label{eq:ga314}
 \eea
 where
 $ \cI (\wt{\o}) \df \underset{r \in [0,T]}{\sup}  \big|   (\o \oti_t \wt{\o})   \big( r \ld s_1 (\wt{\o}) \big)
 \-  (\o  \oti_t \wt{\o})  \big(r \ld s_2 (\wt{\o}) \big)  \big|
 \= \underset{r \in [s_1 (\wt{\o})   ,  s_2 (\wt{\o}) ]}{\sup} \big|
   (\o \oti_t \wt{\o})  (r   )  \-  (\o \oti_t \wt{\o})  \big( s_1 (\wt{\o})    \big) \big| $.

 \no {\bf 2)} Next, we discuss by three cases on values of $\wp (\o)$  and  $ \wp (\o')$:

 \no (i) When $ \wp (\o) \ld \wp (\o') \ge t    $,   Lemma \ref{lem_shift_stopping_time} shows
 that $  \wp^{t,\o} $ and $  \wp^{t,\o'} $ belong  to $ \cT^t$,
 so does $\z  \n : = \n  \ga  \n \land \n  \wp^{t,\o}  \n \land \n  \wp^{t,\o'} $.
 For any $ \wt{\o} \ins \O^t $, as  $ s_1 (\wt{\o})
 \= \z (\wt{\o}) \gs t $, \eqref{eq:ga311} implies that
 $   \cI (\o)  \n = \n   \underset{r \in [s_1 (\wt{\o}),s_2 (\wt{\o}) ]}{\sup} \big|     \wt{\o}   (r   )   \n  -  \n  \wt{\o}   (s_1 (\wt{\o}) ) \big|
 \n \le \n 
  \underset{r \in [\z (\wt{\o}),(\z (\wt{\o})+\d) \land T  ]}{\sup}
  \big|    B^t_r ( \wt{\o}   )   \n  -  \n   B^t_\z (  \wt{\o} )   \big|  $.
  Putting it back into \eqref{eq:ga314}
 and   taking expectation $\hE_\hP [~]$, we see from \eqref{eq:ex015} that
 \bea \label{eq:a103}
   \hE_\hP \Big[ \Big| \wh{Y}^{t,\o}_\ga   \n - \n     \wh{Y}^{t,\o'}_\ga \Big| \Big]   \ls
   \hE_\hP \bigg[ \rY \Big( \d  \n + \n   \underset{r \in [\z  ,(\z  +\d) \land T  ]}{\sup}
  \big|    B^t_r     \n  -  \n   B^t_\z     \big|  \Big) \bigg]    \n \le  \n   \wrY (\d)    .
 \eea

     \no \(ii\) When $   \wp (\o) \ld \wp (\o')  \n < \n  t  \n \le \n  \wp (\o) \ve \wp (\o') $,
  let $(\ul{\o}, \ol{\o})$ be a possible permutation of $(\o  , \o')$ such that
  $ \wp(\ul{\o}) \=  \wp (\o) \ld \wp (\o') \< t $ and
  $ \wp(\ol{\o}) \=  \wp (\o) \ve \wp (\o') \gs t $. By Lemma \ref{lem_shift_stopping_time},
  $ \wp(\ul{\o} \otimes_t \O^t )  \n \equiv \n  \wp(\ul{\o} )    $
  and $\wp (\ol{\o}  \n \otimes_t \n  \O^t  )  \sb [ t , T ] $. For any $ \wt{\o} \ins \O^t $,
  one has  $ s_1 (\wt{\o}) \= \ga (\wt{\o})  \n \land \n  \wp(\ul{\o}  \n \otimes_t \n  \wt{\o})
  \n \land \n  \wp(\ol{\o}  \n \otimes_t \n  \wt{\o}) \=  \wp(\ul{\o} )  \=  t_1 \< t  $
  and $ s_2 (\wt{\o}) \= \ga (\wt{\o})  \n \land \n \big(  \wp(\ul{\o}  \n \otimes_t \n  \wt{\o})
  \n \vee \n  \wp(\ol{\o}  \n \otimes_t \n  \wt{\o}) \big)
  \= \ga (\wt{\o})  \n \land \n \wp(\ol{\o}  \n \otimes_t \n  \wt{\o}) \gs t  $.
  Since $s_2 (\wt{\o}) \<  s_1 (\wt{\o}) \n + \n  \d  \n < \n  t \n + \n \d$ by \eqref{eq:ga311}
  and since $t_2 \= \wp(\ol{\o}) \ld t \= t$,   we can deduce that
  \beas
  \hspace{-3mm}
  \cI (\wt{\o}) & \tn \dn =  & \tn \dn \Big( \,   \underset{r \in [s_1 (\wt{\o}), t ]}{\sup} \big|    \o    (r  )   \n - \n    \o    \big(s_1 (\wt{\o}) \big) \big| \Big)
  \n \vee \n  \Big(  \, \underset{r \in [t,s_2 (\wt{\o})]}{\sup} \big|     \wt{\o}   (r   ) \n + \n
  \o(t)  \n - \n     \o    \big(s_1 (\wt{\o}) \big) \big|  \Big) \\
   & \tn  \dn  \le  & \tn  \dn   \Big( \,  \underset{r \in [t_1, t ]}{\sup}
  \big|    \o    (r  )   \n - \n    \o    (t_1) \big| \Big)
  \n \vee \n  \Big(  \big|\o(t)  \n - \n     \o    (t_1)\big|
  \n + \n \underset{r \in [t,s_2 (\wt{\o})]}{\sup} \big|     \wt{\o}   (r   ) \n - \n \wt{\o} (t)    \big|  \Big)
      \ls   \underset{r \in [t_1 , t_2  ]}{\sup}
   \big|    \o    (r  )   \n - \n    \o    ( t_1  ) \big|
   \n + \n  \underset{   r \in [t, (t+\d) \land T ] }{\sup}
 \big| B^t_r ( \wt{\o}   )  \n - \n  B^t_t ( \wt{\o}   ) \big|     .
   \eeas
 An analogy to \eqref{eq:a103} shows that
 \bea \label{eq:a105}
   \hE_\hP \Big[ \Big| \wh{Y}^{t,\o}_\ga   \n - \n     \wh{Y}^{t,\o'}_\ga \Big| \Big]   \ls
   \hE_\hP \bigg[ \rY \Big( \d  \n + \n   \underset{r \in [t  ,(t  +\d) \land T  ]}{\sup}
  \big|    B^t_r     \n  -  \n   B^t_t      \big|  \Big) \bigg]    \n \le  \n   \wrY (\d)    .
 \eea

     \no \(iii\) When $  \wp (\o) \ve \wp (\o') \n < \n  t   $,
     we see from Lemma \ref{lem_shift_stopping_time} again that
  $  \wp(\o  \n \otimes_t \n  \O^t) \n \equiv  \n \wp(\o) \n < \n  t $ and
  $  \wp(\o'  \n \otimes_t \n  \O^t) \n \equiv  \n \wp(\o') \n < \n  t $.
  For any $ \wt{\o} \ins \O^t $, as $\ga(\wt{\o})  \n \ge \n  t$, one has
 $s_1 (\wt{\o})  \n  = \n   \wp(\o)  \n \land  \n  \wp(\o') \=   t_1 \< t $
   and $ s_2 (\wt{\o}) \=   \wp(\o)  \n \vee   \n  \wp(\o')  \n  = \n  t_2 \< t $.
   It follows that
    $\cI (\o) 
  \n = \n  \underset{r \in [ t_1 , t_2 ]}{\sup} \big|    \o    (r  )   \n - \n    \o    ( t_1  ) \big|   $,
  then \eqref{eq:a105} still holds for this case.

  Therefore, we have proved the first inequality of the lemma. 
  Since  $ t_2 \- t_1 \=    \big|\wp(\o)    \n  \land  \n    t \- \wp(\o')    \n  \land  \n    t \big|
  \n \le \n |\wp(\o)     \- \wp(\o')     |
  \n \le   \n   \k_\wp   \|\o   \- \o'  \|_{0,t}  $ by \eqref{eq:ec017},
  the second inequality easily follows. \qed

\begin{lemm} \label{lem_P_supermg}
  Assume \(P2\)$-$\(P4\) and let $(Y,\wp) \n  \in \n  \fS$. Given  $\hP  \n \in \n  \cP$,
  $Z$ is a $\hP-$supermartingale  and    $\hE_\hP [Z_{\tau} ]  \n \ge \n  \hE_\hP [Z_{\ga} ]$ holds
  for any $\tau, \ga  \n \in \n  \cT$ with $\tau  \n \le \n  \ga $, \pas

\end{lemm}

 \no {\bf Proof:} Fix $(Y,\wp) \n  \in \n  \fS$ and $\hP  \n \in \n  \cP$.

  \no {\bf 1)} Let $t  \n \in \n  [0,T]$ and $\ga  \n \in \n  \cT$.
   Proposition \ref{prop_Z_conti_in_o} and \eqref{eq:ef221} show that $Z_\ga$ is an $\cF_T-$measurable
 bounded random variable.
  By  Proposition \ref{prop_shift1}, we can find   a $\hP-$null set $\cN$ such that
  $   \hE_\hP [ Z_\ga | \cF_t ] (\o)  \n = \n   \hE_{\hP^{t, \o}} \big[ (Z_\ga )^{t, \o}  \big] $,
   $ \fa \o  \n \in \n   \cN^c $.
  Also,   (P3) shows that for
  some extension $(\O,\cF',\hP')$ of $(\O,\cF_T,\hP)$
  and some $\O' \n \in \n  \cF' $ with $\hP'(\O')  \n = \n  1$,
    $\hP^{ t,   \o}  \n \in \n  \cP_t  $ for any $\o  \n \in \n  \O'$.
  Then   Proposition \ref{prop_Z_martingale} implies that
  $ \hE_\hP [ Z_\ga | \cF_t ] (\o)  \n = \n   \hE_{\hP^{t, \o}} \big[ (Z_\ga )^{t, \o}   \big]
        \n \le \n  \ol{\sE}_{\n t} \big[  Z_\ga \big] (\o)
        \n \le \n  Z_{\ga \land t} (\o) $, $  \fa \o  \n \in \n  \O' \cap    \cN^c   $.
  Using similar arguments that lead to \eqref{eq:ef121}, we can obtain that
     \bea \label{eq:et421}
     \hE_\hP [ Z_\ga | \cF_t ] \le Z_{\ga \land t} , \q  \pas
     \eea

   \no {\bf 2)} Let  $\tau , \ga   \n \in \n  \cT$ with   $\tau  \n \le \n  \ga $, \pas ~
 Also, let $n \n \in \n  \hN$ and $i \n = \n 1,\cds \n ,2^n$. We set $ t^n_i  \n : = \n   i 2^{-n}T $ and
 $A^n_i \n : = \n \{ t^n_{i-1}  \n < \n  \tau  \n \le \n  t^n_i  \} \n \in \n \cF_{t^n_i} $ with $t^n_0 \n := \n 0$.
 Applying \eqref{eq:et421} with $t \n = \n t^n_i$ yields that
  $   \hE_\hP [ Z_\ga | \cF_{t^n_i} ]  \n \le \n  Z_{\ga \land t^n_i}$, \pas ~
   Multiplying $\b1_{A^n_i}$ and taking summation over $i \in \{1,\cds,2^n\}$, we obtain
 $   \hE_\hP [ Z_\ga | \cF_{\tau_n} ]  \n \le \n  Z_{\ga \land \tau_n}$, \pas,
 where $\dis \tau_n  \n : = \n     \sum^{2^n}_{i=1}
 \b1_{ A^n_i } t^n_i   \n \in \n  \cT   $. Then taking the expectation
 $\hE_\hP[~]$ yields that $ \hE_\hP [ Z_\ga ]  \n \le \n  \hE_\hP [ Z_{\ga \land \tau_n} ] $.
 Since $\lmtd{n \to \infty} \tau_n  \n = \n  \tau$
 and since Proposition \ref{prop_DPP_Z} shows that  $Z$ is a bounded
 process with  all continuous paths, an application of   the bounded convergence theorem leads to that
 $ \hE_\hP [ Z_\ga ]  \n \le \n  \hE_\hP [ Z_{\ga \land \tau } ]  \n = \n  \hE_\hP [ Z_\tau  ] $.   \qed

 We need  the following extension of  Lemma 4.5 of    \cite{ETZ_2014} to prove  Theorem \ref{thm_cst}
 and Theorem \ref{thm_RDOSRT}.

 \begin{lemm} \label{lem_sandwich}
 Assume \(P1\).
 Let $\O_0 \subset \O$ and let $ \ul{\th} $, $\th$, $\ol{\th}$ be three real-valued random variables on $\O$
 taking values in a compact interval $I  \n \subset \n  \hR$ with  length $|I |  \n > \n  0$.
 If for any $\o \in \O_0$  there exists a $  \d(\o) > 0$ such that
     \bea \label{eq:gf011}
     \ul{\th} (\o' )  \n \le \n    \th (\o )  \n \le \n
       \ol{\th} (\o' ), \q \fa  \o'  \n \in \n  \ol{O}_{\d(\o)} (\o)
       \n = \n  \big\{\o'  \n \in \n  \O \n : \|\o' \n - \n \o\|_{0,T}  \n \le \n  \d(\o) \big\} ,
       \eea
 then for any $\e  \n > \n  0$ one can find an open subset $\wh{\O} $ of $\O$ and
 a Lipschitz  continuous random variable $\wh{\th}  \n : \O  \n \to \n  I$
  such that $ \underset{\hP \in \cP}{\sup} \,  \hP \big( \wh{\O}^c \big)   \n \le \n  \e $ and that
  $ \ul{\th}  \n - \n  \e   \n < \n  \wh{\th}  \n < \n  \ol{\th}  \n + \n  \e $ on  $ \wh{\O}    \cap    \O_0 $.

 \end{lemm}

  \no {\bf Proof:}
   Since the canonical space $\O $ is a separable complete metric space
   and thus  Lindel\"of,    there exists    a sequence $\{\o_j\}_{j \in \hN}$ of $\O$ such that
   $ \dis \underset{j \in \hN}{\cup} O_j  \n = \n  \O $
   with $ O_j  \n := \n  O_{ \frac12 \d(\o_j)  }(\o_j)  \n = \n
  \{\o  \n \in \n  \O  \n : \|\o  \n - \n  \o_j\|_{0,T}  \n < \n  \frac12 \d(\o_j)  \}$.

 Let $n  \n \in \n  \hN$ with $ n^2  \n > \n  |I|^{-1}$.
 By \eqref{eq:bb237},    $\O_n  \n : = \n  \underset{j = 1}{\overset{n}{\cup}} O_j $ is an open subset of $\O$.
 For   $j  \n = \n  1, \cds \n , n$, we define function $f_{n,j}  \n : [0,\infty)  \n \to \n  [0, 1]$ by:
 $f_{n,j} (x)  \n := \n  1$ for $x  \n \in \n  \big[0, \frac12 \d(\o_j)\big]$,
 $\dis f_{n,j}(x)  \n := \n   n^{-2}|I|^{-1}$  for $x  \n \ge \n  \d(\o_j)$,
 and $f_{n,j}$ is linear in $\big[\frac12\d(\o_j)  , \d(\o_j)\big]$. Clearly,
  $g_{n,j}(\o)  \n := \n  f_{n,j}(\|\o  \n - \n  \o_j\|_{0,T})$,
  $\o  \n \in \n  \O$ is a Lipschitz  continuous random variable on $\O$ with coefficient $< 2/\d (\o_j)$.
  It follows that   $\fg_n  \n := \n    \sum^n_{j=1} g_{n,j} $ is a
  Lipschitz  continuous random variable on $\O$ with values in $  \big[ n^{-1}|I|^{-1},n\big] $
  and that $\sum^n_{j=1}  \th(\o_j ) g_{n,j}$ is a
  Lipschitz  continuous random variable on $\O$ whose absolute values  $\le \sum^n_{j=1}  |\th(\o_j )| $.
  Then one can deduce that
 \beas
  \th_n(\o)   :=   \frac{1}{\fg_n(\o)} \sum^n_{j=1}  \th(\o_j )g_{n,j}(\o) , \q \fa \o   \in   \O
 \eeas
 defines another Lipschitz   continuous random variable on $\O$ with values in $I$.

   Given $\o \n \in \n  \O_n    \cap    \O_0$, as $\o $ belongs   $O_j$ for some $j \n = \n  1, \cds \n  , n$,
 we see that   the index set
 $J_n(\o)  \n := \n  \{1  \n \le \n  j  \n \le \n  n  \n : \|\o   -   \o_j\|_{0,T}
  \n \le \n  \d (\o_j)\}    $ is not empty
 and that $\fg_n(\o) > 1$. Then one can deduce from \eqref{eq:gf011} that
 \beas
 \th_n(\o) - \ol{\th}(\o) & = & \frac{1}{\fg_n(\o)}
 \left( \,  \sum_{j \in J_n(\o)} [\th(\o_j ) - \ol{\th}(\o)]g_{n,j}(\o) +
\sum_{j \notin J_n(\o)} [\th(\o_j ) - \ol{\th}(\o)]g_{n,j}(\o) \right) \\
& \le & \frac{1}{\fg_n(\o)}   \sum_{j \notin J_n(\o)}
|I |g_{n,j}(\o) = \frac{1}{\fg_n(\o)} \sum_{j \notin J_n(\o)} \frac{1}{n^2} < \frac{1}{n} ,
 \eeas
 and similarly, $\th_n(\o) \n - \n  \ul{\th}(\o)  \n > \n  - \frac{1}{n}
    $. Since $\cP$ is a weakly compact subset of $\fP_0$ by (P1)
    and since $\underset{n \in \hN}{\cup} \O_n  \n = \n  \O$, Lemma 8 of \cite{LHP_2011} shows that
  $\lmtd{n \to \infty} \underset{\hP \in \cP}{\sup} \,  \hP \big( \O^c_n \big)  \n = \n  0$.
  Hence, for any $\e  \n > \n  0$, there exists an integer $N  \n > \n  1/\e$  such that for any $n  \n \ge \n  N$,
  $ \underset{\hP \in \cP}{\sup} \, \hP \big( \O^c_n \big)  \n \le \n  \e $. Then
   we  take $\big(\wh{\O} , \wh{\th}\,\big)  \n = \n  \big( \O_N , \th_N \big)$.   \qed

 One can find $\bF-$stopping times that are locally Lipschitz continuous as follows. This result
 and its consequence, Lemma \ref{lem_sandwich2},
 are crucial for our approximating $\tau_0$ by Lipschitz continuous stopping times in Proposition \ref{prop_wp_n}.

     \begin{lemm} \label{lem_conti_st}
 Let $(T_0,\o_0) \n \ins \n (0,T] \ti \O$
 and   $  \fR ,  \k   \n > \n  0 $.
 There exists an $\bF-$stopping time $\z$ valued in $(0,T_0]$   such that
 $ \z    \n \equiv \n  T_0 $ on $  \ol{O}^{T_0}_{\fR} (\o_0)
  \n = \n  \{\o  \n \in \n  \O \n : \|\o  \n - \n  \o_0\|_{0,T_0}  \n \le \n  \fR \}$
 and that given   $\o_1, \o_2  \n \in \n  \O$,
  \bea \label{eq:a109}
   | \z (\o_1)  \n - \n  \z (\o_2) |   \n \le \n  \k\| \o_1 \n - \n  \o_2 \|_{0,t_0}
   \eea
      holds for any
  $t_0  \n \in \n  [b,T_0]  \n \cup \n  \big\{ t  \n \in \n  [a,b) \n :
  t  \n \ge \n  a  \n + \n  \k \| \o_1   \n - \n   \o_2 \|_{0,t} \big\} $,
 where $ a  \n := \n  \z (\o_1)  \n \land \n  \z (\o_2) $ and  $ b  \n := \n  \z (\o_1)  \n \vee \n  \z (\o_2) $.

 \end{lemm}

  \no {\bf Proof:} Given $(t,\o)  \n \in \n  [0,T] \n \times \n \O$,
 the continuity of   paths $\o(\cd)$, $\o_0 (\cd)$ implies that
  \beas
 X_t(\o) \n : = \n  \|\o \n - \n \o_0 \|_{0,t}
 \n = \n \underset{r \in   [0,t]}{\sup} |B_r (\o)  \n - \n \o_0(r)|
 \n = \n \underset{r \in \hQ \cap [0,t]}{\sup} |B_r (\o) \n - \n \o_0(r)|
  \n \in \n  [0,\infty)  .
  \eeas
  As the random variable $\underset{r \in \hQ \cap [0,t]}{\sup} |B_r   \n - \n \o_0(r)|$ is
 $   \cF_t -$measurable, we see that   $X$ is an $\bF-$adapted process with all
 continuous paths.

 Define   $f(x)   \n : = \n
    \n -        x /\k     \+ T_0/\k    \n + \n     \fR     $, $ \fa x  \n \in \n  [ 0 , T_0 ]  $.
 Since $ \z_0  \n := \n  \inf\{t  \n \in \n  [0,T] \n :    f(t \land T_0)  \n - \n  X_t  \n \le \n  0\} \ld T$
 is an $\bF-$stopping time,
  $ \z  \n : = \n  \z_0 \ld T_0  \n = \n  \inf\{t  \n \in \n  [0,T_0] \n :   X_t  \n \ge \n  f(t  )  \}
  \n  \land \n  T_0  $
 is also an $\bF-$stopping time taking values in $(0,T_0]$:
 Given $\o  \n \in \n  \O$, since $ X_0 (\o) \- f(0) \n = \n  0 \- (T_0/\k    \n + \n     \fR) \< 0  $
 and since the path  $X_\cd (\o) \- f(\cd)$ is continuous,
 there exists some $t_\o  \n \in \n  (0,T_0) $ such that
 $X_t (\o)  \- f(t)  \n \le \n  - \frac12 (T_0/\k    \n + \n     \fR) \< 0 $, $\fa t  \n \in \n  [0,t_\o]$. Thus
 $ \z(\o)  \n >  \n  t_\o  \n > \n  0 $.

    Let $\o \n \in \n \O$. If  $\|\o \n - \n \o_0\|_{0,T_0}  \n \le \n  \fR   $,
      one can deduce that
   $   X_t (\o)  \n = \n  \|\o \n - \n \o_0\|_{0,t}  \n \le \n  \|\o \n - \n \o_0\|_{0,T_0}  \n \le \n  \fR
    \n = \n   f(T_0)  \n < \n  f(t) $, $ \fa  t  \n \in \n  [0,T_0) $,
    thus,  $\z(\o) \n = \n  T_0$.

 Next,   let  $\o_1, \o_2 \in \O$. If $ \z(\o_1)  \=   \z(\o_2)   $,   \eqref{eq:a109} holds automatically.
 So let us  assume without loss of generality that
 $a  \n   := \n  \z(\o_1)  \n < \n   \z(\o_2)  \n  := \n  b    $.
 We claim that
 \bea \label{eq:es091}
 \hb{if   $t_0 \n \in \n  [a  , b  ]  $ satisfies
 $t_0  \n - \n a   \n \ge \n   \k \| \o_1   \n - \n   \o_2 \|_{0,t_0} $, then
 $ | \z (\o_1)  \n - \n  \z (\o_2) |  \n = \n b  \n - \n  a  \n \le \n
      \k\| \o_1 \n - \n  \o_2 \|_{0,t_0}  $.  }
 \eea
 To see this, we let  $t_0 \n \in \n  [a  , b  ]  $ satisfying
 $t_0  \n - \n a   \n \ge \n   \k \| \o_1   \n - \n   \o_2 \|_{0,t_0} $,
 and set $\d  \n : = \n  \| \o_1   \n - \n   \o_2 \|_{0,t_0} $,
   $\wh{t}  \n : = \n  a  \n + \n  \k \d  \n \le \n  t_0  $.
 As $ \z(\o_1) \< T_0 $,
 the continuity  of process $X$ and function $f$ implies that
 $ \|\o_1  \n - \n  \o_0\|_{0,a}  \n = \n  \|\o_1  \n - \n  \o_0\|_{0,\z(\o_1)}  \n = \n   f (\z(\o_1))  \n = \n   f (a) $.
 Then one can deduce   that
   \beas
  \q \|\o_2 \n - \n  \o_0\|_{0,\wh{t}}     \ge \n   \|\o_1  \n - \n  \o_0\|_{0,\wh{t} }
    - \n   \|\o_1  \n - \n  \o_2 \|_{0,\wh{t}}  \n \ge \n  \|\o_1  \n - \n  \o_0\|_{0,a}
   \n - \n  \|\o_1  \n - \n  \o_2 \|_{0, t_0}
  \n = \n  f (a)  \n - \n  \d
   \= f \big( \wh{t} \big) .
 \eeas
 So $b = \z (\o_2) \le \wh{t}  $. It follows that
 $  | \z (\o_1)  \n - \n  \z (\o_2) |  \n = \n b  \n - \n  a  \n \le \n   \wh{t}  \n - \n  a  \n = \n
     \k \d  \n = \n  \k\| \o_1 \n - \n  \o_2 \|_{0,t_0}  $, proving the claim.

\if{0}
     On the other hand, for any $ t_0  \n  \in \n   [b, T_0] $,
     if $ b \n - \n a  \n \ge \n  \k \| \o_1   \n - \n   \o_2 \|_{0,b}$,
     applying \eqref{eq:es091} with $t_0 = b$ shows  that
     $ | \z (\o_1)  \n - \n  \z (\o_2) |   \n  \le \n   \k\| \o_1 \n - \n  \o_2 \|_{0,b}
      \n  \le \n   \k\| \o_1 \n - \n  \o_2 \|_{0,t_0} $; Otherwise,
     if $ b \n - \n a  \n \le \n  \k \| \o_1   \n - \n   \o_2 \|_{0,b}$, we directly have
     $ | \z (\o_1)  \n - \n  \z (\o_2) |  \n  = \n  b \n - \n a  \n \le \n  \k\| \o_1 \n - \n  \o_2 \|_{0,b}
     \le \k\| \o_1 \n - \n  \o_2 \|_{0,t_0} $.  \qed
\fi

If   $ b \n - \n a  \n > \n  \k \| \o_1   \n - \n   \o_2 \|_{0,b}$ held,
applying    \eqref{eq:es091} with $t_0 \= b$ would yield  that
     $ b \- a  \= | \z (\o_1)  \n - \n  \z (\o_2) |  \n  \le \dn   \k\| \o_1 \n - \n  \o_2 \|_{0,b}
       $, a contradiction appears. Hence, we must have
       $  | \z (\o_1)  \n - \n  \z (\o_2) |  \n = \n b \n - \n a  \n \le \n  \k \| \o_1   \n - \n   \o_2 \|_{0,b}
       \ls \k\| \o_1 \n - \n  \o_2 \|_{0,t_0} $,  $ \fa t_0  \n  \in \n   [b, T_0] $. \qed


 \begin{lemm} \label{lem_sandwich2}
 Let  $ \th_1 $, $\th_2$, $\th_3$ be three real-valued random variables on $\O$ satisfying:
     for some   $\d  \n  > \n  0$, it holds for $i=1,2$ and any $\o  \n \in \n  \O$ that
     \bea \label{eq:eb033}
     \th_i (\o')  \n \le \n     \th_{i+1} (\o), \q \fa  \o'  \n \in \n  \ol{O}^{\th_{i+1}(\o)}_\d (\o)
       \n = \n  \big\{\o'  \n \in \n  \O \n : \|\o' \n - \n \o\|_{0,\th_{i+1}(\o)}  \n \le \n  \d \big\} .
       \eea
 If $\th_2$ takes values in $(0,T]$,
 then for any   $\k  \n > \n  T/\d$, there exists    an $\bF-$stopping time $\wp $
  such that   $ \th_1     \n \le \n  \wp   \n \le \n  \th_3   $ on   $   \O $.
    Moreover, given  $\o_1, \o_2  \n \in \n  \O$,
      \bea \label{eq:a111}
      \big| \wp (\o_1)  \n - \n  \wp (\o_2) \big|   \n \le \n  \k\| \o_1 \n - \n  \o_2 \|_{0,t_0}
      \eea
      holds for any
  $t_0  \n \in \n  [b,T]  \n \cup \n  \big\{ t  \n \in \n  [a,b) \n : t  \n \ge \n  a
   \n + \n  \k \| \o_1   \n - \n   \o_2 \|_{0,t} \big\} $,
 where $ a  \n := \n  \wp (\o_1)  \n \land \n  \wp (\o_2) $
 and  $ b  \n := \n  \wp (\o_1)  \n \vee \n  \wp (\o_2) $.

 \end{lemm}

  \no {\bf Proof:} We fix $ \k \n > \n  T/\d $ and set $\d_0  \n := \n  \d  \n - \n  T/\k   $.
   Since the canonical space $\O $ is a separable complete metric space
   and thus  Lindel\"of,    there exists
   a countable  dense subset $\big\{ \o_j   \big\}_{j \in \hN}$ of $\O$ under   norm $\|~\|_{0,T}$.
 Given $j \n \in \n  \hN$, we   set $t_j  \n : = \n  \th_2 (\o_j)  \n \in \n  ( 0 , T ]  $
 and $\k_j \df \frac{t_j  }{\d    -    \d_0}$.
 Applying Lemma \ref{lem_conti_st} with
  $(\o_0,   T_0, \fR , \k )  \n = \n  (\o_j,    t_j , \d_0 , \k_j   ) $
  yields an $\bF-$stopping time $\z_j$ valued in   $ (0, t_j ]  $ such that
  \bea \label{eq:es201}
   \z_j (\o)   \n \equiv \n  t_j  , \q
    \fa \o  \n \in \n  \ol{O}^{t_j}_{\d_0} (\o_j) .
  \eea
  Given   $\o_1, \o_2  \n \in \n  \O$, it holds for any
  $t_0  \n \in \n  [b_j,t_j]  \n \cup \n  \big\{ t  \n \in \n  [a_j,b_j) \n :
  t  \n \ge \n  a_j  \n + \n  \k_j \| \o_1   \n - \n   \o_2 \|_{0,t} \big\} $ that
 \bea \label{eq:eb045}
  \big| \z_j (\o_1) \n - \n  \z_j (\o_2) \big|
  \n \le \n  \k_j \| \o_1   \n - \n   \o_2 \|_{0,t_0}
  \n \le \n  \k \| \o_1   \n - \n   \o_2 \|_{0,t_0}   ,
 \eea
  where $ a_j := \z_j (\o_1)  \n \land \n  \z_j (\o_2) $
 and  $ b_j := \z_j (\o_1)  \n \vee \n  \z_j (\o_2) $.

  Clearly,   $\wp \n : = \n  \underset{j \in \hN}{\sup} \, \z_j $ defines
  an $\bF-$stopping time taking values in $(0,T]$.
   Let  $\o_1, \o_2  \n \in \n  \O$.
  If $ \wp(\o_1)  \=   \wp(\o_2)    $,  one has \eqref{eq:a111}   automatically.
  So let us  assume without loss of generality that
  $a  \n   := \n  \wp (\o_1)  \n < \n   \wp (\o_2)  \n  := \n  b    $.
 We claim that
 \bea \label{eq:ec043}
 \hb{if $t_0 \n \in \n   [ a  , b  ]  $ satisfies
 $t_0  \n - \n a  \n \ge \n  \k \| \o_1   \n - \n   \o_2 \|_{0,t_0} $,
 then $\big| \wp (\o_1) \- \wp (\o_2) \big|   \n \le \n  \k \| \o_1   \n - \n   \o_2 \|_{0,t_0}   $. }
 \eea
 To see this, we let $t_0 \n \in \n   [ a  , b  ]  $ satisfying
 $t_0  \n - \n a  \n \ge \n  \k \| \o_1   \n - \n   \o_2 \|_{0,t_0} $,
 and let $\l  \n \in \n (0, b  \n - \n  a ]  $.
 There exists a $ j \= j(\l ) \n \in \n  \hN $ such that $   \z_j (\o_2)  \n \ge \n    b   \n - \n  \l  $. As
 $  \z_j (\o_2)  \n \ge \n  a  \n = \n  \wp (\o_1)  \n \ge \n  \z_j (\o_1)    $,
 we see that $a_j  \n = \n  \z_j (\o_1) $ and  $b_j  \n = \n  \z_j (\o_2)$. Then
 $ t_0 $ is in $  [a_j  , T  ] $ and satisfies
 $t_0  \n - \n a_j  \n \ge \n t_0  \n - \n a  \n \ge \n  \k \| \o_1   \n - \n   \o_2 \|_{0,t_0}
 \n \ge \n  \k_j \| \o_1   \n - \n   \o_2 \|_{0,t_0} $.
 So by \eqref{eq:eb045},
 $ \big| \wp (\o_1)  \n - \n  \wp (\o_2) \big|  \n = \n  b  \n - \n  a
  \n \le \n   \z_j (\o_2) \n + \n \l \n - \n  \z_j (\o_1)
  \n \le \n  \k \| \o_1   \n - \n   \o_2 \|_{0,t_0} \n + \n \l$.
  Letting $\l \to 0$ yields that
 $  \big| \wp (\o_1) - \wp (\o_2) \big|   \n \le \n  \k \| \o_1   \n - \n   \o_2 \|_{0,t_0}   $,
 proving the claim.

\if{0}
      On the other hand, let $ t_0  \n  \in \n    [ b, T  ] $.
     If $ b \n - \n a  \n \ge \n  \k \| \o_1   \n - \n   \o_2 \|_{0,b}$,
     applying \eqref{eq:ec043} with $t_0  \n = \n  b $  shows   that
     $ | \wp (\o_1)  \n - \n  \wp (\o_2) |   \n  \le \n   \k\| \o_1 \n - \n  \o_2 \|_{0,b}
      \n  \le \n   \k\| \o_1 \n - \n  \o_2 \|_{0,t_0} $; Otherwise,
     if $ b \n - \n a  \n \le \n  \k \| \o_1   \n - \n   \o_2 \|_{0,b}$, it holds automatically that
     $ | \wp (\o_1)  \n - \n  \wp (\o_2) |  \n  = \n  b \n - \n a  \n \le \n  \k\| \o_1 \n - \n  \o_2 \|_{0,b}
     \le \k\| \o_1 \n - \n  \o_2 \|_{0,t_0} $.
\fi

If   $ b \n - \n a  \n > \n  \k \| \o_1   \n - \n   \o_2 \|_{0,b}$ held,
applying claim \eqref{eq:ec043} with $t_0 = b$ would yield  that
     $ b \- a  \= | \wp (\o_1)  \n - \n  \wp (\o_2) |  \n  \le \n   \k\| \o_1 \n - \n  \o_2 \|_{0,b}
       $, a contradiction appears. Hence, we must have
       $ | \wp (\o_1)  \n - \n  \wp (\o_2) |  \n  = \n b \n - \n a  \n \le \n  \k \| \o_1   \n - \n   \o_2 \|_{0,b} \ls \k\| \o_1 \n - \n  \o_2 \|_{0,t_0} $,   $ \fa  t_0  \n  \in \n   [b, T ] $.

    Now,  let us fix $\o \n \in \n  \O$.
  Since $O_{\d_0} \big(\o_j\big)  \n \subset \n  O^{t_j}_{\d_0} \big(\o_j\big) $ for any
    $  j  \n \in \n  \hN$, one has
  $ \O  \n = \n  \underset{j \in \hN}{\cup} O_{\d_0} \big(\o_j\big)
   \n \subset  \n  \underset{j \in \hN}{\cup}  O^{t_j}_{\d_0} \big(\o_j\big)  \n \subset \n  \O $.
  So  $\o  \n \in \n  O^{t_j}_{\d_0} \big( \o_j \big)$
  for some $j  \n \in \n  \hN$ and it follows from \eqref{eq:es201} that
  $  \wp (\o)  \n \ge \n  \z_j (\o)  \n = \n  t_j  \n > \n 0 $.
  Since $\|\o \n - \n \o_j\|_{0,\th_2 (\o_j)}  \n = \n  \|\o \n - \n \o_j\|_{0, t_j}   \n < \n   \d_0  \n < \n  \d  $,
  taking $(i,\o,\o')  \n = \n  (1, \o_j , \o )$ in \eqref{eq:eb033} shows that
  $  \th_1 (\o)  \n \le \n  \th_2 (\o_j)  \n = \n  t_j    \n = \n  \z_j (\o)  \n \le \n  \wp (\o)  $.

  We claim   that $ \z_\ell (\o)  \n \le \n  \th_3 (\o) $, $ \fa \ell  \n \in \n \hN $: Assume not,
  i.e. $ \z_\ell (\o)  \n > \n  \th_3 (\o) $ for some $   \ell  \n \in \n \hN $.
  From the proof of Lemma \ref{lem_conti_st}, we see that
  $ \z_\ell (\o) = \inf\{t  \n \in \n  [0,t_\ell] \n :   \|\o - \o_\ell  \|_{0,t}  \n \ge \n  f_\ell (t  )  \}
  \n  \land \n  t_\ell   $, where $ f_\ell (x)  \n : = \n  -   x / \k_\ell \n + \n  t_j/ \k_\ell  \n + \n  \d_0 $,
  $\fa x  \n \in \n  [0,t_\ell]$.
  Since $ \|\o_\ell  \n - \n  \o\|_{0,\th_3(\o)} \n \le \n  \|\o  \n - \n  \o_\ell  \|_{0, \z_\ell (\o)}
   \n \le \n  f_\ell \big( \z_\ell (\o) \big)  \n < \n   f_\ell  ( 0  )  \n = \n
    t_\ell / \k_\ell  \n + \n  \d_0 \= \d   $,
  taking $(i,\o,\o')  \n = \n  \big(2, \o , \o_\ell \big)$ in \eqref{eq:eb033} leads to a contradiction:
  $   \th_3(\o)  \n \ge \n  \th_2 \big( \o_\ell \big)  \n = \n  t_\ell \ge \z_\ell (\o) $\,!
  Hence, $   \z_\ell (\o)  \n \le \n   \th_3 (\o) $, $ \fa   \ell  \n \in \n \hN $.
  It follows that
  $\wp (\o)  \n = \n  \underset{\ell    \in   \hN}{\sup} \z_\ell (\o)  \n \le \n   \th_3 (\o) $. \qed

\subsection{Proofs of Starred Inequalities in Section \ref{sec:proofs}}

    \no {\bf Proof of \eqref{eq:a107}:}
    Let $r \n \in \n  [t,T]$.
    If $r  \n < \n  t_1$,
     as $ \{  \ga    \n  < \n  \nu \}  \n \in \n  \cF^t_{ \ga \land \nu}    \n \subset \n  \cF^t_{ \ga } $,
    one has
    $ \{ \wh{\ga}_\l   \n \le \n  r \}  \n = \n  \{  \ga   \n   < \n  \nu \}   \cap    \{ \ga  \n \le \n  r\}
       \n \in \n  \cF^t_r $. Otherwise, if $r  \n \ge \n  t_1$, let
       $k$ be the largest integer such that $t_k  \n \le \n  r  $.
      Since
       $ \{  \ga  \n  \ge  \n   \nu \} \n \cap \n \{\ga  \n \le \n  r\} \n \subset \n \{ \nu  \n \le \n  r \}
       \n \subset \n \{ \nu  \n \ne \n  t_i \}  \n \subset \n \cA^i_0 $
       for $i \n = \n k \n + \n 1 \cds \n, m$
       and since $  \{ \ga  \n \ge \n  \nu \}   \n \cap \n     \cA^i_j
     \n = \n \{ \ga  \n \ge \n  t_i \}   \n  \cap \n  \{ \nu \n = \n t_i \}  \n \cap \n
     \Big(  O^{t_i}_{\d_j }(\wt{\o}_j)   \backslash \underset{j' < j}{\cup}  O^{t_i}_{\d_{j'} }(\wt{\o}_{j'}) \Big)
      \n \in \n  \cF^t_{t_i}  \n \subset \n  \cF^t_r $ for $i \n = \n 1,\cds \n , k$ and $j \n = \n 1, \cds \n , \l$,
      one can deduce that
        \beas
    \{  \wh{\ga}_\l  \n \le \n  r \}
   \n = \n \big( \{  \ga    \n  < \n  \nu \}  \n \cap \n  \{ \ga  \n \le \n  r\} \big)  \n \cup \n
   \bigg[ \{  \ga    \n  \ge  \n  \nu \}  \n \cap \n  \{ \ga  \n \le \n  r\}
    \n \cap \n  \Big( \underset{i=1}{\overset{k}{\cap}} \cA^i_0 \Big) \bigg]
    \n \cup \n  \bigg[
    \underset{i=1}{\overset{k}{\cup}} \underset{j=1}{\overset{\l}{\cup}}
    \n  \Big(  \{  \ga   \n  \ge  \n  \nu \}  \n \cap  \n    \cA^i_j
     \n \cap  \n  \{ \ga^i_j  (\Pi^t_{t_i})  \n \le \n  r \}     \Big) \bigg]  \n \in \n  \cF^t_r .
       \eeas
  Hence,  $\wh{\ga}_\l   \n \in \n  \cT^t$.

 \no {\bf Proof of \eqref{eq:a114}:}
 We let $\wh{\k}_n$ be the Lipschitz coefficient of $\wh{\d}_n$.
 Given  $\o  \n \in \n  \O$ and  $\e  \n > \n  0$,
   set $ \wh{\l}_n  \n = \n  \wh{\l}_n (\o,\e)  \n := \n   \frac{\e}{3}  \n \land \n
 \frac{(\phi^\o_T)^{-1} (\e/3)}{\wh{\k}_n} $ and
 let $\o'  \n \in \n  O_{\wh{\l}_n} (\o)$.

 Let $0  \n \le \n  r  \n \le \n  r'  \n \le \n  T$ with $r' \n - \n r  \n \le \n  \wh{\d}_n (\o)$.
 If $\wh{\d}_n (\o)  \n \le \n  \wh{\d}_n (\o')$, then
 \bea  \label{eq:es031}
 |\o(r') \n - \n \o(r)|  \n \le \n  |\o(r') \n - \n \o'(r')| \n + \n  |\o'(r') \n - \n \o'(r)|
  \n + \n  |\o'(r ) \n - \n \o(r )|
  \n \le \n  \phi^{\o'}_T \n \big(   \wh{\d}_n (\o) \big)  \n + \n  2 \|\o'  \n - \n \o \|_{0,T}
  \n < \n  \phi^{\o'}_T \n \big(   \wh{\d}_n (\o') \big)  \n + \n  \frac23 \e . \q
 \eea
 Otherwise if $\wh{\d}_n (\o') \n < \n  \wh{\d}_n (\o)$, we set $s'  \n : = \n  r'  \n \land \n  (r \n + \n \wh{\d}_n (\o'))$.
 Since \eqref{eq:ec017} shows that $r' \n - \n s'
  \n = \n  r'  \n \land \n  (r \n + \n \wh{\d}_n (\o))  \n - \n  r'  \n \land \n  (r \n + \n \wh{\d}_n (\o'))
  \n \le \n  \wh{\d}_n (\o)  \n - \n  \wh{\d}_n (\o') \n \le \n \wh{\k}_n \|\o  \n - \n \o' \|_{0,T} $
 and that $s' \n - \n r  \n = \n  r'  \n \land \n  (r \n + \n \wh{\d}_n (\o'))
  \n - \n  r'  \n \land \n  r  \n \le \n  \wh{\d}_n (\o') $,
 we can deduce that
  \beas
 |\o(r') \n - \n \o(r)| & \tn \le & \tn  |\o (r') \n - \n \o (s')|  \n + \n  |\o (s') \n - \n \o(r )|
  \n \le \n  \phi^\o_T \big( \wh{\k}_n \|\o  \n - \n \o' \|_{0,T} \big)
  \n + \n  |\o (s') \n - \n \o' (s')|  \n + \n  |\o'(s') \n - \n \o'(r)|  \n + \n  |\o'(r ) \n - \n \o(r )| \\
  & \tn \le & \tn  \phi^\o_T \big( \wh{\k}_n \|\o  \n - \n \o' \|_{0,T} \big)
  \n + \n  \phi^{\o'}_T \n \big(   \wh{\d}_n (\o') \big)  \n + \n  2 \|\o'  \n - \n \o \|_{0,T}
  \n < \n  \phi^{\o'}_T \n \big(   \wh{\d}_n (\o') \big)  \n + \n  \e .
 \eeas
 Combining it with \eqref{eq:es031}   and taking supremum over the pair $(r,r')$ yields that
 $ \phi^\o_T \big( \wh{\d}_n (\o) \big) \le \phi^{\o'}_T \n \big(   \wh{\d}_n (\o') \big) + \e $.

 On the other hand,   let $0 \ls \wt{r} \ls \wt{r}' \le T$ with $\wt{r}' \- \wt{r} \ls  \wh{\d}_n (\o')$.
 If $\wh{\d}_n (\o') \ls  \wh{\d}_n (\o)$, an analogy to \eqref{eq:es031} shows that
   \bea \label{eq:es033}
 |\o'(\wt{r}')-\o'(\wt{r})| 
 < \phi^\o_T \big(   \wh{\d}_n (\o) \big) + \frac23 \e .
 \eea
 Otherwise if $\wh{\d}_n (\o) < \wh{\d}_n (\o')$, one can deduce that
 \beas
  |\o'(\wt{r}') \n - \n \o'(\wt{r})| & \tn \dn \le & \tn \dn
  |\o'(\wt{r}') \n - \n \o(\wt{r}')| \n + \n  |\o(\wt{r}') \n - \n \o(\wt{r})|
   \n + \n   |\o(\wt{r} ) \n - \n \o'(\wt{r} )|
  \n \le \n  \phi^\o_T \big( \wh{\d}_n (\o')   \big)   \n + \n  2 \|\o  \n - \n \o' \|_{0,T}  \\
  & \tn \dn \le & \tn \dn \phi^\o_T \big( \wh{\d}_n (\o')  \n - \n  \wh{\d}_n (\o) \big)
  \n + \n  \phi^\o_T \big(   \wh{\d}_n (\o) \big)  \n + \n  2 \|\o  \n - \n \o' \|_{0,T}
  \n <  \n  \phi^\o_T \big( \wh{\k}_n \|\o  \n - \n \o' \|_{0,T} \big)
  \n + \n  \phi^\o_T \big(   \wh{\d}_n (\o) \big)  \n + \n  \frac23 \e
  \n \le \n  \phi^\o_T \big(   \wh{\d}_n (\o) \big)  \n + \n  \e .
 \eeas
 Combining it with \eqref{eq:es033}   and taking supremum over the pair $\big(\wt{r},\wt{r}'\big)$ yields that
 $ \phi^{\o'}_T \n \big( \wh{\d}_n (\o') \big) \le \phi^\o_T \big(   \wh{\d}_n (\o) \big) + \e $.

 Hence $\o \to \phi^\o_T \big(\wh{\d}_n (\o) \big)$ is   a continuous random variable on $\O$. \qed

 \no {\bf Proof of \eqref{eq:a117}:}    Let $\o, \o' \n \in \n \O$ and set
 $t : = \wh{\th}_n (\o )$, $s : = \wh{\th}_n (\o' )$.  We see from \eqref{eq:gd011} and \eqref{eq:es017} that
       \beas
     \big|  Z_s (\o  ) \n - \n  Z_s (\o') \big|
  & \tn \dn \le   & \tn \dn       \wrY \Big( (1 \n + \n \k_\wp) \|\o  \n - \n \o' \|_{0,T}
  \n + \n \phi^\o_T \big( \k_\wp   \|\o   \n  - \n  \o'  \|_{0,T} \big)   \Big)  , \q \hb{and} \\
  \big| Z_t (\o) \n - \n  Z_s (\o) \big|
    & \tn \dn  = & \tn \dn  \big| Z_{t \land s} (\o) \n - \n  Z_{t \vee s} (\o) \big|
    \n  \le  \n  2 C_{\n \vr} M_Y   \Big( |s \n - \n t|^{ \frac{q_1}{2}}
     \n \vee  \n  |s \n - \n t|^{q_2-\frac{q_1}{2}} \Big)  \n + \n \wrY \big(|s \n - \n t|\big)
   \n + \n    \wrY  \big(   \d'_{t,s} (\o)  \big)
    \n \vee \n  \wh{\wh{\rho}\,}_Y \big(   \d'_{t,s} (\o)  \big)       ,
  \eeas
  where  $\d'_{t,s} (\o) \n := \n   (1 \n + \n \k_\wp)
 \Big(  |s \n - \n t|^{\frac{q_1}{2}}  \n + \n \phi^\o_T \big(|s \n - \n t|\big) \Big)       $.
 Adding them up, one can deduce from  the Lipschitz continuity of random variable $\wh{\th}_n$  that
  $Z_{\wh{\th}_n}$ is a continuous   random variable on $\O$. \qed

   \no {\bf Proof of \eqref{eq:a119}:}
   If   $\wp_n (\o_1)  \n \land \n   \wp_n (\o_2)   \n  + \n  2^{-k}    \n  > \n  t_1 $,
   one has  $H_{t_1} (\o_1) \n = \n  H_{t_1} (\o_2)  \n = \n 0 $. 
   On the other hand, suppose that  $\wp_n (\o_1)   \n   \land  \n     \wp_n (\o_2)   \n  + \n  2^{-k}
      \n \le \n  t_1    $.
     When $\|\o_1  \n - \n  \o_2\|_{0,t_1}  \n \ge \n  2^{-k} \k^{-1}_n $,
     we automatically have
     $      \big|H_{t_1} (\o_1) \n - \n  H_{t_1} (\o_2) \big|  \n \le \n  1
     \n \le \n    2^k \k_n   \|\o_1  \n - \n  \o_2\|_{0,t_1}   $;
     When
   $\|\o_1  \n - \n  \o_2\|_{0,t_1}  \n < \n  2^{-k} \k^{-1}_n $, since
   $  \wp_n (\o_1)  \n \land \n   \wp_n (\o_2)  \n + \n  \k_n \|\o_1  \n - \n  \o_2\|_{0,t_1} 
    \n  < \n  \wp_n (\o_1)  \n \land \n   \wp_n (\o_2)  \n + \n  2^{-k}  \n \le \n  t_1 $,
   applying Proposition \ref{prop_wp_n} (2) with $t_0  \n = \n  t_1$ yields that
   $ \big| \wp_n (\o_1) \n - \n \wp_n (\o_2) \big|
    \n \le \n  \k_n  \|\o_1  \n - \n  \o_2\|_{0,t_1} $.
   Then  \eqref{eq:ec017} implies  that
   \beas
   \q  \big|H_{t_1} \n  (\o_1) \n - \n  H_{t_1} \n  (\o_2) \big|
      \ls  \big| \big( 2^k  (  t_1  \n - \n  \wp_n (\o_1)) \n - \n 1 \big)^+
    \n \-   \big( 2^k  (  t_1  \n - \n \wp_n (\o_2)) \n - \n 1   \big)^+ \big|
    \ls     2^k \big|  \wp_n (\o_1) \n - \n  \wp_n (\o_2) \big|
        \ls   2^k \k_n   \|\o_1  \n - \n  \o_2\|_{0,t_1}  . \q \hb{\qed}
   \eeas

  \no {\bf Proof of \eqref{eq:et415}:}
  Let $\o' \n \in \n \O$. If the set
  $\{t'  \n \in \n  [0,T] \n : \sX (t' , \o')  \n \le \n  0\}$  is not empty,
    Proposition \ref{prop_wp_n} (1) implies that
    $ \lmtu{n \to \infty} \wp_n (\o')  \n = \n  \tau_0 (\o') $,
  however, $ \wp_n (\o')  \n < \n  \tau_0 (\o') $ for any $n \in \hN$.
  Then one can deduce that $ \lmt{n \to \infty} \b1_{ [   0, \wp_n (\o') ]   } (t') \n = \n
   \b1_{ [   0, \tau_0 (\o')  )   }  (t')    $, $ \fa t'  \n \in \n  [0,T]$, and the continuity of the path
   $U_\cd (\o')$ implies that
   \beas
   \lmt{n \to \infty} \sY^n_{t'} (\o') & \tn =  & \tn   \lmt{n \to \infty}
   \big(   \b1_{\{ t' \le \wp_n (\o') \}} L (t',\o')   \n + \n
 \b1_{\{ t' > \wp_n (\o')    \}} U (\wp_n (\o') , \o')  \big) \\
  & \tn =  & \tn   \b1_{\{ t' < \tau_0 (\o') \}} L (t',\o')   \n + \n
 \b1_{\{ t' \ge \tau_0 (\o')    \}} U (\tau_0 (\o') , \o')
 = \sY (\tau_0 (\o') \land t' , \o')
  =  \wh{\sY}_{t'} (\o')  , \q \fa t'  \n \in \n  [0,T] .
   \eeas
 On the other hand, if the set
  $\{t'  \n \in \n  [0,T] \n : \sX (t' ,   \o')  \n \le \n  0\}$  is  empty, the continuity of
  path $\sX_\cd  (\o')$ implies that $ \underset{t' \in [0,T]}{\inf} \sX (t' ,   \o') > 0 $.
  For large enough $n \in \hN$,  the set  $\big\{t'  \n \in \n  [0,T] \n : \sX (t',\o')  \n \le \n
    \big(  \lceil \log_2(n \n + \n 2) \rceil  \n + \n \lfloor  \sX^{-1}_0  \rfloor \n - \n 1 \big)^{-1}  \big\}
     $ is also empty, thus  $ T = \tau_n (\o')  \n = \n  \wp_n (\o')  \n = \n  \tau_0 (\o') $
     by Proposition \ref{prop_wp_n} (1).
    Then   (A2) shows that for any  $t'  \n \in \n  [0,T]$
   \beas
 \qq   \lmt{n \to \infty} \sY^n_{t'} (\o') & \tn \dn  =  & \tn  \dn   \lmt{n \to \infty}
   \big(   \b1_{\{ t' \le \wp_n (\o') \}} L (t',\o')   \n + \n
 \b1_{\{ t' > \wp_n (\o')    \}} U (\wp_n (\o') , \o')  \big)
 \n =  \n  \b1_{\{ t' \le T \}} L (t',\o')     \\
  & \tn  \dn =  & \tn  \dn   \b1_{\{ t' < T \}} L (t',\o')   \n + \n
 \b1_{\{ t' = T    \}} U (T , \o')
   \n = \n    \b1_{\{ t' < \tau_0 (\o') \}} L (t',\o')   \n + \n
 \b1_{\{ t' \ge \tau_0 (\o')    \}} U (\tau_0 (\o') , \o')
   \n = \n     \wh{\sY}_{t'} (\o')  . \qq \hb{\qed}
   \eeas

 \no {\bf Proof of \eqref{eq:a125}:}
 Let $\o \n \in \n \O$.  Since $\wh{Y}^{\ell,\ell}_t (\o)  \n = \n  L_t (\o)$ over $ [0, \wp_\ell (\o) \n + \n  2^{- \ell} ]
  \n \supset \n   [0, \wp_n ) $, one has
  \bea
  \wh{\z}^\a_{i,\ell} \n \land \n  \wp_n
  \=  \inf\big\{t  \n \in \n  [0,\wp_n) \n : Z^{\ell,\ell}_t
     \n  \le  \n   \wh{Y}^{\ell,\ell}_t  \n + \n  1/i  \n + \n  1/\a \big\}   \n \land \n  \wp_n
     \=  \inf\big\{t  \n \in \n  [0,\wp_n) \n : Z^{\ell,\ell}_t  \n \le \n  L_t  \n + \n  1/i  \n + \n  1/\a \big\}
   \n \land \n  \wp_n
  \=  \z^\a_{i,\ell} \n \land \n  \wp_n  . \q    \label{eq:et363}
  \eea

 If
 $ Z^{m_j,m_j}_t (\o) \n = \n  L_t (\o) $ for some $t  \n \in \n  \big[0,\wp_n (\o)\big)$, applying \eqref{eq:et341}
 with $k \n = \n m_j$ and $k \n = \n \ell$ respectively shows that
 $  \sZ_t(\o)  \n \le \n     Z^{m_j,m_j}_t (\o)   \n + \n     \ol{\e}_{m_j}
   \n \le \n   L_t (\o)   \n + \n     \ol{\e}_\ell
    \n < \n    L_t (\o)  \n + \n  \frac{1}{2i}  \n + \n  \frac{1}{\a} $ and thus
 $  Z^{\ell,\ell}_t (\o)  \n \le \n   \sZ_t(\o)     \n + \n    \ol{\e}_\ell
        \n < \n   L_t (\o)  \n + \n  \frac{1}{ i}  \n + \n  \frac{1}{\a} $.
 So $ \inf\{t  \n \in \n  [0,\wp_n(\o)) \n : Z^{\ell,\ell}_t (\o) \n \le \n  L_t (\o) \n + \n  1/i  \n + \n  1/\a \}
  \n \le \n   \inf\{t  \n \in \n  [0,\wp_n(\o)) \n : Z^{m_j,m_j}_t (\o) \n = \n     L_t(\o) \}  $.
 As $\wh{Y}^{m_j,m_j}_t (\o) \n = \n  L_t (\o) $
 over $ 
   [ 0, \wp_n (\o) ) $,  one can deduce that
 \beas
  \z^\a_{i,\ell} (\o) \n \land \n  \wp_n (\o) & \tn \dn = & \tn \dn
    \inf\big\{t  \n \in \n  [0,\wp_n (\o) ) \n : Z^{\ell,\ell}_t  (\o)
     \n \le \n  L_t  (\o)  \n + \n  1/i  \n + \n  1/\a \big\}
   \n \land \n  \wp_n  (\o)
   \n \le \n \inf\{t  \n \in \n  [0,\wp_n (\o) ) \n : Z^{m_j,m_j}_t  (\o)
    \n = \n     L_t (\o)  \}  \n \land \n  \wp_n (\o)     \\
   & \tn \dn = & \tn \dn
   \inf\{t  \n \in \n  [0,\wp_n (\o) ) \n : Z^{m_j,m_j}_t  (\o)
    \n = \n     \wh{Y}^{m_j,m_j}_t (\o)  \}  \n \land \n  \wp_n (\o)
   \n = \n   \nu_{m_j}  (\o)  \n \land \n  \wp_n  (\o)       .
 \eeas
 On the other hand, if the set $\big\{t  \n \in \n  [0,\wp_n (\o)) \n : Z^{m_j,m_j}_t (\o) \n = \n L_t (\o)   \big\}$
 is empty, we can deduce that $  \nu_{m_j} (\o)  \n \ge \n  \wp_n (\o) $. Then
 $\nu_{m_j} (\o) \n \land \n  \wp_n (\o)  \n = \n \wp_n (\o)
 \n \ge \n  \z^\a_{i,\ell} (\o) \n \land \n  \wp_n (\o) $ holds automatically. \qed

 \no {\bf Proof of \eqref{eq:a127}:}
 Set  $\wh{\z}'_{i,\ell} \df \lmtu{\a \to \infty} \wh{\z}^\a_{i,\ell}   \le \wh{\z}_{i,\ell}    $.
 As   the continuity of $Z^{\ell,\ell}  \n - \n \wh{Y}^{\ell,\ell}$    shows that
 $  Z^{\ell,\ell}_{\wh{\z}^\a_{i,\ell}}  \n - \n \wh{Y}^{\ell,\ell}_{\wh{\z}^\a_{i,\ell}}     \n \le \n
 \frac{1}{i} + \frac{1}{\a } $, $\fa \a \n > \n \ell $, letting $\a  \n \to \n  \infty$, we see from
 the continuity of $Z^{\ell,\ell}  \n - \n \wh{Y}^{\ell,\ell}$ again that
 $    Z^{\ell,\ell}_{\wh{\z}'_{i,\ell}}  \n - \n \wh{Y}^{\ell,\ell}_{\wh{\z}'_{i,\ell}}   \n \le \n  1/i $,
 which implies that
 $ \wh{\z}_{i,\ell} = \wh{\z}'_{i,\ell} = \lmtu{\a \to \infty} \wh{\z}^\a_{i,\ell}    $. \qed

  \no {\bf Proof of \eqref{eq:et367}:}
 Let  $\o \n \in \n  \O$. If the set
 $ \cI (\o)  \n := \n  \big\{t  \n \in \n  [0,T] \n : \sZ_t (\o)  \n \le \n  L_t (\o) \n + \n  1/i \big\} $
 is empty,  the continuity of   path $\sZ_\cd (\o) \n - \n  L_\cd (\o) $ implies that
 $\eta (\o)  \n : = \n  \underset{t \in [0,T]}{\inf} \big(\sZ_t (\o)  \n - \n  L_t (\o)\big)  \n > \n  1/i$.
 For any integer $h \n > \n    ( \eta (\o)  \n - \n  1 / i )^{-1}     $,
 since $   \underset{t \in [0,T]}{\inf} (\sZ_t (\o)  \n - \n  L_t (\o))
  \n = \n  \eta (\o)  \n > \n  1/i  \n + \n  1/h$,
 the set $\cI_h (\o) : = \big\{t  \n \in \n  [0,T] \n :
   \sZ_t (\o) \n \le \n  L_t (\o) \n + \n  1/i  \n + \n  1/h  \big\}$ is also empty and thus
     $\ga^h_i (\o)  \n = \n  T$.
 It follows that $\lmtu{h \to \infty} \ga^h_i (\o)  \n = \n  T  \n = \n  \ga_i (\o) $.

 On the other hand, if   $ \cI (\o) $ is not empty,
 we set  $\ga'_i (\o)  \n := \n  \lmtu{h \to \infty} \ga^h_i  (\o)    \n \le \n  \ga_i  (\o)
 \n = \n \inf \cI (\o) $.
 For any $ h  \n \in \n  \hN $,   $ \cI_h (\o) $ contains $\cI (\o)$ and is thus not empty.
 The continuity of   path $\sZ_\cd (\o)  \n - \n L_\cd (\o)$  then implies that
 $  \sZ \big(\ga^h_i (\o) , \o \big)  \n - \n L \big(\ga^h_i (\o) , \o \big)    \n \le \n
 \frac{1}{i}  \n + \n  \frac{1}{h } $.
 Letting $h  \n \to \n  \infty$, we see from
 the continuity of   path $\sZ_\cd (\o)  \n - \n L_\cd (\o)$ again that
 $    \sZ \big(\ga'_i (\o) , \o \big)    \n - \n L \big(\ga'_i (\o) , \o \big)   \n \le \n  1/i $,
 which shows  that $  \ga_i  (\o)
 \n = \n \inf \cI (\o) \n \le \n \ga'_i (\o) $. Thus
 $ \ga_i (\o)  \n = \n  \ga'_i (\o)  \n = \n  \lmtu{h \to \infty} \ga^h_i (\o)    $. \qed

  \no {\bf Proof of \eqref{eq:et383}:}
  Set $s_n  \n = \n  s_n ( \o)  \n := \n ( \ga_*    \n \land \n   \wp_n ) (\o) $
  and  let $\e  \n > \n  0$. By the continuity of   path $\sZ_\cd (\o)$,
  there exists a $\d_n  \n = \n  \d_n ( \o) \n > \n 0 $ such that
  $\big|\sZ_t (\o)  \n - \n  \sZ \big(  s_n , \o \big)\big| \le \e $,
  $\fa t  \n \in \n  \big[ (s_n  \n - \n \d_n)^+ , s_n \big]$.
  We see from \eqref{eq:et381} that for  large enough $i  \n \in \n  \hN$,
  both $(\ga_i    \n \land \n  \wp_n) (\o) $ and
  $(\ga_{2i}    \n \land \n  \wp_n) (\o) $ are in $\big[ (s_n  \n - \n \d_n)^+ , s_n \big]$,
  so $\cJ_{n,i} (\o)  \n \subset \n  \big[ (s_n  \n - \n \d_n)^+ , s_n \big]$.
  It follows that
  $  \sZ \big(  s_n , \o \big)  \n - \n  \e  \n \le \n
  \underset{t \in \cJ_{n,i} (\o)}{\inf}   \sZ(t,\o)  \n \le \n  \underset{t \in \cJ_{n,i} (\o)}{\sup}   \sZ(t,\o)
   \n \le \n    \sZ \big(  s_n , \o \big)  \n + \n  \e $.
  As $i  \n \to \n  \infty$, we obtain
  $  \sZ \big(  s_n , \o \big)  \n - \n  \e
   \n \le \n  \linf{i \to \infty}   \underset{t \in \cJ_{n,i} (\o)}{\inf}   \sZ(t,\o)
   \n \le \n  \lsup{i \to \infty} \underset{t \in \cJ_{n,i} (\o)}{\sup}   \sZ(t,\o)
   \n \le \n    \sZ \big(  s_n , \o \big)  \n + \n  \e $.
  Letting $\e  \n \to \n  0$ then yields to \eqref{eq:et383}. \qed

\bibliographystyle{siam}
\bibliography{RDOSRT_bib}

\end{document}